\documentclass[10pt,reqno]{article}

\usepackage{geometry}
\usepackage{amssymb}
\usepackage{amsmath}
\usepackage{amsthm}

\usepackage{tikz}
\usepackage{tikz-cd}
\usepackage{float}
\usepackage{caption}
\usepackage{subcaption}

\usepackage{scalerel}

\usepackage{enumitem}
\usepackage{hyperref}
\hypersetup{
  colorlinks   = true,
  urlcolor     = blue,
  linkcolor    = blue,
  citecolor   = red
}

\numberwithin{equation}{section}
\allowdisplaybreaks

\newtheorem{theorem}{Theorem}[section]
\newtheorem{definition}[theorem]{Definition}

\newtheorem{remark}[theorem]{Remark}
\newtheorem{lemma}[theorem]{Lemma}
\newtheorem{proposition}[theorem]{Proposition}

\newtheorem{assumption}[theorem]{Assumption}
\newtheorem{problem}[theorem]{Problem}

\newtheorem{construction}[theorem]{Construction}

\newtheorem*{assum}{\it{Assumption}}
\newtheorem*{step i}{\it{\bfseries{Step i}}}
\newtheorem*{step ii}{\it{\bfseries{Step ii}}}

\newcommand{\ed}{\mathrm{d} \makebox[0ex]{}}
\newcommand{\us}{\underline{s} \makebox[0ex]{}}

\newcommand{\circg}{\mathring{g} \makebox[0ex]{}}
\newcommand{\subcircg}{\mathring{g} \makebox[0ex]{}}

\newcommand{\uC}{\underline{C} \makebox[0ex]{}}
\newcommand{\uL}{\underline{L} \makebox[0ex]{}}

\newcommand{\uchi}{\underline{\chi} \makebox[0ex]{}}
\newcommand{\ueta}{\underline{\eta} \makebox[0ex]{}}
\newcommand{\uomega}{\underline{\omega} \makebox[0ex]{}}
\newcommand{\hatchi}{\hat{\chi} \makebox[0ex]{}}
\newcommand{\hatuchi}{\underline{\hat{\chi}} \makebox[0ex]{}}
\newcommand{\tr}{\mathrm{tr} \makebox[0ex]{}}

\newcommand{\slashg}{g  \mkern-7.2mu \scaleto{\pmb{\slash}}{1.6ex} \mkern+1mu \makebox[0ex]{}}  
\newcommand{\subslashg}{g \mkern-7.7mu \scaleto{\pmb{\slash}}{1ex} \mkern+1mu \makebox[0ex]{}}  

\newcommand{\ualpha}{\underline{\alpha} \makebox[0ex]{}}
\newcommand{\ubeta}{\underline{\beta} \makebox[0ex]{}}
\newcommand{\slashepsilon}{\epsilon \mkern-5.8mu \raisebox{0.12ex}[1.4ex][0ex]{$\scaleto{\pmb{\slash}}{1.4ex}$} \mkern+0.5mu \makebox[0ex]{}}
\newcommand{\Ric}{\mathrm{Ric} \makebox[0ex]{}}

\newcommand{\circnabla}{\mathring{\nabla} \makebox[0ex]{}}
\newcommand{\dvol}{\mathrm{dvol} \makebox[0ex]{}}

\newcommand{\ucalH}{\underline{\mathcal{H}} \makebox[0ex]{}}
\newcommand{\uh}{\underline{h} \makebox[0ex]{}}
\newcommand{\uf}{\underline{f} \makebox[0ex]{}}
\newcommand{\ufl}[1]{\mkern+1mu{}^{#1}\mkern-4mu \underline{f} \makebox[0ex]{}}
\newcommand{\slashd}{\mathrm{d} \mkern-8.9mu \raisebox{+0.3ex}[0ex][0ex]{$\scaleto{\pmb{\slash}}{2ex}$} \mkern+2.3mu \makebox[0ex]{}}
\newcommand{\udelta}{\underline{\delta} \makebox[0ex]{}}

\newcommand{\bSigma}{\bar{\Sigma} \makebox[0ex]{}}
\newcommand{\fl}[1]{\mkern+1mu{}^{#1}\mkern-4mu f \makebox[0ex]{}}
\newcommand{\os}{\overline{s} \makebox[0ex]{}}
\newcommand{\flt}{\mkern+1mu{}^{t}\mkern-4mu f \makebox[0ex]{}}
\newcommand{\uflt}{\mkern+1mu{}^{t}\mkern-4mu \underline{f} \makebox[0ex]{}}
\newcommand{\ue}{\underline{e} \makebox[0ex]{}}
\newcommand{\uvarepsilon}{\underline{\varepsilon} \makebox[0ex]{}}
\newcommand{\Xlt}{\mkern+1mu{}^{t}\mkern-4mu X \makebox[0ex]{}}
\newcommand{\relt}{\mkern+1mu{}^{t}\mkern-1mu \mathrm{re} \makebox[0ex]{}}
\newcommand{\ddcircnabla}{\ddot{\raisebox{0ex}[2.1ex][0ex]{$\mathring{\nabla}$}} \makebox[0ex]{}}
\newcommand{\circDelta}{\mathring{\Delta} \makebox[0ex]{}}
\newcommand{\ddcircDelta}{\ddot{\raisebox{0ex}[2.1ex][0ex]{$\mathring{\Delta}$}} \makebox[0ex]{}}
\newcommand{\ud}{\underline{d} \makebox[0ex]{}}
\newcommand{\ddcircdiv}{\ddot{\raisebox{0ex}[2.1ex][0ex]{$\mathring{\mathrm{div}}$}} \makebox[0ex]{}}

\newcommand{\bL}{\bar{L} \makebox[0ex]{}}
\newcommand{\buL}{\bar{\underline{L}} \makebox[0ex]{}}
\newcommand{\dpartial}{\dot{\partial} \makebox[0ex]{}}
\newcommand{\dL}{\dot{L} \makebox[0ex]{}}
\newcommand{\duL}{\dot{\underline{L}} \makebox[0ex]{}}

\newcommand{\dslashd}{\dot{\mathrm{d}} \mkern-8.9mu \raisebox{+0.3ex}[0ex][0ex]{$\scaleto{\pmb{\slash}}{2ex}$} \mkern+2.3mu \makebox[0ex]{}}
\newcommand{\db}{\dot{b} \makebox[0ex]{}}
\newcommand{\dslashg}{\dot{g}  \mkern-8.2mu \scaleto{\pmb{\slash}}{1.6ex} \mkern+1mu \makebox[0ex]{}}
\newcommand{\dslashepsilon}{\dot{\epsilon} \mkern-5.8mu \raisebox{0.12ex}[1.4ex][0ex]{$\scaleto{\pmb{\slash}}{1.4ex}$} \mkern+0.5mu \makebox[0ex]{}}

\newcommand{\dchi}{\dot{\chi} \makebox[0ex]{}} 
\newcommand{\dtr}{\dot{\mathrm{tr}} \makebox[0ex]{}}
\newcommand{\duchi}{\dot{\underline{\chi}} \makebox[0ex]{}}
\newcommand{\deta}{\dot{\eta} \makebox[0ex]{}}
\newcommand{\duomega}{\dot{\underline{\omega}} \makebox[0ex]{}}

\newcommand{\slashnabla}{\nabla \mkern-13mu \raisebox{0.3ex}[0ex][0ex]{$\scaleto{\pmb{\slash}}{1.7ex}$} \mkern+5mu \makebox[0ex]{}}
\newcommand{\sym}{\mathrm{sym} \makebox[0ex]{}}
\newcommand{\slashDelta}{\Delta \mkern-10mu \raisebox{0.4ex}[0ex][0ex]{$\scaleto{\pmb{\slash}}{1.7ex}$} \mkern+2mu \makebox[0ex]{}}
\newcommand{\slashGamma}{\Gamma \mkern-9.5mu \raisebox{0.4ex}[0ex][0ex]{$\scaleto{\pmb{\slash}}{1.7ex}$} \mkern+1mu \makebox[0ex]{}}
\newcommand{\circGamma}{\mathring{\Gamma} \makebox[0ex]{}}
\newcommand{\circtriangle}{\mathring{\triangle} \makebox[0ex]{}}

\newcommand{\dualpha}{\dot{\underline{\alpha}} \makebox[0ex]{}}
\newcommand{\dubeta}{\dot{\underline{\beta}} \makebox[0ex]{}}
\newcommand{\dsigma}{\dot{\sigma} \makebox[0ex]{}}
\newcommand{\drho}{\dot{\rho} \makebox[0ex]{}}
\newcommand{\dbeta}{\dot{\beta} \makebox[0ex]{}}
\newcommand{\dalpha}{\dot{\alpha} \makebox[0ex]{}}

\newcommand{\dslashnabla}{\dot{\nabla} \mkern-13mu \raisebox{0.3ex}[0ex][0ex]{$\scaleto{\pmb{\slash}}{1.7ex}$} \mkern+5mu \makebox[0ex]{}}
\newcommand{\dslashGamma}{\dot{\Gamma} \mkern-9.5mu \raisebox{0.4ex}[0ex][0ex]{$\scaleto{\pmb{\slash}}{1.7ex}$} \mkern+1mu \makebox[0ex]{}}

\newcommand{\flu}{\mkern+1mu{}^{u}\mkern-4mu f \makebox[0ex]{}}
\newcommand{\ddL}{\ddot{L} \makebox[0ex]{}}
\newcommand{\dduL}{\underline{\ddot{L}} \makebox[0ex]{}}
\newcommand{\ddpartial}{\ddot{\partial} \makebox[0ex]{}}
\newcommand{\dB}{\dot{B} \makebox[0ex]{}}
\newcommand{\ddslashg}{\ddot{g} \mkern-8.2mu \scaleto{\pmb{\slash}}{1.6ex} \mkern+1mu \makebox[0ex]{}}
\newcommand{\subddslashg}{\ddot{g} \mkern-7.7mu \scaleto{\pmb{\slash}}{1ex} \mkern+1mu \makebox[0ex]{}} 
\newcommand{\dvarepsilon}{\dot{\varepsilon} \makebox[0ex]{}}
\newcommand{\ddslashd}{\ddot{\raisebox{0ex}[1.5ex][0ex]{$\mathrm{d}$}} \mkern-8.9mu \raisebox{+0.3ex}[0ex][0ex]{$\scaleto{\pmb{\slash}}{2ex}$} \mkern+2.3mu \makebox[0ex]{}}
\newcommand{\ddslashepsilon}{\ddot{\epsilon} \mkern-5.8mu  \raisebox{0.12ex}[0ex][0ex]{$\scaleto{\pmb{\slash}}{1.4ex}$} \mkern+0.5mu \makebox[0ex]{}}

\newcommand{\dduchi}{\underline{\ddot{\chi}} \makebox[0ex]{}}
\newcommand{\ddchi}{\ddot{\chi} \makebox[0ex]{}}
\newcommand{\ddeta}{\ddot{\eta} \makebox[0ex]{}}
\newcommand{\circdot}{\makebox[2ex]{$\circ\mkern-7mu\cdot$} \makebox[0ex]{}}

\newcommand{\ddalpha}{\ddot{\alpha} \makebox[0ex]{}}
\newcommand{\ddbeta}{\ddot{\beta} \makebox[0ex]{}}
\newcommand{\ddrho}{\ddot{\rho} \makebox[0ex]{}}
\newcommand{\ddsigma}{\ddot{\sigma} \makebox[0ex]{}}
\newcommand{\ddubeta}{\underline{\ddot{\beta}} \makebox[0ex]{}}
\newcommand{\ddualpha}{\underline{\ddot{\alpha}} \makebox[0ex]{}}

\newcommand{\bal}[1]{{}^{#1} \mkern-1.5mu \bar{a} \makebox[0ex]{}}
\newcommand{\balu}{{}^{u} \mkern-1.5mu \bar{a} \makebox[0ex]{}}
\newcommand{\bslashg}{\bar{g} \mkern-8.2mu \scaleto{\pmb{\slash}}{1.6ex} \mkern+1mu \makebox[0ex]{}}
\newcommand{\subbslashg}{\bar{g} \mkern-7.7mu \scaleto{\pmb{\slash}}{1ex} \mkern+1mu \makebox[0ex]{}}
\newcommand{\bslashgl}[1]{{}^{#1} \mkern-1.5mu \bar{g} \mkern-8.2mu \scaleto{\pmb{\slash}}{1.6ex} \mkern+1mu \makebox[0ex]{}}
\newcommand{\bslashglu}{{}^{u} \mkern-1.5mu \bar{g} \mkern-8.2mu \scaleto{\pmb{\slash}}{1.6ex} \mkern+1mu \makebox[0ex]{}}
\newcommand{\subbslashgl}[1]{{}^{#1} \mkern-1.5mu \bar{g} \mkern-7.7mu \scaleto{\pmb{\slash}}{1ex} \mkern+1mu \makebox[0ex]{}}
\newcommand{\subbslashglu}{{}^{u} \mkern-1.5mu \bar{g} \mkern-7.7mu \scaleto{\pmb{\slash}}{1ex} \mkern+1mu \makebox[0ex]{}}
\newcommand{\bslashd}{\bar{\mathrm{d}} \mkern-8.9mu \raisebox{+0.2ex}[0ex][0ex]{$\scaleto{\pmb{\slash}}{2ex}$} \mkern+2.3mu \makebox[0ex]{}}

\newcommand{\buchi}{\underline{\bar{\chi}} \makebox[0ex]{}}
\newcommand{\buchil}[1]{{}^{#1}\mkern-1.5mu \underline{\bar{\chi}} \makebox[0ex]{}}
\newcommand{\buchilu}{{}^{u}\mkern-1.5mu \underline{\bar{\chi}} \makebox[0ex]{}}

\newcommand{\bchil}[1]{{}^{#1}\mkern-1.5mu \bar{\chi} \makebox[0ex]{}}
\newcommand{\bchilu}{{}^{u}\mkern-1.5mu \bar{\chi} \makebox[0ex]{}}

\newcommand{\btal}[1]{{}^{#1}\mkern-1.5mu \bar{\eta} \makebox[0ex]{}}
\newcommand{\btalu}{{}^{u}\mkern-1.5mu \bar{\eta} \makebox[0ex]{}}

\newcommand{\buomegal}[1]{{}^{#1}\mkern-1.5mu \underline{\bar{\omega}} \makebox[0ex]{}}
\newcommand{\buomegalu}{{}^{u}\mkern-1.5mu \underline{\bar{\omega}} \makebox[0ex]{}}

\newcommand{\balphal}[1]{{}^{#1}\mkern-1.5mu \bar{\alpha} \makebox[0ex]{}}
\newcommand{\bbetal}[1]{{}^{#1}\mkern-1.5mu \bar{\beta} \makebox[0ex]{}}
\newcommand{\bsigmal}[1]{{}^{#1}\mkern-1.5mu \bar{\sigma} \makebox[0ex]{}}
\newcommand{\brhol}[1]{{}^{#1}\mkern-1.5mu \bar{\rho} \makebox[0ex]{}}
\newcommand{\bubetal}[1]{{}^{#1}\mkern-1.5mu \bar{\underline{\beta}} \makebox[0ex]{}}
\newcommand{\bualphal}[1]{{}^{#1}\mkern-1.5mu \bar{\underline{\alpha}} \makebox[0ex]{}}

\newcommand{\balphalu}{{}^{u}\mkern-1.5mu \bar{\alpha} \makebox[0ex]{}}
\newcommand{\bbetalu}{{}^{u}\mkern-1.5mu \bar{\beta} \makebox[0ex]{}}
\newcommand{\bsigmalu}{{}^{u}\mkern-1.5mu \bar{\sigma} \makebox[0ex]{}}
\newcommand{\brholu}{{}^{u}\mkern-1.5mu \bar{\rho} \makebox[0ex]{}}
\newcommand{\bubetalu}{{}^{u}\mkern-1.5mu \bar{\underline{\beta}} \makebox[0ex]{}}
\newcommand{\bualphalu}{{}^{u}\mkern-1.5mu \bar{\underline{\alpha}} \makebox[0ex]{}}

\newcommand{\dR}{\dot{R} \makebox[0ex]{}}
\newcommand{\lie}{\mathcal{L} \makebox[0ex]{}}
\newcommand{\uhl}[1]{\mkern+1mu{}^{#1}\mkern-2mu \underline{h} \makebox[0ex]{}}
\newcommand{\uhlt}{\mkern+1mu{}^{t}\mkern-2mu \underline{h} \makebox[0ex]{}}
\newcommand{\Xl}[1]{\mkern+1mu{}^{#1}\mkern-4mu X \makebox[0ex]{}}
\newcommand{\rel}[1]{\mkern+1mu{}^{#1}\mkern-1mu \mathrm{re} \makebox[0ex]{}}
\newcommand{\lo}[1]{l.\{#1\} \makebox[0ex]{} \makebox[0ex]{}}
\newcommand{\hi}[1]{h.\{#1\} \makebox[0ex]{} \makebox[0ex]{}}

\newcommand{\circepsilon}{\mathring{\epsilon} \makebox[0ex]{}}
\newcommand{\hatdchi}{\hat{\dot{\chi}} \makebox[0ex]{}}
\newcommand{\hatduchi}{\hat{\dot{\underline{\chi}}} \makebox[0ex]{}}
\newcommand{\ddtr}{\ddot{\mathrm{tr}} \makebox[0ex]{}}
\newcommand{\hatddchi}{\hat{\ddot{\chi}} \makebox[0ex]{}}
\newcommand{\hatdduchi}{\hat{\ddot{\underline{\chi}}} \makebox[0ex]{}}
\newcommand{\uslashd}{\underline{d} \mkern-7.5mu \raisebox{+0.3ex}[0ex][0ex]{$\scaleto{\pmb{\slash}}{1.6ex}$} \mkern+2.3mu \makebox[0ex]{}}
\newcommand{\uflu}{\mkern+1mu{}^{u}\mkern-4mu \underline{f} \makebox[0ex]{}}

\newcommand{\bmul}[1]{{}^{#1}\mkern-1.5mu \bar{\mu} \makebox[0ex]{}}
\newcommand{\bmulu}{{}^{u}\mkern-1.5mu \bar{\mu} \makebox[0ex]{}}

\newcommand{\bKl}[1]{{}^{#1}\mkern-1.5mu \bar{K} \makebox[0ex]{}}
\newcommand{\bKlu}{{}^{u}\mkern-1.5mu \bar{K} \makebox[0ex]{}}

\newcommand{\btr}{\bar{\mathrm{tr}} \makebox[0ex]{}}

\newcommand{\ddKlu}{{}^{u}\mkern-1.5mu \ddot{K} \makebox[0ex]{}}

\newcommand{\ddchil}[1]{{}^{#1}\mkern-1.5mu \ddot{\chi} \makebox[0ex]{}}
\newcommand{\ddchilu}{{}^{u}\mkern-1.5mu \ddot{\chi} \makebox[0ex]{}}

\newcommand{\hatdduchil}[1]{{}^{#1}\mkern-1.5mu \hat{\ddot{\underline{\chi}}} \makebox[0ex]{}}
\newcommand{\hatdduchilu}{{}^{u}\mkern-1.5mu \hat{\ddot{\underline{\chi}}} \makebox[0ex]{}}

\newcommand{\hatddchil}[1]{{}^{#1}\mkern-1.5mu \hat{\ddot{\chi}} \makebox[0ex]{}}
\newcommand{\hatddchilu}{{}^{u}\mkern-1.5mu \hat{\ddot{\chi}} \makebox[0ex]{}}

\newcommand{\bslashdiv}{\bar{\mathrm{div}} \mkern-16.5mu \raisebox{+0.2ex}[0ex][0ex]{$\scaleto{\pmb{\slash}}{2ex}$} \mkern+9mu \makebox[0ex]{}}

\newcommand{\ddslashdiv}{\ddot{\mathrm{div}} \mkern-16.5mu \raisebox{+0.2ex}[0ex][0ex]{$\scaleto{\pmb{\slash}}{2ex}$} \mkern+9mu \makebox[0ex]{}}

\newcommand{\ddrhol}[1]{{}^{#1}\mkern-1.5mu \ddot{\rho} \makebox[0ex]{}}
\newcommand{\ddrholu}{{}^{u}\mkern-1.5mu \ddot{\rho} \makebox[0ex]{}}

\newcommand{\br}{\bar{r} \makebox[0ex]{}}

\newcommand{\bslashcurl}{\bar{\mathrm{curl}} \mkern-16.5mu \raisebox{+0.2ex}[0ex][0ex]{$\scaleto{\pmb{\slash}}{2ex}$} \mkern+9mu \makebox[0ex]{}}

\newcommand{\bslashDelta}{\bar{\Delta} \mkern-10mu \raisebox{0.4ex}[0ex][0ex]{$\scaleto{\pmb{\slash}}{1.7ex}$} \mkern+2mu \makebox[0ex]{}}

\newcommand{\bslashnabla}{\bar{\nabla} \mkern-13mu \raisebox{0.3ex}[0ex][0ex]{$\scaleto{\pmb{\slash}}{1.7ex}$} \mkern+5mu \makebox[0ex]{}}

\newcommand{\hatbuchil}[1]{{}^{#1}\mkern-1.5mu \hat{\bar{\underline{\chi}}} \makebox[0ex]{}}
\newcommand{\hatbuchilu}{{}^{u}\mkern-1.5mu \hat{\bar{\underline{\chi}}} \makebox[0ex]{}}

\newcommand{\hatbchil}[1]{{}^{#1}\mkern-1.5mu \hat{\bar{\chi}} \makebox[0ex]{}}
\newcommand{\hatbchilu}{{}^{u}\mkern-1.5mu \hat{\bar{\chi}} \makebox[0ex]{}}

\newcommand{\bslashepsilonl}[1]{{}^{#1}\mkern-1.5mu \bar{\epsilon} \mkern-5.8mu  \raisebox{0.12ex}[0ex][0ex]{$\scaleto{\pmb{\slash}}{1.4ex}$} \mkern+0.5mu \makebox[0ex]{}}
\newcommand{\bslashepsilonlu}{{}^{u}\mkern-1.5mu \bar{\epsilon} \mkern-5.8mu  \raisebox{0.12ex}[0ex][0ex]{$\scaleto{\pmb{\slash}}{1.4ex}$} \mkern+0.5mu \makebox[0ex]{}}

\newcommand{\ddmulu}{{}^{u}\mkern-1.5mu \ddot{\mu} \makebox[0ex]{}}

\newcommand{\dduchil}[1]{{}^{#1}\mkern-1.5mu \underline{\ddot{\chi}} \makebox[0ex]{}}
\newcommand{\dduchilu}{{}^{u}\mkern-1.5mu \underline{\ddot{\chi}} \makebox[0ex]{}}

\newcommand{\ddetal}[1]{{}^{#1}\mkern-1.5mu \ddot{\eta} \makebox[0ex]{}}
\newcommand{\ddetalu}{{}^{u}\mkern-1.5mu \ddot{\eta} \makebox[0ex]{}}

\newcommand{\duomegalu}{{}^{u}\mkern-1.5mu \underline{\dot{\omega}} \makebox[0ex]{}}

\newcommand{\ous}{\overline{\underline{s}} \makebox[0ex]{}}
\newcommand{\bcircnabla}{\bar{\raisebox{0ex}[2.1ex][0ex]{$\mathring{\nabla}$}} \makebox[0ex]{}}

\newcommand{\ddualphal}[1]{{}^{#1}\mkern-1.5mu \ddot{\underline{\alpha}} \makebox[0ex]{}}
\newcommand{\ddubetal}[1]{{}^{#1}\mkern-1.5mu \ddot{\underline{\beta}} \makebox[0ex]{}}
\newcommand{\ddsigmal}[1]{{}^{#1}\mkern-1.5mu \ddot{\sigma} \makebox[0ex]{}}
\newcommand{\ddbetal}[1]{{}^{#1}\mkern-1.5mu \ddot{\beta} \makebox[0ex]{}}
\newcommand{\ddalphal}[1]{{}^{#1}\mkern-1.5mu \ddot{\alpha} \makebox[0ex]{}}

\newcommand{\osc}[1]{\mathrm{osc} \{ #1 \} \makebox[0ex]{}}
\newcommand{\barR}{\bar{R} \makebox[0ex]{}}

\newcommand{\bP}{\bar{P} \makebox[0ex]{}}
\newcommand{\bN}{\bar{N} \makebox[0ex]{}}
\newcommand{\itbSigma}{\mathit{\bar{\Sigma}} \makebox[0ex]{}}
\newcommand{\itbXi}{\mathit{\bar{\Xi}} \makebox[0ex]{}}

\newcommand{\Xlu}{\mkern+1mu{}^{u}\mkern-4mu X \makebox[0ex]{}}
\newcommand{\Mlu}{\mkern+1mu{}^{u}\mkern-1mu \mathrm{M} \makebox[0ex]{}}
\newcommand{\relu}{\mkern+1mu{}^{u}\mkern-1mu \mathrm{re} \makebox[0ex]{}}

\newcommand{\philu}{\mkern+1mu{}^{u}\mkern-4mu \phi \makebox[0ex]{}}
\newcommand{\circdiv}{\mathring{\mathrm{div}} \makebox[0ex]{}}

\newcommand{\calA}{\mathcal{A} \makebox[0ex]{}}
\newcommand{\dmu}{\dot{\mu} \makebox[0ex]{}}
 
\newcommand{\Vvert}{\vert \mkern-1.5mu \vert \mkern-1.5mu \vert \makebox[0ex]{}}

\title{\textsc{Global Existence and Geometry of Constant Mass Aspect Function Foliation in Perturbed Schwarzschild Spacetime}}
\author{Pengyu Le}
\newcommand{\Address}{{
  \bigskip
  \footnotesize
  \textsc{Beijing Institute of Mathematical Sciences and Applications, Beijing, China}
  
  \textit{E-mail address}: \texttt{pengyu.le@bimsa.cn}
}}
\date{}

\begin{document}

\maketitle

\begin{abstract}
The constant mass function foliation has been shown useful for studying the null Penrose inequality on a null hypersurface, because of the monotonicity formula of Hawking mass along such a foliation. In this paper, we show the global existence of the constant mass aspect function foliation on a nearly spherically symmetric incoming null hypersurface, emanating from a spacelike surface near the apparent horizon to the past null infinity in a vacuum perturbed Schwarzschild spacetime. Moreover, we study the geometry of the constant mass aspect function foliation, by comparing with the spherically symmetric foliation in the Schwarzschild spacetime. The knowledge about the geometry of the foliation is essential for investigating the perturbation of the constant mass aspect function foliation, which is the core in the application to the null Penrose inequality for a vacuum perturbed Schwarzschild spacetime.
\end{abstract}

\tableofcontents

\section{Introduction}\label{sec 1}
The constant mass aspect function foliation was introduced first in the monumental proof of the global nonlinear stability of the Minkowski spacetime \cite{CK1993} by Christodoulou and Klainerman on a maximal spacelike hypersurface. It plays an important role in the continuity argument in the proof, which resolves the difficulty of constructing the optical functions one order smoother than the metric. We refer to the well written exposition on this issue in \cite{C2008}.

Another motivation to consider the constant mass aspect function foliation is due to an elegant variation formula of the Hawking mass along an arbitrary foliation in a spacelike hypersurface appearing first in \cite{C2008}. The variation formula in the case of an inverse mean curvature flow was already well known, obtained in \cite{G1973} by Geroch, leading to the Geroch monotonicity of the Hawking mass along an inverse mean curvature flow as an argument in support of the positive energy conjecture. The positive energy conjecture was proved by Schoen and Yau in \cite{SY1979} using minimal surfaces and later by Witten in \cite{W1981} using spinors. Geroch's approach was modified by Jang and Wald in \cite{JW1977} for an application to the Penrose inequality on a time symmetric initial data. This approach was finally accomplished in the proof of the Riemannian Penrose inequality by Huisken and Ilmanen in \cite{HI2001}. At the same time, Bray in \cite{Br2001} proved the Riemannian Penrose inequality independently by a different method.

In this paper, we concern the constant mass aspect function foliation in a null hypersurface defined first in \cite{C2003}. Similarly as its analogy in a spacelike surface, its counterpart in a null hypersurface also has the property of the smoothing effect as in \cite{CK1993} and the similar variation formula of the Hawking mass as in \cite{C2008}. The variation formula along the null direction was first derived by Hawking in \cite{H1968} for a general spacetime, and in particular he considered the special foliation parameterised by a luminosity parameter, which could be viewed as the analogy of the inverse mean curvature flow on a null hypersurface. The variation formula for the constant mass aspect function foliation was first derived by Christodoulou in \cite{C2003}.

The formula leads to the monotonicity of Hawking mass along both the foliation parameterised by a luminosity parameter and the constant mass aspect function foliation on a nearly spherically symmetric null hypersurface. Motivated by this property, Sauter in his thesis \cite{S2008} initiated the project to employ both foliations to prove the null Penrose inequality. This approach is similar to the inverse mean curvature flow approach to the Riemannian Penrose inequality in \cite{JW1977} and \cite{HI2001}. Although \cite{S2008} didnot prove the null Penrose inequality for a nearly spherically symmetric null hypersurface, it recognised the main obstacle of this approach towards the proof being the asymptotic geometry of the foliation at null infinity, which determines whether the Hawking mass converges to the Bondi mass along the foliation.

In order to overcome the above obstacle, Christodoulou and Sauter proposed the project to investigate the change of the asymptotic geometry at null infinity when varying the two foliations. The author's works \cite{L2020} \cite{L2022} \cite{L2023} are preliminary steps to carry out their proposal for the constant mass aspect function foliation. In \cite{L2020}, we investigated the structure of the set of marginally trapped surface in a perturbed Schwarzschild spacetime, which prepares us to study the constant mass aspect function foliation emanating from a marginally trapped surface. In \cite{L2022}, we show the global existence of null hypersurfaces from a spacelike surface near the apparent horizon to the past null infinity in a perturbed Schwarzschild exterior, based on which we can construct the constant mass aspect function foliation on such null hypersurfaces. In \cite{L2023}, we obtained a detailed result on the linearised perturbation of the constant mass aspect function foliation for the spherically symmetric foliation in a Schwarzschild spacetime, which serves as a model to illustrate the structure of the linearised perturbation of the foliation in a perturbed Schwarzschild spacetime.

In this paper, we proceed with Christodoulou and Sauter's proposal by one step further. We shall prove the global existence of the constant mass aspect function foliation emanating from an initial leaf near the apparent horizon to the past null infinity in a vacuum perturbed Schwarzschild spacetime. Moreover we shall investigate in details the geometry of the foliation by comparing with the spherically symmetric foliation in the Schwarzschild spacetime.

We briefly describe the vacuum perturbed Schwarzschild metric considered in this paper.
\begin{definition}[Rough version of vacuum perturbed Schwarzschild metrics, see definition \ref{def 2.4}]
Let $g_S$ be the Schwarzschild metric taking the following form in the double null coordinate system $\{s, \us, \theta^1, \theta^2 \}$
\begin{align*}
	g_S
	= 
	2 \Omega_S^2 ( \ed s\otimes \ed \us + \ed \us \otimes \ed s  ) + r_S^2 \cdot \circg_{ab} \ed \theta^a \otimes \ed \theta^b,
\end{align*}
where
\begin{align*}
\Omega_S^2 = \frac{s+r_0}{r_S}\exp\frac{\us+s+r_0-r_S}{r_0},
\quad
(r_S-r_0)\exp\frac{r_S}{r_0} = s\exp\frac{\us+s+r_0}{r_0}.
\end{align*}
Let $M_{\kappa}$ be a neighbourhood of the null hypersurface $\uC_{\us=0,s\geq 0}$. Let $g_{\epsilon}$ be a vacuum metric on $M_{\kappa}$ taking the form
\begin{align*}
	g_{\epsilon}
	=
	2\Omega_{\epsilon}^2  ( \ed s  \otimes \ed \us  +\ed \us \otimes \ed s ) 
	+ 
	( \slashg_{\epsilon})_{ab} 
	( \ed \theta^a - b_{\epsilon}^{a} \ed s) \otimes ( \ed \theta^b - b_{\epsilon}^{b} \ed s).
\end{align*}
Let $r_{\epsilon}$ be the area radius of the surface $(\Sigma_{s,\us}, \slashg_{\epsilon})$. We say that a $(N+2)$-th differentiable metric $g_{\epsilon}$ is $(n+2)$-th order $\epsilon$-close to the Schwarzschild metric $g_S$ in $M_{\kappa}$ if the difference between $g_{\epsilon}$ and $g_S$, and the derivatives of the difference up to the $(n+2)$-th order are controlled by certain bound $\epsilon r_0^{\tau} r_{\epsilon}^{\lambda}$ with suitable decaying rate $\lambda$ at the past null infinity. Moreover we assume that the differences between the connection coefficients, the curvature components are also bounded by suitable $\epsilon r_0^{\tau} r_{\epsilon}^{\lambda}$. The decaying rates are taken as the same as the ones in the nonlinear stability of Minkowski spacetime \cite{CK1993}.
\end{definition}

We sketch the main theorem on the global existence of the constant mass aspect function foliation.
\begin{theorem}[Sketch of main theorem \ref{thm 8.2}]
Let $(M_{\kappa}, g_{\epsilon})$ be a vacuum perturbed Schwarzschild metric in definition \ref{def 2.4}. Let $\bSigma_{u=0}$ be a spacelike surface in an incoming null hypersurface $\ucalH$. If $\ucalH$ is nearly spherically symmetric and sufficiently close to a null hypersurface $\uC_{\us}$, and $\bSigma_{u=0}$ is sufficiently close to the Schwarzschild horizon $C_{s=0}$, then there exists a constant mass aspect function foliation $\{\bSigma_u\}$  on $\ucalH$ emanating from $\bSigma_{u=0}$ to the past null infinity. The geometry along the foliation is controlled by $\epsilon$, the size of the perturbation of the metric $g_{\epsilon}$ from $g_S$, and the size of the perturbation of $\bSigma_{u=0}$ from some $\Sigma_{s=0,\us}$ in the Schwarzschild horizon.
\end{theorem}

In the following, we briefly explain the important aspects and ideas to prove the main theorem.
\begin{enumerate}[label=\alph*.]
\item
Find the appropriate method to parameterise a foliation $\{\bSigma_u\}$ in an incoming null hypersurface $\ucalH$. We present two methods in section \ref{sec 3}. 

The first one is to parameterise the incoming null hypersurface $\ucalH$ by a function $\uh(s,\vartheta)$ as its graph in the double null coordinate system, then parameterise each leaf $\bSigma_u$ by a function $\flu(\vartheta)$ as its graph in the coordinate system $\{s, \vartheta\}$ of $\ucalH$. The parameterisation function $\uh$ of $\ucalH$ is determined by its restriction $\uh(s=0,\vartheta)$ at $s=0$, which is defined as $\ufl{s=0}(\vartheta)$. Then we use $(\ufl{s=0}, \flu)$ to parameterise $\{ \bSigma_u \}$.

The second one is to parameterise each leaf $\bSigma_u$ by a pair of functions $(\us, s) =(\uflu(\vartheta), \flu(\vartheta))$ as its graph of $(\us,s)$ over the $\vartheta$ coordinate.

The transformation between two kinds of parameterisation is given in section \ref{sec 3.3}.

These two methods of parameterisation are also used to parameterise spacelike surfaces in \cite{L2020}.

\item
Derive the formulae of the geometric quantities associated with a foliation $\{ \bSigma_u \}$. We obtain these formulae in section \ref{sec 4}.

The derivation consists of two steps. The first step is to derive the formulae of the geometric quantities associated with the special foliation $\{\Sigma_s = C_s \cap \ucalH \}$ in terms of the background quantities on $\Sigma_s$ and parameterisation function $\uh$ of $\ucalH$. Then the second step is to derive the formulae of the geometric quantities associated with $\{\bSigma_u\}$ in terms of the geometric quantities associated with $\{\Sigma_s\}$ and the parameterisation function $\flu$ of $\bSigma_u$.

The formulae are also derived by the above method in the special case of the Schwarzschild spacetime in \cite{L2023}.

\item
With the help of the formulae derived in section \ref{sec 4}, we obtain the estimates of the geometric quantities associated with $\{\bSigma_u\}$ in terms of $\epsilon$, the size of the perturbation of the metric $g_{\epsilon}$ from $g_S$, and the sizes of the parameterisation functions $\ufl{s=0}$ and $\fl{u=0}$.

In the process of obtaining the estimates, we decompose the geometric quantities into two parts: the first order main part and the high order remainder in section \ref{sec 5.2}. The decompositions reveal the structure of the formulae of the geometric quantities.

\item
To study the geometry of the constant mass aspect function foliation, we employ a system collecting the basic equations satisfied by the parameterisation function $\flu$, the lapse $\balu$ and the geometric quantities associated with the foliation in section \ref{sec 6.2}.

The system consists of propagation equations and elliptic equations \eqref{eqn 6.3} - \eqref{eqn 6.13}. This propagation-elliptic system is used in \cite{L2023} to calculate the linearised perturbations of the geometric quantities of the foliation.

\item
In section \ref{sec 7}, we obtain the geometry of the initial leaf $\bSigma_{u=0}$ near a surface $\Sigma_{s=0,\us}$ in the Schwarzschild horizon $C_{s=0}$.

The estimates of the geometric quantities in section \ref{sec 4} give necessary information on the geometry of $\bSigma_{u=0}$, except for the lapse function $\bal{u=0}$. However, the estimates in section \ref{sec 4} donot give the optimal regularities for the torsion $\btal{u=0}$ and the acceleration $\buomegal{u=0}$. This difficulty is overcome through the elliptic equations in section \ref{sec 6}.

The lapse function $\bal{u=0}$ is estimated through the inverse lapse equation \eqref{eqn 6.5}
\begin{align}
	&
	\left\{
	\begin{aligned}
		&
		\bslashDelta \log \balu  
		= 
		- (\ddrholu - \overline{\ddrholu})
		- \frac{1}{2} [ (\hatdduchilu, \hatddchilu') - \overline{( \hatdduchilu, \hatddchilu')}^u ]
		- \bslashdiv \ddetalu, 
	\\
		&
		\overline{\balu\, \ddtr \dduchilu}^u
		=\overline{\btr \buchilu}^u 
		= \frac{2}{\br_u}.
	\end{aligned}
	\right.
	\tag{\ref{eqn 6.5}\ensuremath{'}}
\end{align}
The torsion $\btal{u=0}$ is estimated with the optimal regularity through the elliptic system \eqref{eqn 6.12}
\begin{align}
	&
	\left\{
	\begin{aligned}
		&
		\bslashcurl \btalu
		=
		\frac{1}{2} \hatbchilu' \wedge \hatbuchilu 
		+ \bsigmalu,
	\\
		&
		\bslashdiv \btalu
		= 
		- \brholu
		- \frac{1}{2} ( \hatbuchilu, \hatbchilu' ) 
		- \bmulu,
	\end{aligned}
	\tag{\ref{eqn 6.12}\ensuremath{'}}
	\right.
\end{align}
The acceleration $\buomegal{u=0}$ is estimated with the optimal regularity through the elliptic equation \eqref{eqn 6.13}
\begin{align}
	&
	\left\{
	\begin{aligned}
		2 \bslashDelta \buomegalu
		&=
		-\frac{3}{2} ( \bmulu\, \btr \buchilu - \overline{\bmulu\, \btr \buchilu}^u )
		+ \frac{1}{2} ( \btr \buchilu | \btalu |^2 - \overline{\btr \buchilu | \btalu |^2}^u )
	\\
		&
		\phantom{=}
		+ \frac{1}{4} ( \btr \bchilu' | \hatbuchilu |^2 - \overline{ \btr \bchilu' | \hatbuchilu|^2}^u )
		+ 4 (\bslashdiv \hatbuchilu, \btalu )
		+ 4 ( \hatbuchilu, \bslashnabla \btalu )
		- 2 \bslashdiv \bubetalu,
	\\
		\overline{\buomegalu}^u
		&=
		-\frac{r_u}{2} \overline{(\buomegalu - \overline{\buomegalu}^u)(\btr\buchilu - \overline{\btr\buchilu}^u)}^u
		- \frac{r_u}{8} \overline{(\btr \buchilu - \overline{\btr \buchilu}^u)^2}^u
		+ \frac{r_u}{4} \overline{|\hatbuchilu|^2}^u,
	\end{aligned}
	\tag{\ref{eqn 6.13}\ensuremath{'}}
	\right.
\end{align}
The estimates of the shears $\hatbuchil{u=0}$, $\hatbchil{u=0}'$ are also improved through the elliptic equations \eqref{eqn 6.10} \eqref{eqn 6.11}
\begin{align}
	&
	\bslashdiv \hatbuchilu
	- \frac{1}{2} \bslashd \btr \buchilu
	- \hatbuchilu \cdot \btalu 
	+ \frac{1}{2} \btr \buchilu \, \btalu
	=
	- \bubetalu,
	\tag{\ref{eqn 6.10}\ensuremath{'}}
\\
	&
	\bslashdiv \hatbchilu'
	- \frac{1}{2} \bslashd \btr \bchilu'
	+ \hatbchilu' \cdot \btalu 
	- \frac{1}{2} \btr \bchilu' \, \btalu
	=
	- \bbetalu,
	\tag{\ref{eqn 6.11}\ensuremath{'}}
\end{align}

\item
We adopt the bootstrap argument to prove the global existence of the foliation $\{\bSigma_u\}$. One way to study the existence of the foliation is to consider the inverse lapse system consisting of
\begin{align}
	&
	\uL \flu = \balu.
	\tag{\ref{eqn 6.3}\ensuremath{'}}
\end{align}
and the inverse lapse equation \eqref{eqn 6.5} of $\balu$. However the inverse lapse system has the problem of regularity loss as follows,
\begin{align*}
	\flu \in \mathrm{W}^{n+2,p}(\mathbb{S}^2)
	\underset{i.}{\rightarrow}
	\ddetalu \in \mathrm{W}^{n,p}(\mathbb{S}^2)
	\underset{ii.}{\rightarrow}
	\balu \in \mathrm{W}^{n+1,p}(\mathbb{S}^2)
	\underset{iii.}{\rightarrow}
	\flu \in \mathrm{W}^{n+1,p}(\mathbb{S}^2).
\end{align*}
\begin{enumerate}[label=\textit{\roman*}.]
\item
By propositions \ref{prop 5.9}, \ref{prop 5.10} on the estimates of the geometric quantities associated with the foliation.

\item
By the inverse lapse equation \eqref{eqn 6.5}.

\item
Propagation equation \eqref{eqn 6.3}.
\end{enumerate}

To overcome the above regularity loss, we shall consider another elliptic equation \eqref{eqn 8.1} for $\flu$
\begin{align}
	\balu^{-1} \btr \bchilu' 
	= 
	- 2 (r_S|_{\bSigma_u})^{-2} \circDelta \flu 
	+ \tr \chi'_S|_{\bSigma_u}
	+ \hi{\ddtr \ddchilu'}.
	\tag{\ref{eqn 8.1}\ensuremath{'}}
\end{align}
Although it is sufficient to obtain the parameterisation function $\flu$ for $u\in[0,+\infty)$ to describe the foliation $\{\bSigma_u\}$, we need to consider the geometry of the foliation simultaneously, since the above equations \eqref{eqn 6.5}, \eqref{eqn 8.1} are coupled with the geometric quantities associated with the foliation. Therefore to adopt the bootstrap argument, we shall consider equation \eqref{eqn 8.1} and the propagation-elliptic system in section \ref{sec 6} together.

\item
One major difficulty of the bootstrap argument is to come up with the appropriate bootstrap assumption. The bootstrap assumption involves the assumptions not only on the parameterisation function $\flu$, but also on the lapse $\balu$, the metric components, the connection coefficients and the curvature components. We successfully find the appropriate assumption \ref{assum 8.1} and show that the estimates in assumption \ref{assum 8.1} can be improved on any closed interval $u\in[0,u_a]$, which is the main task of sections \ref{sec 8.4} - \ref{sec 8.9}.

\item
Another difficulty of the bootstrap argument is to prove that one can always extend the bootstrap assumption to a slightly larger interval. Given that the bootstrap assumption is satisfied in $[0,u_a]$, then the stronger estimates than the ones in the bootstrap assumption are satisfied at $u=u_a$. To extend the bootstrap assumption beyond $u_a$, we need the local existence result of the foliation proved in \cite{S2008}.

However the regularities of $\ucalH$ and $\bSigma_{u_a}$ are not subject to the condition of the local existence result, thus we shall approximate $\ucalH$, $\bSigma_{u_a}$ by more regular ones. Moreover we need the more differentiabilities for the metric $g_{\epsilon}$. The additional differentiability assumption for the metric $g_{\epsilon}$ is a technical assumption since there are no bounds required for the additional derivatives of the metric.

In section \ref{sec 8.10}, we construct the regularised approximation for the foliation $\{\bSigma_u \}$ near $u=u_a$, and show that there exists a uniform existence interval $[u_a - \tau_u, u_a+\tau_u]$ for the regularised approximation where assumption \ref{assum 8.1} holds. Then the limit of the regularised approximation gives the local existence of the foliation $\{\bSigma_u \}$ near $u=u_a$, and the assumption \ref{assum 8.1} is preserved in the limit by the weak convergence.
\end{enumerate}

In the last section \ref{sec 9}, we use the geometry of the foliation to study the asymptotic geometry at the past null infinity. In particular, we discuss the asymptotic reference frame defined by an asymptotically round foliation and the associated energy-momentum vector. These concepts are important for the application to the null Penrose inequality.

\section{Vacuum perturbed Schwarzschild spacetime}\label{sec 2}

In this section, we introduce the class of vacuum perturbed Schwarzschild spacetimes considered in the paper. We first review the Schwarzschild metric in a double null coordinate system, then give the definition of the vacuum perturbed Schwarzschild metric by comparing the geometric quantities in the double null coordinate system.

\subsection{Review of the Schwarzschild metric}\label{sec 2.1}
In the coordinate system $\{t,r,\theta, \phi\}$, the Schwarzschild metric $g_S$ takes the form
\begin{align*}
g_S=-( 1-\frac{2m}{r} ) \ed t^2 + ( 1-\frac{2m}{r} )^{-1} \ed r ^2 + r^{2} ( \ed \theta^2 + \sin^2 \theta \ed \phi^2 ).
\end{align*}
Define $r_0=2m$, where $m$ is the mass, and $r_0$ is the area radius of the event horizon of the black hole. We use $(\mathcal{S}, g_S)$ to denote the Schwarzschild spacetime. Introduce the coordinate transformation between $\{t, r\}$ and $\{\us, s\}$ as in \cite{L2018}\cite{L2022}\cite{L2023}
\begin{align*}
\left\{
\begin{aligned}
&
(r-r_0)^{\frac{1}{2}} \exp \frac{t+r}{2r_0}
=
\exp\frac{\us}{r_0},
\\
&
(r-r_0)^{\frac{1}{2}} \exp \frac{-t+r}{2r_0}
=
s \exp \frac{s+r_0}{r_0},
\end{aligned}
\right.
\end{align*}
then $g_S$ takes the form
\begin{align*}
g_S
= 
2 \Omega_S^2 ( \ed s\otimes \ed \us + \ed \us \otimes \ed s  ) + r_S^2 ( \ed \theta^2 + \sin^2 \theta \ed \phi^2 ),
\end{align*}
where
\begin{align*}
\Omega_S^2 = \frac{s+r_0}{r_S}\exp\frac{\us+s+r_0-r_S}{r_0},
\quad
(r_S-r_0)\exp\frac{r_S}{r_0} = s\exp\frac{\us+s+r_0}{r_0}.
\end{align*}
We shall use $\circg$ to denote the standard metric on the sphere of radius $1$ that $\circg = \ed \theta^2 + \sin^2 \theta \ed \phi^2$. Let $\{\theta^1, \theta^2\}$ be a coordinate system of the sphere, then we also use $\{\us, s, \theta^1, \theta^2\}$ as a coordinate system of the Schwarzschild spacetime.

The coordinate system $\{\us, s, \theta^1, \theta^2\}$ is a double null coordinate system, as $\us$, $s$ are both optical functions, i.e. the level sets of $\us$, $s$ are both null hypersurfaces. Introduce the following notations:
\begin{enumerate}[label=\footnotesize\textbullet]
\item $C_s$: the level set of $s$. $\uC_{\us}$: the level set of $\us$. $\Sigma_{s, \us}$: the intersection of $C_s$ and $\uC_{\us}$.
\item $L$: tangential null vector field on $C_s$ with $L \us=0$. $\uL$: tangential null vector field on $\uC_{\us}$ with $\uL s=1$.
\item $L'$: tangential null vector field on $C_s$ conjugate to $\uL$ that $g(L', \uL)=2$. $\uL'$: tangential null vector field on $\uC_{\us}$ conjugate to $L$ that $g(\uL', L)=2$.
\end{enumerate}
The above notations are valid for an arbitrary double null coordinate system in a general $4$-dimensional Lorentzian manifold. For a general $4$-dimensional Lorentzian manifold $(M,g)$, we can define the double null coordinate system $\{\us, s, \theta^1, \theta^2\}$ similarly by requiring that $\us$, $s$ are both optical functions. The metric $g$ in an arbitrary double null coordinate system takes the form
\begin{align*}
g=2\Omega^2  ( \ed s  \otimes \ed \us  +\ed \us \otimes \ed s ) + \slashg_{ab} ( \ed \theta^a - b^{a} \ed s - \underline{b}^{a} \ed \us ) \otimes ( \ed \theta^b - b^{b} \ed s - \underline{b}^{b} \ed \us ),
\end{align*}
where $\slashg$ is the intrinsic metric on the spacelike surface $\Sigma_{s,\us}$. Up to coordinate transformations of $\{\theta^1, \theta^2\}$ on each $\Sigma_{s,\us}$, we can find a double null coordinate system in which the metric $g$ takes the form
\begin{align*}
g=2\Omega^2  ( \ed s  \otimes \ed \us  +\ed \us \otimes \ed s ) + \slashg_{ab} ( \ed \theta^a - b^{a} \ed s ) \otimes ( \ed \theta^b - b^{b} \ed s ).
\end{align*}
In this paper, we always choose the above particular double null coordinate system.

Associated with the double null coordinate system $\{\us, s, \theta^1, \theta^2\}$, we can define the corresponding connection coefficients and curvature components.
\begin{definition}[Connection coefficients]\label{def 2.1}
Let $\{s, \us, \theta^1, \theta^2\}$ be a double null coordinate system of a $4$-dimensional Lorentzian manifold $(M,g)$. We define that
\begin{align*}
&
\uchi_{ab} = g(\nabla_{\partial_a} \uL, \partial_b),
&&
\chi_{ab} = g(\nabla_{\partial_a} L, \partial_b),
\\
&
\uchi'_{ab} = g(\nabla_{\partial_a} \uL', \partial_b),
&&
\chi'_{ab} = g(\nabla_{\partial_a} L', \partial_b),
\\
&
\eta_a = \frac{1}{2} g(\nabla_{\partial_a} \uL, L'),
&&
\ueta_a = \frac{1}{2} g(\nabla_{\partial_a} L, \uL'),
\\
&
\omega = \frac{1}{4} g(\nabla_L L, \uL') =  L \log \Omega
&&
\uomega = \frac{1}{4} g(\nabla_{\uL}\uL, L') =  \uL \log \Omega,
\end{align*}
where subscripts $a,b$ denote indices in $\{1,2\}$.
We decompose $\chi, \uchi$ into their traces and trace-free parts relative to $\slashg$:
\begin{align*}
\chi = \hatchi + \frac{1}{2} \tr \chi \slashg,
\quad
\uchi = \hatuchi + \frac{1}{2} \tr \uchi \slashg,
\end{align*}
and similarly for $\chi', \uchi'$. $\hatchi$, $\hatuchi$ are called shears and $\tr\chi$, $\tr\uchi$ are called null expansions of $\Sigma_{s,\us}$ relative to $L$, $\uL$ respectively.
\end{definition}
\begin{definition}[Curvature components]\label{def 2.2}
Let $\{s, \us, \theta^1, \theta^2\}$ be a double null coordinate system of a $4$-dimensional vacuum spacetime $(M,g)$, i.e. $\Ric=0$. 
\begin{align*}
&
\ualpha_{ab} = \mathrm{R}_{\uL a b \uL},
&&
\alpha_{ab}= \mathrm{R}_{L' a b L'},
\\
&
\ubeta_a = \frac{1}{2} \mathrm{R}_{a \uL L' \uL},
&&
\beta_a = \frac{1}{2} \mathrm{R}_{a L' \uL L'},
\\
&
\rho = \frac{1}{4} \mathrm{R}_{L' \uL \uL L'},
&&
\sigma \slashepsilon_{ab} = \frac{1}{2} \mathrm{R}_{a b \uL L'},
\end{align*}
where subscripts $a,b$ denote indices in $\{1,2\}$ and $\slashepsilon$ is the volume form of $(\Sigma_{s,\us}, \slashg)$.\footnote{Convention for the curvature tensor: $\mathrm{R}(X,Y)Z = \nabla_X \nabla_Y Z - \nabla_Y \nabla_X Z - \nabla_{[X,Y]} Z$, $\mathrm{R}(X,Y,Z,W) = g(\mathrm{R}(X,Y)Z, W)$.}
\end{definition}
The above curvature components contain all the information of the curvature tensor, because
\begin{align*}
&
\mathrm{R}_{a b c \uL} = \slashepsilon_{ab} \slashepsilon_{cd} \ubeta^d,
&&
\mathrm{R}_{a b c L'} = \slashepsilon_{ab} \slashepsilon_{cd} \beta^d,
\\
&
\mathrm{R}_{a b c d} = \rho \slashepsilon_{ab} \slashepsilon_{cd},
&&
\mathrm{R}_{\uL a b L'} = \rho \slashg_{ab} - \sigma \slashepsilon_{ab}.
\end{align*}

Adopting the above definitions to the double null coordinate system $\{s, \us, \theta^1, \theta^2\}$ of the Schwarzschild spacetime, we obtain the corresponding connection coefficients and curvature components: for the sake of brevity, we omit the subscript $S$ denoting the quantities of the Schwarzschild spacetime in the following formulae,
\begin{align}
	\begin{aligned}
		&
		\begin{aligned}
		\partial_{\us} r= \frac{r-r_0}{r}, 
		\quad
		\partial_s r= \frac{s+r_0}{r} \cdot \frac{r-r_0}{s},
		\end{aligned}
	\\
		&
		\left\{
		\begin{aligned}
			&
			\tr \chi'= \frac{2s}{r(s + r_0)},
		\\	
			&
			\tr\uchi = \frac{2(s+r_0)}{r^2} \cdot \frac{r-r_0}{s} = \frac{2(s+r_0)}{r^2} \exp\frac{\us+s+r_0-r}{r_0},
		\\
			&
			\hatchi' = \hatuchi =0,
		\\
			&
			\eta = \ueta= 0,
		\\
			&
			\omega =\frac{r_0}{2r^2}, 
		\\	
			&
			\uomega =\frac{1}{2(s+r_0)} +\frac{1}{2r_0} - \Big( \frac{1}{2r} + \frac{1}{2r_0} \Big)  \frac{s+r_0}{r} \exp \frac{\us+s+r_0-r}{r_0},
		\end{aligned}
		\right.
	\end{aligned}
	\label{eqn 2.1}
\end{align}
and
\begin{align}
\alpha = \ualpha =0,
\quad
\beta = \ubeta = 0,
\quad
\sigma=0,
\quad
\rho = - \frac{ r_0}{r^3}.
\label{eqn 2.2}
\end{align}

\subsection{Vacuum perturbed Schwarzschild metric}\label{sec 2.2}
We introduce first a neighbourhood of a spherically symmetric incoming null hypersurface in the Schwarzschild spacetime, then define the class of vacuum perturbed Schwarzschild metrics considered in this paper.
\begin{definition}[$\kappa$-neighbourhood $M_{\kappa}$ of $\uC_{\us=0, s\geq 0}$]\label{def 2.3}
Let $\{\us,s\}$ be the double null coordinates of the Schwarzschild spacetime $(\mathcal{S},g_S )$ introduced in section \ref{sec 2.1}. We define the $\kappa$-neighbourhood $M_{\kappa}$ of the truncated null hypersurface $\uC_{\us=0, s\geq 0}$ as the following open set of the Schwarzschild spacetime:
\begin{align*}
M_{\kappa} = \{ p\in \mathcal{S}: s(p)>-\kappa r_0, \vert \us \vert < \kappa r_0  \}.
\end{align*}
See figure \ref{fig 1}. In this paper, we assume that $\kappa\leq 0.1$.
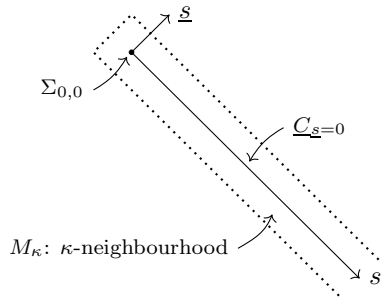
\begin{figure}[H]
\centering
\begin{tikzpicture}
\draw[->] (0,0) -- (3,-3) node[right] {$s$};
\draw[->] (0,0) -- (0.5,0.5) node[right] {$\us$};
\draw[dotted,thick] (3.25, -2.75) -- (0,0.5) -- (-0.5,0) -- (2.75,-3.25);
\draw[fill] (0,0)  circle [radius=0.03];
\draw[->] (-0.5,-0.5) node[left] {\footnotesize $\Sigma_{0,0}$} to [out=30,in=-120] (-0.06,-0.06);
\draw[->] (2,-1) node[right] {\footnotesize $\uC_{\us=0}$} to [out=-150,in=60]   (1.55,-1.45);
\draw[->] (1.4,-2.6) node[left] {\footnotesize$M_{\kappa}$: $\kappa$-neighbourhood }to [out=30,in=-110] (1.85,-2.15);
\end{tikzpicture}
\caption{$\kappa$-neighbourhood $M_{\kappa}$}
\label{fig 1}
\end{figure}
\end{definition}
In the following, we introduce the class of vacuum perturbed Schwarzschild metrics in the $\kappa$-neighbourhood $M_{\kappa}$. Denote the vacuum perturbed Schwarzschild metric by $g$. Let $\{s, \us\}$ be the coordinate functions on $M_{\kappa}$ inherited from the Schwarzschild spacetime. We assume that $\{s, \us\}$ remains to be a double null coordinate system for $g$, and choose $\{\theta^1, \theta^2\}$ such that $g$ takes the form
\begin{align*}
g
=
2\Omega^2  ( \ed s  \otimes \ed \us  +\ed \us \otimes \ed s ) + \slashg_{ab} ( \ed \theta^a - b^{a} \ed s ) \otimes ( \ed \theta^b - b^{b} \ed s ).
\end{align*}
Adopting definitions \ref{def 2.1} and \ref{def 2.2} to $(M_{\kappa}, g)$ and the double null coordinate system $\{s, \us, \theta^1, \theta^2\}$, we can obtain the corresponding connection coefficients and curvature components. We shall use the subscript $S$ to denote the quantities in the Schwarzschild spacetime, to distinguish with the ones in $(M_{\kappa}, g)$. The vacuum perturbed Schwarzschild metric is characterised by comparing the corresponding metric components, connection coefficients and curvature components with the Schwarzschild metric.
\begin{definition}\label{def 2.4}
Let $\epsilon$ be a positive number and $g_{\epsilon}$ be a Ricci-flat Lorentzian metric on $M_{\kappa}$ that in coordinates $\{\us,s,\theta^1,\theta^2\}$
\begin{align}
	g_{\epsilon}
	=
	2\Omega_{\epsilon}^2  ( \ed s  \otimes \ed \us  +\ed \us \otimes \ed s ) 
	+ 
	( \slashg_{\epsilon})_{ab} 
	( \ed \theta^a - b_{\epsilon}^{a} \ed s) \otimes ( \ed \theta^b - b_{\epsilon}^{b} \ed s).
\end{align}
We define the area radius function $r_{\epsilon}(\us,s)$ by
\begin{align}
4\pi r_{\epsilon}^2(\us,s) =\int_{\Sigma_{\us,s}} 1 \cdot \dvol_{\subslashg_{\epsilon}}.
\end{align}
Let $N\geq n$ be two positive integers. A $(N+2)$-th order differentiable metric $g_{\epsilon}$ is called $(n+2)$-th order $\epsilon$-close to the Schwarzschild metric $g_S$ on $M_{\kappa}$, if the following listed assumptions on the difference between the area radius, the metric components, connection coefficients, curvature components of $g_{\epsilon}$ and $g_S$ hold.
\begin{enumerate}[label=\Alph*.]

\item Assumption on the area radius:
\begin{align*}
1-\epsilon \leq \vert \frac{r_{\epsilon}}{r_S} \vert \leq 1+ \epsilon.
\end{align*}

\item Assumptions on the metric components: 

\begin{enumerate}
\item[$\Omega_{\epsilon}$:] for $k\geq 0$, $l\geq 1$, $m\geq 2$, $k +l +m \leq n+2$,
\begin{align*}
\begin{aligned}
&  \vert \log \Omega_{\epsilon} - \log \Omega_{S} \vert \leq \frac{\epsilon r_0}{r_\epsilon},
&& \vert \circnabla^k ( \log \Omega_{\epsilon} - \log \Omega_{S} ) \vert \leq \frac{\epsilon r_0}{r_{\epsilon}},
\\
&  \vert \circnabla^k \partial_s^l ( \log \Omega_{\epsilon} - \log \Omega_{S} ) \vert \leq \frac{\epsilon r_0}{r_{\epsilon}^{1+l}},
\\
& \vert \circnabla^k \partial_{\us} ( \log \Omega_{\epsilon} - \log \Omega_{S} ) \vert \leq \frac{\epsilon r_0}{r_{\epsilon}^2},
&&  \vert \circnabla^k \partial_{\us}^m ( \log \Omega_{\epsilon} - \log \Omega_{S} ) \vert \underset{m\geq 2}{\leq} \frac{\epsilon r_0}{r_{\epsilon}^3 r_0^{m-2}},
\\
& \vert \circnabla^k \partial_s^l \partial_{\us} ( \log \Omega_{\epsilon} - \log \Omega_{S} ) \vert \leq \frac{\epsilon r_0}{r_{\epsilon}^{2+l}},
&&  \vert \circnabla^k \partial_s^l \partial_{\us}^m ( \log \Omega_{\epsilon} - \log \Omega_{S} ) \vert \underset{m\geq 2}{\leq} \frac{\epsilon r_0}{r_{\epsilon}^{3+l} r_0^{m-2}}.
\end{aligned}
\end{align*}

\item[$\vec{b}_{\epsilon}$:] for $k\geq 0$, $l\geq 1$, $m\geq 1$, $k +l +m \leq n+2$,
\begin{align*}
\begin{aligned}
& \vert \vec{b}_{\epsilon} \vert_{\subcircg} \leq \frac{\epsilon r_0 \vert\us\vert}{r^3},  
&& \vert \circnabla^{k} \vec{b}_{\epsilon} \vert_{\subcircg} \leq \frac{\epsilon r_0 \vert\us\vert}{r_{\epsilon}^3}, 
\\
& \vert \circnabla^k \partial_{s}^l \vec{b}_{\epsilon} \vert_{\subcircg} \leq \frac{\epsilon r_0 \vert\us\vert}{r_{\epsilon}^{3+l}},
&& \vert \circnabla^{k} \partial_{\us}^{m} \vec{b}_{\epsilon} \vert_{\subcircg} \leq \frac{\epsilon}{r_{\epsilon}^{3}r_0^{m-2}}  ,
&& \vert \circnabla^k \partial_s^l \partial_{\us}^m \vec{b}_{\epsilon} \vert_{\subcircg} \leq \frac{\epsilon}{r_{\epsilon}^{3+l}r_0^{m-2}}.
\end{aligned}
\end{align*}

\item[$\slashg_{\epsilon}$:] for $k\geq 0$, $l\geq 1$, $m\geq 1$, $k +l +m \leq n+2$,
\begin{align*}
\begin{aligned}
&
\begin{aligned}
&  \vert \slashg_{\epsilon} - \slashg_S \vert_{\subcircg} \leq \epsilon r_{\epsilon}^2, 
&&  \vert \circnabla^k ( \slashg_{\epsilon} - \slashg_S ) \vert_{\subcircg} \leq \epsilon r_{\epsilon}^2, 
\\
&  \vert \circnabla^k \partial_s^l ( \slashg_{\epsilon} - \slashg_S ) \vert_{\subcircg} \leq \frac{\epsilon r_0 }{r_{\epsilon}^{l-1}},
&& \vert \circnabla^k \partial_{\us}^m ( \slashg_{\epsilon} - \slashg_S ) \vert_{\subcircg} \leq \frac{\epsilon r_{\epsilon}}{r_0^{m-1}},
\end{aligned}
\\
& \vert \circnabla^k \partial_s^l \partial_{\us}^m ( \slashg_{\epsilon} - \slashg_S ) \vert_{\subcircg} \leq \frac{\epsilon}{ r_0^{m-1} r_{\epsilon}^{l-1}}.
\end{aligned}
\end{align*}
\end{enumerate}

\item Assumptions on the connection coefficients : 

\begin{enumerate}
\item[$\tr \chi_{\epsilon}$:] for $k\geq 0$, $l\geq 1$, $m\geq 1$, $k+l +m \leq n+1$,
\begin{align*}
\begin{aligned}
&
\begin{aligned}
& \vert \tr\chi_{\epsilon} - \tr\chi_S \vert_{\subcircg} \leq \frac{\epsilon}{ r_{\epsilon} },
&& \vert \circnabla^k  (  \tr\chi_{\epsilon} - \tr\chi_S ) \vert_{\subcircg} \leq \frac{\epsilon }{ r_{\epsilon} } , 
\\
& \vert  \circnabla^k \partial_s^l  (  \tr\chi_{\epsilon} - \tr\chi_S ) \vert_{\subcircg} \leq \frac{\epsilon }{r_{\epsilon}^{1+l}},
&& \vert  \circnabla^k \partial_{\us}^m  (  \tr\chi_{\epsilon} - \tr\chi_S )  \vert_{\subcircg} \leq \frac{\epsilon }{ r_{\epsilon}^{2}r_0^{m-1}},
\end{aligned}
\\
& \vert \circnabla^k \partial_s^l \partial_{\us}^m   (  \tr\chi_{\epsilon} - \tr\chi_S ) \vert_{\subcircg} \leq \frac{\epsilon }{r_{\epsilon}^{2+l}r_0^{m-1}}.
\end{aligned}
\end{align*}

\item[$\hatchi_{\epsilon}$:] for $k\geq 0$, $l\geq 1$, $m\geq 1$, $k+l +m \leq n+1$,
\begin{align*}
\begin{aligned}
& \vert \hatchi_{\epsilon} \vert_{\subcircg} \leq \epsilon r_{\epsilon},
&& \vert \circnabla^k \hatchi_{\epsilon} \vert_{\subcircg} \leq \epsilon r_{\epsilon},
\\
& \vert \circnabla^k \partial_{s}^l \hatchi_{\epsilon} \vert_{\subcircg} \leq \frac{\epsilon}{r_{\epsilon}^{l-1}},
&& \vert  \circnabla^k \partial_{\us}^m \hatchi_{\epsilon}  \vert_{\subcircg} \leq \frac{\epsilon r_{\epsilon}}{r_0^{m}},
&& \vert \circnabla^k \partial_{s}^l \partial_{\us}^m  \hatchi_{\epsilon} \vert_{\subcircg} \leq \frac{\epsilon}{ r_{\epsilon}^{l-1} r_0^{m}}.
\end{aligned}
\end{align*}

\item[$\tr\uchi_{\epsilon}$:] for $k\geq 0$, $l\geq 1$, $m\geq 1$, $k+l +m \leq n+1$,
\begin{align*}
\begin{aligned}
&
\begin{aligned}
& \vert \tr\uchi_{\epsilon} - \tr\uchi_{S}\vert_{\subcircg} \leq \frac{\epsilon r_0 }{ r_{\epsilon}^2},
&& \vert \circnabla^{k}  ( \tr\uchi_{\epsilon} - \tr\uchi_{S} ) \vert_{\subcircg} \leq \frac{\epsilon r_0}{r_{\epsilon}^2},
\\
& \vert  \circnabla^k \partial_{s}^l  ( \tr\uchi_{\epsilon} - \tr\uchi_{S} ) \vert_{\subcircg} \leq \frac{\epsilon r_0}{r_{\epsilon}^{2+l}},
&& \vert  \circnabla^k \partial_{\us}^m ( \tr\uchi_{\epsilon} - \tr\uchi_{S} )  \vert_{\subcircg} \leq \frac{\epsilon}{ r_{\epsilon}^{3}r_0^{m-2}},
\end{aligned}
\\
& \vert \circnabla^k \partial_{s}^l \partial_{\us}^m  ( \tr\uchi_{\epsilon} - \tr\uchi_{S} ) \vert_{\subcircg} \leq \frac{\epsilon}{r_{\epsilon}^{3+l} r_0^{m-2}}.
\end{aligned}
\end{align*}

\item[$\hatuchi_{\epsilon}$:] for $k\geq 0$, $l\geq 1$, $m\geq 2$, $k+l +m \leq n+1$,
\begin{align*}
&
\begin{aligned}
& \vert \hatuchi_{\epsilon} \vert_{\subcircg} \leq \epsilon r_0,
&& \vert \circnabla^k \hatuchi_{\epsilon} \vert_{\subcircg} \leq \epsilon r_0,
&& \vert \circnabla^k \partial_{s}^l \hatuchi_{\epsilon} \vert_{\subcircg} \leq \frac{\epsilon r_0^{\frac{3}{2}}}{ r_{\epsilon}^{\frac{1}{2}+l}},
\end{aligned}
\\
&
\begin{aligned}
& \vert  \circnabla^k \partial_{\us}^m \hatuchi_{\epsilon}  \vert_{\subcircg} \leq \frac{\epsilon}{ r_0^{m-1}},
&& \vert \circnabla^k \partial_{s}^l \partial_{\us}  \hatuchi_{\epsilon} \vert_{\subcircg} \leq \frac{\epsilon r_0^{\frac{3}{2}}}{ r_{\epsilon}^{\frac{3}{2}+l}},
&& \vert \circnabla^k \partial_{s}^l \partial_{\us}^m  \hatuchi_{\epsilon} \vert_{\subcircg} \underset{m\geq 2}{\leq} \frac{\epsilon}{ r_0^{m-3} r_{\epsilon}^{2+l}}.
\end{aligned}
\end{align*}

\item[$\eta_{\epsilon}$:] for $k\geq 0$, $l\geq 1$, $m\geq 1$, $k+l +m \leq n+1$,
\begin{align*}
\begin{aligned}
& \vert \eta_{\epsilon} \vert_{\subcircg} \leq \frac{\epsilon r_0}{ r_{\epsilon}},
&& \vert \circnabla^{k}  \eta_{\epsilon} \vert_{\subcircg} \leq \frac{\epsilon r_0}{r_{\epsilon}},
\\
& \vert  \circnabla^k \partial_{s}^l  \eta_{\epsilon} \vert_{\subcircg} \leq \frac{\epsilon r_0}{r_{\epsilon}^{1+l}},
&& \vert  \circnabla^k \partial_{\us}^m \eta_{\epsilon}  \vert_{\subcircg} \leq \frac{\epsilon }{ r_{\epsilon} r_0^{m-1}},
&& \vert \circnabla^k \partial_{s}^l \partial_{\us}^m  \eta_{\epsilon} \vert_{\subcircg} \leq \frac{\epsilon }{r_{\epsilon}^{1+l} r_0^{m-1}}.
\end{aligned}
\end{align*}

\item[$\ueta_{\epsilon}$:] for $k\geq 0$, $l\geq 1$, $m\geq 1$, $k+l +m \leq n+1$,
\begin{align*}
\begin{aligned}
& \vert \ueta_{\epsilon} \vert_{\subcircg} \leq \frac{\epsilon r_0}{ r_{\epsilon}},
&& \vert \circnabla^{k}  \ueta_{\epsilon} \vert_{\subcircg} \leq \frac{\epsilon r_0}{r_{\epsilon}},
\\
& \vert  \circnabla^k \partial_{s}^l  \ueta_{\epsilon} \vert_{\subcircg} \leq \frac{\epsilon r_0}{r_{\epsilon}^{1+l}},
&& \vert  \circnabla^k \partial_{\us}^m \ueta_{\epsilon}  \vert_{\subcircg} \leq \frac{\epsilon}{ r_{\epsilon} r_0^{m-1}},
&& \vert \circnabla^k \partial_{s}^l \partial_{\us}^m  \ueta_{\epsilon} \vert_{\subcircg} \leq \frac{\epsilon}{r_{\epsilon}^{1+l} r_0^{m-1}}.
\end{aligned}
\end{align*}

\item[$\omega_{\epsilon}$:] for $k\geq 0$, $l\geq 1$, $m\geq 1$, $k+l +m \leq n+1$,
\begin{align*}
\begin{aligned}
&
\begin{aligned}
& \vert \omega_{\epsilon} -\omega_{S} \vert_{\subcircg} \leq \frac{\epsilon r_0 }{r_{\epsilon}^2},
&& \vert \circnabla^{k} (\omega_{\epsilon} -\omega_{S})  \vert_{\subcircg} \leq \frac{\epsilon  r_0 }{r_{\epsilon}^2},
\\
& \vert  \circnabla^k \partial_{s}^l  (\omega_{\epsilon} -\omega_{S}) \vert_{\subcircg} \leq \frac{\epsilon r_0}{r_{\epsilon}^{2+l}},
&& \vert  \circnabla^k \partial_{\us}^m (\omega_{\epsilon} -\omega_{S}) \vert_{\subcircg} \leq \frac{\epsilon r_0 }{r_{\epsilon}^{3} r_0^{m-1} },
\end{aligned}
\\
& \vert \circnabla^k \partial_{s}^l \partial_{\us}^m (\omega_{\epsilon} -\omega_{S}) \vert_{\subcircg} \leq \frac{\epsilon r_0}{r_{\epsilon}^{3+l} r_0^{m-1}}.
\end{aligned}
\end{align*}

\item[$\uomega_{\epsilon}$:] for $k\geq 0$, $l\geq 1$, $m\geq 1$, $k+l +m \leq n+1$,
\begin{align*}
\begin{aligned}
&
\begin{aligned}
& \vert \uomega_{\epsilon} -  \uomega_{S} \vert_{\subcircg} \leq \frac{\epsilon r_0}{r_{\epsilon}^2},
&& \vert \circnabla^{k} (\uomega_{\epsilon} -  \uomega_{S}) \vert_{\subcircg} \leq \frac{\epsilon r_0^{\frac{3}{2}}}{r_{\epsilon}^{\frac{5}{2}}},
\\
& \vert  \circnabla^k \partial_{s}^l (\uomega_{\epsilon} -  \uomega_{S}) \vert_{\subcircg} \leq \frac{\epsilon r_0}{r_{\epsilon}^{2+l}},
&& \vert  \circnabla^k \partial_{\us}^m (\uomega_{\epsilon} -  \uomega_{S}) \vert_{\subcircg} \leq \frac{\epsilon r_0}{r_{\epsilon}^{3} r_0^{m-1}},
\end{aligned}
\\
& \vert \circnabla^k \partial_{s}^l \partial_{\us}^m  (\uomega_{\epsilon} -  \uomega_{S}) \vert_{\subcircg} \leq \frac{\epsilon r_0}{r_{\epsilon}^{3+l} r_0^{m-1}}.
\end{aligned}
\end{align*}

\end{enumerate}

\item Assumptions on the curvature components: 

\begin{enumerate}
\item[$\alpha_{\epsilon}$:] for $k\geq 0$, $l\geq 1$, $m\geq 1$, $k+ l+ m \leq n$,
\begin{align*}
\begin{aligned}
&
\begin{aligned}
& \vert \alpha_{\epsilon} \vert_{\subcircg} \leq \frac{\epsilon r_\epsilon}{r_0},
\quad
&& \vert \circnabla^k  \alpha_{\epsilon}  \vert_{\subcircg} \leq \frac{\epsilon r_\epsilon}{r_0},
\quad
&& \vert \circnabla^k \partial_{s}^l \alpha_{\epsilon}  \vert_{\subcircg} \leq \frac{\epsilon }{r_{\epsilon}^{l-1} r_0},
\end{aligned}
\\
&
\begin{aligned}
& \vert  \circnabla^k \partial_{\us}^m  \alpha_{\epsilon}   \vert_{\subcircg} \leq \frac{\epsilon r_{\epsilon} }{r_0^{m+1}},
\quad
&& \vert \circnabla^k \partial_{s}^l \partial_{\us}^m  \alpha_{\epsilon}   \vert_{\subcircg} \leq \frac{\epsilon}{r_{\epsilon}^{l-1}r_0^{m+1}}.
\end{aligned}
\end{aligned}
\end{align*}

\item[$\beta_{\epsilon}$:] for $k\geq 0$, $l\geq 1$, $m\geq 1$, $k+ l+ m \leq n$,
\begin{align*}
\begin{aligned}
&
\begin{aligned}
& \vert \beta_{\epsilon} \vert_{\subcircg} \leq\frac{\epsilon }{r_{\epsilon}},
\quad
&& \vert \circnabla^k  \beta_{\epsilon} \vert_{\subcircg} \leq \frac{\epsilon }{r_{\epsilon}},
\quad
&& \vert \circnabla^k \partial_{s}^l  \beta_{\epsilon} \vert_{\subcircg} \leq \frac{\epsilon }{r_{\epsilon}^{1+l}},
\end{aligned}
\\
&
\begin{aligned}
& \vert  \circnabla^k \partial_{\us}^m   \beta_{\epsilon}  \vert_{\subcircg} \leq \frac{\epsilon }{r_{\epsilon} r_0^m},
\quad
&& \vert \circnabla^k \partial_{s}^l \partial_{\us}^m   \beta_{\epsilon}  \vert_{\subcircg} \leq\frac{\epsilon }{r_{\epsilon}^{1+l} r_0^m}.
\end{aligned}
\end{aligned}
\end{align*}

\item[$\rho_{\epsilon}$:] for $k\geq 0$, $l\geq 1$, $m\geq 1$, $k+ l+ m \leq n$,
\begin{align*}
\begin{aligned}
&
\begin{aligned}
& \vert \rho_{\epsilon} - \rho_{S} \vert_{\subcircg} \leq\frac{\epsilon r_0}{r_{\epsilon}^{3}},
\quad
&& \vert \circnabla^k (\rho_{\epsilon} - \rho_{S})  \vert_{\subcircg} \leq \frac{\epsilon r_0}{r_{\epsilon}^{3}},
\quad
&& \vert \circnabla^k \partial_{s}^l (\rho_{\epsilon} - \rho_{S}) \vert_{\subcircg} \leq \frac{\epsilon r_0}{r_{\epsilon}^{3+l}},
\end{aligned}
\\
&
\begin{aligned}
& \vert  \circnabla^k \partial_{\us}^m  (\rho_{\epsilon} - \rho_{S}) \vert_{\subcircg} \leq \frac{\epsilon r_0}{r_{\epsilon}^{3} r_0^m},
\quad
&& \vert \circnabla^k \partial_{s}^l \partial_{\us}^m  (\rho_{\epsilon} - \rho_{S}) \vert_{\subcircg} \leq \frac{\epsilon r_0}{r_{\epsilon}^{3+l} r_0^m}.
\end{aligned}
\end{aligned}
\end{align*}

\item[$\sigma_{\epsilon}$:] for $k\geq 0$, $l\geq 1$, $m\geq 1$, $k+ l+ m \leq n$,
\begin{align*}
\begin{aligned}
&
\begin{aligned}
& \vert \sigma_{\epsilon} \vert_{\subcircg} \leq \frac{\epsilon r_0}{r_{\epsilon}^{3}},
\quad
&& \vert \circnabla^k \sigma_{\epsilon} \vert_{\subcircg} \leq \frac{\epsilon r_0}{r_{\epsilon}^{3}},
\quad
&& \vert \circnabla^k \partial_{s}^l \sigma_{\epsilon} \vert_{\subcircg} \leq \frac{\epsilon r_0}{r_{\epsilon}^{3+l}},
\end{aligned}
\\
&
\begin{aligned}
& \vert  \circnabla^k \partial_{\us}^m  \sigma_{\epsilon}  \vert_{\subcircg} \leq \frac{\epsilon r_0}{r_{\epsilon}^{3} r_0^m},
\quad
&& \vert \circnabla^k \partial_{s}^l \partial_{\us}^m  \sigma_{\epsilon}  \vert_{\subcircg} \leq \frac{\epsilon r_0}{r_{\epsilon}^{3+l} r_0^m}.
\end{aligned}
\end{aligned}
\end{align*}

\item[$\ubeta_{\epsilon}$:] for $k\geq 0$, $l\geq 1$, $m\geq 1$, $k+ l+ m \leq n$,
\begin{align*}
\begin{aligned}
&
\begin{aligned}
& \vert \ubeta_{\epsilon} \vert_{\subcircg} \leq \frac{\epsilon r_0^{\frac{3}{2}}}{r_{\epsilon}^{\frac{5}{2}}},
\quad
&& \vert \circnabla^k  \ubeta_{\epsilon} \vert_{\subcircg} \leq \frac{\epsilon r_0^{\frac{3}{2}}}{r_{\epsilon}^{\frac{5}{2}}},
\quad
&& \vert \circnabla^k \partial_{s}^l \ubeta_{\epsilon} \vert_{\subcircg} \leq \frac{\epsilon r_0^{\frac{3}{2}}}{r_{\epsilon}^{\frac{5}{2}+l}},
\end{aligned}
\\
&
\begin{aligned}
& \vert  \circnabla^k \partial_{\us}^m  \ubeta_{\epsilon}  \vert_{\subcircg} \leq \frac{\epsilon r_0}{r_{\epsilon}^{3} r_0^{m-1}},
\quad
&& \vert \circnabla^k \partial_{s}^l \partial_{\us}^m  \ubeta_{\epsilon} \vert_{\subcircg} \leq \frac{\epsilon r_0}{r_{\epsilon}^{3+l} r_0^{m-1}}.
\end{aligned}
\end{aligned}
\end{align*}

\item[$\ualpha_{\epsilon}$:] for $k\geq 0$, $l\geq 1$, $m\geq 2$, $k+ l+ m \leq n$,
\begin{align*}
	\begin{aligned}
		&
		\begin{aligned}
			& \vert \ualpha_{\epsilon} \vert_{\subcircg} \leq \frac{\epsilon r_0^{\frac{3}{2}}}{r_{\epsilon}^{\frac{3}{2}}},
			\quad
			&& \vert \circnabla^k \ualpha_{\epsilon} \vert_{\subcircg} \leq \frac{\epsilon r_0^{\frac{3}{2}}}{r_{\epsilon}^{\frac{3}{2}}},
			\quad
			&& \vert \circnabla^k \partial_s^l \ualpha_{\epsilon}  \vert_{\subcircg} \leq \frac{\epsilon r_0^{\frac{3}{2}}}{r_{\epsilon}^{\frac{3}{2}+l}},
		\end{aligned}
	\\
		&
		\begin{aligned}
			& \vert  \circnabla^k \partial_{\us}  \ualpha_{\epsilon}  \vert_{\subcircg} \leq \frac{\epsilon r_0^{\frac{3}{2}}}{r_{\epsilon}^{\frac{5}{2}}},
			&&\vert  \circnabla^k \partial_{\us}^{m}  \ualpha_{\epsilon}  \vert_{\subcircg} \underset{m\geq 2}{\leq} \frac{\epsilon r_0}{r_{\epsilon}^{3} r_0^{m-2}},
		\\
			&  \vert \circnabla^k \partial_s^l \partial_{\us} \ualpha_{\epsilon}   \vert_{\subcircg} \leq \frac{\epsilon r_0^{\frac{3}{2}}}{r_{\epsilon}^{\frac{5}{2}+l}},
			&& \vert \circnabla^k \partial_s^l \partial_{\us}^{m} \ualpha_{\epsilon}   \vert_{\subcircg} \underset{m\geq 2}{\leq} \frac{\epsilon r_0}{r_{\epsilon}^{3+l} r_0^{m-2}}.
		\end{aligned}
	\end{aligned}
\end{align*}

\end{enumerate}

\end{enumerate}
\end{definition}

\begin{remark}
The decay assumptions in the above definition of $g_{\epsilon}$ are motived from the results in the nonlinear stability of Minkowski spacetime \cite{CK1993}. For example in the decay assumptions for the curvature components $\ualpha_{\epsilon}$ and $\ubeta_{\epsilon}$, since
\begin{align*}
	\slashg_{\epsilon} \sim r_{\epsilon}^2 \circg,
\end{align*}
we have
\begin{align*}
	| \alpha_{\epsilon} |_{\subslashg_{\epsilon}} 
	= 
	( \slashg_{\epsilon}^{ac} \slashg_{\epsilon}^{bd}  (\alpha_{\epsilon})_{ab} (\alpha_{\epsilon})_{cd})^{\frac{1}{2}} 
	\simeq
	r_{\epsilon}^{-2} ( \circg^{ac} \circg^{bd}  (\alpha_{\epsilon})_{ab} (\alpha_{\epsilon})_{cd})^{\frac{1}{2}}
	=
	r_{\epsilon}^{-2} | \alpha_{\epsilon} |_{\subcircg},
\end{align*}
and similarly for $\ualpha$. Also
\begin{align*}
	| \beta_{\epsilon} |_{\subslashg_{\epsilon}} 
	= 
	( \slashg_{\epsilon}^{ab} (\beta_{\epsilon})_{a} (\beta_{\epsilon})_{b})^{\frac{1}{2}} 
	\simeq
	r_{\epsilon}^{-1} ( \circg^{ab}  (\beta_{\epsilon})_{a} (\beta_{\epsilon})_{b})^{\frac{1}{2}}
	=
	r_{\epsilon}^{-1} | \beta_{\epsilon} |_{\subcircg},
\end{align*}
and similarly for $\ubeta$. Therefore the first assumptions for $\alpha_{\epsilon}$, $\beta_{\epsilon}$, $\ualpha_{\epsilon}$ are equivalent to
\begin{align*}
	&
	\begin{aligned}
	|\alpha_{\epsilon}|_{\subslashg_{\epsilon}} \lesssim \frac{\epsilon}{r_0} \cdot \frac{1}{r_{\epsilon}},
	\end{aligned}
\\
	&
	\begin{aligned}
	| \beta_{\epsilon} |_{\subslashg_{\epsilon}} \lesssim \epsilon \cdot \frac{1}{r_{\epsilon}^2},
	\end{aligned}
\\
	&
	\begin{aligned}
	| \ualpha_{\epsilon} |_{\subslashg_{\epsilon}} \lesssim \epsilon r_0^{\frac{3}{2}} \cdot \frac{1}{r_{\epsilon}^{\frac{7}{2}}}.
	\end{aligned}
\end{align*}
These decay assumptions are exactly the ones proved in \cite{CK1993}.
\end{remark}

The rest of the paper mainly works in the Lorentzian manifold $(M_{\kappa},g_{\epsilon})$ where the parameters $N$, $n$ may vary depending on the context. For the sake of brevity, we simply use $(M,g)$ to denote $(M_{\kappa},g_{\epsilon})$ and omit the subscript $\epsilon$ in the geometric quantities with respect to $g_{\epsilon}$. Also for the notations of norms, we shall omit the subscript $\circg$ if the norm is taken w.r.t. $\circg$, i.e. $\vert \cdot \vert$ and $\Vert \cdot \Vert$ means $\vert \cdot \vert_{\subcircg}$ and $\Vert \cdot \Vert_{\subcircg}$.

\section{Parameterisation of foliation on null hypersurface}\label{sec 3}

In this section, we review the methods to parameterise a foliation on a null hypersurface in a double null coordinate system. These methods are the same as in \cite{L2018}\cite{L2023}.

\subsection{Parameterisation of null hypersurface}\label{sec 3.1}
Let $\ucalH$ be an incoming null hypersurface in $(M,g)$. We shall present the method to parameterise $\ucalH$ same as in \cite{L2018}\cite{L2022}. Let $\{\us, s, \vartheta\}$ be a double null coordinate system of $(M,g)$ as in definition \ref{def 2.4}, here we use $\vartheta$ to denote $(\theta^1, \theta^2)$ for the sake of brevity. Suppose that $\ucalH$ is parameterised by a function $\uh$ as its graph of $\us$ over the $\{s,\vartheta\}$ domain in the double null coordinate system $\{\us, s, \vartheta\}$, i.e.
\begin{align*}
	\ucalH = \{ (\us, s, \vartheta): \us = \uh(s,\vartheta)\},
\end{align*}
then the parameterisation function $h$ satisfies the following equation
\begin{align}
	\partial_s \uh = -b^i \partial_i \uh + \Omega^2 (\slashg^{-1})^{ij} \partial_i \uh \partial_j \uh,
	\label{eqn 3.1}
\end{align}
where $b^i$, $\Omega^2$, $(\slashg^{-1})^{ij}$ take values at the point $(\us = \uh(s,\vartheta), s, \vartheta)$.

In \cite{L2022}, it was shown that when $\uh(s=0, \cdot)$ is nearly constant, then $\uh$ exists globally for $s \in (-\kappa r_0, +\infty)$. Before stating this result, we introduce the background foliation $\{\Sigma_s\}$ and the parameterisation of each $\Sigma_s$ first. Let $\Sigma_s$ be the level set of $s$ in $\ucalH$. Geometrically, $\Sigma_s=\ucalH \cap \uC_s$ the intersection of $\ucalH$ and the outgoing null hypersurface $\uC_s$. Suppose that each $\Sigma_s$ is parameterised by a pair of functions $(\ufl{s}, s)$ as its graph of $(\us,s)$ over the $\vartheta$ domain in the double null coordinate system. We have that $\ufl{s} = \uh(s,\cdot)$.
\begin{proposition}\label{prop 3.1}
Let $\ucalH$ be an incoming null hypersurface in $(M,g)$ which is parameterised by $\us = \uh(s, \vartheta)$. Let $\ufl{s=0}$ be the parameterisation function of $\Sigma_0$, the intersection of $\ucalH$ with $C_{s=0}$,
\begin{align*}
	\ufl{s=0}(\vartheta) = \uh(0,\vartheta).
\end{align*}
Suppose that $\ufl{s=0}$ satisfies that
\begin{align*}
	\Vert \slashd \ufl{s=0} \Vert^{n+1,p}
	\leq
	\udelta_o r_0,
	\quad
	\vert \overline{\ufl{s=0}}^{\subcircg} \vert
	\leq
	\udelta_m r_0,
\end{align*}
where 
$
	n
	\geq 
	n_p
	= 	
	\left\{
	\begin{aligned}
		&
		1, 
		&&
		p>2,
	\\
		&
		2, 
		&&
		2\geq p >1.
	\end{aligned}
	\right.
$
Then there exist a small positive constant $\delta$, and constants $c_o, c_{m,m}, c_{m,o}$, all depending on $n,p$, such that if $\epsilon, \udelta_o, \udelta_m$ are suitably bounded that $\epsilon, \udelta_o, \epsilon \udelta_m \leq \delta$, then $\uh$ satisfies that
\begin{align*}
	&
	\Vert \slashd \uh \Vert^{n+1,p} 
	\leq
	c_o \Vert \slashd \ufl{s=0} \Vert^{n+1,p} 
	\leq
	c_o \udelta_o r_0,
\\
	&
	\vert \overline{\uh(s,\cdot)}^{\subcircg} - \overline{\ufl{s=0}}^{\subcircg} \vert 
	\leq
	\frac{c_{m,m}}{r_0} \epsilon \vert \overline{\ufl{s=0}}^{\subcircg} \vert \cdot \Vert \slashd \ufl{s=0} \Vert^{n+1,p}  
	+ 
	\frac{c_{m,o}}{r_0} ( \Vert \slashd \ufl{s=0} \Vert^{n+1,p} )^2
\\
	&\phantom{\vert \overline{\uh(s,\cdot)}^{\subcircg} - \overline{\ufl{s=0}}^{\subcircg} \vert }
	\leq 
	( c_{m,m} \epsilon \udelta_m \udelta_o + c_{m,o} \udelta_o^2 ) r_0,
\end{align*}
for all $s\in (-\kappa r_0, +\infty)$.
\end{proposition}
\begin{remark}\label{rem 3.2}
As shown in \cite{L2022}, proposition \ref{prop 3.1} requires the $L^{\infty}$ bounds of $b, \Omega, \slashg$ up to $(n+2)$-th order tangential derivatives on each $C_s$, which are exactly the top order derivatives with $L^{\infty}$ bounds in definition \ref{def 2.4}.
\end{remark}

\subsection{Parameterisation of foliation}\label{sec 3.2}
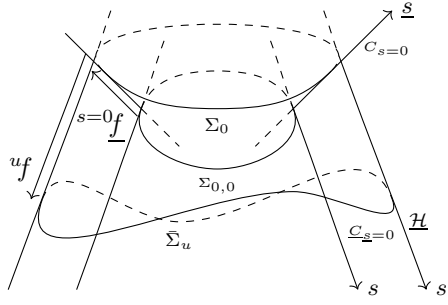
\begin{figure}[H]
\begin{center}
\begin{tikzpicture}
\draw[dashed] (-1,0)
to [out=70, in=180] (0,0.5)
to [out=0,in=110] (1,0);
\draw (1,0)
to [out=-70,in=0] (0,-0.8) node[below]{\tiny $\Sigma_{0,0}$}
to [out=180,in=-110] (-1,0); 
\draw[dashed] (-1,0) to [out=70,in=-110] (-0.7,0.9);
\draw (-1,0) to [out=-110,in=70] (-1.85,-2.4);
\draw[dashed] (1,0) to [out=110,in=-70] (0.7,0.9);
\draw[->] (1,0) to [out=-70,in=110] (1.85,-2.4) node[right] {\small $s$}; 
\node[above right] at (1.6,-1.9) {\tiny $\uC_{\us=0}$};
\draw[dashed] (-1,0) to [out=-45,in=135] (-0.5,-0.5);
\draw (-1,0) to [out=135, in= -45] (-2,1);
\draw[dashed] (1,0) to [out=-135,in=45] (0.5,-0.5);
\draw[->] (1,0) to [out=45, in= -135] (1.8,0.8) node[right] {\tiny $C_{s=0}$}to [out=45,in=-135] (2.3,1.3) node[right] {\small $\us$}; 
\draw[dashed] (-1.5,0.5) to [out=135,in=45] (1.5,0.5);
\draw (1.5,0.5) to [out=-135,in=0] (0,0) node[below]{\scriptsize $\Sigma_0$} to [out=180,in=-45] (-1.5,0.5); 
\draw[dashed] (-1.4,1.3) to [out=-110,in=70] (-1.6,0.7);
\draw (-1.6,0.7) to [out=-110,in=70] (-2.75,-2.4);
\draw[dashed] (1.4,1.3) to [out=-70,in=110] (1.6,0.7);
\draw[->] (1.6,0.7) to [out=-70,in=110] (2.75,-2.4) node[right]{\small $s$}; 
\node[right] at (2.4,-1.5) {\scriptsize $\ucalH$}; 
\draw[dashed] (-2.3,-2.2+1) to [out=70,in=180] (-0.5,-2.5+1) node[below] {\scriptsize $\bSigma_u$} to [out=0,in=110] (2.3,-2.2+1);
\draw (2.3,-2.2+1) to [out=-70,in=0] (1,-2+0.9) to [out=180,in=-110] (-2.3,-2.2+1); 
\draw[->] (-1.73,0.73) to [out=-110,in=70] (-2.45,-1.2);
\node[left] at (-2.3,-0.8) {$\fl{u}$};
\draw[->] (-1.05,-0.1) to [out=135,in=-45] (-1.65,0.5);
\node[below] at (-1.55,0.1) {$\ufl{s=0}$};
\end{tikzpicture}
\end{center}
\caption{The second kind of parameterisation of $\Sigma$.}
\label{fig 2}
\end{figure}

Let $\{\bSigma_u\}$ be a foliation on the incoming null hypersurface $\ucalH$. We parameterise $\{\bSigma_u\}$ in the coordinate system $\{s, \vartheta\}$ in $\ucalH$. Suppose that $\bSigma_u$ is parameterised by a function $\fl{u}$ as its graph of $s$ over the $\vartheta$ domain in $\{s, \vartheta\}$. Then combining with the parameterisation function $\ufl{s=0}$ of the surface $\Sigma_{s=0} = \ucalH \cap C_{s=0}$ which determines $\ucalH$, the pair of functions $(\ufl{s=0}, \fl{u})$ determines the location of the leaf $\bSigma_u$. Thus we use $(\ufl{s=0}, \fl{u})$ as the parameterisation of $\bSigma_u$. See figure \ref{fig 2}.

Another method to parameterise $\bSigma_u$ is via the double null coordinate system $\{\us, s, \vartheta\}$ directly. We assume that $\bSigma_u$ is parameterised by a pair of functions $(\ufl{u}, \fl{u})$ as its graph of $(\us, s)$ over the $\vartheta$ domain in the double null coordinate system $\{\us, s, \vartheta\}$, i.e. $\bSigma_u = \{ (\us, s, \vartheta): \us = \ufl{u}(\vartheta), s = \fl{u}(\vartheta) \}$.

We summarise that there are two methods to parameterise a foliation $\{\bSigma_u\}$ in a null hypersurface $\ucalH$: one is the family of pairs of functions $\{(\ufl{s=0}, \fl{u})\}_u$ and the other is $\{ (\ufl{u}, \fl{u}) \}_u$. 
\begin{enumerate}[label=\roman*.]
\item
The first kind parameterisation $\{(\ufl{s=0}, \fl{u})\}_u$: the parameterisation function $\ufl{s=0}$ of $\Sigma_{s=0} = \ucalH\cap C_{s=0}$ which determines $\ucalH$, and the parameterisation function $\fl{u}$ of $\bSigma_u$ in the $\{s, \vartheta\}$ coordinate system in $\ucalH$.
\item
The second kind parameterisation $\{ (\ufl{u}, \fl{u}) \}_u$: the parameterisation of $\bSigma_u$ in the double null coordinate system $\{\us, s, \vartheta\}$.
\end{enumerate}

\subsection{Transformation between two kinds of parameterisations}\label{sec 3.3}
In \cite{L2018}\cite{L2020}, we study the transformation between the above two kinds of parameterisations. We state the result in the setting of this paper in the following.
\begin{proposition}\label{prop 3.3}
Let $\Sigma$ be a spacelike surface embedded in $(M,g)$ with the parameterisations $(\ufl{s=0}, f)$ and $(\uf,f)$. Suppose that the parameterisation functions $\ufl{s=0}$, $f$ satisfy that
\begin{align*}
	\Vert \slashd \ufl{s=0} \Vert^{n+1,p} \leq \udelta_o r_0,
	\quad
	\vert \overline{\ufl{s=0}}^{\circg} \vert \leq \udelta_m r_0,
	\quad
	\Vert \slashd f \Vert^{n+1,p} \leq \delta_o (r_0+\os),
	\quad
	\overline{f}^{\circg} = \os > -\frac{r_0}{2},
\end{align*}
where $n\geq 1+n_p$, $p >1$. Then there exist a small positive constant $\delta$, and constants $c_o, c_{m,m}, c_{m,o}$ all depending on $n,p$, such that if $\epsilon, \udelta_o, \udelta_m, \delta_o$ suitably bounded that $\epsilon, \udelta_o, \epsilon \udelta_m, \delta_o \leq \delta$, then the parameterisation function $\uf$ satisfies that
\begin{align}
	\begin{aligned}
		&
		\Vert \slashd \uf \Vert^{n,p} 
		\leq c_o \udelta_o r_0,
	\\
		&
		\vert \overline{\uf}^{\circg} - \overline{\ufl{s=0}}^{\circg} \vert
		\leq 
		( c_{m,m} \epsilon \udelta_m \udelta_o + c_{m,o} \udelta_o^2 ) r_0.
	\end{aligned}
\label{eqn 3.2}
\end{align}
\end{proposition}
The proof of the above proposition follows the same route as in the proof of similar propositions in \cite{L2018}\cite{L2020}. We sketch its proof in the following and complete it in appendix \ref{appen prop 3.3}.
\begin{proof}[\it{Proof sketch}]
Let $\ucalH$ be the incoming null hypersurface where $\Sigma$ is embedded. Let $\Sigma_{s=0}= \ucalH \cap C_{s=0}$. 

Firstly we construct a family of surfaces $\{S_t \}_{t\in[0,1]}$ in $\ucalH$ deforming from $\Sigma_{s=0}$ to $\Sigma$ as follows: let $S_t$ be the surface with the parameterisation $(\ufl{s=0}, \flt)$ where $\flt = tf$, then $S_{t=0} = \Sigma_{s=0}$ and $S_{t=1} = \Sigma$.

Secondly, suppose that $S_t$ has the parameterisation $(\uflt, \flt)$ by the other method. We derived the following equation for $\uflt$ in \cite{L2018}\cite{L2020}:
\begin{align}
	\begin{aligned}
		&
		\partial_t \uflt 
		= 
		F( f,\ t b^i f_i,\ t e^i f_i,\ t \ue^i f_i,\ \uvarepsilon,\  b^i (\uflt)_i,\  e^i (\uflt)_i,\  \ue^i (\uflt)_i  ),
	\\
		&
		F= f \cdot [ 1-tb^i f_i - t\ue^i f_i - t e^i f_i \cdot \uvarepsilon ]^{-1} \cdot [ \uvarepsilon - b^i (\uflt)_i - \ue^i (\uflt)_i - e^i (\uflt)_i \cdot \uvarepsilon ],
	\end{aligned}
	\label{eqn 3.3}
\end{align}
where
\begin{align*}
	\begin{aligned}
		&
		\underline{\varepsilon} = \frac{ -|\underline{e}|^2}{(2\Omega^2 + e\cdot \underline{e}) + \sqrt{(2\Omega^2 + e\cdot \underline{e})^2 -|e|^2 |\underline{e}|^2}},
	\\
		&
		|e|^2 = \slashg_{ij}e^ie^j,
		\quad
		|\underline{e}|^2 = \slashg_{ij} \underline{e}^i \underline{e}^j, 
		\quad
		e\cdot \underline{e} =\slashg_{ij} e^i \underline{e}^j,
	\\
		&
		e^k =-2\Omega^2 (\flt)_i (B^{-1})_j^i (\slashg^{-1})^{jk}, 
		\quad
		\underline{e}^k = -2\Omega^2 (\uflt)_i (B^{-1})_j^i (\slashg^{-1})^{jk},
	\\
		&
		B_i^j= \delta_i^j - \flt_i b^j.
	\end{aligned}
\end{align*}
The initial condition of equation \eqref{eqn 3.3} is $\ufl{t=0} = \ufl{s=0}$.

Thirdly, to estimate $\uflt$ from equation \eqref{eqn 3.3}, we differentiate it to obtain the following equation for $\ddcircDelta \uflt$, where $\ddcircDelta$ is the Laplacian operator of the metric $\circg$ on $S_t$:
\begin{align}
	\label{eqn 3.4}
	\partial_t ( \ddcircDelta \uflt )
	=
	\Xlt^i \partial_i ( \ddcircDelta \uflt ) + \relt,
\end{align}
where
\begin{align*}
	\Xlt^i
	=
	&
	\partial_{t \ue^i f_i}F \cdot [-2\Omega^2 (B^{-1})_j^i ( \slashg^{-1} )^{jk} \cdot t f_k  ]
	+
	\partial_{b^i  (\uflt)_i} F \cdot b^i
\\
	&
	+
	\partial_{\ue^i  (\uflt)_i} F \cdot [ \ue^i -2\Omega^2 (B^{-1})_j^i ( \slashg^{-1} )^{jk} (\uflt)_k ]
\\
	&
	+
	\partial_{\uvarepsilon} F \cdot 
	\partial_{\vert \ue \vert^2} \uvarepsilon \cdot [ 8\Omega^4 ( B^{-1} )_j^i  ( \slashg^{-1} )^{jk} (B^{-1})_k^l (\uflt)_l ]
\\
	&
	+
	\partial_{\uvarepsilon} F \cdot \partial_{e\cdot \ue} \uvarepsilon \cdot [  4\Omega^4 (B^{-1})_k^l ( \slashg^{-1} )^{jk} ( B^{-1} )_j^i \cdot t f_l   ],
\end{align*}
and
\begin{align*}
	\relt
	&=
	\partial_{f} F \cdot \ddcircDelta f 
	+
	\partial_{t f_i \ue^i} F \cdot \{ \ddcircDelta ( t f_i \ue^i ) - [-2\Omega^2 (B^{-1})_j^i ( \slashg^{-1} )^{jk} \cdot t f_k ] (\ddcircDelta \uflt )_i  \}
\\
	&\phantom{=}
	+
	\partial_{t e^i f_i} F \cdot \ddcircDelta ( te^i f_i )
	+
	\partial_{t b^i f_i} F \cdot \ddcircDelta( t b^i f_i )
\\
	&\phantom{=}
	+ 
	\partial_{\uvarepsilon} F \cdot \partial_{\vert \ue \vert^2} \uvarepsilon \cdot \{ \ddcircDelta \vert \ue \vert^2 - 8\Omega^4\ (\uflt)_l (B^{-1})_k^l ( \slashg^{-1} )^{jk} ( B^{-1} )_j^i ( \ddcircDelta \uflt )_i \}
\\
	&\phantom{=}
	+
	\partial_{\uvarepsilon} F \cdot \partial_{e \cdot \ue} \uvarepsilon \cdot \{ \ddcircDelta ( e \cdot \ue ) -4\Omega^4\ t f_l (B^{-1})_k^l ( \slashg^{-1} )^{jk} ( B^{-1} )_j^i ( \ddcircDelta \uflt )_i  \}
\\
	&\phantom{=}
	+
	\partial_{e^i (\uflt)_i} F \cdot \{ \ddcircDelta ( e^i  (\uflt)_i ) - e^i  ( \ddcircDelta\uflt )_i \}
	+
	\partial_{b^i (\uflt)_i} F \cdot \{ \ddcircDelta ( b^i  (\uflt)_i ) - b^i  ( \ddcircDelta\uflt )_i \}
\\
	&\phantom{=}
	+
	\partial_{\ue^i (\uflt)_i }F \cdot \{ \ddcircDelta( \ue^i  (\uflt)_i ) - [ \ue^i -2\Omega^2 (B^{-1})_j^i ( \slashg^{-1} )^{jk} (\uflt)_k ] ( \ddcircDelta \uflt )_i \}
\\
	&\phantom{=}
	+
	\partial^2_{\pmb{a} \pmb{b}} F \cdot \ddcircnabla_i \pmb{a} \ddcircnabla^i \pmb{b},
\\
	\pmb{a}, \pmb{b}:
	&\phantom{=}
	f,\ t b^i f_i,\ t e^i f_i,\ t \ue^i f_i,\ \uvarepsilon,\  b^i  (\uflt)_i,\  e^i  (\uflt)_i,\  \ue^i  (\uflt)_i.
\end{align*}
For the heuristic understanding of equation \eqref{eqn 3.4}, we refer to \cite{L2020}.

Fourthly, we introduce the following notations of bounds similarly as in \cite{L2020} that
\begin{align}
	\ud_o = c_o \udelta_o, 
	\quad 
	\ud_m = [ 1+ c_{m,m} \epsilon \udelta_o ] \udelta_m + c_{m,o} \udelta_o^2,
	\label{eqn 3.5}
\end{align}
where the constants $c_o$, $c_{m,m}$, $c_{m,o}$ are determined later in the proof.

After the above preparations, we use the bootstrap argument to prove the proposition. Introduce the following bootstrap assumption.
\begin{assum}
Estimates \eqref{eqn 3.2} hold for $\uflt$ in the closed interval $[0,t_a]\subset [0,1]$.
\end{assum}
We show that estimates \eqref{eqn 3.2} can be strengthened to strict inequalities in the interval $[0,t_a]$, thus they can be extended to a slightly larger interval. Then the bootstrap argument implies that the assumption holds on the maximal interval $[0,1]$.

By the continuity argument, as long as $c_o>1$, $c_{m,m},c_{m,o}>0$, there exists a small neighbourhood $[0,t_a]$ of $0$ such that the assumption is true.

We shall use the assumption to estimate $F$, $\Xlt$, $\relt$ on the interval $[0,t_a]$. Assume that $\delta$, $c_o$, $c_{m,m}$, $c_{m,o}$ satisfy the following primary bounds
\begin{align*}
	\delta \leq \frac{1}{2},
	\quad
	(c_o+c_{m,m} +c_{m,o})\delta \leq 1.
\end{align*}
Then we have the following estimates
\begin{align}
	\begin{aligned}
		&
		\vert F \vert
		\leq
		c(n,p) \frac{(|\os|/r_0+ \delta_o)r_0^2}{(r_0 + t \os)^2} \Big[ \ud_o +  \frac{\epsilon \ud_m r_0}{r_0+ t \os} \Big] \ud_o r_0,
	\\
		&
		\Vert \Xlt \Vert^{n,p}
		\leq
		c(n,p) \frac{(|\os|/r_0+ \delta_o)r_0^2}{(r_0 + t \os)^2} \Big[ \ud_o +  \frac{\epsilon \ud_m r_0}{r_0+ t \os} \Big]
		\leq
		c(n,p) \frac{(|\os|/r_0+ \delta_o)r_0^2}{(r_0 + t \os)^2}  ,
	\\
		&
		\Vert \relt \Vert^{n-1,p}
		\leq
		c(n,p) \frac{(|\os|/r_0+ \delta_o)r_0^2}{(r_0 + t \os)^2} \Big[ \ud_o +  \frac{\epsilon \ud_m r_0}{r_0+ t \os} \Big] \ud_o r_0.
	\end{aligned}
	\label{eqn 3.6}
\end{align}
The proof of above estimates leaves to appendix \ref{appen prop 3.3}. By Sobolev inequality, we have that
\begin{align*}
	\vert \ddcircdiv \Xlt \vert 
	\leq
	\Vert \Xlt \Vert^{n,p}
	\leq
	c(n,p) \frac{(|\os|/r_0+ \delta_o)r_0^2}{(r_0 + t \os)^2},
	\quad
	\int_0^t \vert \ddcircdiv \Xlt \vert  \ed t
	\leq
	c(n,p).
\end{align*}
Then by equations \eqref{eqn 3.2}\eqref{eqn 3.3} and Gronwall's inequality, we obtain that
\begin{align*}
	&
	\begin{aligned}
		\vert \overline{\uflt}^{\circg} - \overline{\ufl{s=0}}^{\circg} \vert
		\leq
		\int_0^t \vert \overline{F}^{\circg} \vert \ed t
		&
		\leq
		c(n,p) \frac{t(|\os|/r_0+ \delta_o)r_0}{r_0 + t \os} ( \ud_o +  \epsilon \ud_m ) \ud_o r_0
		\leq
		c(n,p) ( \ud_o +  \epsilon \ud_m ) \ud_o r_0
	\\
		&
		\leq
		c(n,p) c_o ( c_o + c_{m,o} \epsilon \udelta_o) \udelta_o^2 r_0
		+
		c(n,p) c_o (1 + c_{m,m} \epsilon \udelta_o) \epsilon \udelta_m \udelta_o r_0,
	\end{aligned}
\\
	&
	\begin{aligned}
		\Vert \slashd \uflt \Vert^{n,p}
		&
		\leq
		c(n,p) \Vert \ddcircDelta \uflt \Vert^{n-1,p}
		\leq
		c(n,p)\{
		\Vert \ddcircDelta \ufl{s=0} \Vert^{n-1,p}
		+
		\int_0^t \Vert \relt \Vert^{n-1,p} \ed t \}
	\\
		&
		\leq
		c(n,p) \udelta_o r_0
		+
		c(n,p) ( \ud_o +  \epsilon \ud_m ) \ud_o r_0
	\\
		&
		\leq
		c(n,p) 
		\{ 1 + c_o ( c_o + c_{m,o} \epsilon \udelta_o)\udelta_o
		+
		c_o (1 + c_{m,m} \epsilon \udelta_o) \epsilon \udelta_m  \} \udelta_o r_0.
	\end{aligned}
\end{align*}
We choose suitable constants $\delta$, $c_o$, $c_{m,m}$, $c_{m,o}$ to close the bootstrap argument. Considering the following system of inequalities
\begin{align*}
	c_o >1,
	\quad
	\delta \leq \frac{1}{2},
	\quad
	(c_o+c_{m,m} +c_{m,o})\delta \leq 1,
\\
	c(n,p) c_o ( c_o + c_{m,o} \delta^2 ) 
	<
	c_{m,o},
\\
	c(n,p) c_o (1 + c_{m,m} \delta^2)
	<
	c_{m,m},
\\
	c(n,p) 
	[ 1 + c_o ( c_o + c_{m,o} \delta^2 ) \delta
	+
	c_o (1 + c_{m,m} \delta^2) \delta ]
	<
	c_o,
\end{align*}
the solution of the above system is nonempty, for example let $c_o = \max\{ 2, 2c(n,p) \}$, $c_{m,o} = 2 c(n,p) c_o^2$, $c_{m,m} = 2 c(n,p) c_o$, and $\delta$ sufficiently small, then the inequalities are satisfied.
\end{proof}
\begin{remark}\label{rem 3.4}
Proposition \ref{prop 3.3} requires the $L^{\infty}$ bounds of the metric components $\vec{b}, \Omega, \slashg$ up to $(n+1)$-th order derivatives given in definition \ref{def 2.4}, see appendix \ref{appen prop 3.3}.
\end{remark}

\section{Formulae of geometric quantity associated with foliation}\label{sec 4}

In this section, we present the formulae of geometric quantities associated with a foliation. The method to derive these formulae follows the same method employed in \cite{L2018}\cite{L2020}\cite{L2023}.

\subsection{Two-step procedure to derive formulae}
Let $\{ \bSigma_u \}$ be a foliation of an incoming null hypersurface $\ucalH$. Construct the conjugate null frame $\{ \buL^u, \bL'^u \}$ associated with $\{ \bSigma_u\}$, where
\begin{align*}
\buL^u u =1,
\quad
g(\buL^u, \bL'^u) =2.
\end{align*}
We follow the following two-step procedure to derive the formulae of geometric quantities relative to $\{ \buL^u, \bL'^u\}$ associated with $\{ \bSigma_u \}$.
\begin{step i}\label{step i}
Construct the foliation $\{ \Sigma_s \}$ of $\ucalH$ where $\Sigma_s= C_s \cap \ucalH$, i.e. $s$ is the restriction of the double null coordinate function $s$ to $\ucalH$. Then we derive the formulae of geometric quantities associated with this special foliation $\{\Sigma_s\}$. 
\end{step i}
\begin{step ii}\label{step ii}
Then we view $\{ \Sigma_s\}$ as the background foliation and $\{s,\vartheta\}$ as the background coordinate system of $\ucalH$. Suppose that $\{ \bSigma_u \}$ is parameterised by the family of functions $\{ \fl{u} \}$ as the graphs $s=\fl{u}(\vartheta)$ in $\{s,\vartheta\}$ coordinate system. Then we derive the formulae of geometric quantities relative to $\{\buL^u, \bL'^u\}$ associated with $\{\bSigma_u \}$, in terms of the parametrisation functions $\fl{u}$ and the geometric quantities associated with $\{\Sigma_s\}$ obtained in \hyperref[step i]{Step i}.
\end{step ii}

\subsection{\hyperref[step i]{Step i}: formulae for special foliation by coordinate $s$}\label{sec 4.2}

$\{s, \vartheta\}$ is the coordinate system of $\ucalH$ and $\{\Sigma_s\}$ is the coordinate foliation of $\ucalH$. Use $\cdot$ on the top to indicate the corresponding notation being associated with $\{\Sigma_s\}$ or $\{s, \vartheta\}$ on $\ucalH$. We list the formulae of the geometric quantities associated with $\{\Sigma_s\}$ in the following. Suppose that $\ucalH$ is parameterised by a function $\uh$ as described in section \ref{sec 3.1}.

\begin{enumerate}[label=\alph*.,leftmargin=15pt]
\item
Let $\dpartial_s$, $\dpartial_i$ be the coordinate vector of $\{s, \vartheta\}$ on $\ucalH$, then
\begin{align*}
\dpartial_s = \partial_s + \partial_s \uh\cdot \partial_{\us},
\quad
\dpartial_i = \partial_i + \partial_i \uh\cdot \partial_{\us}.
\end{align*}

\item 
Introduce the conjugate null frame $\{ \duL, \dL'\}$ associated with $\{\Sigma_s\}$,
\begin{align*}
\left\{
\begin{aligned}
&
\dL'=L',
\\
&
\duL=\uL+\uvarepsilon L + \uvarepsilon^i \partial_i,
\end{aligned}
\right.
\end{align*}
where
\begin{align*}
\uvarepsilon=-\Omega^2 (\slashg^{-1} )^{ij} \uh_i \uh_j = -\Omega^2 \vert \dslashd \uh \vert_{\subslashg}^2,
\quad
\uvarepsilon^i = -2\Omega^2 ( \slashg^{-1} )^{ij} \uh_j.
\end{align*}

\item
The shifting vector $\db$ between $\dpartial_s$ and $\duL$ is given by
\begin{align*}
\duL = \dpartial_s + \db^i \dpartial_i
\quad
\Rightarrow
\quad
\db^i 
= 
b^i - 2 \Omega^2 ( \slashg^{-1} )^{ij} \uh_j.
\end{align*}

Let $\dslashg$ be the intrinsic metric on $\Sigma_s$ and $\dslashepsilon$ be the volume form of $\dslashg$, then 
\begin{align*}
\dslashg_{ij} = \slashg_{ij},
\quad
( \dslashg^{-1} )^{ij} = ( \slashg^{-1} )^{ij},
\quad
\dslashepsilon_{ij} = \slashepsilon_{ij}.
\end{align*}
The degenerated metric on $\ucalH$ in $\{s, \vartheta\}$ coordinate system is
\begin{align*}
g\vert_{\ucalH} =  \dslashg_{ij} ( \ed \theta^i - \db^i \ed s  ) \otimes ( \ed \theta^j - \db^j \ed s ).
\end{align*}

\item
The connection coefficients on $\ucalH$ relative to $\{\Sigma_s\}$ and $\{\duL, \dL' \}$ are given by the following formulae:
\begin{align}
\begin{aligned}
&
\begin{aligned}
\dchi'_{ij} 
=
\chi'_{ij},  
\quad 
\dtr \dchi' 
= 
\tr \chi',
\end{aligned}
\\
&
\begin{aligned}
\duchi_{ij}
&=
\uchi_{ij} 
-\Omega^2 \vert \dslashd \uh \vert^2_{\subslashg} \chi_{ij}  
-2\Omega^2 \slashnabla_{ij}^2\uh -4\Omega^2 \sym \{\ueta \otimes \dslashd \uh \}_{ij} 
- 4\omega\Omega^2 (\dslashd \uh \otimes \dslashd \uh )_{ij} 
\\
&\phantom{=}
+ 4\Omega^2 \sym \{\dslashd\uh \otimes ( \chi \cdot \dslashd \uh ) \}_{ij},
\\
\dtr \duchi 
&=
( \dslashg^{-1} )^{ij} \duchi_{ij}
=  
\tr \uchi 
- 2\Omega^2 \slashDelta \uh 
- \Omega^2 \vert \dslashd \uh \vert_{\subslashg}^2 \tr \chi 
-4\Omega^2 \ueta \cdot \dslashd \uh 
- 4 \Omega^2 \omega \vert \dslashd \uh \vert^2_{\subslashg} 
\\
&\phantom{= ( \dslashg^{-1} )^{ij} \duchi_{ij}=}
+ 4 \Omega^2 \chi ( \slashnabla \uh, \slashnabla \uh),
\\
\deta_i 
&=
\eta_i 
+ ( \chi \cdot \dslashd \uh )_i,
\\
\duomega 
&= 
\uomega 
- 2\Omega^2 \eta \cdot \dslashd \uh 
- \Omega^2 \chi ( \slashnabla \uh, \slashnabla \uh).
\end{aligned}
\end{aligned}
\label{eqn 4.1}
\end{align}
In the above formulae, we use $\cdot$ to denote the inner product w.r.t. $\slashg$, $\dtr$ to denote the trace w.r.t. $\dslashg=\slashg$. $\dslashd$ is differential operator on $\Sigma_s$. $\slashnabla$ in $\slashnabla \uh, \slashnabla^2_{ij} \uh$ is the pull back of the covariant derivative $\slashnabla$ of $(\Sigma_{\us,s}, \slashg)$ to $\Sigma_s$.\footnote{Here we abuse the notation $\slashnabla$ to denote both the covariant derivative of $(\Sigma_{\us,s}, \slashg)$ and its pull back to $\Sigma_s$. Which meaning $\slashnabla$ represents in a concrete formula depends on where the differentiated function, vector field or tensor field is defined. For example, if a vector field $V$ is defined on $\Sigma_s$, then $\slashnabla$ in $\slashnabla V$ is the pull back of the covariant derivative of $(\Sigma_{\us,s}, \slashg)$ on $\Sigma_s$. If $\slashnabla$ can be interpreted in both way in a formula, we will state the precise meaning of $\slashnabla$ in that formula to avoid ambiguity.  \label{footnote 2}}
$\slashDelta$ in $\slashDelta \uh$ is the operator $ ( \slashg^{-1}   )^{ij} \slashnabla^2_{ij}$.

The precise meaning of the pull back $\slashnabla$ on $\Sigma_s$ is as follows: let $\slashGamma_{ij}^k$ be the Christoffel symbol of the covariant connection $\slashnabla$ of $(\Sigma_{\us,s}, \slashg)$, then
\begin{align*}
&
\text{$\phi$: a function on $\Sigma_s$}
&&
\slashnabla_i \phi = \dpartial_i \phi, \quad \slashnabla^i \phi= ( \slashg^{-1} )^{ij} \dpartial_j \phi,
\\
&
\text{$V$: a vector field on $\Sigma_s$}
&&
\slashnabla_i V^k =  \dpartial_i V^k + \slashGamma_{ij}^k V^j,
\\
&
\text{$T$: a tensor field on $\Sigma_s$}
&&
\slashnabla_{i} T_{i_1 \cdots i_k}^{j_1 \cdots j_l}
=
\dpartial_i T_{i_1 \cdots i_k}^{j_1 \cdots j_l} 
-  \slashGamma_{i i_m}^{r}  T_{i_1 \cdots \underset{\hat{i_m}}{r}\cdots i_k}^{j_1 \cdots j_l}
+ \slashGamma_{i s}^{j_n} T_{i_1 \cdots i_k}^{j_1 \cdots \overset{\hat{j_n}}{s}\cdots  j_l}.
\end{align*}

Note the Christoffel symbol $\slashGamma$ has the following decomposition by the Christoffel symbol $\circGamma$ of the round sphere $(\Sigma_{\us, s} , \circg)$,
\begin{align*}
\slashGamma_{ij}^k
=
\circGamma_{ij}^k
+
\circtriangle_{ij}^k,
\quad
\circtriangle_{ij}^k
=
\frac{1}{2} ( \slashg^{-1} )^{kl} ( \circnabla_i \slashg_{jl} + \circnabla_j \slashg_{il} -\circnabla_l \slashg_{ij} ).
\end{align*}

\item
The curvature components on $\ucalH$ relative to $\{\Sigma_s\}$ and $\{ \duL, \dL' \}$ are given by the following formulae:
\begin{align} 
\begin{aligned}
\dualpha_{ij} 
&=
\ualpha_{ij} 
+ ( \slashepsilon_{ik} \slashepsilon_{lj} + \slashepsilon_{jk} \slashepsilon_{li}) \uvarepsilon^l \ubeta^k
-2 \Omega^2( \uh_j \ubeta_i +  \uh_i \ubeta_j)
+4\Omega^4\uh_i \uh_j  \rho 
+2\Omega^2 \uvarepsilon \slashg_{ij}  \rho
\\
&\phantom{=}
+\uvarepsilon^k \uvarepsilon^l \slashepsilon_{ki} \slashepsilon_{jl} \rho 
- \Omega^2\uvarepsilon^k( \uh_i \slashg_{jk} + \uh_j \slashg_{ik})  \rho 
+\Omega^2 \uvarepsilon^k ( \uh_i \slashepsilon_{jk} + \uh_j \slashepsilon_{ik}) \sigma 
\\
&\phantom{=}
+ 2\Omega^2 \uvarepsilon^l( \uh_i \slashepsilon_{jl} + \uh_j \slashepsilon_{il}) \sigma
- 4 \Omega^4 \uh_i \uh_j\uvarepsilon^l \beta_l 
+2\Omega^4\uvarepsilon ( \uh_i \beta_j + \uh_j \beta_i) 
\\
&\phantom{=}
+\Omega^2 \uvarepsilon \uvarepsilon^l( \slashepsilon_{ik} \slashepsilon_{lj} + \slashepsilon_{jk} \slashepsilon_{li})  \beta^k 
+\Omega^2 \uvarepsilon^k \uvarepsilon^l( \uh_i \slashepsilon_{jl} + \uh_j \slashepsilon_{il}) \slashepsilon_{kh} \beta^h 
+ \Omega^4 \uvarepsilon^2 \alpha_{ij}
\\
&\phantom{=}
+ \Omega^4 \uvarepsilon^k \uvarepsilon^l \uh_i \uh_j \alpha_{kl} 
- \Omega^4 \uvarepsilon \uvarepsilon^l ( \uh_j \alpha_{il} + \uh_i \alpha_{jl} ),
\\
\dubeta_i
&=
\ubeta_i 
-3 \Omega^2 \uh_i \rho  -\frac{3}{2}  \uvarepsilon^l \sigma 
- \Omega^2 \uvarepsilon \beta_i 
+ 2 \Omega^2 \uh_i \uvarepsilon^l \beta_l 
+\frac{1}{2}  \uvarepsilon^k \uvarepsilon^l \slashepsilon_{ki} \slashepsilon_{lh} \beta^h 
\\
&\phantom{=}
- \frac{1}{2} \Omega^2 \uh_i \uvarepsilon^k \uvarepsilon^l \alpha_{kl}
+ \frac{1}{2} \Omega^2 \uvarepsilon \uvarepsilon^l \alpha_{il},
\\
\dsigma \cdot \slashepsilon_{ij} 
&=
\sigma \cdot \slashepsilon_{ij} 
+ \Omega^2(\uh_j  \beta_i -\uh_i \beta_j ) 
+ \frac{1}{2}  \uvarepsilon^k \slashepsilon_{kl} \beta^l \slashepsilon_{ij}
+\frac{1}{2} \Omega^2  \uvarepsilon^k ( \uh_i \alpha_{jk} - \uh_j \alpha_{ik}) ,
\\
\drho 
&=
\rho 
- \uvarepsilon^l \beta_l 
+ \frac{1}{4} \uvarepsilon^k \uvarepsilon^l \alpha_{kl} ,
\\ 
\dbeta_i 
&=
\beta_i - \frac{1}{2} \uvarepsilon^j \alpha_{ij},
\\
\dalpha_{ij} 
&= 
\alpha_{ij}.
\end{aligned}
\label{eqn 4.2}
\end{align}

\item Let $\dslashnabla$ be the covariant derivative of $(\Sigma_s, \dslashg)$ and $\dslashGamma_{ij}^k$ be the Christoffel symbol of $\dslashnabla$. Introduce the notation $\triangle_{ij}^k$ as
\begin{align*}
\triangle_{ij}^k
=
( \slashg^{-1} )^{kl} ( \dpartial_i \uh \cdot \chi_{jl} + \dpartial_j \uh \cdot \chi_{il} - \dpartial_l \uh \cdot \chi_{ij}  ),
\end{align*}
then $\dslashGamma_{ij}^k$ is given by
\begin{align*}
\dslashGamma_{ij}^k 
& = 
\slashGamma_{ij}^k 
+ \triangle_{ij}^k.
\end{align*}
\end{enumerate}

\subsection{\hyperref[step ii]{Step ii}: formulae for general foliation}\label{sec 4.3}
$\{ \vartheta \}$ is the coordinate system of each $\bSigma_u$. Suppose that $\bSigma_u$ is parameterised by $\flu$ as its graph of $s$ over $\vartheta$ domain in the $\{s, \vartheta\}$ coordinate system of $\ucalH$, i.e.
\begin{align*}
\bSigma_u = \{(s,\vartheta): s= \flu(\vartheta) \}. 
\end{align*}
 
In the following we shall introduce another conjugate null frame $\{ \dduL, \ddL' \}$ of $\bSigma_u$, which is intermediate between $\{\duL, \dL'\}$ of $\{\Sigma_s\}$ and $\{\buL^u, \bL'^u \}$ of $\bSigma_u$. Use $\cdot\cdot$ on the top to indicated the corresponding notation being associated with $\bSigma_u$ and $\{\dduL, \ddL'\}$. We describe the geometric quantities associated with $\bSigma_u$ and $\{\dduL, \ddL'\}$ first.
\begin{enumerate}[label=\alph*.,leftmargin=15pt]
\item
Let $\ddpartial_i$ be the coordinate vector of $\{ \vartheta \}$ on $\bSigma_u$, then
\begin{align*}
\ddpartial_i = \dpartial_i + (\flu)_i \dpartial_s = \dB_i^j \dpartial_j + (\flu)_i \duL,
\quad
\dB_i^j = \delta_i^j - (\flu)_i \db^j.
\end{align*}

\item
Let $\ddslashg$ as the intrinsic metric on $\bSigma_u$, then
\begin{align*}
\ddslashg_{ij} 
= 
\dB_i^k \dB_j^l \dslashg_{kl} 
=
\slashg_{ij} 
- (\slashg_{ik} \db^k (\flu)_j + \slashg_{jl} \db^l (\flu)_i)
+ (\flu)_i (\flu)_j \slashg_{kl} \db^k \db^l.
\end{align*}
Let $\ddslashepsilon$ be the volume form of $(\bSigma_u, \ddslashg)$.

\item
Introduce the conjugate null frame $\{\dduL, \ddL'\}$ of $\bSigma_u$ by
\begin{align*}
\left\{
\begin{aligned}
&
\ddL' = \dL' + \dvarepsilon' \duL + \dvarepsilon'^i \dpartial_i,
\\
&
\dduL=\duL,
\end{aligned}
\right.
\end{align*}
where
\begin{align*}
\begin{aligned}
\dvarepsilon'^i
=&
-2 ( \dslashg^{-1} )^{ik} ( \dB^{-1} )_k^j (\flu)_j,
\\
\dvarepsilon'
=&
-\vert \ddslashd \flu \vert_{\subddslashg}^2
=
-( \ddslashg^{-1} )^{ij} (\flu)_i (\flu)_j
= 
- ( \dslashg^{-1} )^{kl} ( \dB^{-1} )_k^i ( \dB^{-1} )_l^j (\flu)_i (\flu)_j,
\end{aligned}
\end{align*}
where $\ddslashd$ is the differential operator on $\bSigma_u$.

\item
$\dslashnabla$ is the covariant derivative of $(\Sigma_s, \dslashg)$, and we also use it to denote the pull back of $\dslashnabla$ to $\bSigma_u$. The precise interpretation of $\dslashnabla$ shall be understood in the context.\footnote{This is similar to the pull back of $\slashnabla$. See footnote \ref{footnote 2}. Which meaning $\dslashnabla$ is interpreted as depends on where the differentiated object is defined. For example, if $V$ is a vector field defined on $\bSigma$, then $\dslashnabla$ in $\dslashnabla V$ should be interpreted as the pull back to $\bSigma$. If it can be interpreted in both way, we will point out the precise meaning of $\dslashnabla$ to avoid ambiguity.} 
We have the following formulae for the pull back of $\dslashnabla$ to $\bSigma_u$
\begin{align*}
&
\text{$\phi$: a function on $\bSigma_u$,}
&&
\dslashnabla_i \phi = \ddpartial_i \phi, \quad \dslashnabla^i \phi= ( \dslashg^{-1} )^{ij} \ddpartial_j \phi = ( \slashg^{-1} )^{ij} \ddpartial_j \phi ,
\\
&
\text{$V$: a vector field on $\bSigma_u$,}
&&
\dslashnabla_i V^k =  \ddpartial_i V^k + \dslashGamma_{ij}^k V^j,
\\
&
\text{$T$: a tensor field on $\bSigma_u$,}
&&
\dslashnabla_{i} T_{i_1 \cdots i_k}^{j_1 \cdots j_l}
=
\ddpartial_i T_{i_1 \cdots i_k}^{j_1 \cdots j_l} 
-  \dslashGamma_{i i_m}^{r}  T_{i_1 \cdots \underset{\hat{i_m}}{r}\cdots i_k}^{j_1 \cdots j_l}
+ \dslashGamma_{i s}^{j_n} T_{i_1 \cdots i_k}^{j_1 \cdots \overset{\hat{j_n}}{s}\cdots  j_l}.
\end{align*}

\item
The connection coefficients on $\bSigma_u$ relative to $\{\dduL, \ddL' \}$ are given by the following formulae:
\begin{align}
\begin{aligned}
\dduchi_{ij}
&=
\duchi_{ij} 
+ 2 \sym  \{   (-\duchi (\db  )  )  \otimes \ddslashd \flu   \}_{ij} 
+ \duchi (\db,\db  )    (\flu)_i    (\flu)_j,
\\
\ddchi'_{ij} 
&=
\dchi'_{ij} 
+ \dvarepsilon' \duchi_{ij} 
+  ( \db \circdot \vec{\dvarepsilon}' -2  ) \dslashnabla^2_{ij} \flu
\\
&\phantom{=}
+2\sym  
\{  
\ddslashd \flu 
\otimes 
[ \dslashnabla \db \circdot \vec{\dvarepsilon}' 
- \duchi (\vec{\dvarepsilon}'  ) 
-\dvarepsilon' \duchi (\db  ) 
- \dchi' (\db  ) -2 \deta ]
\}_{ij}
\\
&\phantom{=}
+ 
[ 2\duchi (\db,\vec{\dvarepsilon}'  ) 
+\dvarepsilon' \duchi (\db,\db  ) 
+\dchi' (\db,\db  ) 
+4\deta (\db  )
-\dslashnabla_{\db} \db \circdot \vec{\dvarepsilon}' 
-\dpartial_s \db \circdot \vec{\dvarepsilon}' 
-4\duomega
]    (\flu)_i    (\flu)_j,
\\
\ddeta_i 
&= 
\deta_i 
+ \frac{1}{2} \duchi (\vec{\dvarepsilon}'  )_i 
+  
[ 2\duomega
-\deta (\db  ) 
-\frac{1}{2} \duchi (\db,\vec{\dvarepsilon}' )
]  (\flu)_i,
\end{aligned}
\end{align}
where
\begin{align*}
( \dslashnabla \db \circdot \vec{\dvarepsilon}'  )_i 
= 
\dslashg_{kl} \dslashnabla_i \db^k \dvarepsilon'^l,
\quad
\dslashnabla_{\db} \db \circdot \vec{\dvarepsilon}' 
= 
\db^i \dslashnabla_i \db^k \cdot \dvarepsilon'^l \dslashg_{kl},
\quad
\dpartial_s \db \circdot \vec{\dvarepsilon}' 
= 
\dslashg_{ij} \dpartial_s \db^i \dvarepsilon'^j.
\end{align*}
We use the circled dot $\circdot$ to denote the inner product w.r.t. the metric $\dslashg$. However since $\dslashg=\slashg$, numerically $\circdot$ is equal to the inner product $\cdot$ w.r.t. $\slashg$, thus we simply use the dot $\cdot$ to replace $\circdot$ in the following.

\item
The curvature components on $\bSigma_u$ relative to $\{\dduL, \ddL'\}$ are given by the following formulae:
\begin{align}
\begin{aligned}
	\ddualpha_{ij}&=\dB_i^k \dB_j^l \dualpha_{kl},
\\
	\ddubeta_i &= \dB_i^j \dubeta_j -\frac{1}{2} \dB_i^j \dvarepsilon^k \dualpha_{jk},
\\ 	
	\ddrho = &\drho - \dvarepsilon^j\dubeta_j + \frac{1}{4} \dvarepsilon^i \dvarepsilon^j \dualpha_{ij},
\\ 
	\ddsigma \cdot \ddslashepsilon_{ij}
	&=
	\dB_i^k \dB_j^l \dsigma \dslashepsilon_{kl} 
	+ ( \dB_j^k f_i -\dB_i^k f_j) \dubeta_k  
	+\frac{1}{2}\dB_i^k \dB_j^l \dvarepsilon^m \dslashepsilon_{lk} \dslashepsilon_{mn} \dubeta^n 
\\
	&\phantom{=}
	+ \frac{1}{2} (f_j \dB_i^k - f_i \dB_j^k) \dvarepsilon^m \dualpha_{km},
\\
	\ddbeta_i
	&=
	\dB_i^j\dbeta_j 
	- 3 f_i \drho 
	+ \frac{3}{2}\dB_i^j \dvarepsilon^k \dslashepsilon_{jk}  \dsigma
	- \dB_i^j \dvarepsilon \dubeta_j 
	+ 2 f_i \dvarepsilon^l \dubeta_l 
	+ \frac{1}{2} \dB_i^j \dvarepsilon^k \dvarepsilon^l \dslashepsilon_{kj} \dslashepsilon_{lm} \dubeta^m
\\
	&\phantom{=}
	+ \frac{1}{2} \dB_i^j \dvarepsilon^l  \dvarepsilon \dualpha_{jl}
	-\frac{1}{2} f_i \dvarepsilon^k \dvarepsilon^l \dualpha_{kl},
\\
	\ddalpha_{ij}
	&=
	\dB_i^k \dB_j^l \dalpha_{kl} 
	-2 ( \dB_i^k f_j + \dB_j^k f_i  )\dbeta_k 
	+ \dB_i^k \dB_j^l \dvarepsilon^n \dslashepsilon_{nl} \dslashepsilon_{km} \dbeta^m
\\
	&\phantom{=}
	+ \dvarepsilon^m \dB_i^k \dB_j^l \dslashepsilon_{mk} \dslashepsilon_{ln} \dbeta^n 
	-  (f_j \dB_i^k\dvarepsilon^n + f_i \dB_j^l \dvarepsilon^m) ( \drho  \dslashg_{nk} -\dsigma \dslashepsilon_{nk}) 
	+4 f_i f_j \drho 
\\
	&\phantom{=}
	+ 2\dvarepsilon \dB_i^k \dB_j^l \dslashg_{kl} \drho 
	+ \dB_i^k \dB_j^l \dvarepsilon^m \dvarepsilon^n \dslashepsilon_{mk} \dslashepsilon_{ln} \drho
	+ 2 (\dB_i^k f_j +\dB_j^k f_i  )\dvarepsilon^n \dslashepsilon_{nk} \dsigma
\\
	&\phantom{=}
	- 4 f_i f_j \dvarepsilon^n \dubeta_n
	+ 2 (\dB_j^k f_i +\dB_i^k  f_j ) \dvarepsilon \dubeta_k 
	+ (\dB_i^k \dB_j^l + \dB_j^k \dB_i^l) \dvarepsilon \dvarepsilon^n \dslashepsilon_{nl} \dslashepsilon_{km} \dubeta^m
\\
	&\phantom{=}
	+ (\dB_i^k f_j + \dB_j^k f_i ) \dvarepsilon^m \dvarepsilon^n \dslashepsilon_{km} \dslashepsilon_{nl} \dubeta^l
	+ \dB_i^k \dB_j^l \dvarepsilon^2 \dualpha_{kl}
	+ f_i f_j \dvarepsilon^m \dvarepsilon^n \dualpha_{mn}
\\
	&\phantom{=}
	- (\dB_i^k f_j + \dB_j^k  f_i) \dvarepsilon \dvarepsilon^n \dualpha_{nk}.
\end{aligned}
\label{eqn 4.4}
\end{align}
\end{enumerate}

Now we apply the above to derive the formulae of the geometric quantities associated with $\{\bSigma_u\}$ and $\{\buL^u, \bL'^u\}$.
\begin{enumerate}[label=\alph*.,leftmargin=15pt]
\item
Let $\balu$ be the lapse function of $\{\bSigma_u\}$ relative to $\{ \Sigma_s \}$, i.e.
\begin{align*}
\balu \dduL u =1.
\end{align*}
Then the conjugate null frame $\{ \buL^u, \bL'^u\}$ is given by
\begin{align*}
\buL^u = \balu \dduL,
\quad
\bL'^u = \balu^{-1} \ddL'.
\end{align*}
Then the coordinate vector $\ddpartial_u$ of the coordinate system $\{ u, \vartheta \}$ on $\ucalH$ is given by
\begin{align*}
	\left\{
	\begin{aligned}
		&
		\ddpartial_i = \dpartial_i + (\flu)_i \dpartial_s,
	\\
		&
		\buL = \balu \dduL = \balu \dpartial_s + \balu \db^i \dpartial_i,
	\end{aligned}
	\right.
	\quad
	\Rightarrow
	\quad
	\left\{
	\begin{aligned}
		&
		\buL = \ddpartial_u + \balu \db^i \ddpartial_i,
	\\
		&
		\ddpartial_u = \balu [ 1 - \db^i (\flu)_i ] \dpartial_s.
	\end{aligned}
	\right.
\end{align*}

\item
Denote the intrinsic metric on $\bSigma_u$ by $\bslashglu$. $\bslashglu = \ddslashg$. 

\item
The connection coefficients on $\bSigma_u$ relative to $\{\buL^u, \bL'^u\}$ are given by
\begin{align*}
\buchilu
&=
\balu \cdot \dduchi,
\\
\bchilu'
&=
\balu^{-1} \cdot \ddchi',
\\
\btalu
&=
\ddeta + \bslashd \log \balu,
\\
\buomegalu
&=
\balu \cdot \duomega + \frac{1}{2} \frac{\ed}{\ed u} \balu.
\end{align*}

\item
The curvature components on $\bSigma_u$  relative to $\{\buL^u, \bL'^u\}$ are given by
\begin{align*}
\bualphalu = \balu^2 \cdot \ddualpha,
\quad
\bubetalu = \balu \cdot \ddubeta,
\quad
\brholu = \ddrho,
\quad
\bsigmalu = \ddsigma,
\quad
\bbetalu = \balu^{-1} \cdot \ddbeta,
\quad
\balphalu = \balu^{-2} \cdot \ddalpha.
\end{align*}

\end{enumerate}

\section{Estimate of geometric quantity associated with foliation}

In this section, we apply the formulae of geometric quantities obtained in section \ref{sec 4} to obtain their estimates.

\subsection{Estimate of differential and Hessian of parameterisation function $\uh$}
In order to obtain the estimates of the geometric quantities associated with the foliation $\{\bSigma_u \}$, we see that the estimates of the differential $\dslashd \uh$ and $\circnabla^2 \uh$ on each leaf $\bSigma_u$ are necessary, since they naturally appear in formulae \eqref{eqn 4.1}, \eqref{eqn 4.2}. A general method to estimate them was described in \cite{L2018}\cite{L2020}, thus we adopt the method for our case. Before we start to estimate $\dslashd \uh$ and $\circnabla^2 \uh$, we review the basic knowledge of rotation vector fields on the sphere.

\subsubsection{Basic of rotational vector field}\label{sec 5.1.1}
We introduce the rotational vector field on a surface $\Sigma$ as follows.
\begin{enumerate}[label=\roman*.]
\item
There exists a natural mapping from the Schwarzschild spacetime $(M,g_S)$ to $(\mathbb{S}^2, \circg)$, such that $g_S|_{\Sigma_{0,0}} = r_0^2 \circg$. 

\item
Restricting the mapping on $\Sigma_s$ gives the mapping between $\Sigma$ and $(\mathbb{S}^2, \circg)$, then one can pull back the rotational vector field on $(\mathbb{S}^2,\circg)$ to $\Sigma$. 

\item
Fix an isometric embedding from $(\mathbb{S}^2,\circg)$ to $(\mathbb{R}^3, \delta_{ij} \ed x^i \ed x^j)$, then define the rotational vector field $R_i$ as $\sum_{j,k=1,2,3} \epsilon_{ijk} x^j \partial_k$ where $\epsilon_{ijk}$ is the permutation symbol of $\{1,2,3\}$. Thus we define the rotational vector field on $\Sigma$ as the pull back of $R_i$ as in step ii. We also pull back the coordinate functions $x^1, x^2, x^3$ to $\Sigma$ which are useful when we state the properties of the rotational vector field.
\end{enumerate}
Then we introduce two rotational vector fields in the following.
\begin{enumerate}
\item
Applying the above construction to the double null foliation $\{\Sigma_{s,\us}\}$ of $(M,g)$, then we obtain the rotational vector fields on each sphere $\Sigma_{s,\us}$, which we also denote by $R_i$. We have that
\begin{align*}
[\partial_s, R_i] = [\partial_{\us}, R_i] =0.
\end{align*}
\item
Applying the above construction to the foliation $\{\Sigma_s\}$ of the incoming null hypersurface $\ucalH$, we obtain the rotational vector field on each sphere $\Sigma_s$ of the foliation, denoted by $\dR_i$. We have that
\begin{align*}
\dR_i = R_i + (\dR_i \uh) \partial_{\us},
\quad
\lie_{\dpartial_s} \dR_i = [\dpartial_s, \dR_i] = 0.
\end{align*}
\end{enumerate}

We introduce the notation of the rotational vector component of a tensor field on $(\Sigma,\circg)$, similar to the notation of the coordinate component.
\begin{enumerate}[label=\alph*.]
\item
For a vector field $v$ and a one-form $\omega$ on $\Sigma$, the corresponding $R_i$ rotational vector components are
\begin{align*}
v^{R,i} = \circg(v, R_i),
\quad
\omega_{R,i} = \omega(R_i).
\end{align*}
We use the letter $R$ in the subscript to indicate that the corresponding component is the rotational vector component.

\item
The rotational vector component of a tensor $T$ can be defined inductively by above, or explicitly
\begin{align*}
T^{R,j_1 \cdots j_l}_{R,i_1 \cdots i_k} = \circg(T(R_{i_1}, \cdots, R_{i_k}), R_{j_1} \otimes \cdots \otimes R_{j_l}).
\end{align*}
\end{enumerate}
Furthermore we introduce the notation of the mixed component with rotational vector and coordinate of a tensor field $(\Sigma,\circg)$ as follows.
\begin{enumerate}[label=\alph*.]\addtocounter{enumi}{2}
\item
For a $(2,0)$-tensor field $v$, its mixed component with the first slot being rotational vector component and second slot being coordinate component is defined by
\begin{align*}
v^{R,i \overline{j}}
=
\circg(v, R_i \otimes \partial_l) \circg^{jl}.
\end{align*}
We use the overlined indices to denote the coordinate component. Similarly for a $(0,2)$-tensor field $\omega$, its mixed component is defined by
\begin{align*}
\omega_{R,i \overline{j}} = \omega(R_i, \partial_j).
\end{align*}

\item
For a general tensor field $T$, its mixed component is defined by
\begin{align*}
T^{R,j_1 \cdots \overline{j_r} \cdots j_l}_{R,i_1 \cdots \overline{i_s} \cdots i_k} = \circg(T(R_{i_1}, \cdots, \partial_{i_s}\cdots, R_{i_k}), R_{j_1} \otimes \cdots \otimes \partial_{j_r'} \otimes \cdots \otimes R_{j_l}) \circg^{j_r' j_r}.
\end{align*}
\end{enumerate}

We list several useful properties of the rotational vector component and the rotational vector derivative in appendix \ref{appen r.v.}.

\subsubsection{Propagation equations of $1$st and $2$nd rotational vector derivatives of $\uh$}
Recall the propagation equation \eqref{eqn 3.1} of the parameterisation function $\uh$,
\begin{align}
\dpartial_s \uh = -b^i \dpartial_i \uh + \Omega^2 (\slashg^{-1})^{ij} \dpartial_i \uh \dpartial_j \uh.
\tag{\ref{eqn 3.1}\ensuremath{'}}
\end{align}
Differentiate the above equation with respect to $\dR_k$ and write in terms of the rotational vector component and the mixed component, we obtain that
\begin{align}
\begin{aligned}
\dpartial_s (\dR_k \uh)
&=
[- b^i + 2 \Omega^2 ( \slashg^{-1} )^{R,\overline{i}j} (\dR_j \uh)] \dpartial_i ( \dR_k \uh )
\\
&\phantom{=}
- [ R_k, b]^{R,i} (\dR_i \uh)
- ( \dR_k \uh ) \partial_{\us} b^{R,i} (\dR_i \uh)
+ ( \lie_{R_k} (\Omega^2 \slashg^{-1} ) )^{R,ij} (\dR_i \uh) (\dR_j \uh)
\\
&\phantom{=}
+ (\dR_k \uh ) ( \partial_{\us} (\Omega^2 \slashg^{-1} ) )^{R,ij} (\dR_i \uh) (\dR_j \uh).
\end{aligned}
\label{eqn 5.1}
\end{align}
Introduce the simplified notation $\uh_{\dR,k}$ to denote $\dR_k \uh$ and use the square bracket $[R,\xi]$ to denote the Lie derivative of a tensor field $\xi$ with respect to a vector field $R$, then equation \eqref{eqn 5.1} is rewritten as
\begin{align}
\begin{aligned}
\dpartial_s \uh_{\dR,k}
&=
[- b^i + 2 \Omega^2 ( \slashg^{-1} )^{R,\overline{i}j} (\dR_j \uh)] \dpartial_i \uh_{\dR,k}
\\
Q_k&
\left\{
\begin{aligned}
&
- [ R_k, b]^{R,i} \uh_{\dR,i}
- (\partial_{\us} b)^{R,i} \uh_{\dR,i} \uh_{\dR,k}
+ [ R_k, \Omega^2 \slashg^{-1} ]^{R,ij} \uh_{\dR,i} \uh_{\dR,j}
\\
&
+ ( \partial_{\us} (\Omega^2 \slashg^{-1} ) )^{R,ij} \uh_{\dR,i} \uh_{\dR,j} \uh_{\dR,k}
\end{aligned}
\right.
\end{aligned}
\tag{\ref{eqn 5.1}\ensuremath{'}}
\end{align}
Differentiate the above equation with respect to the rotational vector field again and use the notation $\uh_{\dR,lk}$ to denote $\dR_l \dR_k \uh$, we obtain that
\begin{align}
\begin{aligned}
\dpartial_s \uh_{\dR,lk} 
&=
[- b^i + 2 \Omega^2 ( \slashg^{-1} )^{R,\overline{i}j} \uh_{\dR,j}] 
\dpartial_i ( \uh_{\dR,lk} )
\\
Q_{lk}
&
\left\{
\begin{aligned}
&
-  [\dR_l, b ]^{\dR,i} \uh_{\dR,ik}  
\\
&
+ 
\{ 
2  [ \dR_l, \Omega^2 ( \slashg^{-1} )  ]^{\dR,ij}   \uh_{\dR,j}  \uh_{\dR,ik} 
+ 
2 \Omega^2 ( \slashg^{-1} )^{R,ij} \uh_{\dR,ik} \uh_{\dR,jl}
\}
\\
&
-
\{
 [ \dR_l, [ R_k, b ] ]^{\dR,i} \uh_{\dR,i} 
+
[ R_k, b ]^{R,i} \uh_{\dR,il}
\}
\\
&
- 
\{
 [ \dR_l,  \partial_{\us} b  ]^{\dR,i} \uh_{\dR,i} \uh_{\dR,k}
+ ( \partial_{\us} b )^{R,i} \uh_{\dR,il} \uh_{\dR,k}
+ ( \partial_{\us} b )^{R,i} \uh_{\dR,i} \uh_{\dR,lk}
\}
\\
&
+
\{
 [\dR_l, [ R_k, \Omega^2 \slashg^{-1} ]  ]^{\dR, ij} \uh_{\dR,i} \uh_{\dR,j}
+ [ R_k, \Omega^2 \slashg^{-1} ]^{R, ij} \uh_{\dR,il} \uh_{\dR,j}
\\
&\phantom{+  \{[}
+ [ R_k, \Omega^2 \slashg^{-1} ]^{R, ij} \uh_{\dR,i} \uh_{\dR,jl}
\}
\\
&
+ 
\{
 [\dR_l, \partial_{\us} (\Omega^2 \slashg^{-1} )  ]^{\dR,ij} \uh_{\dR,i} \uh_{\dR,j} \uh_{\dR,k}
+ ( \partial_{\us} (\Omega^2 \slashg^{-1} ) )^{R,ij} \uh_{\dR,il} \uh_{\dR,j} \uh_{\dR,k}
\\
&\phantom{+  \{[}
+ ( \partial_{\us} (\Omega^2 \slashg^{-1} ) )^{R,ij} \uh_{\dR,i} \uh_{\dR,jl} \uh_{\dR,k}
+ ( \partial_{\us} (\Omega^2 \slashg^{-1} ) )^{R,ij} \uh_{\dR,i} \uh_{\dR,j} \uh_{\dR,lk}
\},
\end{aligned}
\right.
\end{aligned}
\label{eqn 5.2}
\end{align}
where $\dR_l = R_l + \uh_{\dR,l} \partial_{\us}$,
\begin{align*}
[\dR_l, \xi]^{\dR,i_1\cdots i_k} = [ R_l, \xi]^{R, i_1\cdots i_k} + \uh_{\dR,l} \big( \partial_{\us} \xi \big)^{R,i_1\cdots i_k},
\end{align*}
with $\xi$ being $b$, $\Omega^2 \slashg^{-1}$, $[R_k, b]$, $\partial_{\us} b$, $[ R_k, \Omega^2 \slashg^{-1} ]$, $\partial_{\us} (\Omega^2 \slashg^{-1} )$.

We shall use equations \eqref{eqn 5.1}, \eqref{eqn 5.2} to derive estimate the differential $\dslashd \uh$ and the Hessian $\dslashnabla^2 \uh$ on a spacelike surface $\Sigma$ in the incoming null hypersurface $\ucalH$. The method is similar to the proof of proposition \ref{prop 3.3}, which we summarise in the following lemma.
\begin{lemma}\label{lem 5.1}
Let $\Sigma$ be a spacelike surface in the incoming null hypersurface $\ucalH$. Suppose that $\ucalH$ is parameterised by $\us = \uh(s,\vartheta)$, and $\Sigma$ is parameterised by $(\ufl{s=0}, f)$ and $(\uf,f)$ by the two methods in section \ref{sec 3.2}. Construct the family of surfaces $\{S_t\}_{t\in[0,1]}$ in $\ucalH$ deforming from $\Sigma_{s=0}$ to $\Sigma$ that $S_t$ is the surface with the parameterisation $(\ufl{s=0}, tf)$. Introduce the notations $\uhlt_{\dR,k}$ and $\uhlt_{\dR,lk}$ to denote the restrictions of $\uh_{\dR,k}$ and $\uh_{\dR,lk}$ on $S_t$, i.e.
\begin{align*}
\uhlt_{\dR,k} = \uh_{\dR,k} |_{S_t},
\quad
\uhlt_{\dR,lk} = \uh_{\dR,lk} |_{S_t}.
\end{align*}
We have the following propagation equations of $\uhlt_{\dR,k}$ and $\uhlt_{\dR,lk}$ along $\{S_t\}_{t\in[0,1]}$,
\begin{align}
\begin{aligned}
&
\partial_t\, \uhlt_{R,k}
=
\Xlt_{\uh}^i \ddpartial_i\, \uhlt_{\dR,k} + \relt_{\uh,k},
\quad
\uhl{t=0}_{\dR,k} = \ufl{s=0}_{\dR,k}= \dR_k ( \ufl{s=0} ),
\\
&
\partial_t \uhlt_{\dR,lk}
=
\Xlt_{\uh}^i \ddpartial_i \uhlt_{\dR,lk}
+
\relt_{\uh,lk},
\quad
\uhl{t=0}_{\dR,lk} = \ufl{s=0}_{\dR,lk} = \dR_l \dR_k (\ufl{s=0}),
\end{aligned}
\label{eqn 5.3}
\end{align}
where $\ddpartial_i$ is the tangential coordinate vector field on $S_t$ and
\begin{align*}
&
\Xlt_{\uh} 
=
f [ 1- t b^m f_m - 2 t \Omega^2 ( \slashg^{-1} )^{R,\overline{m}j} f_m (\uhlt_{\dR,j})  ]^{-1}
\cdot
[ - b^i + 2 \Omega^2 ( \slashg^{-1} )^{R,\overline{i}j} (\uhlt_{\dR,j})  ]
\ddpartial_i,
\\
&
\relt_{\uh,k}
=
f [ 1- t b^m f_m - 2 t \Omega^2 ( \slashg^{-1} )^{R,\overline{m}j} f_m (\uhlt_{\dR,j}) ]^{-1}
\cdot
Q_k|_{S_t},
\\
&
\relt_{\uh,lk}
=
f [ 1- t b^m f_m - 2 t \Omega^2 ( \slashg^{-1} )^{R,\overline{m}j} f_m (\uhlt_{\dR,j}) ]^{-1}
\cdot
Q_{lk}|_{S_t}.
\end{align*}
\end{lemma}

\subsubsection{Estimate of differential of $\uh$}
We state the estimate of the differential $\dslashd \uh$ on a surface $\Sigma$ in the following proposition.
\begin{proposition}\label{prop 5.2}
Let $\Sigma$ be a spacelike surface in the incoming null hypersurface $\ucalH$. Suppose that $\ucalH$ is parameterised by $\uh$ and $\Sigma$ is parameterised by $(\ufl{s=0}, f)$ and $(\uf,f)$ with the two methods described in section \ref{sec 3.2}. Suppose that the parameterisation functions $\ufl{s=0}$, $f$ satisfies that
\begin{align*}
\Vert \slashd \ufl{s=0} \Vert^{n+1,p} \leq \udelta_o r_0,
\quad
\vert \overline{\ufl{s=0}}^{\circg} \vert \leq \udelta_m r_0,
\quad
\Vert \slashd f \Vert^{n+1,p} \leq \delta_o (r_0+\os),
\quad
\vert \overline{f}^{\circg} \vert = \os > -\frac{r_0}{2},
\end{align*}
where $n\geq 1+n_p$, $p >1$. Then there exist a small positive constant $\delta$ and a constant $c_{\uh}$ both depending on $n,p$, such that if $\epsilon, \udelta_o, \udelta_m, \delta_o$ are suitably bounded that $\epsilon, \udelta_o, \epsilon \udelta_m, \delta_o \leq \delta$, then the differential $(\dslashd \uh)|_{\Sigma}$ satisfies the following estimate
\begin{align}\label{eqn 5.4}
\Vert (\dslashd \uh)|_{\Sigma} \Vert^{n+1,p} 
\leq
c_{\uh} \Vert \dslashd \ufl{s=0} \Vert^{n+1,p}
\leq
c_{\uh} \udelta_o r_0.
\end{align}
\end{proposition}
\begin{proof}
We adopt the construction in lemma \ref{lem 5.1}. Assume that $\delta \leq \frac{1}{2}$ is sufficiently small such that proposition \ref{prop 3.3} holds. Recalling the notations $\ud_o$, $\ud_m$ in formulae \eqref{eqn 3.5}, we have that
\begin{align*}
\Vert \slashd \uflt \Vert^{n,p} \leq \ud_o r_0,
\quad
\vert \overline{\uflt}^{\circg} \vert \leq \ud_m r_0.
\end{align*}
We prove the proposition by the bootstrap argument. Introduce the following bootstrap assumption.
\begin{assum}
Estimate \eqref{eqn 5.4} holds for $(\dslashd \uh)|_{S_t}$ in the closed interval $[0,t_a] \subset [0,1]$.
\end{assum}
Firstly as long as $c_{\uh}>1$, there exists a small number $t_a>0$ such that the bootstrap assumption holds by the continuity argument.

We shall show that if the bootstrap assumption holds, then estimate \eqref{eqn 5.4} can be strengthen to a strict inequality at $t=t_a$ for suitably chosen $\delta$ and $c_{\uh}$ independent of $t_a$ in the following. We estimate $(\dslashd \uh)|_{\Sigma_t}$ by integrating equation \eqref{eqn 5.3} for $\uhlt_{R,k}$. Thus it is necessary to estimate the terms $\Xlt_{\uh}$, $\relt_{\uh,k}$ in the interval $[0, t_a]$.

Assume that $C_{\uh} \delta \leq 1$. We have the following estimates for $\Xlt_{\uh}$, $\relt_{\uh,k}$ in $[0,t_a]$,
\begin{align*}
&
\Vert \Xlt_{\uh} \Vert^{n+1,p}
\leq
c(n,p) 
\frac{(|\os|/r_0 + \delta_o) r_0^2}{(r_0 + t \os)^2} 
\Big[ \frac{\epsilon (\ud_o + \ud_m)r_0}{r_0 + t \os} + c_{\uh} \udelta_o \Big]
\leq
c(n,p) 
\frac{(|\os|/r_0 + \delta_o) r_0^2}{(r_0 + t \os)^2},
\\
&
\Vert \relt_{\uh,k} \Vert^{n+1,p}
\leq
c(n,p) 
\frac{(|\os|/r_0 + \delta_o) r_0^2}{(r_0 + t \os)^2} 
\Big[ \frac{\epsilon (\ud_o + \ud_m) \cdot c_{\uh} \udelta_o r_0}{r_0 + t \os} + \epsilon c_{\uh}^2 \udelta_o^2 + \frac{c_{\uh}^3 \udelta_o^3 r_0}{r_0 + t \os} \Big] r_0,
\end{align*}
which follows from
\begin{align}
\begin{aligned}
\Vert b \Vert_{S_t}^{n+1,p}, \Vert [R_k, b] \Vert_{S_t}^{n+1,p} 
\leq 
c(n,p) \frac{\epsilon ( \ud_m + \ud_o ) r_0^2}{(r_0 + t \os)^3},
\\
\Vert \partial_{\us} b \Vert_{S_t}^{n+1,p} 
\leq 
c(n,p) \frac{\epsilon r_0}{(r_0 + t \os)^3},
\\
\Vert [ R_k, \Omega^2 \slashg^{-1}] \Vert_{S_t}^{n+1,p} 
\leq 
c(n,p) \frac{\epsilon }{(r_0 + t \os)^2},
\\
\Vert \partial_{\us} (\Omega^2 \slashg^{-1}) \Vert_{S_t}^{n+1,p} 
\leq
\frac{c(n,p)}{(r_0 + t \os)^3},
\end{aligned}
\label{eqn 5.5}
\end{align}
essentially obtained in appendix \ref{appen prop 3.3}. Then by Gronwall's inequality, we obtain that
\begin{align*}
\Vert (\dslashd \uh)|_{S_t} \Vert^{n+1,p}
&\leq
c(n,p) \sum_{k=1}^3 \Vert \uhlt_{\dR,k} \Vert^{n+1,p}
\\
&\leq
c(n,p) \sum_{k=1}^3 \{ \Vert \uhl{t=0}_{\dR,k} \Vert^{n+1,p}  + \int_0^t \Vert \relt_{\uh,k} \Vert^{n+1,p} \ed t\}
\\
&\leq
c(n,p) [1 +   \epsilon (\ud_o + \ud_m) c_{\uh} + \epsilon \udelta_o c_{\uh}^2+ \udelta_o^2 c_{\uh}^3]  \udelta_o r_0 
\\
&\leq
c(n,p) (1 +   \delta c_{\uh} + \delta^2 c_{\uh}^2+ \delta^2 c_{\uh}^3) \udelta_o r_0.
\end{align*}
We choose suitable $\delta$ and $c_{\uh}$ such that estimate \eqref{eqn 5.4} can be strengthen to a strict inequality by the above estimate. It is sufficient to require that
\begin{align*}
c(n,p) (1 +   \delta c_{\uh} + \delta^2 c_{\uh}^2+ \delta^2 c_{\uh}^3) < c_{\uh}.
\end{align*}
We choose that $c_{\uh} = 2 c(n,p)$ and $\delta$ suitably small such that the above inequality holds. Then for such $\delta$ and $c_{\uh}$, the bootstrap argument implies the proposition.
\end{proof}
\begin{remark}\label{rem 5.3}
Proposition \ref{prop 5.2} requires the $L^{\infty}$ bounds of $b, \Omega, \slashg$ up to $(n+2)$-th order derivatives given in definition \ref{def 2.4}, which originates from estimates \eqref{eqn 5.5}.

Comparing the orders of differentiability in the estimate of $(\dslashd \uh)|_{\Sigma}$ in proposition \ref{prop 5.2} and in the estimate of $\slashd \uf$ in proposition \ref{prop 3.3}, the former one is one order higher and is optimal. We show that the estimate \eqref{eqn 5.4} of $(\dslashd \uh)|_{\Sigma}$ will imply the estimate of $\slashd \uf$ of the optimal order of differentiability, improving estimate \eqref{eqn 3.2} of proposition \ref{prop 3.3}. See the following proposition. 
\end{remark}
\begin{proposition}\label{prop 5.4}
Adopt the assumptions of propositions \ref{prop 3.3} and \ref{prop 5.2}. Let $\delta$ be suitably small such that propositions \ref{prop 3.3} and \ref{prop 5.2} hold. Then there exists a constant $c_{o'}$ depending on $n,p$ such that
\begin{align*}
\Vert \slashd \uf \Vert^{n+1,p}
\leq
c_{o'} \udelta_o r_0.
\end{align*}
\end{proposition}
\begin{proof}
We have that
\begin{align*}
\slashd \uf 
= 
\ddslashd ( \uh|_{\Sigma} ) = (\dslashd \uh)|_{\Sigma} + \slashd f \cdot (\partial_s \uh)|_{\Sigma},
\end{align*}
then substitute equation \eqref{eqn 3.1},
\begin{align*}
\slashd \uf 
= 
(\dslashd \uh)|_{\Sigma} + \slashd f \cdot (-b^i \dpartial_i \uh + \Omega^2 (\slashg^{-1})^{ij} \dpartial_i \uh \dpartial_j \uh)|_{\Sigma}.
\end{align*}
Thus applying estimates \eqref{eqn 5.4} and \eqref{eqn 5.5}, we obtain that
\begin{align*}
\Vert \ddslashd \uf \Vert^{n+1,p}
\leq
c_{\uh} \udelta_o r_0 + c(n,p) \delta_o (\epsilon (\ud_m+\ud_o) c_{\uh} \udelta_o  + c_{\uh}^2 \udelta_o^2 ) r_0
\leq
(c_{\uh} + c(n,p)) \udelta_o r_0.
\end{align*}
Then choose $c_{o'} = c_{\uh} + c(n,p)$ and the proposition follows.
\end{proof}

\subsubsection{Estimate of Hessian of $\uh$}
We state the estimate of the Hessian $\circnabla^2 \uh$ on a surface $\Sigma$ in the following proposition.
\begin{proposition}\label{prop 5.5}
Under the same assumptions of proposition \ref{prop 5.2}, there exist a small positive constant $\delta$ and a constant $c_{\uh,2}$ both depending on $n,p$, such that if $\epsilon, \udelta_o, \udelta_m, \delta_o$ are suitably bounded that $\epsilon, \udelta_o, \epsilon \udelta_m, \delta_o \leq \delta$, then the Hessian $( \circnabla^2 \uh )|_{\Sigma}$ satisfies the following estimate
\begin{align}
	\Vert ( \circnabla^2 \uh )|_{\Sigma} \Vert^{n,p} 
	\leq 
	c_{\uh,2} \Vert \dslashd \ufl{s=0} \Vert^{n+1,p} 
	\leq 
	c_{\uh,2} \udelta_o r_0.
	\label{eqn 5.6}
\end{align}
\end{proposition}
\begin{proof}
By appendix \ref{appen r.v.} and proposition \ref{prop 5.2}, it is sufficient to estimate $\dR_l \dR_k \uh$ on $\Sigma$. We follow the same strategy as in the proof of proposition \ref{prop 5.2} by the bootstrap argument.

Assume that $\delta \leq \frac{1}{2}$ is sufficiently small such that propositions \ref{prop 3.3}, \ref{prop 5.2} hold. Introduce the following bootstrap assumption.
\begin{assum}
Estimate \ref{eqn 5.6} holds for $(\circnabla^2 \uh)|_{S_t}$ in the closed interval $[0,t_a] \subset [0,1]$.
\end{assum}
Choose that $c_{\uh,2} >1$, then there exists a small number $t_a>0$ such that the bootstrap assumption holds by the continuity argument. Moreover we choose $c_{\uh,2} \geq c_{\uh}$ and $c_{\uh,2} \delta \leq 1$.

As shown in the proof of proposition \ref{prop 5.2}, $\Xlt_{\uh}$ satisfies the estimate in $[0,1]$
\begin{align*}
	&
	\Vert \Xlt_{\uh} \Vert^{n+1,p}
	\leq
	c(n,p) 
	\frac{(|\os|/r_0 + \delta_o) r_0^2}{(r_0 + t \os)^2} 
	\Big[ \frac{\epsilon (\ud_o + \ud_m)r_0}{r_0 + t \os} + c_{\uh} \udelta_o \Big]
	\leq
	c(n,p) 
	\frac{(|\os|/r_0 + \delta_o) r_0^2}{(r_0 + t \os)^2}.
\end{align*}
We have the following estimate for $\relt_{\uh,l,k}$ in $[0,t_a]$
\begin{align*}
	\Vert \relt_{\uh,lk} \Vert^{n,p}
	\leq
	c(n,p) 
	\frac{(|\os|/r_0 + \delta_o) r_0^2}{(r_0 + t \os)^2} 
	\Big[ \frac{\epsilon (\ud_o + \ud_m) \cdot c_{\uh,2} \udelta_o r_0}{r_0 + t \os} + c_{\uh,2}^2 \udelta_o^2 \Big] r_0,
\end{align*}
which follows from the $\Vert \cdot \Vert_{S_t}^{n,p}$ estimates of $b$, $\Omega$, $\slashg$ and their derivatives up to $2$nd order in appendix \ref{appen prop 3.3}.

By Gronwall's inequality and lemma \ref{lem 5.1}, we obtain that
\begin{align*}
	\Vert \uhlt_{\dR, lk} \Vert^{n,p}
	&\leq
	c(n,p) \{ \Vert \uhl{t=0}_{\dR,lk} \Vert^{n+1,p}  + \int_0^t \Vert \relt_{\uh,lk} \Vert^{n+1,p} \ed t\}
\\
	&\leq
	c(n,p) [1 +   \delta c_{\uh,2 } + \delta c_{\uh,2}^2]  \delta r_0.
\\
	\Vert (\circnabla^2 \uh)_{S_t} \Vert^{n,p}
	&\leq
	c(n,p) \{ \sum_{l,k=1,2,3} \Vert \uhlt_{\dR, lk} \Vert^{n,p} + \sum_{k=1,2,3} \Vert \uhlt_{\dR, k} \Vert^{n,p} \}
\\
	&\leq
	c(n,p) [1 +   \delta c_{\uh,2 } + \delta c_{\uh,2}^2 + c_{\uh}]  \delta r_0.
\end{align*}
We choose $c_{\uh,2} = c(n,p)(2+c_{\uh})$ and $\delta$ sufficiently small, then
\begin{align*}
	c(n,p) [1 +   \delta c_{\uh,2 } + \delta c_{\uh,2}^2 + c_{\uh}] < c_{\uh,2},
\end{align*}
thus the bootstrap argument implies the proposition.
\end{proof}
\begin{remark}\label{rem 5.6}
Proposition \ref{prop 5.5} requires the $L^{\infty}$ bounds of $b, \Omega, \slashg$ up to $(n+2)$-th order derivatives given in definition \ref{def 2.4}, same as in proposition \ref{prop 5.2}.
\end{remark}

\subsection{Decomposition of geometric quantity}\label{sec 5.2}
We decompose the geometric quantities associated with the foliation $\{ \bSigma_u \}$ of an incoming null hypersurface $\ucalH$ for the preparations to estimate them. The idea of the decomposition is subtracting a part of the higher order smallness comparing to the geometric quantity to obtain the main part being as simple as possible. We shall call the main part as the first order main part, and call the subtracted quadratic small part as the high order remainder. Denote the first order main part by $\lo{\cdot}$ and the high order remainder by $\hi{\cdot}$.

\subsubsection{Decomposition of geometric quantity associated with $s$ foliation}\label{sec 5.2.1}
In section \ref{sec 4.2}, we obtain the formulae for the geometric quantities associated with the special foliation $\{\Sigma_s\}$ by coordinate $s$. In the following we give the decomposition of these geometric quantities by their formulae.
\begin{enumerate}[label= \alph*.]
\item
For the shifting vector $\db$, the metric $\dslashg_{ij}$ and its volume form $\dslashepsilon$ on the foliation $\{\Sigma_s\}$ of $\ucalH$:
\begin{subequations}
\begin{align}
	\db:
	&
	\left\{
	\begin{aligned}
		\lo{\db}
		&=
		- 2\Omega_S^2 (\slashg_S^{-1})^{ij} \uh_j,
	\\
		\hi{\db^i}
		&=
		b^i 
		- 2(\Omega^2 - \Omega_S^2) (\slashg^{-1})^{ij} \uh_j 
		- 2 \Omega_S^2 ( \slashg^{-1} - \slashg_S^{-1})^{ij} \uh_j,
	\end{aligned}
	\right.
\\
	\dslashg:
	&
	\left\{
	\begin{aligned}
		\lo{\dslashg} 
		&=
		( r_S|_{\Sigma_s} )^2 \circg, 
	\\
		\hi{\dslashg_{ij}}
		&=
		\slashg_{ij} - ( r_S|_{\Sigma_s} )^2 \circg_{ij}.
	\end{aligned}
	\right.
\\
	\dslashg^{-1}:
	&
	\left\{
	\begin{aligned}
		\lo{\dslashg^{-1}} 
		&=
		( r_S|_{\Sigma_s} )^{-2}\cdot  \circg^{-1} , 
	\\
		\hi{(\dslashg^{-1})^{ij}}
		&=
		- ( r_S|_{\Sigma_s} )^{-2} ( \circg^{-1} )^{ik} \cdot \hi{ \dslashg_{kl} }  \cdot (\dslashg^{-1})^{lj}.
	\end{aligned}
	\right.
\\
	\dslashepsilon:
	&
	\left\{
	\begin{aligned}
		\lo{\dslashepsilon}
		&=
		(r_S|_{\Sigma_s})^2 \circepsilon,
	\\
		\hi{\dslashepsilon}
		&=
		\frac{
		(\slashg_S^{-1})^{kl} \hi{\dslashg_{kl}} 
		+ \det \{ (\slashg_S^{-1})^{kr} \hi{\dslashg_{rl}} \}
		}{
		\sqrt{ \det \{ (\slashg_S^{-1})^{kr} \dslashg_{rl} \} }
		+1
		}
		(r_S|_{\bSigma_u})^2 \circepsilon.
	\end{aligned}
	\right.
\end{align}
\end{subequations}

\item 
For the connection coefficients on $\ucalH$ relative to $\{\Sigma_s\}$ and $\{\duL, \dL'\}$:
\begin{subequations}
\begin{align}
	&
	\begin{aligned}
		\dchi':
		&
		\left\{
		\begin{aligned}
			\lo{\dchi'}
			&=
			\chi'_S|_{\Sigma_s},
			\\
			\hi{\dchi'}
			&=
			\chi'-\chi'_S|_{\Sigma_s},
		\end{aligned}
		\right.
	\\
		\dtr \dchi':
		&
		\left\{
		\begin{aligned}
			\lo{\dtr \dchi'}
			&=
			\tr \chi'_S|_{\Sigma_s},
			\\
			\hi{\dtr \dchi'}
			&=
			(\dslashg^{-1})^{ij} \hi{\dchi'_{ij}}
			+ \hi{(\dslashg^{-1})^{ij}} \dchi'_{ij},
		\end{aligned}
		\right.
	\\
		\hatdchi':
		&
		\left\{
		\begin{aligned}
			&
			\lo{\hatdchi'}
			=
			0,
		\\
			&
			\hi{\hatdchi'}
			=
			- \frac{1}{2} 
			(\hi{(\dslashg^{-1})^{kl}} \lo{\dchi'_{kl}} \dslashg
			+ \lo{\dtr \dchi'} \hi{\dslashg})
			+\widehat{\hi{\dchi'_{ij}}},
		\end{aligned}
		\right.
	\end{aligned}
\\
	&
	\begin{aligned}
		\duchi:
		&
		\left\{
		\begin{aligned}
			\lo{\duchi}
			&=
			\uchi_S|_{\Sigma_s} - 2\Omega_S^2 \circnabla^2 \uh,
		\\
			\hi{\duchi}
			&=
			(\uchi - \uchi_S)|_{\Sigma_s}
			- 2 (\Omega^2 - \Omega_S^2) \slashnabla^2 \uh
			- 2 \Omega_S^2 (\slashnabla^2 \uh - \circnabla^2 \uh)
		\\
			&\phantom{=}
			-\Omega^2 \vert \dslashd \uh \vert^2_{\subslashg} \chi
			-4\Omega^2 \sym \{\ueta \otimes \dslashd \uh \}
			- 4\omega\Omega^2 (\dslashd \uh \otimes \dslashd \uh )
		\\
			&\phantom{=}
			+ 4\Omega^2 \sym \{\dslashd\uh \otimes ( \chi \cdot \dslashd \uh ) \},
		\end{aligned}
		\right.
	\\
		\dtr \duchi:
		&
		\left\{
		\begin{aligned}
			\lo{\dtr \duchi}
			&=
			\lo{(\dslashg^{-1})^{ij}} \lo{\duchi}
			=
			\tr\uchi_S|_{\Sigma_s} - 2\Omega_S^2 \circDelta \uh,
		\\
			\hi{\dtr \duchi}
			&=
			(\dslashg^{-1})^{ij} \cdot \hi{\duchi_{ij}}
			+\hi{(\dslashg^{-1})^{ij}} \cdot \duchi_{ij},
		\end{aligned}
		\right.	
	\\
		\hatduchi:
		&
		\left\{
		\begin{aligned}
			&
			\lo{\hatduchi}
			=
			-2 \Omega_S^2 ( \circnabla^2 - \frac{1}{2} \circDelta) \uh,
		\\
			&
			\hi{\hatduchi}
			=
			- \frac{1}{2} 
			(\hi{(\dslashg^{-1})^{kl}} \lo{\duchi_{kl}} \dslashg
			+ \lo{\dtr \duchi} \hi{\dslashg})
			+\widehat{\hi{\duchi_{ij}}},
		\end{aligned}
		\right.
	\end{aligned}
\\
	&
	\begin{aligned}
		\deta:
		&
		\left\{
		\begin{aligned}
			\lo{\deta_i}
			&=
			(\chi_S|_{\Sigma_s})_{ij} \cdot (\slashg_S^{-1})^{jk} \cdot \uh_k,
		\\
			\hi{\deta_i}
			&=
			\eta_i
			+(\chi - \chi_S|_{\Sigma_s})_{ij} \cdot (\slashg^{-1})^{jk} \cdot \uh_k
			+(\chi_S|_{\Sigma_s})_{ij} \cdot (\slashg^{-1} - \slashg_S^{-1})^{jk} \cdot \uh_k,
		\end{aligned}
		\right.
	\end{aligned}
\\
	&
	\begin{aligned}
		\duomega:
		&
		\left\{
		\begin{aligned}
			\lo{\duomega}
			&=
			\uomega_S|_{\Sigma_s},
		\\
			\hi{\duomega}
			&=
			(\uomega - \uomega_S|_{\Sigma_s})
			- 2\Omega^2 \eta \cdot \dslashd \uh 
			- \Omega^2 \chi ( \slashnabla \uh, \slashnabla \uh).
		\end{aligned}
		\right.
	\end{aligned}
	\label{eqn 5.8d}
\end{align}
\end{subequations}

\item
For the curvature components on $\ucalH$ relative to $\{\Sigma_s\}$ and $\{ \duL, \dL'\}$:
\begin{subequations}
\begin{align}
	&
	\dualpha:
	\left\{
	\begin{aligned}
		\lo{\dualpha}
		&=
		0,
	\\
		\hi{\dualpha}
		&=
		\dualpha,
	\end{aligned}
	\right.
\\
	&
	\dubeta:
	\left\{
	\begin{aligned}
		\lo{\dubeta}
		&=
		-3\Omega_S^2 \rho_S \dslashd \uh,
	\\
		\hi{\dubeta}
		&=
		\ubeta_i 
		-3 (\Omega^2 - \Omega_S^2) \rho \dslashd \uh  
		-3 \Omega_S^2 (\rho- \rho_S) \dslashd \uh  
		-\frac{3}{2}  \uvarepsilon^l \sigma 
		- \Omega^2 \uvarepsilon \beta_i  
	\\
		&\phantom{=}
		+ 2 \Omega^2 \uh_i \uvarepsilon^l \beta_l 
		+\frac{1}{2}  \uvarepsilon^k \uvarepsilon^l \slashepsilon_{ki} \slashepsilon_{lh} \beta^h
		- \frac{1}{2} \Omega^2 \uh_i \uvarepsilon^k \uvarepsilon^l \alpha_{kl}
		+ \frac{1}{2} \Omega^2 \uvarepsilon \uvarepsilon^l \alpha_{il},
	\end{aligned}
	\right.
\\
	&
	\dsigma:
	\left\{
	\begin{aligned}
		\lo{\dsigma}
		&=
		0,
	\\
		\hi{\dsigma}
		&=
		\dsigma.
	\end{aligned}
	\right.
\\
	&
	\drho:
	\left\{
	\begin{aligned}
		\lo{\drho}
		&=
		\rho_S|_{\Sigma_s},
	\\
		\hi{\drho}
		&=
		(\rho - \rho_S)|_{\Sigma_s}
		- \uvarepsilon^l \beta_l 
		+ \frac{1}{4} \uvarepsilon^k \uvarepsilon^l \alpha_{kl},
	\end{aligned}
	\right.
\\
	&
	\dbeta:
	\left\{
	\begin{aligned}
		\lo{\dbeta}
		&=
		0,
	\\
		\hi{\dbeta}
		&=
		\dbeta,
	\end{aligned}
	\right.
\\
	&
	\dalpha:
	\left\{
	\begin{aligned}
		\lo{\dalpha}
		&=
		0,
	\\
		\hi{\dalpha}
		&=
		\alpha.
	\end{aligned}
	\right.
\end{align}
\end{subequations}

\end{enumerate}

\subsubsection{Decomposition of geometric quantity associated with general foliation}\label{sec 5.2.2}
We give the decomposition of the geometric quantities associated with a general foliation $\{\bSigma_u\}$ by the formulae obtained in section \ref{sec 4.3}.
\begin{enumerate}[label=\alph*.]

\item
For the preliminary quantity $\vec{\dvarepsilon}'$:
\begin{align}
\begin{aligned}
	\vec{\dvarepsilon}':
	\left\{
	\begin{aligned}
		\lo{\dvarepsilon'^i}
		&=
		-2 (\slashg_{S}^{-1}|_{\bSigma_u})^{ij} (\flu)_j,
	\\
		\hi{\dvarepsilon'^i}
		&=
		-2 \hi{(\dslashg^{-1})^{ik}} (\dB^{-1})_k^j (\flu)_j
		- 2 (\slashg_{S}^{-1}|_{\bSigma_u})^{ik} (\dB^{-1})_k^l (\flu)_l \db^j (\flu)_j. 
	\end{aligned}
	\right.
\end{aligned}
\end{align}

\item
For the metric $\ddslashg$ and the volume form $\ddslashepsilon$ on a general foliation $\{\bSigma_u\}$ of $\ucalH$:
\begin{subequations}
\begin{align}
	&
	\begin{aligned}
		&
		\ddslashg:
		\left\{
		\begin{aligned}
			\lo{\ddslashg}
			&=
			\lo{\dslashg}
			=
			(r_S|_{\bSigma_u})^2 \circg,
		\\
			\hi{\ddslashg}
			&=
			\hi{\dslashg}
			- (\slashg_{ik} \db^k (\flu)_j + \slashg_{jl} \db^l (\flu)_i)
			+ (\flu)_i (\flu)_j \slashg_{kl} \db^k \db^l,
		\end{aligned}
		\right.
	\end{aligned}
\\
	&
	\begin{aligned}
		&
		\ddslashg^{-1}:
		\left\{
		\begin{aligned}
			\lo{\ddslashg^{-1}}
			&=
			\lo{\dslashg^{-1}}
			=
			(r_S|_{\bSigma_u})^{-2}\circg^{-1},
		\\
			\hi{(\ddslashg^{-1})^{ij}}
			&=
			-(r_S|_{\bSigma_u})^{-2} (\circg^{-1})^{ik} \cdot \hi{\ddslashg_{kl}} \cdot (\ddslashg^{-1})^{lj},
		\end{aligned}
		\right.
	\end{aligned}
\\
	&
	\begin{aligned}
		\ddslashepsilon:
		\left\{
		\begin{aligned}
			\lo{\ddslashepsilon}
			&=
			\lo{\dslashepsilon}
			=
			(r_S|_{\bSigma_u})^2 \circepsilon,
		\\
			\hi{\ddslashepsilon}
			&=
			\frac{
			(\slashg_S^{-1})^{kl} \hi{\dslashg_{kl}} 
			+ \det \{ (\slashg_S^{-1})^{kr} \hi{\dslashg_{rl}} \}
			}{
			\sqrt{ \det \{ (\slashg_S^{-1})^{kr} \dslashg_{rl} \} }
			+1
			}
			(r_S|_{\bSigma_u})^2 \circepsilon.
		\end{aligned}
		\right.
	\end{aligned}
\end{align}
\end{subequations}

\item
For the connection coefficients on $\ucalH$ relative to $\{\bSigma_u\}$ and $\{\dduL, \ddL'\}$:
\begin{subequations}
\begin{align}	
	&
	\begin{aligned}
		\ddchi':
		&
		\left\{
		\begin{aligned}
			\lo{\ddchi'}
			&=
			\lo{\dchi'}
			- 2 \circnabla^2 \flu
			=
			\chi'_S|_{\bSigma_u} 
			- 2 \circnabla^2 \flu,
		\\	
			\hi{\ddchi'_{ij}}
			&=
			\hi{\dchi'_{ij}} 
			+ \dvarepsilon' \duchi_{ij} 
			+ \db \circdot \vec{\dvarepsilon}' \dslashnabla^2_{ij} \flu
			- 2 (\dslashnabla^2_{ij} \flu - \circnabla^2_{ij} \flu)
		\\	
			&\phantom{=}
			+2\sym  
			\{  
				\ddslashd \flu 
				\otimes 
				[ 
					\dslashnabla \db \circdot \vec{\dvarepsilon}' 
					- \duchi (\vec{\dvarepsilon}'  ) 
					-\dvarepsilon' \duchi (\db  ) 
					- \dchi' (\db  ) -2 \deta 
				]
			\}_{ij}
		\\	
			&\phantom{=}
			+ 
			[ 
				2\duchi (\db,\vec{\dvarepsilon}'  ) 
				+\dvarepsilon' \duchi (\db,\db  ) 
				+\dchi' (\db,\db  ) 
				+4\deta (\db  )
		\\	
			&\phantom{=
			+ 
			[ 
				2\duchi (\db,\vec{\dvarepsilon}'  ) 
				+\dvarepsilon' \duchi (\db,\db  ) }
				-\dslashnabla_{\db} \db \circdot \vec{\dvarepsilon}' 
				-\dpartial_s \db \circdot \vec{\dvarepsilon}' 
				-4\duomega
			]    (\flu)_i    (\flu)_j,
		\end{aligned}
		\right.
		\label{eqn 5.12a}
	\\	
		\ddtr \ddchi':
		&
		\left\{
		\begin{aligned}
			\lo{\ddtr \ddchi'}
			&=
			\lo{\ddtr \duchi}
			=
			\tr \chi'_S|_{\bSigma_u} - 2 (r_S|_{\bSigma_u})^{-2} \circDelta \flu,
		\\	
			\hi{\ddtr \ddchi'}
			&=
			(\ddslashg^{-1})^{ij} \hi{\ddchi'_{ij}}
			+ \hi{(\ddslashg^{-1})^{ij}} \cdot \ddchi'_{ij},
		\end{aligned}
		\right.
	\\
		\hatddchi':
		&
		\left\{
		\begin{aligned}
			&
			\lo{\hatddchi'}
			=
			-2[\circnabla^2 \flu- \frac{1}{2} (\circDelta \flu) \circg],
		\\
			&
			\hi{\hatdchi'}
			=
			- \frac{1}{2} 
			(\hi{(\ddslashg^{-1})^{kl}} \lo{\ddchi'_{kl}} \ddslashg
			+ \lo{\ddtr \ddchi'} \hi{\ddslashg})
			+\widehat{\hi{\ddchi'_{ij}}},
		\end{aligned}
		\right.
	\end{aligned}
\\
	&
	\begin{aligned}
		\dduchi:
		&
		\left\{
		\begin{aligned}
			\lo{\dduchi}
			&=
			\lo{\duchi}
			=
			\uchi_S|_{\bSigma_u} - 2\Omega_S^2 \circnabla^2 \uh,
		\\	
			\hi{\dduchi}
			&=
			\hi{\duchi}
			+ 2 \sym  \{   (-\duchi (\db  )  )  \otimes \ddslashd \flu   \}_{ij} 
			+ \duchi (\db,\db  )    (\flu)_i    (\flu)_j,
		\end{aligned}
		\right.
	\\	
		\ddtr \dduchi:
		&
		\left\{
		\begin{aligned}
			\lo{\ddtr \dduchi}
			&=
			\lo{\ddtr \duchi}
			=
			\tr \uchi_S|_{\bSigma_u} - 2\Omega_S^2 (r_S|_{\bSigma_u})^{-2}\circDelta \uh,
		\\	
			\hi{\ddtr \dduchi}
			&=
			(\ddslashg^{-1})^{ij} \hi{\dduchi_{ij}}
			+ \hi{(\ddslashg^{-1})^{ij}} \cdot \dduchi_{ij},
		\end{aligned}
		\right.
	\\
		\hatdduchi:
		&
		\left\{
		\begin{aligned}
			&
			\lo{\hatdduchi}
			=
			-2\Omega_S^2 [\circnabla^2 \uh- \frac{1}{2} (\circDelta \uh) \circg],
		\\
			&
			\hi{\hatdduchi}
			=
			- \frac{1}{2} 
			(\hi{(\ddslashg^{-1})^{kl}} \lo{\dduchi_{kl}} \ddslashg
			+ \lo{\ddtr \dduchi} \hi{\ddslashg})
			+\widehat{\hi{\dduchi_{ij}}},
		\end{aligned}
		\right.
	\end{aligned}
\\
	\ddeta:
	&
	\left\{
	\begin{aligned}
		\lo{\ddeta}
		&=
		\lo{\deta} 
		+ \frac{1}{2} \lo{\duchi} ( \lo{\vec{\dvarepsilon}'})
		\\
		&=
		(\chi_S|_{\bSigma_u})_{ij} \cdot (\slashg_S^{-1})^{jk} \cdot \uh_k
		- [(\uchi_S|_{\bSigma_u})_{ij} - 2\Omega_S^2 \circnabla^2_{ij} \uh] \cdot (\slashg_S^{-1})^{jk} \cdot (\flu)_k,
	\\
		\hi{\ddeta}
		&=
		\hi{\deta}
		+ \frac{1}{2} [\duchi(\hi{\vec{\dvarepsilon}'}) + \hi{\duchi} ( \lo{\vec{\dvarepsilon}'})]
		+ 
		[ 
			2\duomega
			-\deta (\db  ) 
			-\frac{1}{2} \duchi (\db,\vec{\dvarepsilon}' )
		]  \ddslashd \flu.
	\end{aligned}
	\right.
\end{align}
\end{subequations}

\item
For the curvature components on $\ucalH$ relative to $\{\bSigma_u\}$ and $\{\dduL, \ddL' \}$:
\begin{subequations}
\begin{align}
	\ddualpha:
	&
	\left\{
	\begin{aligned}
		\lo{\ddualpha}
		&=
		\lo{\dualpha}
		=0,
	\\
		\hi{\ddualpha}
		&=
		\ddualpha,
	\end{aligned}
	\right.
\\
	\ddubeta:
	&
	\left\{
	\begin{aligned}
		\lo{\ddubeta}
		&=
		\lo{\ddubeta}
		=
		-3\Omega_S^2 \rho_S \dslashd \uh,
	\\
		\hi{\ddubeta}
		&=
		-(\flu)_i \db^j \dubeta_j
		+ \dB_i^j \hi{\dubeta_j}
		-\frac{1}{2} \dB_i^j \dvarepsilon^k \dualpha_{jk},
	\end{aligned}
	\right.
\\
	\ddrho:
	&
	\left\{
	\begin{aligned}
		\lo{\ddrho}
		&=
		\lo{\drho}
		=
		\rho_S|_{\bSigma_u},
	\\
		\hi{\ddrho}
		&=
		\hi{\drho}
		- \dvarepsilon^j \dubeta_j
		+ \frac{1}{4} \dvarepsilon^i \dvarepsilon^j \dualpha_{ij},
	\end{aligned}
	\right.
\\
	\ddsigma:
	&
	\left\{
	\begin{aligned}
		\lo{\ddsigma}
		&=
		0,
	\\
		\hi{\ddsigma}
		&=
		\ddsigma,
	\end{aligned}
	\right.
\\
	\ddbeta:
	&
	\left\{
	\begin{aligned}
		\lo{\ddbeta}
		&=
		-3\lo{\drho} \ddslashd \flu
		=
		-3 \rho_S|_{\bSigma_u} \ddslashd \flu,
	\\
		\hi{\ddbeta_i}
		&=
		\dB_i^j\dbeta_j 
		- 3 \hi{\drho} f_i
		+ \frac{3}{2}\dB_i^j \dvarepsilon^k \dslashepsilon_{jk}  \dsigma
		- \dB_i^j \dvarepsilon \dubeta_j 
		+ 2 f_i \dvarepsilon^l \dubeta_l 
	\\
		&\phantom{=}
		+ \frac{1}{2} \dB_i^j \dvarepsilon^k \dvarepsilon^l \dslashepsilon_{kj} \dslashepsilon_{lm} \dubeta^m
		+ \frac{1}{2} \dB_i^j \dvarepsilon^l  \dvarepsilon \dualpha_{jl}
		-\frac{1}{2} f_i \dvarepsilon^k \dvarepsilon^l \dualpha_{kl},
	\end{aligned}
	\right.
\\
	\ddalpha:
	&
	\left\{
	\begin{aligned}
		\lo{\ddalpha}
		&=
		0,
	\\
		\hi{\ddalpha}
		&=
		\ddalpha.
	\end{aligned}
	\right.
\end{align}
\end{subequations}

\end{enumerate}

\subsection{Estimate of first order main part of geometric quantity}\label{sec 5.3}
We estimate the first order main part of the geometric quantities in the decompositions constructed in section \ref{sec 5.2}. Let $\{\bSigma_u\}$ be a general foliation of an incoming null hypersurface $\ucalH$. We shall obtain the estimates on each leaf $\bSigma_u$ of the foliation. For the sake of simplification of the statement of the estimates, we recall the notations $\ud_o$, $\ud_m$ in formulae \eqref{eqn 3.5}
\begin{align}
	\ud_o = c_o \udelta_o, 
	\quad 
	\ud_m = [ 1+ c_{m,m} \epsilon \udelta_o ] \udelta_m + c_{m,o} \udelta_o^2,
	\tag{\ref{eqn 3.5}\ensuremath{'}}
\end{align}
and introduce the similar notations $\uslashd_m$, $\ud_{\uh}$ by
\begin{align}
	\uslashd_m 
	= 
	c_{m,m} \epsilon \udelta_o \udelta_m 
	+ c_{m,o} \udelta_o^2,
	\quad
	\ud_{\uh} 
	=
	\max\{c_{\uh}, c_{\uh,2} \} \udelta_o.
	\label{eqn 5.14}
\end{align}
Let $\uh$ be the parameterisation of $\ucalH$, $(\ufl{s=0}, \flu)$ be the parametrisation of $\bSigma_u$ and $(\uflu, \flu)$ be the parametrisation of $\bSigma_u$ in the double null coordinate system. Then the conclusions of propositions \ref{prop 3.3}, \ref{prop 5.2}, \ref{prop 5.5} can be written as
\begin{align}
\begin{aligned}
	&
	\Vert \ddslashd \uflu \Vert^{n,p} 
	\leq 
	\ud_o r_0,
	\quad
	\vert \overline{\uflu}^{\circg} - \overline{\ufl{s=0}}^{\circg} \vert
	\leq 
	\uslashd_m r_0,
	\quad
	\vert \overline{\uflu}^{\circg} \vert \leq \ud_m r_0,
\\
	&
	\Vert (\dslashd \uh)|_{\bSigma_u} \Vert^{n+1,p},
	\Vert (\circnabla^2 \uh)|_{\bSigma_u} \Vert^{n,p}
	\leq
	\ud_{\uh} r_0.
\end{aligned}
\label{eqn 5.15}
\end{align}

\subsubsection{Estimate of first order main part associated with $s$ foliation}
We state the estimates of the first order main parts of geometric quantities associated with the $s$ foliation $\{\Sigma_s\}$ on each leaf $\bSigma_u$ in the following proposition.
\begin{proposition}\label{prop 5.7}
Suppose that the parameterisation functions $\ufl{s=0}$, $\flu$ of the leaf $\bSigma_u$ satisfy the following
\begin{align*}
	\Vert \slashd \ufl{s=0} \Vert^{n+1,p} \leq \udelta_o r_0,
	\quad
	\vert \overline{\ufl{s=0}}^{\circg} \vert \leq \udelta_m r_0,
	\quad
	\Vert \slashd \flu \Vert^{n+1,p} \leq \delta_o (r_0+\os_u),
	\quad
	\overline{\flu}^{\circg} = \os_u > -\frac{r_0}{2},
\end{align*}
where $n\geq 1+n_p$, $p >1$. Then there exists a small positive constant $\delta$ depending on $n,p$, such that if $\epsilon, \udelta_o, \udelta_m, \delta_o$ suitably bounded that $\epsilon, \udelta_o, \epsilon \udelta_m, \delta_o \leq \delta$, then the first order main parts of geometric quantities associated with the $s$ foliation $\{ \Sigma_s \}$ satisfies the estimates, which we collect as follows.
\begin{enumerate}[label=\alph*.]
\item 
For the shifting vector $\db$,the metric $\dslashg$ and its volume form $\dslashepsilon$ of $\{\Sigma_s\}$ on $\bSigma_u$:
\begin{align*}
	\begin{aligned}
		&
		\Vert \lo{\db}|_{\bSigma_u} \Vert^{n+1,p}
		\leq
		\frac{c(n,p) \ud_{\uh} r_0}{(r_0+\os_u)^2},
	\\
		&
		\Vert \lo{\dslashg}|_{\bSigma_u} \Vert^{n+1,p}
		\leq
		c(n,p) (r_0+\os_u)^2,
	\\
		&
		\Vert \lo{\dslashg^{-1}}|_{\bSigma_u} \Vert^{n+1,p}
		\leq
		\frac{c(n,p)}{(r_0+\os_u)^2},
	\\
		&
		\Vert \lo{\dslashepsilon}|_{\bSigma_u} \Vert^{n+1,p}
		\leq
		c(n,p) (r_0+\os_u)^2.
	\end{aligned}
\end{align*}

\item
For the connection coefficients on $\ucalH$ relative to $\{\Sigma_s\}$ and $\{\duL, \dL'\}$ on $\bSigma_u$:
\begin{align*}
	\begin{aligned}
		&
		\Vert \lo{\dchi'}|_{\bSigma_u} \Vert^{n+2,p}
		\leq
		c(n,p)  (\os_u + \delta_o r_0),
	\\
		&
		\left\{
		\begin{aligned}
			&
			\Vert \lo{\dtr \dchi'}|_{\bSigma_u} \Vert^{n+2,p}
			\leq
			\frac{c(n,p) (\os_u + \delta_o r_0)}{ (r_0 + \os_u)^2},
		\\
			&
			\Vert \ddslashd (\lo{\dtr \dchi'}|_{\bSigma_u}) \Vert^{n+1,p}
			\leq
			\frac{c(n,p) \delta_o}{r_0 + \os_u}
			+ \frac{c(n,p) (\os_u + \delta_o r_0)\ud_o r_0}{ (r_0 + \os_u)^3},
		\end{aligned}
		\right.
	\\
		&
		\Vert \lo{\hatdchi'}|_{\bSigma_u} \Vert^{n+2,p}
		=
		0,
	\\
		&
		\Vert \lo{\duchi}|_{\bSigma_u} \Vert^{n,p}
		\leq
		c(n,p) (r_0+\os_u),
	\\
		&
		\left\{
		\begin{aligned}
			&
			\Vert \lo{\dtr \duchi}|_{\bSigma_u} \Vert^{n,p}
			\leq
			\frac{c(n,p)}{r+\os_u},
		\\
			&
			\Vert \ddslashd \lo{\dtr \duchi}|_{\bSigma_u} \Vert^{n-1,p}
			\leq
			\frac{c(n,p) \delta_o}{r+\os_u}
			+ \frac{c(n,p)(\ud_o + \ud_h) r_0}{(r+\os_u)^2},
		\end{aligned}
		\right.
	\\
		&
		\Vert \lo{\hatduchi}|_{\bSigma_u} \Vert^{n,p}
		\leq
		c(n,p) \ud_{\uh} r_0,
	\\
		&
		\Vert \lo{\deta}|_{\bSigma_u} \Vert^{n+1,p}
		\leq
		\frac{c(n,p) \ud_{\uh} (\os_u + \delta_o r_0) r_0}{(r+\os_u)^2},
	\\
		&
		\left\{
		\begin{aligned}
			&
			\Vert \lo{\duomega}|_{\bSigma_u} \Vert^{n+1,p}
			\leq
			\frac{c(n,p) (\ud_m + \ud_o) r_0^2}{(r+\os_u)^3},
		\\
			&
			\Vert \ddslashd \lo{\duomega}|_{\bSigma_u} \Vert^{n,p}
			\leq
			\frac{c(n,p) [(\ud_m + \ud_o) \delta_o + \ud_o] r_0^2}{(r+\os_u)^3}.
		\end{aligned}
		\right.
	\end{aligned}
\end{align*}
	
\item
For the curvature components on $\ucalH$ relative to $\{\Sigma_s\}$ and $\{ \duL, \dL'\}$ on $\bSigma_u$:
\begin{align*}
		&
		\Vert \lo{\dualpha}|_{\bSigma_u} \Vert^{n+1,p}
		=
		0,
	\\
		&
		\Vert \lo{\dubeta}|_{\bSigma_u} \Vert^{n+1,p}
		\leq
		\frac{c(n,p) \ud_{\uh} r_0^2}{(r_0 + \os_u)^3},
	\\
		&
		\Vert \lo{\dsigma}|_{\bSigma_u} \Vert^{n+1,p}
		=
		0,
	\\
		&
		\Vert \lo{\drho}|_{\bSigma_u} \Vert^{n+1,p}
		\leq
		\frac{c(n,p) r_0}{(r_0 + \os_u)^3},
	\\
		&
		\Vert \lo{\dbeta}|_{\bSigma_u} \Vert^{n+1,p}
		=
		0,
	\\
		&
		\Vert \hi{\dalpha}|_{\bSigma_u} \Vert^{n+1,p}
		=
		0.
\end{align*}
\end{enumerate}
\end{proposition}
\begin{proof}
Choose $\delta$ sufficiently small such that propositions \ref{prop 3.3}, \ref{prop 5.2}, \ref{prop 5.5} are true. Their conclusions, the estimates of $\uflu$ and $(\dslashd \uh)|_{\bSigma_u}$, $(\circnabla^2 \uh)|_{\bSigma_u}$ listed in \eqref{eqn 5.15}, hold. By the choice of $\delta$ in the proofs of propositions \ref{prop 3.3}, \ref{prop 5.2}, \ref{prop 5.5}, we have that
\begin{align*}
\ud_o \leq 1,
\quad
\uslashd_m \leq 1,
\quad
\ud_m \leq 2,
\quad
\ud_{\uh} \leq 1.
\end{align*}
Then the proposition simply follows from direct checks by using the assumptions of the vacuum perturbed Schwarzschild metric in definition \ref{def 2.4}, the decomposition formulae in section \ref{sec 5.2.1}, and the estimates of $\uflu$ and $(\dslashd \uh)|_{\bSigma_u}$, $(\circnabla^2 \uh)|_{\bSigma_u}$ in \eqref{eqn 5.15} from propositions \ref{prop 3.3}, \ref{prop 5.2}, \ref{prop 5.5}.
\end{proof}

\subsubsection{Estimate of first order main part associated with general foliation}

We state the estimates of the first order main parts of geometric quantities associated with a general foliation $\{\bSigma_u\}$ in the following proposition.
\begin{proposition}\label{prop 5.8}
Suppose that the parameterisation functions $\ufl{s=0}$, $\flu$ of the leaf $\bSigma_u$ satisfy the following
\begin{align*}
	\Vert \slashd \ufl{s=0} \Vert^{n+1,p} \leq \udelta_o r_0,
	\quad
	\vert \overline{\ufl{s=0}}^{\circg} \vert \leq \udelta_m r_0,
	\quad
	\Vert \slashd \flu \Vert^{n+1,p} \leq \delta_o (r_0+\os_u),
	\quad
	\overline{\flu}^{\circg}= \os_u > -\frac{r_0}{2},
\end{align*}
where $n\geq 1+n_p$, $p >1$. Then there exists a small positive constant $\delta$ depending on $n,p$, such that if $\epsilon, \udelta_o, \udelta_m, \delta_o$ suitably bounded that $\epsilon, \udelta_o, \epsilon \udelta_m, \delta_o \leq \delta$, then the first order main parts of geometric quantities associated with the foliation $\{ \bSigma_u \}$ satisfy the estimates, which we collect as follows.
\begin{enumerate}[label=\alph*.]
\item
For the preliminary quantity $\vec{\dvarepsilon}'$ on $\bSigma_u$:
\begin{align*}
	\Vert \lo{\vec{\dvarepsilon}'} \Vert^{n+1,p}
	\leq
	\frac{c(n,p) \delta_o}{r_0+\os_u}.
\end{align*}

\item 
For the metric $\ddslashg$ and its volume form $\ddslashepsilon$ of $\{\bSigma_u\}$:
\begin{align*}
	\begin{aligned}
		&
		\Vert \lo{\ddslashg}|_{\bSigma_u} \Vert^{n+1,p}
		\leq
		c(n,p) (r_0+\os_u)^2,
	\\
		&
		\Vert \lo{\ddslashg^{-1}}|_{\bSigma_u} \Vert^{n+1,p}
		\leq
		\frac{c(n,p) }{(r_0+\os_u)^2},
	\\
		&
		\Vert \lo{\ddslashepsilon}|_{\bSigma_u} \Vert^{n+1,p}
		\leq
		c(n,p) (r_0+\os_u)^2.
	\end{aligned}
\end{align*}

\item
For the connection coefficients on $\ucalH$ relative to $\{\Sigma_s\}$ and $\{\duL, \dL'\}$ on $\bSigma_u$:
\begin{align*}
	\begin{aligned}
		&
		\Vert \lo{\ddchi'}|_{\bSigma_u} \Vert^{n,p}
		\leq
		c(n,p) (\os_u + \delta_o r_o),
	\\
		&
		\Vert \lo{\ddtr \ddchi'}|_{\bSigma_u} \Vert^{n,p}
		\leq
		\frac{c(n,p)(\os_u + \delta_o r_o)}{(r_0 + \os_u)^2},
	\\
		&
		\Vert \lo{\hatddchi'}|_{\bSigma_u} \Vert^{n,p}
		\leq
		c(n,p) \delta_o (r_0 + \os_u),
	\\
		&
		\Vert \lo{\dduchi}|_{\bSigma_u} \Vert^{n,p}
		\leq
		c(n,p) (r_0+\os_u),
	\\
		&
		\Vert \lo{\ddtr \dduchi}|_{\bSigma_u} \Vert^{n,p}
		\leq
		\frac{c(n,p)}{r+\os_u},
	\\
		&
		\Vert \lo{\hatdduchi}|_{\bSigma_u} \Vert^{n,p}
		\leq
		c(n,p) \ud_{\uh} r_0,
	\\
		&
		\Vert \lo{\ddeta}|_{\bSigma_u} \Vert^{n,p}
		\leq
		c(n,p) \delta_o
		+ \frac{c(n,p) \ud_{\uh} \os_u r_0}{(r+\os_u)^2}.
	\end{aligned}
\end{align*}
	
\item
For the curvature components on $\ucalH$ relative to $\{\Sigma_s\}$ and $\{ \duL, \dL'\}$ on $\bSigma_u$:
\begin{align*}
		&
		\Vert \lo{\ddualpha}|_{\bSigma_u} \Vert^{n+1,p}
		=
		0,
	\\
		&
		\Vert \lo{\ddubeta}|_{\bSigma_u} \Vert^{n+1,p}
		\leq
		\frac{c(n,p) \ud_{\uh} r_0^2}{(r_0 + \os_u)^3},
	\\
		&
		\Vert \lo{\ddsigma}|_{\bSigma_u} \Vert^{n+1,p}
		=
		0,
	\\
		&
		\Vert \lo{\ddrho}|_{\bSigma_u} \Vert^{n+1,p}
		\leq
		\frac{c(n,p) r_0}{(r_0 + \os_u)^3},
	\\
		&
		\Vert \lo{\ddbeta}|_{\bSigma_u} \Vert^{n+1,p}
		\leq
		\frac{c(n,p) \delta_o r_0 }{(r_0 + \os_u)^2},
	\\
		&
		\Vert \lo{\ddalpha}|_{\bSigma_u} \Vert^{n+1,p}
		=
		0.
\end{align*}
\end{enumerate}
\end{proposition}
\begin{proof}
The proof follows the same route as the proof of proposition \ref{prop 5.7}.
\end{proof}

\subsection{Estimate of high order remainder of geometric quantity}\label{sec 5.4}
We estimate the high order remainders of the geometric quantities in the decompositions constructed in section \ref{sec 5.2}. Follow the notations in section \ref{sec 5.3}.

\subsubsection{Estimate of high order remainder associated with $s$ foliation}
We state the estimates of the high order remainders of geometric quantities associated with the $s$ foliation $\{\Sigma_s\}$ on each leaf $\bSigma_u$ in the following proposition.
\begin{proposition}\label{prop 5.9}
Suppose that the parameterisation functions $\ufl{s=0}$, $\flu$ of the leaf $\bSigma_u$ satisfy the following
\begin{align*}
	\Vert \slashd \ufl{s=0} \Vert^{n+1,p} \leq \udelta_o r_0,
	\quad
	\vert \overline{\ufl{s=0}}^{\circg} \vert \leq \udelta_m r_0,
	\quad
	\Vert \slashd \flu \Vert^{n+1,p} \leq \delta_o (r_0+\os_u),
	\quad
	\overline{\flu}^{\circg} = \os_u > -\frac{r_0}{2},
\end{align*}
where $n\geq 1+n_p$, $p >1$. Then there exists a small positive constant $\delta$ depending on $n,p$, such that if $\epsilon, \udelta_o, \udelta_m, \delta_o$ suitably bounded that $\epsilon, \udelta_o, \epsilon \udelta_m, \delta_o \leq \delta$, then the high order remainders of geometric quantities associated with the $s$ foliation $\{ \Sigma_s \}$ satisfied the estimates, which we collect as follows.
\begin{enumerate}[label=\alph*.]
\item For the preliminary quantity $\circtriangle_{ij}^k$ on $\bSigma_u$:
\begin{align*}
	\Vert \circDelta \Vert^{n+1,p}
	\leq
	c(n,p) \epsilon.
\end{align*}

\item 
For the shifting vector $\db$,the metric $\dslashg$ and its volume form $\dslashepsilon$ of $\{\Sigma_s\}$ on $\bSigma_u$:
\begin{align*}
	\begin{aligned}
		&
		\Vert \hi{\db}|_{\bSigma_u} \Vert^{n+1,p}
		\leq
		\frac{c(n,p) \epsilon (\ud_m + \ud_o)r_0^2}{(r_0+\os_u)^3}
		+ \frac{c(n,p) \epsilon \ud_{\uh} r_0}{(r_0+\os_u)^2},
	\\
		&
		\Vert \hi{\dslashg}|_{\bSigma_u} \Vert^{n+1,p}
		\leq
		c(n,p) \epsilon (r_0+\os_u)^2,
	\\
		&
		\Vert \hi{\dslashg^{-1}}|_{\bSigma_u} \Vert^{n+1,p}
		\leq
		\frac{c(n,p) \epsilon}{(r_0+\os_u)^2},
	\\
		&
		\Vert \hi{\dslashepsilon}|_{\bSigma_u} \Vert^{n+1,p}
		\leq
		c(n,p) \epsilon (r_0+\os_u)^2.
	\end{aligned}
\end{align*}

\item
For the connection coefficients on $\ucalH$ relative to $\{\Sigma_s\}$ and $\{\duL, \dL'\}$ on $\bSigma_u$:
\begin{align*}
	\begin{aligned}
		&
		\Vert \hi{\dchi'}|_{\bSigma_u} \Vert^{n+1,p}
		\leq
		c(n,p) \epsilon (r_0 + \os_u),
	\\
		&
		\Vert \hi{\dtr \dchi'}|_{\bSigma_u} \Vert^{n+1,p}
		\leq
		\frac{c(n,p) \epsilon}{ r_0 + \os_u},
	\\
		&
		\Vert \hi{\hatdchi'}|_{\bSigma_u} \Vert^{n+1,p}
		\leq
		c(n,p) \epsilon (r_0 + \os_u),
	\\
		&
		\Vert \hi{\duchi}|_{\bSigma_u} \Vert^{n,p}
		\leq
		c(n,p) \epsilon r_0 + \frac{c(n,p) \ud_{\uh}^2 r_0^2}{r_0+\os_u},
	\\
		&
		\Vert \hi{\dtr \duchi}|_{\bSigma_u} \Vert^{n,p}
		\leq
		\frac{c(n,p) \epsilon r_0}{(r+\os_u)^2} + \frac{c(n,p) \ud_{\uh}^2 r_0^2}{(r_0+\os_u)^3},
	\\
		&
		\Vert \hi{\hatduchi}|_{\bSigma_u} \Vert^{n,p}
		\leq
		c(n,p) \epsilon r_0 + \frac{c(n,p) \ud_{\uh}^2 r_0^2}{r_0+\os_u},
	\\
		&
		\Vert \hi{\deta}|_{\bSigma_u} \Vert^{n+1,p}
		\leq
		\frac{c(n,p) \epsilon r_0}{r+\os_u},
	\\
		&
		\Vert \hi{\duomega}|_{\bSigma_u} \Vert^{n+1,p}
		\leq
		\frac{c(n,p) (\epsilon + \ud_{\uh}^2) r_0^2}{(r+\os_u)^3},
	\end{aligned}
\end{align*}
	
\item
For the curvature components on $\ucalH$ relative to $\{\Sigma_s\}$ and $\{ \duL, \dL'\}$ on $\bSigma_u$:
\begin{align*}
		&
		\Vert \hi{\dualpha}|_{\bSigma_u} \Vert^{n,p}
		\leq
		\frac{c(n,p) \epsilon r_0^{\frac{3}{2}}}{(r_0 + \os_u)^{\frac{3}{2}}}
		+
		\frac{c(n,p)\ud_{\uh}^2 r_0^3}{(r_0 + \os_u)^3},
	\\
		&
		\Vert \hi{\dubeta}|_{\bSigma_u} \Vert^{n,p}
		\leq
		\frac{c(n,p) \epsilon r_0^{\frac{3}{2}}}{(r_0 + \os_u)^{\frac{5}{2}}},
	\\
		&
		\Vert \hi{\dsigma}|_{\bSigma_u} \Vert^{n,p}
		\leq
		\frac{c(n,p) \epsilon r_0}{(r_0 + \os_u)^3},
	\\
		&
		\Vert \hi{\drho}|_{\bSigma_u} \Vert^{n,p}
		\leq
		\frac{c(n,p) \epsilon r_0}{(r_0 + \os_u)^3},
	\\
		&
		\Vert \hi{\dbeta}|_{\bSigma_u} \Vert^{n,p}
		\leq
		\frac{c(n,p) \epsilon}{r_0 + \os_u},
	\\
		&
		\Vert \hi{\dalpha}|_{\bSigma_u} \Vert^{n,p}
		\leq
		\frac{c(n,p) \epsilon (r_0 + \os_u)}{r_0}.
\end{align*}
\end{enumerate}
\end{proposition}
\begin{proof}
The proof follows the same route as the proof of proposition \ref{prop 5.7}.
\end{proof}

\subsubsection{Estimate of high order remainder associated with general foliation}

We state the estimates of the high order remainders of geometric quantities associated with a general foliation $\{\bSigma_u\}$ in the following proposition.
\begin{proposition}\label{prop 5.10}
Suppose that the parameterisation functions $\ufl{s=0}$, $\flu$ of the leaf $\bSigma_u$ satisfy the following
\begin{align*}
	\Vert \slashd \ufl{s=0} \Vert^{n+1,p} \leq \udelta_o r_0,
	\quad
	\vert \overline{\ufl{s=0}}^{\circg} \vert \leq \udelta_m r_0,
	\quad
	\Vert \slashd \flu \Vert^{n+1,p} \leq \delta_o (r_0+\os_u),
	\quad
	\overline{\flu}^{\circg}= \os_u > -\frac{r_0}{2},
\end{align*}
where $n\geq 1+n_p$, $p >1$. Then there exists a small positive constant $\delta$ depending on $n,p$, such that if $\epsilon, \udelta_o, \udelta_m, \delta_o$ suitably bounded that $\epsilon, \udelta_o, \epsilon \udelta_m, \delta_o \leq \delta$, then the high order remainders of geometric quantities associated with the foliation $\{ \bSigma_u \}$ satisfied the estimates, which we collect as follows.
\begin{enumerate}[label=\alph*.]
\item
For the preliminary quantity $\vec{\dvarepsilon}'$ and $\triangle_{ij}^k$ on $\bSigma_u$:
\begin{align*}
	\Vert \hi{\vec{\dvarepsilon}'} \Vert^{n+1,p}
	\leq
	\frac{c(n,p)\epsilon \delta_o}{r_0+\os_u},
\\
	\Vert \triangle_{ij}^k \Vert^{n+1,p}
	\leq
	\frac{c(n,p) \ud_{\uh} r_0}{r_0+\os_u}.
\end{align*}

\item 
For the metric $\ddslashg$ and its volume form $\ddslashepsilon$ of $\{\bSigma_u\}$:
\begin{align*}
	\begin{aligned}
		&
		\Vert \hi{\ddslashg}|_{\bSigma_u} \Vert^{n+1,p}
		\leq
		c(n,p) \epsilon (r_0+\os_u)^2,
	\\
		&
		\Vert \hi{\ddslashg^{-1}}|_{\bSigma_u} \Vert^{n+1,p}
		\leq
		\frac{c(n,p) \epsilon}{(r_0+\os_u)^2},
	\\
		&
		\Vert \hi{\ddslashepsilon}|_{\bSigma_u} \Vert^{n+1,p}
		\leq
		c(n,p) \epsilon (r_0+\os_u)^2.
	\end{aligned}
\end{align*}

\item
For the connection coefficients on $\ucalH$ relative to $\{\Sigma_s\}$ and $\{\duL, \dL'\}$ on $\bSigma_u$:
\begin{align*}
	\begin{aligned}
		&
		\Vert \hi{\ddchi'}|_{\bSigma_u} \Vert^{n,p}
		\leq
		c(n,p) \epsilon (r_0 + \os_u)
		+ c(n,p) \delta_o^2 r_0
		+ c(n,p) \ud_{\uh} \delta_o r_0,
	\\
		&
		\Vert \hi{\ddtr \ddchi'}|_{\bSigma_u} \Vert^{n,p}
		\leq
		\frac{c(n,p) \epsilon}{r_0 + \os_u}
		+ \frac{c(n,p) \delta_o^2 r_0}{ (r_0+ \os_u)^2}
		+ \frac{c(n,p) \ud_{\uh} \delta_o r_0}{ (r_0+ \os_u)^2},
	\\
		&
		\Vert \hi{\hatddchi'}|_{\bSigma_u} \Vert^{n,p}
		\leq
		c(n,p) \epsilon (r_0 + \os_u)
		+ c(n,p) \delta_o^2 r_0
		+ c(n,p) \ud_{\uh} \delta_o r_0,
	\\
		&
		\Vert \hi{\dduchi}|_{\bSigma_u} \Vert^{n,p}
		\leq
		c(n,p) \epsilon r_0 
		+ \frac{c(n,p) \ud_{\uh}^2 r_0^2}{r_0+\os_u}
		+ \frac{c(n,p) \ud_{\uh} \delta_o r_0^2}{r_0+\os_u},
	\\
		&
		\Vert \hi{\ddtr \dduchi}|_{\bSigma_u} \Vert^{n,p}
		\leq
		\frac{c(n,p) \epsilon r_0}{(r+\os_u)^2} 
		+ \frac{c(n,p) \ud_{\uh}^2 r_0^2}{(r_0+\os_u)^3}
		+ \frac{c(n,p) \ud_{\uh} \delta_o r_0^2}{(r_0+\os_u)^3},
	\\
		&
		\Vert \hi{\hatdduchi}|_{\bSigma_u} \Vert^{n,p}
		\leq
		c(n,p) \epsilon r_0 
		+ \frac{c(n,p) \ud_{\uh}^2 r_0^2}{r_0+\os_u}
		+ \frac{c(n,p) \ud_{\uh} \delta_o r_0^2}{r_0+\os_u},
	\\
		&
		\Vert \hi{\ddeta}|_{\bSigma_u} \Vert^{n,p}
		\leq
		\frac{c(n,p) \epsilon r_0}{r+\os_u}
		+ \frac{c(n,p) \ud_{\uh}^2 \delta_o r_0^2}{(r+\os_u)^2}
		+ \frac{c(n,p) (\ud_m + \ud_o) \delta_o r_0^2}{(r+\os_u)^2}
		+ \frac{c(n,p) \ud_{\uh} \delta_o^2 r_0^2}{(r+\os_u)^2}.
	\end{aligned}
\end{align*}
	
\item
For the curvature components on $\ucalH$ relative to $\{\Sigma_s\}$ and $\{ \duL, \dL'\}$ on $\bSigma_u$:
\begin{align*}
		&
		\Vert \hi{\ddualpha}|_{\bSigma_u} \Vert^{n,p}
		\leq
		\frac{c(n,p) \epsilon r_0^{\frac{3}{2}}}{(r_0 + \os_u)^{\frac{3}{2}}}
		+
		\frac{c(n,p)\ud_{\uh}^2 r_0^3}{(r_0 + \os_u)^3},
	\\
		&
		\Vert \hi{\ddubeta}|_{\bSigma_u} \Vert^{n,p}
		\leq
		\frac{c(n,p) \epsilon r_0^{\frac{3}{2}}}{(r_0 + \os_u)^{\frac{5}{2}}}
		+ \frac{c(n,p) \delta_o \ud_{\uh}^2 r_0^3}{(r_0+\os_u)^4},
	\\
		&
		\Vert \hi{\ddsigma}|_{\bSigma_u} \Vert^{n,p}
		\leq
		\frac{c(n,p) \epsilon r_0}{(r_0 + \os_u)^3}
		+ \frac{c(n,p) \delta_o \ud_{\uh} r_0^2}{(r_0+\os_u)^4},
	\\
		&
		\Vert \hi{\ddrho}|_{\bSigma_u} \Vert^{n,p}
		\leq
		\frac{c(n,p) \epsilon r_0}{(r_0 + \os_u)^3}
		+ \frac{c(n,p) \delta_o \ud_{\uh} r_0^2}{(r_0+\os_u)^4},
	\\
		&
		\Vert \hi{\ddbeta}|_{\bSigma_u} \Vert^{n,p}
		\leq
		\frac{c(n,p) \epsilon }{r_0 + \os_u}
		+ \frac{c(n,p) \delta_o^2 \ud_{\uh} r_0^2}{(r_0+\os_u)^3},
	\\
		&
		\Vert \hi{\ddalpha}|_{\bSigma_u} \Vert^{n,p}
		\leq
		\frac{c(n,p) \epsilon (r_0 + \os_u)}{r_0}
		+ \frac{c(n,p) \delta_o^2 r_0}{(r_0+\os_u)}.
\end{align*}
\end{enumerate}
\end{proposition}
\begin{proof}
The proof follows the same route as the proof of proposition \ref{prop 5.7}.
\end{proof}

\section{Constant mass aspect function foliation}\label{sec 6}
In this section, we review the basics of the constant mass aspect function foliation on an incoming null hypersurface, including its definition and the basic equations associated with the foliation.

\subsection{Definition}\label{sec 6.1}
We recall the definition of the mass aspect function of a foliation on an incoming null hypersurface.
\begin{definition}
Let $\ucalH$ be an incoming null hypersurface and $\{ \bSigma_u \}$ be a foliation on $\ucalH$. Let $\{ \buL^u, \bL'^u \}$ be the conjugate null frame associated with $\{\bSigma_u\}$ where $\buL^u u =1$, $g(\buL^u, \bL'^u)=2$. The mass aspect function of the foliation $\{ \bSigma_u \}$ is defined by
\begin{align}
	\bmulu 
	= 
	\bKlu 
	- \frac{1}{4} \btr \bchilu' \ \btr \buchilu 
	- \bslashdiv \btalu
	=
	- \brholu
	- \frac{1}{2} ( \hatbuchilu, \hatbchilu')
	- \bslashdiv \btalu.
	\label{eqn 6.1}
\end{align}
Similarly we define the mass aspect function $\ddmulu$ of $\bSigma_u$ associated with the conjugate null frame $\{ \dduL^u, \ddL'^u \}$ by
\begin{align}
\begin{aligned}
	\ddmulu 
	&= 
	\ddKlu 
	- \frac{1}{4} \ddtr \ddchilu' \ \ddtr \dduchilu 
	- \ddslashdiv \ddetalu
	= 
	\bKlu 
	- \frac{1}{4} \btr \bchilu' \ \btr \buchilu 
	- \bslashdiv \ddetalu
\\
	&=
	- \ddrholu
	- \frac{1}{2} ( \hatdduchilu, \hatddchilu')
	- \bslashdiv \ddetalu
	=
	- \brholu
	- \frac{1}{2} ( \hatbuchilu, \hatbchilu')
	- \bslashdiv \ddetalu.
\end{aligned}
	\label{eqn 6.2}
\end{align}
\end{definition}
The constant mass aspect function foliation on an incoming null hypersurface is defined by the mass aspect function being  constant.
\begin{definition}
A foliation $\{ \bSigma_u \}$ on an incoming null hypersurface is called a constant mass aspect function foliation if its mass aspect function $\bmulu = \overline{\bmulu}^{\subbslashglu}$.
\end{definition}
We note that a relabelling of the parameter $u$ of a constant mass aspect function foliation gives rise to another constant mass aspect function foliation. To remove this freedom of relabelling, we choose the foliation parameterised by the area radius, i.e. imposing the condition $\br_u = \br_{u_0} + u - u_0$ where $|\bSigma_u| = 4\pi \br_u^2$. It is equivalent to $\overline{\btr \buchilu}^{\subbslashglu} = \frac{2}{\br_u}$. In this following, without specifications, all the constant mass aspect function foliations are parameterised by area radius.

\subsection{Basic equations}\label{sec 6.2}
We collect the basic equations of a constant mass aspect function foliation parametrised by area radius, see also \cite{L2018}\cite{L2023}. Let $\{\bSigma_u\}$ be a constant mass aspect function foliation and $\flu$ is the parameterisation function of $\bSigma_u$, then we have the following equations for the inverse lapse problem associated with $\{\bSigma_u\}$,
\begin{align}
	&
	\uL \flu = \balu,
	\label{eqn 6.3}
\\
	&
	\buL \log \balu = 2 \buomegalu - 2 \balu\, \duomegalu,
	\label{eqn 6.4}
\\
	&
	\left\{
	\begin{aligned}
		&
		\bslashDelta \log \balu 
		= 
		\ddmulu - \bmulu 
		= 
		- (\ddrholu - \overline{\ddrholu})
		- \frac{1}{2} [ (\hatdduchilu, \hatddchilu') - \overline{( \hatdduchilu, \hatddchilu')}^u ]
		- \bslashdiv \ddetalu,
	\\
		&
		\overline{\balu\, \ddtr \dduchilu}^u
		=\overline{\btr \buchilu}^u 
		= \frac{2}{\br_u}.
	\end{aligned}
	\right.
	\label{eqn 6.5}
\end{align}
where we use $\overline{}^{u}$ to denote the mean value with respect to the metric $\bslashglu$. Equations \eqref{eqn 6.4} and \eqref{eqn 6.5} are compatible. Along $\{\bSigma_u\}$, we have that
\begin{align}
	&
	\phantom{\big\{ \}}
	\bmulu = \overline{\bmulu}^u,
	\label{eqn 6.6}
\\
	&
	\phantom{\big\{ \}}
	\buL \overline{\bmulu}^u 
	=
	- \frac{3}{2} \overline{\bmulu}^u\, \overline{\btr \buchilu}^u
	+ \frac{1}{4} \overline{\btr \bchilu' | \hatbuchilu |^2}^u
	+ \frac{1}{2} \overline{\btr \buchilu |\btalu|^2}^u,
	\label{eqn 6.7}
\\
	&
	\phantom{\big\{ \}}
	\buL \btr \buchilu
	= 
	2 \buomegalu\, \btr \buchilu 
	- | \hatbuchilu |^2 
	- \frac{1}{2} ( \btr \buchilu )^2,
	\label{eqn 6.8}
\\
	&
	\phantom{\big\{ \}}
	\buL \btr \bchilu' 
	= 
	-2 \buomegalu\, \btr \bchilu' 
	- \frac{1}{2} \btr \buchilu\, \btr \bchilu' 
	- 2 | \btalu |^2
	+ 2 \bmulu,
	\label{eqn 6.9}
\\
	&
	\phantom{\big\{ \}}
	\bslashdiv \hatbuchilu
	- \frac{1}{2} \bslashd \btr \buchilu
	- \hatbuchilu \cdot \btalu 
	+ \frac{1}{2} \btr \buchilu \, \btalu
	=
	- \bubetalu,
	\label{eqn 6.10}
\\
	&
	\phantom{\big\{ \}}
	\bslashdiv \hatbchilu'
	- \frac{1}{2} \bslashd \btr \bchilu'
	+ \hatbchilu' \cdot \btalu 
	- \frac{1}{2} \btr \bchilu' \, \btalu
	=
	- \bbetalu,
	\label{eqn 6.11}
\\
	&
	\left\{
	\begin{aligned}
		&
		\bslashcurl \btalu
		=
		\frac{1}{2} \hatbchilu' \wedge \hatbuchilu 
		+ \bsigmalu,
	\\
		&
		\bslashdiv \btalu
		= 
		- \brholu
		- \frac{1}{2} ( \hatbuchilu, \hatbchilu' ) 
		- \bmulu,
	\end{aligned}
	\label{eqn 6.12}
	\right.
\\
	&
	\left\{
	\begin{aligned}
		2 \bslashDelta \buomegalu
		&=
		-\frac{3}{2} ( \bmulu\, \btr \buchilu - \overline{\bmulu\, \btr \buchilu}^u )
		+ \frac{1}{2} ( \btr \buchilu | \btalu |^2 - \overline{\btr \buchilu | \btalu |^2}^u )
	\\
		&
		\phantom{=}
		+ \frac{1}{4} ( \btr \bchilu' | \hatbuchilu |^2 - \overline{ \btr \bchilu' | \hatbuchilu|^2}^u )
		+ 4 (\bslashdiv \hatbuchilu, \btalu )
		+ 4 ( \hatbuchilu, \bslashnabla \btalu )
		- 2 \bslashdiv \bubetalu,
	\\
		\overline{\buomegalu}^u
		&=
		-\frac{r_u}{2} \overline{(\buomegalu - \overline{\buomegalu}^u)(\btr\buchilu - \overline{\btr\buchilu}^u)}^u
		- \frac{r_u}{8} \overline{(\btr \buchilu - \overline{\btr \buchilu}^u)^2}^u
		+ \frac{r_u}{4} \overline{|\hatbuchilu|^2}^u,
	\end{aligned}
	\label{eqn 6.13}
	\right.
\end{align}
where the mean value is taken with respect to the metric $\bslashglu$, and
\begin{align*}
	&
	(\hatbuchilu \cdot \btalu)_a = \hatbuchilu_{ab}\, \btalu_c\, (\bslashglu^{-1})^{bc},
	\quad \quad
	(\hatbchilu' \cdot \btalu)_a = \hatbchilu'_{ab}\, \btalu_c\, (\bslashglu^{-1})^{bc},
\\
	&
	( \hatbchilu', \hatbuchilu ) = \hatbchilu'_{ab}\, \hatbuchilu_{cd}\, ( \bslashglu^{-1} )^{ac}\, ( \bslashglu^{-1} )^{bd},
	\quad \quad
	( \hatbuchilu, \bslashnabla \btalu ) = \hatbchilu'_{ab}\, \bslashnabla_c \btalu_d\, ( \bslashglu^{-1} )^{ac}\, ( \bslashglu^{-1} )^{bd},
\\
	&
	(\bslashdiv \hatbuchilu, \btalu )= (\bslashdiv \hatbuchilu)_a\, \btalu_b\, (\bslashglu^{-1})^{ab},
	\quad \quad
	\hatbchilu' \wedge \hatbuchilu = \hatbchilu'_{ab}\, \hatbuchilu_{cd}\, (\bslashglu^{-1})^{ac}\, \bslashepsilonlu^{bd}.
\end{align*}

\subsection{Schwarzschild spherically symmetric constant mass aspect function foliation}\label{sec 6.3}
We present the parametrisation and the geometry of the spherically symmetric constant mass aspect function foliation in the Schwarzschild spacetime. Let $\{ \bSigma_u \}$ be the spherically symmetric constant mass aspect function foliation in the incoming null hypersurface $\uC_{\us}$ and suppose that $\bSigma_{u=0}=\Sigma_{\us,s_0}$. Then we have the following results for $\{ \bSigma_u\}$.
\begin{enumerate}[label=\alph*.]
\item Let $\br_u$ be the area radius of the leaf $\bSigma_u$, then $\br_{u=0}$ is solved implicitly by the equation
\begin{align*}
	(\br_{u=0}-r_0)\exp\frac{\br_{u=0}}{r_0} = s_0\exp\frac{\us+s_0+r_0}{r_0}.
\end{align*}
Then the area radius of a general leaf $\bSigma_u$ is given by the formula $\br_u = \br_{u=0} +u$. For the special case $s_0=0$, we have $\br_{u=0} = r_0$ and $\br_u = r_0+u$.

\item Suppose that $\bSigma_u$ is the surface $\Sigma_{\us, s_u}$, then $s_u$ is given by the equation
\begin{align*}
	(\br_u-r_0)\exp\frac{\br_u}{r_0} = s_u\exp\frac{\us+s_u+r_0}{r_0}.
\end{align*}

\item The lapse function $\balu$ of the leaf $\bSigma_u$ is solved by equation \eqref{eqn 6.5} which reduces to the following simple form that
\begin{align*}
	\balu 
	= 
	\frac{\br_u}{(s_u+r_0)} \exp (- \frac{\us + s_u + r_0 - \br_u}{r_0})
	=
	\frac{\br_u}{s_u+r_0} \cdot \frac{s_u}{\br_u - r_0}.
\end{align*}
Note that when $\br_u=r_0$, i.e. $s_u=0$, we have that $\balu = \exp (-\frac{\us}{r_0})$.

\item The connection coefficients associated with the foliation $\{ \bSigma_u \}$ with respect to the frame $\{\buL^u, \bL'^u\}$ are given by
\begin{align*}
	\btr \buchilu = \frac{2}{\br_u},
	\quad
	\hatbuchilu =0,
	\quad
	\btr \bchilu' = \frac{2}{\br_u}(1-\frac{r_0}{\br_u}),
	\quad
	\hatbchilu' = 0,
	\quad
	\btalu = 0,
	\quad
	\buomegalu = 0.
\end{align*}

\item The curvature components associated with the foliation $\{ \bSigma_u \}$ with respect to the frame $\{\buL^u, \bL'^u\}$ are given by
\begin{align*}
	\balphalu = \bualphalu = 0,
	\quad
	\bbetalu = \bubetalu = 0,
	\quad
	\bsigmalu = 0,
	\quad
	\brholu = - \frac{ r_0}{\br_u^3}.
\end{align*}

\end{enumerate}

We prove a simple bound of $s_u$ which is useful later.
\begin{lemma}\label{lem 6.3}
For $s_u$ coordinate of the leaf $\bSigma$ described in the above item b., we have
\begin{align*}
	\max\{ e^{-0.1} (\br_u - r_0), \br_u - 1.1 r_0 \} 
	\leq 
	s_u 
	\leq 
	\min \{ e^{0.1} (\br_u - r_0), \br_u - 0.9 r_0 \},
\end{align*}
provided that $\br_u - r_0 \geq 0$.
\end{lemma}
\begin{proof}
The proof is elementary. Note that $ \us \in (-0.1 r_0, 0.1 r_0)$ in definition \ref{def 2.3}, the equation in item b. implies that given $\br_u - r_0 \geq 0$,
\begin{align*}
	(\br_u - r_0) \exp \frac{\br_u - r_0}{r_0}
	=
	s_u \exp \frac{\us + s_u}{r_0}
	\leq
	s_u \exp \frac{0.1 r_0 + s_u}{r_0}
	\leq
	\left\{
	\begin{aligned}
		&
		e^{0.1}s_u \exp \frac{e^{0.1} s_u}{r_0},
	\\
		&
		(0.1 r_0+s_u) \exp \frac{0.1 r_0 + s_u}{r_0}.
	\end{aligned}
	\right.
\end{align*}
Then the inequality on the left side follows. The other inequality follows similarly.
\end{proof}

\subsection{Global existence problem and geometry of constant mass aspect function foliation}\label{sec 6.4}
We formulate the main problem investigated in this paper. 
\begin{problem}\label{prob 6.4}
Let $\ucalH$ be a nearly spherically symmetric incoming null hypersurface in the vacuum perturbed Schwarzschild spacetime $(M,g)$ and $\bSigma_{u=0}$ be a nearly spherically symmetric spacelike surface in $\ucalH$ which is close to the Schwarzschild event horizon $C_{s=0}$. Construct the constant mass aspect function foliation $\{\bSigma_u\}$ initiating from $\bSigma_{u=0}$. Show that the foliation exists for $u\in [-\frac{\kappa r_0}{2}, + \infty)$, which implies that the foliation can be extended to the past null infinity. Moreover obtain the estimates of the geometric quantities associated with the foliation.
\end{problem}

The above problem is only intuitively stated and involves several issues required to be clarified, which we roughly explain in the following and leave the precise formulations to later theorems answering the above problem.
\begin{enumerate}[label={\roman*}.]
\item Nearly spherically symmetric incoming null hypersurface $\ucalH$. See the review of the result from \cite{L2022} in section \ref{sec 3.1} and proposition \ref{prop 3.1}. Let $\uh$ be the parameterisation of $\ucalH$, then the nearly spherical symmetry of $\ucalH$ is characterised as the condition that $\ufl{s=0} = \uh|_{s=0}$ is nearly constant where the precise quantitative formulation is given in proposition \ref{prop 3.1}. It is shown that a nearly spherically symmetric incoming null hypersurface $\ucalH$ can be extended regularly to the past null infinity.

\item Nearly spherically symmetric spacelike surface $\bSigma_{u=0}$. Let $\fl{u=0}$ be the parameterisation of $\bSigma_{u=0}$ in the $\{s,\vartheta\}$ coordinate system of $\ucalH$. We formulate the nearly spherical symmetry of $\bSigma_{u=0}$ as the nearly constancy of $\fl{u=0}$, which will be quantitatively characterised as the smallness of the Sobolev norm of the differential $\slashd \ufl{s=0}$.
\end{enumerate}

The method to solve problem \ref{prob 6.4} is to prove the global existence for the solution of the basic equations \eqref{eqn 6.3} - \eqref{eqn 6.13} in the interval $u \in [-\frac{\kappa r_0}{2}, +\infty)$ for the nearly constant initial data $\fl{u=0}$. In the procedure of the proof of the global existence, we shall obtain the estimate of the geometric quantities associated with the foliation.

\section{Geometry of initial leaf of constant mass aspect function foliation}\label{sec 7}

In this section, we obtain the estimate for the geometric quantities on the nearly spherically symmetric initial leaf $\bSigma_{u=0}$ of the constant mass aspect function foliation. 

The form of the estimates given in this section is through the comparison with the geometric quantities on the spherically symmetric constant mass aspect function foliation in the Schwarzschild spacetime. Thus we introduce the following notation to denote the reference objects in the comparison. 

Let $\{ \bSigma_{u, S}^{\us = \ous} \}$ be the spherically symmetric Schwarzschild constant mass aspect function foliation in the incoming null hypersurface $\uC_{\us = \ous}$ with the initial leaf $\bSigma_{u=0,S}^{\us = \ous} = \Sigma_{s=0, \us = \ous}$. The geometric quantities associated with this Schwarzschild foliation $\{ \bSigma_{u, S}^{\us = \ous} \}$ are given in section \ref{sec 6.3}, which we use the subscript $S$ to indicate the corresponding geometric quantity being associated with $\{ \bSigma_{u, S}^{\us = \ous} \}$.

Furthermore, in order to simplify the estimates, we introduce the notation $\ud_{o,\uh}$ by
\begin{align*}
\ud_{o,\uh} = \max \{ \ud_o, \ud_{\uh} \},
\end{align*}
then we can replace the bounds $\ud_o$ and $\ud_{\uh}$ by the single one $\ud_{o,\uh}$, and estimates \eqref{eqn 5.15} imply that
\begin{align}
	\Vert \ddslashd \uflu \Vert^{n,p}, 
	\Vert (\dslashd \uh)|_{\bSigma_u} \Vert^{n+1,p},
	\Vert (\circnabla^2 \uh)|_{\bSigma_u} \Vert^{n,p}
	\leq
	\ud_{o,\uh} r_0.
	\label{eqn 7.1}
\end{align}

\subsection{Assumption on parameterisation functions of initial leaf $\bSigma_{u=0}$}
We summarise the quantitative assumptions of the parameterisation functions of the initial leaf $\bSigma_{u=0}$ which serves as the basic assumption for the estimates on the geometry of the initial leaf.
\begin{assumption}\label{assum 7.1}
Let $\bSigma_{u=0}$ be the initial leaf of the constant mass aspect function foliation in the incoming null hypersurface $\ucalH$. Let $(\ufl{s=0}, \fl{u=0})$ be the parameterisation of $\bSigma_{u=0}$. Suppose that
\begin{align*}
	&
	\Vert \slashd \ufl{s=0} \Vert^{n+1,p} \leq \udelta_o r_0,
	\quad
	\overline{\ufl{s=0}}^{\circg} = \ous,
	\quad
	\vert \ous \vert \leq \udelta_m r_0 \leq \kappa r_0 \leq 0.1 r_0,
\\
	&
	\Vert \slashd \fl{u=0} \Vert^{n+1,p} \leq \delta_o (r_0+\os_{u=0}),
	\quad
	\overline{\fl{u=0}}^{\circg} = \os_{u=0},
	\quad
	\vert \os_{u=0} \vert \leq \delta_m r_0,
\end{align*}
where $n\geq 1+n_p$, $p >1$. Choose a positive constant $\delta$ such that $\epsilon, \udelta_o, \udelta_m, \delta_o, \delta_m $ are suitably bounded by $\epsilon, \udelta_o, \epsilon \udelta_m, \delta_o, \delta_m \leq \delta$. We assume that $\delta$ is sufficiently small depending on $n,p$ such that estimates \eqref{eqn 5.15} hold and
\begin{align*}
	\ud_o \leq 1,
	\quad
	\uslashd_m \leq 1,
	\quad
	\ud_m \leq 2,
	\quad
	\ud_{\uh} \leq 1.
\end{align*}
\end{assumption}

\subsection{Estimate of initial lapse $\bal{u=0}$}
We shall estimate the initial lapse $\bal{u=0}$ for the foliation on the initial leaf $\bSigma_{u=0}$ by equation \eqref{eqn 6.5}.
\begin{proposition}\label{prop 7.2}
Under assumption \ref{assum 7.1}, there exists a small positive constant $\delta$ depending on $n,p$ such that
\begin{align*}
	&
	\Vert \slashd \bal{u=0} \Vert^{n,p}
	\leq
	c(n,p) [\epsilon + \delta_o + \delta_m (\ud_o + \ud_{\uh})],
\\
	&
	\vert \bal{u=0} - \bal{u=0}_S \vert
	\leq
	c(n,p)(\epsilon + \delta_o + \delta_m + \ud_o + \ud_{\uh}).
\end{align*}
\end{proposition}
We estimate $\bal{u=0}$ by equation \eqref{eqn 6.5}. Before proving the proposition, we prepare some preliminarily estimates of the terms appearing in equation \eqref{eqn 6.5} in the following lemma.
\begin{lemma}\label{lem 7.3}
Under assumption \ref{assum 7.1}, we have that
\begin{align*}
	\bslashgl{u=0}
	&:
	\begin{aligned}
		\vert \bslashgl{u=0} - \bslashgl{u=0}_S \vert
		\leq
		c(n,p) (\epsilon + \delta_o + \delta_m) r_0^2,
		\quad
		\Vert \bcircnabla \bslashgl{u=0} \Vert^{n,p}
		\leq
		c(n,p) (\epsilon + \delta_o + \delta_m \ud_o) r_0^2,
	\end{aligned}
\\
	\bslashepsilonl{u=0}
	&:
	\begin{aligned}
		\vert \bslashepsilonl{u=0} - \bslashepsilonl{u=0}_S \vert
		\leq
		c(n,p) (\epsilon + \delta_o + \delta_m) r_0^2,
		\quad
		\Vert \bcircnabla \bslashepsilonl{u=0} \Vert^{n,p}
		\leq
		c(n,p) (\epsilon + \delta_o + \delta_m \ud_o) r_0^2,
	\end{aligned}
\\
	\br_{u=0}
	&:
	\begin{aligned}
		\vert \br_{u=0} - \br_{u=0,S} \vert
		\leq
		c(n,p) (\epsilon + \delta_o + \delta_m) r_0,
	\end{aligned}
\\
	\brhol{u=0}
	&:
	\begin{aligned}
		\vert \brhol{u=0} - \brhol{u=0}_S \vert
		\leq
		\frac{c(n,p) ( \epsilon + \delta_o + \delta_m)}{r_0^2},
		\quad
		\Vert \bslashd \brhol{u=0} \Vert^{n-1,p}
		\leq
		\frac{c(n,p)(\epsilon + \delta_o + \delta_m \ud_o)}{r_0^2},
	\end{aligned}
\\
	\ddtr\dduchil{u=0}
	&:
	\left\{
	\begin{aligned}
		&
		\vert \ddtr \dduchil{u=0} - (\ddtr \dduchil{u=0})_{S} \vert
		\leq
		\frac{c(n,p) (\epsilon + \delta_o + \delta_m + \ud_o + \ud_{\uh})}{r_0},
	\\
		&
		\Vert \bslashd \ddtr \dduchil{u=0} \Vert^{n-1,p}
		\leq
		\frac{c(n,p) (\epsilon + \delta_o + \ud_o + \ud_{\uh})}{r_0}.
	\end{aligned}
	\right.
\end{align*}
\end{lemma}

\subsubsection{Proof of preliminarily estimates}
\begin{proof}[Proof of lemma \ref{lem 7.3}]
We estimate each term in the lemma in the following.
\begin{enumerate}[label={\textbullet}]
\item 
$\pmb{\bslashgl{u=0} - \bslashgl{u=0}_S}$.
Decompose the term as follows
\begin{align*}
	\bslashgl{u=0} - \bslashgl{u=0}_S
	&=
	\lo{\ddslashg|_{\bSigma_{u=0}}} - \bslashgl{u=0}_S 
	+ \hi{\ddslashg|_{\bSigma_{u=0}}}
\\
	&=
	[(r_S|_{\bSigma_{u=0}})^2 - r_0^2] \circg 
	+ \hi{\ddslashg|_{\bSigma_{u=0}}}.
\end{align*}
By the estimates of the parameterisation functions in the assumption and \eqref{eqn 5.15},
\begin{align*}
	\vert (r_S|_{\bSigma_{u=0}}) - r_0 \vert
	\leq
	c(n,p) (\delta_o + \delta_m) r_0,
\end{align*}
thus together with the estimate of $\hi{\ddslashg|_{\bSigma_{u=0}}}$ in proposition \ref{prop 5.10}, we obtain that
\begin{align*}
	\vert \bslashgl{u=0} - \bslashgl{u=0}_S \vert
	\leq
	c(n,p) (\epsilon + \delta_o + \delta_m) r_0^2.
\end{align*}

\item 
$\pmb{\bcircnabla \bslashgl{u=0}}$.
Decompose the term as follows
\begin{align*}
	\bcircnabla \bslashgl{u=0}
	=
	\bcircnabla \lo{\ddslashg|_{\bSigma_{u=0}}} + \bcircnabla \hi{\ddslashg|_{\bSigma_{u=0}}}
	=
	\bslashd [(r_S|_{\bSigma_{u=0}})^2] \circg + \bcircnabla \hi{\ddslashg|_{\bSigma_{u=0}}}.
\end{align*}
By the estimates of the parameterisation functions in the assumption and \eqref{eqn 5.15},
\begin{align*}
	\Vert \bslashd (r_S|_{\bSigma_{u=0}}) \Vert^{n,p}
	\leq
	c(n,p) ( \delta_o + \delta_m \ud_o) r_0,
\end{align*}
thus together with the estimate of $\hi{\ddslashg|_{\bSigma_{u=0}}}$ in proposition \ref{prop 5.10}, we obtain that
\begin{align*}
	\Vert \bcircnabla \bslashgl{u=0} \Vert^{n,p}
	\leq
	c(n,p) (\epsilon + \delta_o + \delta_m \ud_o) r_0^2.
\end{align*}

\item
$\pmb{\bslashepsilonl{u=0} - \bslashepsilonl{u=0}_S, \bcircnabla \bslashepsilonl{u=0}, \br_{u=0} - \br_{u=0,S} }$. Their estimates follow from the above estimates of $\bslashgl{u=0} - \bslashgl{u=0}_S$, $\bcircnabla \bslashgl{u=0}$.

\item
$\pmb{\brhol{u=0} - \brhol{u=0}_S}$.
Decompose the term as follows
\begin{align*}
	\brhol{u=0} - \brhol{u=0}_S
	&=
	\ddrho|_{\bSigma_{u=0}} - \brhol{u=0}_S
	=
	\lo{\ddrho|_{\bSigma_{u=0}}}  - \brhol{u=0}_S
	+ \hi{\ddrho|_{\bSigma_{u=0}}}
\\
	&=
	\rho_S|_{\bSigma_{u=0}} - \rho_S|_{\Sigma_{s=0,\us = \ous}}
	+ \hi{\ddrho|_{\bSigma_{u=0}}}.
\end{align*}
By the estimates of the parameterisation functions in the assumption and \eqref{eqn 5.15},
\begin{align*}
	\vert \rho_S|_{\bSigma_{u=0}} - \rho_S|_{\Sigma_{s=0,\us = \ous}} \vert
	\leq
	\frac{c(n,p) (\delta_o + \delta_m)}{r_0^2},
\end{align*}
thus together with the estimate of $\hi{\ddrho|_{\bSigma_{u=0}}}$ in proposition \ref{prop 5.10}, we obtain that
\begin{align*}
	\vert \brhol{u=0} - \brhol{u=0}_S \vert
	\leq
	\frac{c(n,p) (\epsilon + \delta_o + \delta_m)}{r_0^2}.
\end{align*}

\item
$\pmb{\bslashd \brhol{u=0}}$.
Decompose the term as follows
\begin{align*}
	\bslashd \brhol{u=0}
	=
	\bslashd \lo{\ddrho|_{\bSigma_{u=0}}} + \bslashd \hi{\ddrho|_{\bSigma_{u=0}}}
	=
	\bslashd (\rho_S|_{\bSigma_{u=0}}) + \bcircnabla \hi{\bslashgl{u=0}}.
\end{align*}
By the estimates of the parameterisation functions in the assumption and \eqref{eqn 5.15},
\begin{align*}
	\Vert \bslashd (\rho_S|_{\bSigma_{u=0}}) \Vert^{n-1,p}
	\leq
	\frac{c(n,p) ( \delta_o + \delta_m \ud_o)}{r_0^2},
\end{align*}
thus together with the estimate of $\hi{\ddrho|_{\bSigma_{u=0}}}$ in proposition \ref{prop 5.10}, we obtain that
\begin{align*}
	\Vert \bslashd \brhol{u=0} \Vert^{n-1,p}
	\leq
	\frac{c(n,p) (\epsilon + \delta_o + \delta_m \ud_o)}{r_0^2}.
\end{align*}

\item 
$\pmb{\ddtr \dduchil{u=0} - ( \ddtr \dduchil{u=0})_S}$.
Decompose the term as follows
\begin{align*}
	\ddtr \dduchil{u=0} - (\ddtr \dduchil{u=0})_S
	&
	=
	\lo{\ddtr \dduchil{u=0}}  - (\tr \uchi)_S|_{\bSigma_{u=0,S}}
	+ \hi{\ddtr \dduchil{u=0}}
\\
	&=
	(\tr \uchi)_S|_{\bSigma_{u=0}} - (\tr \uchi)_S|_{\Sigma_{s=0,\us = \ous}}
	- 2 (\Omega_S^2 r_S^{-2})|_{\bSigma_{u=0}} (\circDelta \uh)|_{\bSigma_{u=0}}
\\
	&\phantom{=}
	+ \hi{\ddtr \dduchil{u=0}}.
\end{align*}
By the estimates of the parameterisation functions in the assumption and \eqref{eqn 5.15},
\begin{align*}
	&
	\vert (\tr \uchi)_S|_{\bSigma_{u=0}} - (\tr \uchi)_S|_{\Sigma_{s=0,\us = \ous}} \vert
	\leq
	\frac{c(n,p) (\delta_o + \delta_m + \ud_o)}{r_0},
\\
	&
	\vert 2 (\Omega_S^2 r_S^{-2})|_{\bSigma_{u=0}} (\circDelta \uh)|_{\bSigma_{u=0}} \vert
	\leq
	\frac{c(n,p)\ud_{\uh}}{r_0},
\end{align*}
thus together with the estimate of $\hi{\ddtr \dduchi|_{\bSigma_{u=0}}}$ in proposition \ref{prop 5.10}, we obtain that
\begin{align*}
	\vert \ddtr \dduchil{u=0} - (\ddtr \dduchil{u=0})_S \vert
	\leq
	\frac{c(n,p) (\epsilon + \delta_o + \delta_m + \ud_o + \ud_{\uh})}{r_0}.
\end{align*}

\item
$\pmb{\bslashd \ddtr \dduchil{u=0}}$.
Decompose the term as follows
\begin{align*}
	\bslashd \ddtr \dduchil{u=0} 
	&
	=
	\bslashd\lo{\ddtr \dduchil{u=0}}
	+ \bslashd\hi{\ddtr \dduchil{u=0}}
\\
	&=
	\bslashd [(\tr \uchi)_S|_{\bSigma_{u=0}}] 
	- 2 \bslashd [(\Omega_S^2 r_S^{-2})|_{\bSigma_{u=0}} (\circDelta \uh)|_{\bSigma_{u=0}}]
	+ \bslashd \hi{\ddtr \dduchil{u=0}}.
\end{align*}
By the estimates of the parameterisation functions in the assumption and \eqref{eqn 5.15},
\begin{align*}
	&
	\Vert \bslashd[ (\tr \uchi)_S|_{\bSigma_{u=0}} ] \Vert^{n,p}
	\leq
	\frac{c(n,p) (\delta_o + \ud_o)}{r_0},
\\
	&
	\Vert \bslashd [(\Omega_S^2 r_S^{-2})|_{\bSigma_{u=0}} (\circDelta \uh)|_{\bSigma_{u=0}}] \Vert^{n-1,p}
	\leq
	\frac{c(n,p)\ud_{\uh}}{r_0},
\end{align*}
thus together with the estimate of $\hi{\ddtr \dduchi|_{\bSigma_{u=0}}}$ in proposition \ref{prop 5.10}, we obtain that
\begin{align*}
	\Vert  \bslashd \ddtr \dduchil{u=0} \Vert^{n-1,p}
	\leq
	\frac{c(n,p) (\epsilon + \delta_o + \ud_o + \ud_{\uh})}{r_0}.
\end{align*}
\end{enumerate}
\end{proof}

\subsubsection{Proof of estimate of $\bal{u=0}$}
Now we turn to the proof of the proposition.
\begin{proof}[Proof of proposition \ref{prop 7.2}]
We further require $\delta$ sufficiently small such that
\begin{align*}
c(n,p) (\epsilon + \delta_o + \delta_m + \ud_o + \ud_{\uh} ) \leq \frac{1}{2}, 
\end{align*}
for all constants $c(n,p)$ appearing in the proof.

\begin{enumerate}[label={\textbullet}]
\item
$\pmb{\bslashd \log \bal{u=0}}$.
We estimate each term in the elliptic equation of $\log \bal{u=0}$ in \eqref{eqn 6.5}.
\begin{enumerate}[label=\raisebox{0.1ex}{\scriptsize$\bullet$}]
\item 
$\brhol{u=0} - \overline{\brhol{u=0}}$. By lemma \ref{lem 7.3}, we have 
\begin{align*}
	\Vert \brhol{u=0} - \overline{\brhol{u=0}} \Vert^{n,p} 
	\leq 
	c(n,p) \Vert \dslashd \brhol{u=0} \Vert^{n,p} 
	\leq 
	\frac{c(n,p)(\epsilon + \delta_o + \delta_m \ud_o)}{r_0^2}.
\end{align*}

\item
$(\hatdduchil{u=0}, \hatddchil{u=0}') - \overline{( \hatdduchil{u=0}, \hatddchil{u=0}')}^{u=0}$.
By the estimates in propositions \ref{prop 5.8} and \ref{prop 5.10},
\begin{align*}
	&
	\Vert \hatddchil{u=0}' \Vert^{n,p}
	\leq
	\Vert \lo{\hatddchil{u=0}'} \Vert^{n,p}
	+ \Vert \hi{\hatddchil{u=0}'} \Vert^{n,p}
	\leq
	c(n,p) ( \epsilon + \delta_o ) r_0,
\\
	&
	\Vert \hatdduchil{u=0} \Vert^{n,p}
	\leq
	\Vert \lo{\hatdduchil{u=0}} \Vert^{n,p}
	+ \Vert \hi{\hatdduchil{u=0}} \Vert^{n,p}
	\leq
	c(n,p) ( \epsilon + \ud_{\uh}) r_0,
\end{align*}
which implies that
\begin{align*}
	\Vert (\hatdduchil{u=0}, \hatddchil{u=0}') - \overline{( \hatdduchil{u=0}, \hatddchil{u=0}')}^{u=0} \Vert^{n,p}
	\leq
	\frac{c(n,p) ( \epsilon + \delta_o )( \epsilon + \ud_{\uh})}{r_0^2}.
\end{align*}

\item
$\bslashdiv \ddeta|_{\bSigma_{u=0}}$.
By the estimates in propositions \ref{prop 5.8} and \ref{prop 5.10},
\begin{align*}
	&
	\Vert \ddeta|_{\bSigma_{u=0}} \Vert^{n,p}
	\leq
	\Vert \lo{\ddeta|_{\bSigma_{u=0}}} \Vert^{n,p}
	+ \Vert \hi{\ddeta|_{\bSigma_{u=0}}} \Vert^{n,p}
	\leq
	c(n,p) ( \epsilon + \delta_o + \delta_m \ud_{\uh} ),
\\
	&
	\Vert \bslashnabla \ddeta|_{\bSigma_{u=0}} \Vert^{n-1,p}
	\leq
	\Vert \bcircnabla \ddeta|_{\bSigma_{u=0}} \Vert^{n-1,p}
	+ \Vert \triangle_{ij}^k (\ddeta|_{\bSigma_{u=0}})_k \Vert^{n-1,p}
	\leq
	c(n,p) ( \epsilon + \delta_o + \delta_m \ud_{\uh} ),
\\
	&
	\Vert \bslashdiv \ddeta|_{\bSigma_{u=0}} \Vert^{n-1,p}
	\leq
	c(n,p) \Vert (\bslashgl{u=0})^{-1} \Vert^{n-1,p} \Vert \bslashnabla \ddeta|_{\bSigma_{u=0}} \Vert^{n-1,p}
	\leq
	\frac{c(n,p) ( \epsilon + \delta_o + \delta_m \ud_{\uh} )}{r_0^2}.
\end{align*}
\end{enumerate}
From the above estimates, we obtain that
\begin{align*}
	\Vert \bslashDelta \log \bal{u=0} \Vert^{n-1,p}
	\leq
	\frac{c(n,p) [ \epsilon + \delta_o + \delta_m(\ud_o + \ud_{\uh}) ]}{r_0^2}.
\end{align*}
Then by the estimate of $\bslashgl{u=0}$ in lemma \ref{lem 7.3} and the theory of the elliptic equation on the sphere, we derive that
\begin{align*}
	\Vert \bslashd \log \bal{u=0} \Vert^{n,p}
	\leq
	c(n,p) [ \epsilon + \delta_o + \delta_m (\ud_o + \ud_{\uh}) ].
\end{align*}
As a corollary, we have that by the Sobolev inequality,
\begin{align*}
	\vert \log \bal{u=0}  - \overline{\log \bal{u=0}}^{\circg} \vert 
	\leq 
	c(n,p) [ \epsilon + \delta_o + \delta_m (\ud_o + \ud_{\uh}) ],
\end{align*}
which implies that
\begin{align*}
	\{1- c(n,p) [ \epsilon + \delta_o + \delta_m (\ud_o + \ud_{\uh}) ] \}
	\leq
	\frac{\bal{u=0}}{\exp (\overline{\log \bal{u=0}}^{\circg})}
	\leq
	\{1+ c(n,p) [ \epsilon + \delta_o + \delta_m (\ud_o + \ud_{\uh}) ] \}.
\end{align*}

\item
$\pmb{\overline{\log \bal{u=0}}^{\circg}, \bal{u=0}}$.
We estimate $\overline{\log \bal{u=0}}^{\circg}$ by the second equation in \eqref{eqn 6.5}. By the estimate of $\ddtr \dduchil{u=0} - (\ddtr \dduchil{u=0})_S$ in lemma \ref{lem 7.3}, we have
\begin{align*}
	\vert \ddtr \dduchil{u=0} - \frac{2}{r_0} \exp \frac{\ous}{r_0} \vert
	\leq
	\frac{c(n,p) ( \epsilon + \delta_o + \delta_m + \ud_o + \ud_{\uh})}{r_0}.
\end{align*}
Then by \eqref{eqn 6.5},
\begin{align*}
	\frac{2}{\br_{u=0}}
	&
	\leq
	\exp (\overline{\log \bal{u=0}}^{\circg})
	\cdot
	[ \frac{2}{r_0} \exp \frac{\ous}{r_0} + \frac{c(n,p) ( \epsilon + \delta_o + \delta_m + \ud_o + \ud_{\uh})}{r_0} ]
\\
	&
	\phantom{\leq
	\exp (\overline{\log \bal{u=0}}^{\circg})}
	\cdot
	\{1+ c(n,p) [ \epsilon + \delta_o + \delta_m (\ud_o + \ud_{\uh}) ] \}
\\
	&
	\leq
	\exp (\overline{\log \bal{u=0}}^{\circg}) \cdot \frac{(2 e^{0.1} +1) \cdot 2}{r_0}
\\
	\Rightarrow
	\quad
	\overline{\log \bal{u=0}}^{\circg}
	&
	\geq 
	\log \frac{r_0}{5r_{u=0}}
	\geq
	\log \frac{1}{5[1+ c(n,p) ( \epsilon + \delta_o + \delta_m)]}
	\geq
	\log \frac{1}{10}.
\end{align*}
Similarly, we also obtain the upper bound of $\overline{\log \bal{u=0}}^{\circg}$ that
\begin{align*}
	\frac{2}{\br_{u=0}}
	&
	\geq
	\exp (\overline{\log \bal{u=0}}^{\circg})
	\cdot
	[ \frac{2}{r_0} \exp \frac{\ous}{r_0} - \frac{c(n,p) ( \epsilon + \delta_o + \delta_m + \ud_o + \ud_{\uh})}{r_0} ]
\\
	&
	\phantom{\leq
	\exp (\overline{\log \bal{u=0}}^{\circg})}
	\cdot
	\{1- c(n,p) [ \epsilon + \delta_o + \delta_m (\ud_o + \ud_{\uh}) ] \}
\\
	&
	\geq
	\exp (\overline{\log \bal{u=0}}^{\circg}) 
	\cdot 
	\frac{(2 e^{-0.1} -\frac{1}{2}) \cdot \frac{1}{2}}{r_0}
\\
	\Rightarrow
	\quad
	\overline{\log \bal{u=0}}^{\circg}
	&
	\leq
	\log\frac{4r_0}{\br_{u=0}}
	\leq
	\log \frac{4}{1- c(n,p) (\epsilon + \delta_o + \delta_m)}
	\leq
	\log 8.
\end{align*}
Therefore we obtain that
\begin{align*}
	\log \frac{1}{10}
	\leq
	\overline{\log \bal{u=0}}^{\circg}
	\leq
	\log 8,
	\quad
	\frac{1}{10}
	\leq
	\bal{u=0}
	\leq
	8.
\end{align*}

\item 
$\pmb{\bal{u=0} - \overline{\bal{u=0}}^{u=0}, \bslashd \bal{u=0}}$.
By the estimates of $\frac{\bal{u=0}}{\exp (\overline{\log \bal{u=0}}^{\circg})}$ and $\overline{\log \bal{u=0}}^{\circg}$, we obtain that
\begin{align*}
	&
	\max \{ \bal{u=0} \} - \min \{ \bal{u=0} \}
	\leq
	c(n,p) [\epsilon + \delta_o + \delta_m ( \ud_o + \ud_m)]
\\
	\Rightarrow
	\quad
	&
	\vert \bal{u=0} - \overline{\bal{u=0}}^{u=0} \vert
	\leq
	c(n,p) [\epsilon + \delta_o + \delta_m ( \ud_o + \ud_m)].
\end{align*}
By the estimates of $\bslashd \log \bal{u=0}$, $\bal{u=0}$ and the ranges of $n,p$, we have
\begin{align*}
	\Vert \bslashd \bal{u=0} \Vert^{n,p}
	\leq
	c(n,p) [\epsilon + \delta_o + \delta_m ( \ud_o + \ud_m)].
\end{align*}

\item
$\pmb{\bal{u=0} - \bal{u=0}_S}$.
We decompose $\overline{\bal{u=0}\, \ddtr \dduchil{u=0}}^{u=0}$ as follows
\begin{align*}
	\overline{\bal{u=0}\, \ddtr \dduchil{u=0}}^{u=0}
	=
	\overline{\bal{u=0}}^{u=0} \overline{\ddtr \dduchil{u=0}}^{u=0}
	+ \overline{(\bal{u=0} - \overline{\bal{u=0}}^{u=0})(\ddtr \dduchil{u=0} - \overline{\ddtr \dduchil{u=0}}^{u=0})}^{u=0}.
\end{align*}
Then $\overline{\bal{u=0}\, \ddtr \dduchil{u=0}}^{u=0} - \bal{u=0}_S\, (\ddtr \dduchil{u=0})_S$ can be decomposed as
\begin{align*}
	&\phantom{=}
	\overline{\bal{u=0}\, \ddtr \dduchil{u=0}}^{u=0} - \bal{u=0}_S\, (\ddtr \dduchil{u=0})_S
\\
	&
	=
	\frac{2}{\br_{u=0}} - \frac{2}{\br_{u=0,S}}
\\
	&
	=
	(\overline{\bal{u=0}}^{u=0} - \bal{u=0}_S )\overline{\ddtr \dduchil{u=0}}^{u=0}
	+ \bal{u=0}_S [ \overline{\ddtr \dduchil{u=0}}^{u=0} - (\ddtr \dduchil{u=0})_S]
\\
	&\phantom{=}
	+ \overline{(\bal{u=0} - \overline{\bal{u=0}}^{u=0})(\ddtr \dduchil{u=0} - \overline{\ddtr \dduchil{u=0}}^{u=0})}^{u=0}
\end{align*}
Therefore by the estimates of $\br_{u=0} - \br_{u=0,S}$, $\ddtr \dduchil{u=0} - (\ddtr \dduchil{u=0})_S$, $\bslashd \ddtr \dduchil{u=0}$ in lemma \ref{lem 7.3} and the estimate of $\bslashd \bal{u=0}$ above, we obtain that
\begin{align*}
	\vert \overline{\bal{u=0}}^{u=0} - \bal{u=0}_S \vert
	\leq
	c(n,p) (\epsilon + \delta_o + \delta_m + \ud_o + \ud_{\uh}).
\end{align*}
\end{enumerate}
\end{proof}

\subsection{Estimate of curvature component on initial leaf}
We estimate the curvature components associated with the foliation on the initial leaf by the estimates obtained in propositions \ref{prop 5.8}, \ref{prop 5.10} and the estimate of $\bal{u=0}$ in proposition \ref{prop 7.2}.
\begin{proposition}\label{prop 7.4}
Under assumption \ref{assum 7.1}, and for $\delta$ sufficiently small depending on $n,p$ such that proposition \ref{prop 7.2} holds, we have that
\begin{align*}
	\bualphal{u=0}
	&:
	\begin{aligned}
		\Vert \bualphal{u=0} \Vert^{n,p}
		\leq
		c(n,p) ( \epsilon + \ud_{\uh}^2),
	\end{aligned}
\\
	\bubetal{u=0}
	&:
	\begin{aligned}
		\Vert \bubetal{u=0} \Vert^{n,p}
		\leq
		\frac{c(n,p)(\epsilon + \ud_{\uh})}{r_0},
	\end{aligned}
\\
	\bsigmal{u=0}
	&:
	\begin{aligned}
		\Vert \bsigmal{u=0} \Vert^{n,p}
		\leq
		\frac{c(n,p)(\epsilon + \delta_o \ud_{\uh})}{r_0^2},
	\end{aligned}
\\
	\brhol{u=0}
	&:
	\begin{aligned}
		\vert \brhol{u=0} - \brhol{u=0}_S \vert
		\leq
		\frac{c(n,p) ( \epsilon + \delta_o + \delta_m)}{r_0^2},
		\quad
		\Vert \bslashd \brhol{u=0} \Vert^{n-1,p}
		\leq
		\frac{c(n,p)(\epsilon + \delta_o + \delta_m \ud_o)}{r_0^2},
	\end{aligned}
\\
	\bbetal{u=0}
	&:
	\begin{aligned}
		\Vert \bbetal{u=0} \Vert^{n,p}
		\leq
		\frac{c(n,p) (\epsilon + \delta_o) }{r_0},
	\end{aligned}
\\
	\balphal{u=0}
	&:
	\begin{aligned}
		\Vert \balphal{u=0} \Vert^{n,p}
		\leq
		c(n,p)(\epsilon + \delta_o^2).
	\end{aligned}
\end{align*}
\end{proposition}
\begin{proof}
The estimate of $\brhol{u=0}$ is already obtained in lemma \ref{lem 7.3}. Then the estimates of other curvature components follow from the estimates in propositions \ref{prop 5.8}, \ref{prop 5.10} and proposition \ref{prop 7.2} and the formulae
\begin{align*}
	&
	\bualphal{u=0} = \bal{u=0}^2 \cdot \ddualphal{u=0},
	\quad
	\bubetal{u=0} = \bal{u=0} \cdot \ddubetal{u=0},
	\quad
	\bsigmal{u=0} = \ddsigmal{u=0},
\\
	&
	\bbetal{u=0} = \bal{u=0}^{-1} \cdot \ddbetal{u=0},
	\quad
	\balphal{u=0} = \bal{u=0}^{-2} \cdot \ddalphal{u=0}.
\end{align*}
\end{proof}

\subsection{Estimate of connection coefficient on initial leaf}
\begin{figure}
\centering
\begin{tikzpicture}
\draw (-3,5) rectangle (-1,4); 
\node at (-2,4.5) {$\ufl{s=0}, \fl{u=0}$};
\draw[->] (-1,4.5) -- (-0.5,4.5) 
-- (0,4.5); 
\draw (0,5) rectangle (5+1,4); 
\node at (2.5+.5,4.5) {$\ufl{u=0}, (\dslashd \uh)|_{\bSigma_{u=0}}, (\circnabla^2 \uh)|_{\bSigma_{u=0}},\bslashgl{u=0}$}; 
\draw (-3.2,5.2) rectangle (5.2+1,3.8);
\draw[->] (-2.5-.2, 3.8) -- (-2.5-.2,3.5) 
-- (-2.5-.2, 3); 
\draw (-3-.2,2) rectangle (-2-.2,3); 
\node at (-2.5-.2,2.5) {$\brhol{u=0}$};
\%
\draw[->] (2-.2, 3.8) -- (2-.2,3.5) 
-- (2-.2, 3);
\draw (1-.2,2) rectangle (3+.1-.2,3);
\node at (2+.05-.2,2.5) {$(\hatdduchil{u=0}, \hatddchil{u=0}' )$};
\draw[->] (-2-.2,2.5) -- (-1.5-.2,2.5) 
-- (-1-.2,2.5); 
\draw[->] (1-.2,2.5) -- (0.5-.2,2.5) 
-- (0-.2,2.5);
\draw (-1-.2,2) rectangle (0-.2,3);
\node at (-0.5-.2,2.5) {$\bmul{u=0}$};
\draw[->] (4.5-.2, 3.8) -- (4.5-.2,3.5) 
-- (4.5-.2, 3);
\draw (4-.55-.3,2) rectangle (5-.3,3);
\node at (4.5-.25-.3,2.5) {$\bslashdiv \ddeta|_{\bSigma_{u=0}}$};
\draw[->] (5.2+1, 4.5) -- (6.5,4.5) 
-- (6.5,3);
\draw (5.5-.4-.3,2) rectangle (7.5-.1-.2,3);
\node at (6.5-.2-.3,2.5) {$\ddtr \dduchil{u=0}, \ddtr \ddchil{u=0}'$};
\draw[->] (6.5-.2,2) --  (6.5-.2,1.5) -- (4.5-.2,1.5) -- (4.5-.2,1);
\draw[->] (-3.2, 4.5) -- (-4.5,4.5) 
-- (-4.5,3);
\draw (-5+.1,2) rectangle (-4+.1,3);
\node at (-4.5+.1,2.5) {$\bsigmal{u=0}$};
\draw[->] (-3.2, 4.5) -- (-5.5,4.5) -- (-5.5,3);
\draw (-7-.3+.1,2) rectangle (-5.2+.1,3);
\node at (-6.1-.15+.1,2.5) {$\hatddchil{u=0}' \wedge \hatdduchil{u=0} $};
\draw (-7.4+.1,1.8) rectangle (-3.8+.1,3.2);
\draw (-3.2-.2,1.8) rectangle (3.2-.2,3.2);
\draw[->] (-0.5-.2, 2-.2) -- (-0.5-.2, 1.5) -- (4.5-.2,1.5) -- (4.5-.2,1); 
\draw[->] (4.5-.2, 2) -- (4.5-.2, 1.5) -- (4.5-.2,1.5) 
-- (4.5-.2,1); 
\draw (5.2,0) rectangle (3.8-.2,1); 
\node at (4.5-.1,0.5) {$\slashd \log \bal{u=0}$};
\draw[->] (3.8-.2, 0.5) -- (3, 0.5); 
\draw (2,0) rectangle (3,1); 
\node at (2.5,0.5) {$\bal{u=0}$};
\draw (1.8,-0.2) rectangle (5.4, 1.2);
\draw[->] (-5.5,1.8) -- (-5.5,1.5) 
-- (-5.5,1);
\draw (-6.5,0) rectangle (-4.5,1);
\node at (-5.5,0.5) {$\bslashcurl \btal{u=0}$};
\draw[->] (-3,1.8) -- (-3,1.5) 
-- (-3,1);
\draw (-3.5,0) rectangle (-1.5,1);
\node at (-2.5,0.5) {$\bslashdiv \btal{u=0}$};
\draw (-6.7,-0.2) rectangle (-1.3,1.2);
\draw[->] (-4,-0.2) -- (-4,-0.5) 
-- (-4,-1);
\draw (-4.5,-2) rectangle (-3.5,-1);
\node at (-4,-1.5) {$\btal{u=0}$};
\draw[->] (1.8,0.5) -- (0,0.5) 
-- (0,-1);
\draw (-1.5,-2) rectangle (1.5,-1);
\node at (0,-1.5) {$\btr \buchil{u=0}, \btr \bchil{u=0}'$};
\draw[->] (4,-0.2) -- (4,-0.5) 
-- (4,-1);
\draw (3,-2) rectangle (5,-1);
\node at (4,-1.5) {$\bubetal{u=0}, \bbetal{u=0}$};
\draw (-4.7,-2.2) rectangle (5.2,-0.8);
\draw[->] (-3,-2.2) -- (-3,-2.5) 
-- (-3,-3);
\draw (-4,-4) rectangle (-2,-3);
\node at (-3,-3.5) {$\hatbuchil{u=0}, \hatbchil{u=0}$};
\draw[->] (-2,-3.5) -- (-1.5,-3.5) 
-- (-1,-3.5);
\draw[->] (0,-2.2) -- (0,-2.5) 
-- (0,-3);
\draw (-1,-4) rectangle (1,-3);
\node at (0,-3.5) {$\slashd \buomegal{u=0}$};
\draw[->] (1,-3.5) -- (2,-3.5);
\draw (2,-4) rectangle (4,-3);
\node at (3,-3.5) {$\overline{\buomegal{u=0}}^{\subbslashg}$};
\draw (-1.2,-4.2) rectangle (4.2,-2.8);
\end{tikzpicture}
\caption{Logic order to obtain estimates on the initial leaf.}
\label{fig 3}
\end{figure}
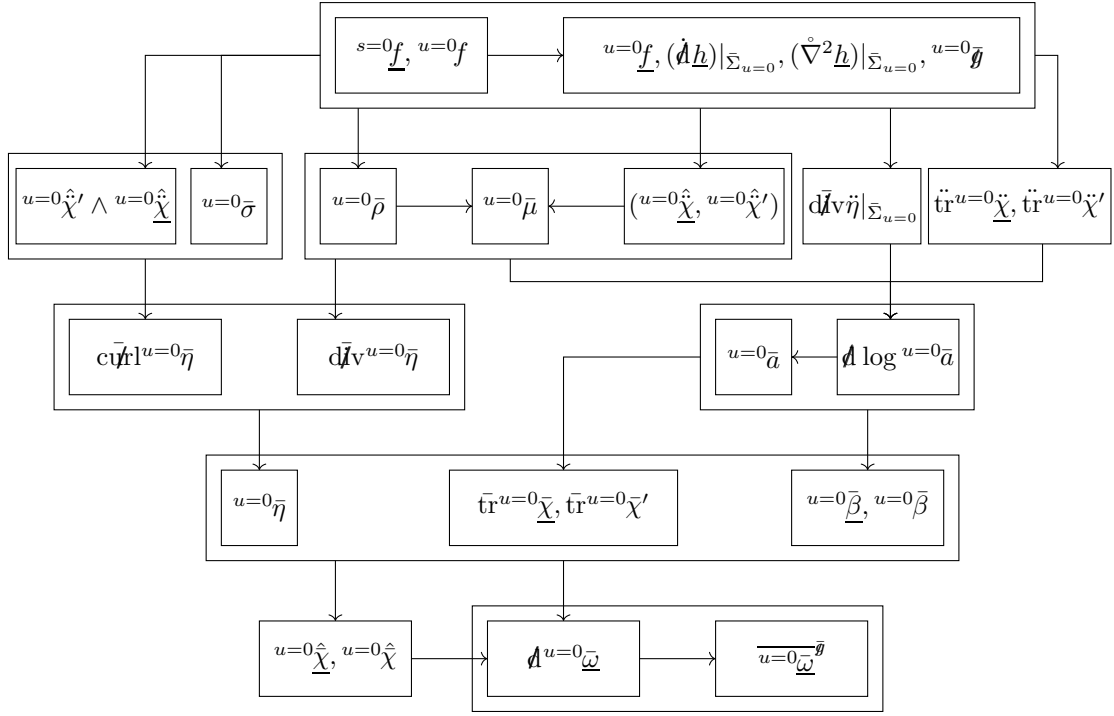

In this section, we estimate the connection coefficients associated with the foliation on the initial leaf. The figure \ref{fig 3} illustrates the logic order to obtain the estimates.

\subsubsection{Estimate of initial null expansions $\btr \buchil{u=0}$, $\btr \bchil{u=0}'$}
\begin{proposition}\label{prop 7.5}
Under assumption \ref{assum 7.1}, and for $\delta$ sufficiently small depending on $n,p$ such that proposition \ref{prop 7.2} holds, we have that
\begin{align*}
	&
	\btr \buchil{u=0}:
	\left\{
	\begin{aligned}
		&
		\vert \overline{\btr \buchil{u=0}}^{u=0} - (\btr \buchil{u=0})_S \vert, 
		\vert \overline{\btr \buchil{u=0}}^{\circg} - (\btr \buchil{u=0})_S \vert
		\leq
		\frac{c(n,p) ( \epsilon + \delta_o + \delta_m)}{r_0},
	\\
		&
		\Vert \bslashd \btr \buchil{u=0} \Vert^{n-1,p}
		\leq
		\frac{c(n,p) (\epsilon + \delta_o + \ud_o + \ud_{\uh})}{r_0},
	\end{aligned}
	\right.
\\
	&
	\btr \bchil{u=0}':
	\left\{
	\begin{aligned}
		&
		\vert \btr \bchil{u=0}' \vert
		\leq
		\frac{c(n,p) ( \epsilon + \delta_o + \delta_m)}{r_0},
	\\
		&
		\Vert \bslashd \btr \bchil{u=0}' \Vert^{n-1,p}
		\leq
		\frac{c(n,p) (\epsilon + \delta_o + \delta_m \ud_o + \delta_m^2 \ud_{\uh})}{r_0}.
	\end{aligned}
	\right.
\end{align*}
\end{proposition}
\begin{proof}
\begin{enumerate}[label={\textbullet}]
\item
$\pmb{\overline{\btr \buchil{u=0}}^{u=0} - (\btr \buchil{u=0})_S, \overline{\btr \buchil{u=0}}^{\circg} - (\btr \buchil{u=0})_S}$.
By equation \eqref{eqn 6.5}, we have
\begin{align*}
	&
	\overline{\btr \buchil{u=0}}^{u=0} - (\btr \buchil{u=0})_S
	=
	\frac{2}{\br_{u=0}} - \frac{2}{\br_{u=0,S}},
\end{align*}
thus the estimate of $\overline{\btr \buchil{u=0}}^{u=0} - (\btr \buchil{u=0})_S$ follows from the estimate of $\br_{u=0}$ in lemma \ref{lem 7.3}. For $\overline{\btr \buchil{u=0}}^{\circg} - (\btr \buchil{u=0})_S$, we have that
\begin{align*}
	&
	\overline{\btr \buchil{u=0}}^{\circg} - \overline{\btr \buchil{u=0}}^{u=0}
	=
	\frac{\int \btr \buchil{u=0}\ \dvol_{\subbslashgl{u=0}_S}}{4\pi \br_{u=0,S}^2}
	-
	\frac{\int \btr \buchil{u=0}\ \dvol_{\subbslashgl{u=0}}}{4\pi \br_{u=0}^2},
\end{align*}
then from the estimates of $\br_{u=0}$ and $\bslashepsilonl{u=0}$ in lemma \ref{lem 7.3}, the estimate of $\overline{\btr \buchil{u=0}}^{\circg} - (\btr \buchil{u=0})_S$ follows.

\item
$\pmb{\bslashd \btr \buchil{u=0}}$.
The estimate follows from the estimates of $\bal{u=0}$ in proposition \ref{prop 7.2}, $\ddtr \dduchil{u=0}$ in lemma \ref{lem 7.3} and the equation
\begin{align*}
	&
	\btr \buchil{u=0}
	=
	\bal{u=0}\cdot \ddtr \dduchil{u=0}.
\end{align*}

\item
$\pmb{\ddtr \ddchil{u=0}', \btr \bchil{u=0}'}$.
By the estimates of $\lo{\ddtr \ddchilu'}$ and $\hi{\ddtr \ddchilu'}$ in propositions \ref{prop 5.8}, \ref{prop 5.10}, we have
\begin{align*}
	&
	\vert \ddtr \ddchil{u=0}' \vert
	\leq
	\frac{c(n,p) ( \epsilon + \delta_o + \delta_m)}{r_0}.
\end{align*}
Combining with the estimate of $\bal{u=0}$ in proposition \ref{prop 7.2}, the estimate of $\btr \bchil{u=0}'$ follows.

\item
$\pmb{\bslashd \ddtr \ddchil{u=0}'}$. By the decomposition of $\ddtr \ddchil{u=0}'$, we have
\begin{align*}
	\bslashd \ddtr \ddchil{u=0}'
 	&
	=
	\bslashd \lo{\ddtr \ddchil{u=0}'} 
	+ \bslashd\hi{\ddtr \ddchil{u=0}'}
\\
	&
	=
	\bslashd (\tr \chi'_S|_{\bSigma_{u=0}}) 
	- 2 \bslashd (r_S^{-2}|_{\bSigma_{u=0}} \circDelta \fl{u=0})
	+ \bslashd \hi{\ddtr \ddchil{u=0}'}
\end{align*}
By estimates \eqref{eqn 5.15},  we have
\begin{align*}
	&
	\Vert \bslashd (\tr \chi'_S|_{\bSigma_{u=0}}) \Vert^{n,p}
	\leq
	\frac{c(n,p)(\delta_o + \udelta_m \ud_o)}{r_0},
\\
	&
	\Vert \bslashd (r_S^{-2}|_{\bSigma_{u=0}} \circDelta \fl{u=0}) \Vert^{n-1,p}
	\leq
	\frac{c(n,p) \delta_o}{r_0}.
\end{align*}
Combining with the estimate of $\hi{\ddtr \ddchil{u=0}'}$ in proposition \ref{prop 5.10},
\begin{align*}
	&
	\Vert \bslashd \ddtr \ddchil{u=0}' \Vert^{n-1,p}
	\leq
	\frac{c(n,p)(\epsilon + \delta_o + \udelta_m \ud_o)}{r_0}
\end{align*}

\item
$\pmb{\bslashd \btr \bchil{u=0}'}$. The estimate follows from the estimates of $\bal{u=0}$ in proposition \ref{prop 7.2} and $\ddtr \ddchil{u=0}'$, $\bslashd \ddtr \ddchil{u=0}'$ above.
\end{enumerate}
\end{proof}

\subsubsection{Estimate of initial torsion $\btal{u=0}$}
\begin{proposition}\label{prop 7.6}
Under assumption \ref{assum 7.1}, we have that
\begin{align*}
	&
	\Vert \btal{u=0} \Vert^{n+1,p}
	\leq
	c(n,p) (\epsilon + \delta_o + \delta_m \ud_o).
\end{align*}
\end{proposition}
\begin{proof}
We estimate $\btal{u=0}$ by equations \eqref{eqn 6.12}.
\begin{enumerate}[label={\textbullet}]
\item
$\pmb{\bslashdiv \btal{u=0}}$.
By equation \eqref{eqn 6.12},
\begin{align*}
	&
	\bslashdiv \btal{u=0}
	=
	- ( \brhol{u=0} - \overline{\brhol{u=0}}^{u=0}) 
	- \frac{1}{2} [(\hatdduchil{u=0}, \hatddchil{u=0}') - \overline{(\hatdduchil{u=0}, \hatddchil{u=0}')}^{u=0}]
\end{align*}
By the estimates of $\brhol{u=0}$ in lemma \ref{lem 7.3} and $(\hatdduchil{u=0}, \hatddchil{u=0}')$ in the proof of proposition \ref{prop 7.2},
\begin{align*}
	&
	\Vert \bslashdiv \btal{u=0} \Vert^{n,p}
	\leq
	\frac{c(n,p)(\epsilon + \delta_o + \delta_m \ud_o)}{r_0^2}.
\end{align*}

\item
$\pmb{\bslashcurl \btal{u=0}}$.
By the estimates of $\hatdduchil{u=0}$, $\hatddchil{u=0}$ in the proof of proposition \ref{prop 7.2}, and $\bsigmal{u=0}$ in proposition \ref{prop 7.4},
\begin{align*}
	&
	\Vert \bslashcurl \btal{u=0} \Vert^{n,p}
	\leq
	\Vert \frac{1}{2} \hatddchil{u=0}' \wedge \hatdduchil{u=0} \Vert^{n,p}
	+
	\Vert \bsigmal{u=0} \Vert^{n,p}
	\leq
	\frac{c(n,p) (\epsilon + \delta_o \ud_{\uh})}{r_0^2}.
\end{align*}

\item
$\pmb{\btal{u=0}}$. By the estimate of $\bslashgl{u=0}$ in lemma \ref{lem 7.3} and the theory of elliptic equations on the sphere, the estimate of $\btal{u=0}$ follows from the estimates of $\bslashdiv \btal{u=0}$ and $\bslashcurl \btal{u=0}$.
\end{enumerate}
\end{proof}

\subsubsection{Estimate of initial shears $\hatbuchil{u=0}$}
\begin{proposition}\label{prop 7.7}
Under assumption \ref{assum 7.1}, there exists a small positive constant $\delta$ depending on $n,p$ such that
\begin{align*}
	&
	\Vert \hatbuchil{u=0} \Vert^{n,p}
	\leq
	c(n,p) (\epsilon + \delta_o + \ud_o + \ud_{\uh}) r_0.
\end{align*}
\end{proposition}
\begin{proof}
Choose $\delta$ sufficiently small such that proposition \ref{prop 7.2} holds. We further require $\delta$ sufficiently small such that
\begin{align*}
c(n,p) (\epsilon + \delta_o + \delta_m + \ud_o + \ud_{\uh} ) \leq \frac{1}{2}, 
\end{align*}
for all constants $c(n,p)$ appearing in the proof. We estimate $\hatbuchil{u=0}$ by equation \eqref{eqn 6.10}.
\begin{enumerate}[label={\textbullet}]
\item
$\pmb{\bslashdiv \hatbuchil{u=0}}$.
By the estimates of $\bubetal{u=0}$, $\bslashd \btr \buchil{u=0}$, $\btal{u=0}$ in propositions \ref{prop 7.4}, \ref{prop 7.5}, \ref{prop 7.6}, we have
\begin{align*}
	&
	\Vert \bslashdiv \hatbuchil{u=0} \Vert^{n-1,p}
	\leq
	\frac{c(n,p)(\epsilon + \delta_o + \ud_o + \ud_{\uh})}{r_0}
	+
	\frac{c(n,p)(\epsilon + \delta_o + \delta_m \ud_o)}{r_0^2} \Vert \hatbuchil{u=0} \Vert^{n-1,p}.
\end{align*}

\item
$\pmb{\hatbuchil{u=0}}$.
By the estimate of $\bslashgl{u=0}$ in lemma \ref{lem 7.3} and the theory of the elliptic equations on the sphere, we have
\begin{align*}
	\Vert \hatbuchil{u=0} \Vert^{n,p}
	&
	\leq
	c(n,p) r_0^2 \Vert \bslashdiv \hatbuchil{u=0} \Vert^{n-1,p}
\\
	&
	\leq
	c(n,p)(\epsilon + \delta_o + \ud_o + \ud_{\uh})
	+
	c(n,p)(\epsilon + \delta_o + \delta_m \ud_o) \Vert \hatbuchil{u=0} \Vert^{n-1,p}
\\
	&
	\leq
	c(n,p)(\epsilon + \delta_o + \ud_o + \ud_{\uh})r_0
	+
	\frac{1}{2} \Vert \hatbuchil{u=0} \Vert^{n-1,p}.
\\
	\Rightarrow \quad
	\Vert \hatbuchil{u=0} \Vert^{n,p}
	&
	\leq
	c(n,p)(\epsilon + \delta_o + \ud_o + \ud_{\uh})r_0.
\end{align*}
\end{enumerate}
\end{proof}

\subsubsection{Estimate of initial shears $\hatbchil{u=0}'$}
\begin{proposition}\label{prop 7.8}
Under assumption \ref{assum 7.1}, there exists a small positive constant $\delta$ depending on $n,p$ such that
\begin{align*}
	&
	\Vert \hatbchil{u=0}' \Vert^{n,p}
	\leq
	c(n,p) (\epsilon + \delta_o + \delta_m \ud_o + \delta_m^2 \ud_{\uh}) r_0.
\end{align*}
\end{proposition}
\begin{proof}
Choose $\delta$ sufficiently small such that proposition \ref{prop 7.2} holds. We further require $\delta$ sufficiently small such that
\begin{align*}
c(n,p) (\epsilon + \delta_o + \delta_m + \ud_o + \ud_{\uh} ) \leq \frac{1}{2}, 
\end{align*}
for all constants $c(n,p)$ appearing in the proof. We estimate $\hatbchil{u=0}'$ by equation \eqref{eqn 6.11}.
\begin{enumerate}[label={\textbullet}]
\item
$\pmb{\bslashdiv \hatbchil{u=0}'}$.
By the estimates of $\bbetal{u=0}$, $\bslashd \btr \bchil{u=0}'$, $\btal{u=0}$ in propositions \ref{prop 7.4}, \ref{prop 7.5}, \ref{prop 7.6}, we have
\begin{align*}
	&
	\Vert \bslashdiv \hatbchil{u=0}' \Vert^{n-1,p}
	\leq
	\frac{c(n,p) (\epsilon + \delta_o + \delta_m \ud_o + \delta_m^2 \ud_{\uh})}{r_0}
	+
	\frac{c(n,p) (\epsilon + \delta_o + \delta_m \ud_o)}{r_0^2} \Vert \hatbuchil{u=0} \Vert^{n-1,p}.
\end{align*}

\item
$\pmb{\hatbchil{u=0}'}$.
By the estimate of $\bslashgl{u=0}$ in lemma \ref{lem 7.3} and the theory of the elliptic equations on the sphere, we have
\begin{align*}
	\Vert \hatbchil{u=0}' \Vert^{n,p}
	&
	\leq
	c(n,p) r_0^2 \Vert \bslashdiv \hatbchil{u=0}' \Vert^{n-1,p}
\\
	&
	\leq
	c(n,p) (\epsilon + \delta_o + \delta_m \ud_o + \delta_m^2 \ud_{\uh})
	+
	c(n,p) (\epsilon + \delta_o + \delta_m \ud_o) \Vert \hatbchil{u=0}' \Vert^{n-1,p}
\\
	&
	\leq
	c(n,p) (\epsilon + \delta_o + \delta_m \ud_o + \delta_m^2 \ud_{\uh}) r_0
	+
	\frac{1}{2} \Vert \hatbchil{u=0}' \Vert^{n-1,p}.
\\
	\Rightarrow \quad
	\Vert \hatbchil{u=0}' \Vert^{n,p}
	&
	\leq
	c(n,p) (\epsilon + \delta_o + \delta_m \ud_o + \delta_m^2 \ud_{\uh}) r_0.
\end{align*}
\end{enumerate}
\end{proof}

\subsubsection{Estimate of initial acceleration $\buomegal{u=0}$}
\begin{proposition}\label{prop 7.9}
Under assumption \ref{assum 7.1}, and for $\delta$ sufficiently small depending on $n,p$ such that propositions \ref{prop 7.2}, \ref{prop 7.7}, \ref{prop 7.8} hold, we have that
\begin{align*}
	&
	\Vert \bslashd \buomegal{u=0} \Vert^{n,p}
	\leq
	\frac{c(n,p) (\epsilon + \delta_o + \ud_o + \ud_{\uh})}{r_0},
\\
	&
	\vert \overline{\buomegal{u=0}}^{u=0} \vert,
	\vert \overline{\buomegal{u=0}}^{\circg} \vert
	\leq
	\frac{c(n,p) (\epsilon + \delta_o + \ud_o + \ud_{\uh})^2}{r_0}.
\end{align*}
\end{proposition}
\begin{proof}
We estimate $\buomegal{u=0}$ by equation \eqref{eqn 6.13}.
\begin{enumerate}[label=\textbullet]
\item
$\pmb{\bslashDelta \buomegal{u=0}}$. We estimate each term in equation \eqref{eqn 6.13}.
\begin{enumerate}[label=\raisebox{0.1ex}{\scriptsize$\bullet$}]
\item
$\bmul{u=0}\, \btr \buchil{u=0} - \overline{\bmul{u=0}\, \btr \buchil{u=0}}^{u=0}$. By the estimate of $\bslashd \btr \buchil{u=0}$ in proposition \ref{prop 7.5},
\begin{align*}
	&
	\Vert \bmul{u=0}\, \btr \buchil{u=0} - \overline{\bmul{u=0}\, \btr \buchil{u=0}}^{u=0} \Vert^{n,p}
	\leq
	\frac{c(n,p)(\epsilon + \delta_o + \ud_o + \ud_{\uh})}{r_0^3}.
\end{align*}

\item
$\btr \buchil{u=0} | \btal{u=0} |^2 - \overline{\btr \buchil{u=0} | \btal{u=0} |^2}^{u=0}$.
By the estimates of $\btr \buchil{u=0}$ and $\btal{u=0}$ in propositions \ref{prop 7.5}, \ref{prop 7.6}, we have
\begin{align*}
	&
	\Vert \btr \buchil{u=0} | \btal{u=0} |^2 - \overline{\btr \buchil{u=0} | \btal{u=0} |^2}^{u=0} \Vert^{n,p}
	\leq
	\frac{c(n,p)(\epsilon + \delta_o + \delta_m \ud_o)^2}{r_0^3}.
\end{align*}

\item
$\btr \bchil{u=0}' | \hatbuchil{u=0} |^2 - \overline{ \btr \bchil{u=0}' | \hatbuchil{u=0}|^2}^{u=0}$.
By the estimates of $\btr \bchil{u=0}'$ and $\hatbuchil{u=0}$ in propositions \ref{prop 7.5}, \ref{prop 7.7}, we have
\begin{align*}
	&
	\Vert \btr \bchil{u=0}' | \hatbuchil{u=0} |^2 - \overline{ \btr \bchil{u=0}' | \hatbuchil{u=0}|^2}^{u=0} \Vert^{n,p}
	\leq
	\frac{c(n,p)(\epsilon + \delta_o + \ud_o + \ud_{\uh} )^2}{r_0^3}.
\end{align*}

\item
$(\bslashdiv \hatbuchil{u=0}, \btal{u=0})$.
By the estimates of $\btal{u=0}$ and $\hatbuchil{u=0}$ in propositions \ref{prop 7.6}, \ref{prop 7.7}, we have
\begin{align*}
	&
	\Vert (\bslashdiv \hatbuchil{u=0}, \btal{u=0}) \Vert^{n-1,p}
	\leq
	\frac{c(n,p) (\epsilon + \delta_o + \ud_o + \ud_{\uh})(\epsilon + \delta_o + \delta_m \ud_o)}{r_0^3}.
\end{align*}

\item
$(\hatbuchil{u=0}, \bslashnabla \btal{u=0})$.
By the estimates of $\btal{u=0}$ and $\hatbuchil{u=0}$ in propositions \ref{prop 7.6}, \ref{prop 7.7}, we have
\begin{align*}
	&
	\Vert (\hatbuchil{u=0}, \bslashnabla \btal{u=0}) \Vert^{n-1,p}
	\leq
	\frac{c(n,p) (\epsilon + \delta_o + \ud_o + \ud_{\uh})(\epsilon + \delta_o + \delta_m \ud_o)}{r_0^3}.
\end{align*}

\item
$\bslashdiv \bubetal{u=0}$.
By the estimate of $\bubetal{u=0}$ in proposition \ref{prop 7.4}, we have
\begin{align*}
	&
	\Vert \bslashdiv \bubetal{u=0} \Vert^{n-1,p}
	\leq
	\frac{c(n,p) (\epsilon + \delta_o)}{r_0^3}.
\end{align*}
\end{enumerate}
Combining the above estimates, we obtain that
\begin{align*}
	&
	\Vert \bslashDelta \buomegal{u=0} \Vert^{n-1,p}
	\leq
	\frac{c(n,p)(\epsilon + \delta_o + \ud_o + \ud_{\uh})}{r_0^3}.
\end{align*}

\item
$\pmb{\bslashd \buomegal{u=0}}$.
By the estimate of $\bslashgl{u=0}$ in lemma \ref{lem 7.3}, and the theory of elliptic equations on the sphere, we have
\begin{align*}
	&
	\Vert \bslashd \buomegal{u=0} \Vert^{n,p}
	\leq
	\frac{c(n,p)(\epsilon + \delta_o + \ud_o + \ud_{\uh})}{r_0}.
\end{align*}

\item
$\pmb{\overline{\buomegal{u=0}}^{u=0}}$.
By equation \eqref{eqn 6.13}, the estimates of $\bslashd \buomegal{u=0}$ above, $\bslashd \btr \buchil{u=0}$, $\hatbuchil{u=0}$ in propositions \ref{prop 7.5}, \ref{prop 7.7}, we have
\begin{align*}
	&
	\vert \overline{\buomegal{u=0}}^{u=0} \vert
	\leq
	\frac{c(n,p)(\epsilon + \delta_o + \ud_o + \ud_{\uh})^2}{r_0}.
\end{align*}

\item
$\pmb{\overline{\buomegal{u=0}}^{\circg}}$.
By the estimate of $\circnabla \bslashepsilonl{u=0}$ in lemma \ref{lem 7.3}, we have
\begin{align*}
	\vert \bslashepsilonl{u=0} - \br_{u=0}^2 \circepsilon \vert
	\leq
	c(n,p) (\epsilon + \delta_o + \delta_m \ud_o) r_0^2.
\end{align*}
Then
\begin{align*}
	\vert \overline{\buomegal{u=0}}^{\circg} \vert
	&
	\leq
	\vert \overline{\buomegal{u=0}}^{u=0} \vert 
	+ \vert \overline{\buomegal{u=0}}^{\circg} - \overline{\buomegal{u=0}}^{u=0} \vert
	\leq
	\frac{c(n,p)(\epsilon + \delta_o + \ud_o + \ud_{\uh})^2}{r_0}.
\end{align*}
\end{enumerate}
\end{proof}

\subsection{Estimate of mass aspect function on initial leaf}
\begin{proposition}\label{prop 7.10}
Under assumption \ref{assum 7.1}, we have that
\begin{align*}
	\vert \bmul{u=0} - \bmul{u=0}_S \vert
	\leq
	\frac{c(n,p) (\epsilon + \delta_o + \delta_m)}{r_0^2}.
\end{align*}
\end{proposition}
\begin{proof}
By the estimate of $\brhol{u=0} - \brhol{u=0}_S$ in lemma \ref{lem 7.3} and the estimates of $\hatdduchil{u=0}$, $\hatddchil{u=0}'$ in propositions \ref{prop 5.8}, \ref{prop 5.10}, we have that
\begin{align*}
	\vert \brhol{u=0} - \brhol{u=0}_S \vert
	\leq
	\frac{c(n,p)(\epsilon + \delta_o + \delta_m)}{r_0^2},
	\quad
	\vert (\hatdduchil{u=0}, \hatddchil{u=0}') \vert
	\leq
	\frac{c(n,p)(\epsilon + \delta_o)(\epsilon + \ud_{\uh})}{r_0^2}.
\end{align*}
Therefore on the initial leaf of the constant mass aspect function foliation,
\begin{align*}
	\bmul{u=0} - \bmul{u=0}_S
	=
	- (\overline{\brhol{u=0}}^{u=0} - \brhol{u=0}_S)
	- \frac{1}{2} \overline{(\hatdduchil{u=0}, \hatddchil{u=0}')}^{u=0}.
\end{align*}
and the proposition follows.
\end{proof}

\section{Global existence and geometry of constant mass aspect function foliation}\label{sec 8}
In this section, we prove the global existence of the constant mass aspect function foliation on a nearly spherically symmetric incoming null hypersurface in $(M,g)$ and study the geometry of the foliation. We shall obtain the estimates of the parameterisation functions and the geometric quantities of the foliation by the bootstrap argument.

\subsection{Bootstrap assumption}

We introduce the bootstrap assumptions on the following quantities associated with the constant mass aspect function foliation:
\begin{enumerate}[label=$\bullet$]
\item
the parameterisation function $\flu$;

\item
the lapse function $\balu$;

\item
the area radius $\br_u$;

\item
the intrinsic metric $\bslashglu$;

\item
the connection coefficients $\btr \buchilu$, $\hatbuchilu$, $\btr \bchilu'$, $\hatbchilu'$, $\btalu$, $\buomegalu$;

\item
the mass aspect function $\bmulu$;

\item
the curvature components $\bubetalu$, $\brholu$, $\bsigmalu$, $\bbetalu$.

\end{enumerate}

\subsubsection{Regularities of quantities in bootstrap assumption}
We list the Sobolev norm of each term in the bootstrap assumptions for which we expect to prove the estimates.
\begin{enumerate}[label=$\bullet$]
\item
The parameterisation function $\flu$: $\vert \flu - \flu_S \vert$, $\Vert \bslashd \flu \Vert^{n+1,p}$.

\item
The lapse function $\balu$: $\vert \balu - \balu_S \vert$, $\Vert \bslashd \balu \Vert^{n,p}$.

\item
The area radius $\br_u$: $\vert \br_u - \br_{u,S} \vert$. 

\item
The intrinsic metric $\bslashglu$: $\vert \bslashglu - \bslashglu_S \vert$, $\Vert \bcircnabla \bslashglu \Vert^{n,p}$.

\item
The connection coefficients $\btr \buchilu$, $\hatbuchilu$, $\btr \bchilu'$, $\hatbchilu'$, $\btalu$, $\buomegalu$:
\begin{align*}
	&
	\vert \overline{\btr \buchilu}^{u} - \frac{2}{\br_u} \vert,
	&&
	\Vert \bslashd \btr \buchilu \Vert^{n-1,p},
	&&
	\Vert \hatbuchilu \Vert^{n,p},
\\
	&
	\vert \overline{\btr \bchilu'}^{u} - \btr \bchilu'_S \vert,
	&&
	\Vert \bslashd \btr \bchilu' \Vert^{n-1,p},	
	&&
	\Vert \hatbchilu' \Vert^{n,p},
\\
	&
	\Vert \btalu \Vert^{n+1,p},
	&&
	\vert \overline{\buomegalu}^u \vert,
	&&
	\Vert \bslashd \buomegalu \Vert^{n,p}.
\end{align*}

\item
The mass aspect function $\bmulu$: $\bmulu - \bmulu_S$.

\item
The curvature components $\bubetalu$, $\brholu$, $\bsigmalu$, $\bbetalu$:
\begin{align*}
	\Vert \bubetalu \Vert^{n,p},
	\quad
	\vert \brholu - \brholu_S \vert,
	\quad
	\Vert \bslashd \brholu \Vert^{n-1,p},
	\quad
	\Vert \bsigmalu \Vert^{n,p},
	\quad
	\Vert \bbetalu \Vert^{n,p}.
\end{align*}
\end{enumerate}
All the above Sobolev norms are consistent with the estimates of the quantities on the initial leaf of the foliation.

\subsubsection{On regularity of parameterisation function $\flu$}\label{sec 8.1.2}
In order to prove the estimate of $\flu$ by the bootstrap argument, equation \eqref{eqn 6.3}
\begin{align}
	\buL \flu = \balu
	\tag{\ref{eqn 6.3}\ensuremath{'}}
\end{align}
is not sufficient, since by integrating it, one can only obtain the Sobolev norm $\Vert \bslashd \flu \Vert^{n,p}$ at most. In order to recover the desired regularity of $\flu$, we shall apply the following elliptic equation of $\flu$, which follows from the formula \eqref{eqn 5.12a} of $\ddtr \ddchilu'$
\begin{align}
	\balu^{-1} \btr \bchilu' 
	= 
	\ddtr \ddchilu' 
	= 
	- 2 (r_S|_{\bSigma_u})^{-2} \circDelta \flu 
	+ \tr \chi'_S|_{\bSigma_u}
	+ \hi{\ddtr \ddchilu'}.
	\label{eqn 8.1}
\end{align}
The Sobolev norm $\Vert \balu^{-1} \btr \bchilu' \Vert^{n,p}$ of the left hand side in the above equation can be estimated by the bootstrap assumption, thus roughly speaking the elliptic theory on the sphere will give the estimate of $\Vert \bslashd \flu \Vert^{n+1,p}$.

The actual approach to obtain the estimate shall also take the regularity of $\hi{\ddchilu'}$ into accounts, since it also a critical term, whose estimate relies on the top order Sobolev norm $\Vert \bslashd \flu \Vert^{n+1,p}$, through the terms $\slashnabla^2 \uh$, $\dslashnabla^2 \flu$, $\dslashnabla^2 \flu - \circnabla^2 \flu$ with the critical regularities.

\subsubsection{Statement of bootstrap assumption}
We state the bootstrap assumption on the estimates of the quantities associated with the foliation in the following.
\begin{assumption}\label{assum 8.1}
The constant mass aspect function foliation $\{ \bSigma_u \}$ exists for $u \in [0,u_a]$ and satisfies the following estimates.
\begin{enumerate}[label=\textbullet]
\item
The estimates of the parameterisation functions.
\begin{enumerate}[label=\raisebox{0.1ex}{\scriptsize$\bullet$}]
\item
The estimates of the parameterisation function $\uflu$ and $\dslashd \uh$, $\circnabla^2 \uh$: estimates \eqref{eqn 5.15}
\begin{align}
\begin{aligned}
	&
	\Vert \ddslashd \uflu \Vert^{n,p} 
	\leq 
	\ud_o r_0,
	\quad
	\vert \overline{\uflu}^{\circg} - \overline{\ufl{s=0}}^{\circg} \vert
	\leq 
	\uslashd_m r_0,
	\quad
	\vert \overline{\uflu}^{\circg} \vert \leq \ud_m r_0,
\\
	&
	\Vert (\dslashd \uh)|_{\bSigma_u} \Vert^{n+1,p},
	\Vert (\circnabla^2 \uh)|_{\bSigma_u} \Vert^{n,p}
	\leq
	\ud_{\uh} r_0.
\end{aligned}
\tag{\ref{eqn 5.15}\ensuremath{'}}
\end{align}

\item
The estimate of the parameterisation function $\flu$:
\begin{align*}
	&
	\vert \flu - \flu_S \vert
	\leq
	c(n,p,\flu) ( \epsilon + \delta_o + \delta_m + \ud_{o,\uh} + \uslashd_m ) \br_u,
\\
	&
	\Vert \bslashd \flu \Vert^{n+1,p}
	\leq
	c(n,p, \bslashd \flu) ( \epsilon + \delta_o + \ud_{o,\uh} ) \br_u.
\end{align*}
\end{enumerate}

\item The estimate of the metric components.
\begin{enumerate}[label=\raisebox{0.1ex}{\scriptsize$\bullet$}]
\item
The estimate of the lapse function $\balu$:
\begin{align*}
	&
	\vert \balu - \balu_S \vert
	\leq
	c(n,p, \balu) ( \epsilon + \delta_o + \delta_m + \ud_{o,\uh} + \uslashd_m ),
\\
	&
	\Vert \bslashd \balu \Vert^{n,p}
	\leq
	c(n,p, \bslashd \balu) ( \epsilon + \delta_o + \ud_{o,\uh} ).
\end{align*}

\item
The estimate of the area radius $\br_u$:
\begin{align*}
	&
	\vert \br_u - \br_{u,S} \vert
	\leq
	c(n,p,\br_u)(\epsilon + \delta_o + \delta_m) r_0.
\end{align*}

\item
The estimate of the intrinsic metric $\bslashglu$:
\begin{align*}
	&
	\vert \bslashglu - \bslashglu_S \vert
	\leq
	c(n,p,\bslashglu)(\epsilon + \delta_o + \delta_m + \ud_{o,\uh} + \uslashd_m) \br_u^2,
\\
	&
	\Vert \bcircnabla \bslashglu \Vert^{n,p}
	\leq
	c(n,p,\bcircnabla \bslashglu) ( \epsilon + \delta_o + \ud_{o,\uh}) \br_u^2.
\end{align*}

\end{enumerate}

\item
The estimates of the connection coefficients.
\begin{enumerate}[label=\raisebox{0.1ex}{\scriptsize$\bullet$}]
\item
The estimate of the null expansion $\btr \buchilu$:
\begin{align*}
	&
	\vert \overline{\btr \buchilu}^{u} - (\btr \buchilu)_S \vert
	\leq
	\frac{c(n,p,\overline{\btr \buchilu}^{u}) ( \epsilon + \delta_o + \delta_m) r_0}{\br_u^2}.
\\
	&
	\Vert \bslashd \btr \buchilu \Vert^{n-1,p}
	\leq
	\frac{c(n,p,\bslashd \btr \buchilu)(\epsilon + \delta_o + \ud_{o,\uh}) r_0}{\br_u^2}
\end{align*}

\item
The estimate of the shear $\hatbuchilu$:
\begin{align*}
	&
	\Vert \hatbuchilu \Vert^{n,p}
	\leq
	c(n,p, \hatbuchilu) (\epsilon + \delta_o + \ud_{o,\uh}) r_0.
\end{align*}

\item
The estimate of the null expansion $\btr \bchilu'$:
\begin{align*}
	&
	\vert \overline{\btr \bchilu'}^{u} - (\btr \bchilu')_S \vert
	\leq
	\frac{c(n,p,\overline{\btr \bchilu'}^{u}) ( \epsilon + \delta_o + \delta_m)}{\br_u}
	+ \frac{c(n,p,\overline{\btr \bchilu'}^{u}) \ud_{o,\uh} u}{\br_u^2},
\\
	&
	\Vert \bslashd \btr \bchilu' \Vert^{n-1,p}
	\leq
	\frac{c(n,p,\bslashd \btr \bchilu') ( \epsilon + \delta_o + \ud_{o,\uh})}{\br_u}.
\end{align*}

\item
The estimate of the shear $\hatbchilu'$:
\begin{align*}
	&
	\Vert \hatbchilu' \Vert^{n,p}
	\leq
	c(n,p,\hatbchilu')(\epsilon + \delta_o + \ud_{o,\uh}) \br_u.
\end{align*}

\item
The estimate of the torsion $\btalu$:
\begin{align*}
	&
	\Vert \btalu \Vert^{n+1,p}
	\leq
	\frac{c(n,p,\btalu)(\epsilon + \delta_o + \ud_{o,\uh})r_0}{\br_u}.
\end{align*}

\item
The estimate of the acceleration $\buomegalu$:
\begin{align*}
	&
	\vert \overline{\buomegalu}^{u} \vert
	\leq
	\frac{c(n,p,\overline{\buomegalu}^{u})(\epsilon + \delta_o + \ud_{o,\uh})^2 r_0^2}{\br_u^3},
\\
	&
	\Vert \bslashd \buomegalu \Vert^{n,p}
	\leq
	\frac{c(n,p,\bslashd \buomegalu) \epsilon r_0^{\frac{3}{2}}}{\br_u^{\frac{5}{2}}}
	+
	\frac{c(n,p,\bslashd \buomegalu) (\delta_o + \ud_{o,\uh}) r_0^2}{\br_u^3}.
\end{align*}
\end{enumerate}

\item
The estimate of the mass aspect function.
\begin{align*}
	\vert \bmulu - \bmulu_S \vert
	\leq
	\frac{c(n,p,\bmulu)(\epsilon + \delta_o + \delta_m) r_0}{\br_u^3}
	+ \frac{c(n,p,\bmulu) \ud_{o,\uh} u r_0}{\br_u^4}.
\end{align*}

\item
The estimates of the curvature components.
\begin{enumerate}[label=\raisebox{0.1ex}{\scriptsize$\bullet$}]
\item
The estimate of $\bubetalu$:
\begin{align*}
	&
	\Vert \bubetalu \Vert^{n,p}
	\leq
	\frac{c(n,p,\bubetalu) \epsilon r_0^{\frac{3}{2}}}{\br_u^{\frac{5}{2}}}
	+
	\frac{c(n,p,\bubetalu) \ud_{o,\uh} r_0^2}{\br_u^3}.
\end{align*}

\item
The estimate of $\brholu$:
\begin{align*}
	&
	\vert \brholu - \brholu_S \vert
	\leq
	\frac{c(n,p,\brholu) (\epsilon + \delta_o + \delta_m + \ud_{o,\uh} + \uslashd_m)r_0}{\br_u^3}.
\\
	&
	\Vert \bslashd \brholu \Vert^{n-1,p}
	\leq
	\frac{c(n,p,\bslashd \brholu) ( \epsilon + \delta_o + \ud_{o,\uh}) r_0}{\br_u^3}.
\end{align*}

\item
The estimate of $\bsigmalu$:
\begin{align*}
	&
	\Vert \bsigmalu \Vert^{n,p}
	\leq
	\frac{c(n,p, \bsigmalu)(\epsilon + \delta_o \ud_{o,\uh} + \ud_{o,\uh}^2) r_0}{\br_u^3}.
\end{align*}

\item
The estimate of $\bbetalu$:
\begin{align*}
	&
	\Vert \bbetalu \Vert^{n,p}
	\leq
	\frac{c(n,p,\bbetalu) \epsilon}{\br_u}
	 +
	 \frac{c(n,p,\bbetalu)(\delta_o + \ud_{o,\uh}) r_0}{\br_u^2}.
\end{align*}

\end{enumerate}

\end{enumerate}
\end{assumption}

\subsection{Statement of main theorem and strategy of proof}\label{sec 8.2}
We state the main theorem on the global existence of the constant mass aspect function foliation and its geometry in the following.

\begin{theorem}\label{thm 8.2}
Let $\bSigma_{u=0}$ be the initial leaf of the constant mass aspect function foliation in the incoming null hypersurface $\ucalH$. Let $(\ufl{s=0}, \fl{u=0})$ be the parameterisation of $\bSigma_{u=0}$. Suppose that 
\begin{align*}
	&
	\Vert \slashd \ufl{s=0} \Vert^{n+1,p} \leq \udelta_o r_0,
	\quad
	\overline{\ufl{s=0}}^{\circg} = \ous,
	\quad
	\vert \ous \vert \leq \udelta_m r_0 \leq \kappa r_0 \leq 0.1 r_0,
\\
	&
	\Vert \slashd \fl{u=0} \Vert^{n+1,p} \leq \delta_o (r_0+\os_{u=0}),
	\quad
	\overline{\fl{u=0}}^{\circg} = \os_{u=0},
	\quad
	\vert \os_{u=0} \vert \leq \delta_m r_0,
\end{align*}
where $n\geq 1+n_p$, $p >1$. Moreover assume that $N \geq n+4+n_p$ in definition \ref{def 2.4}, i.e. the metric $g$ is at least $[(n+4+n_p)+2]$-th order differentiable.

There exists a positive constant $\delta$ depending on $n,p$ such that if $\epsilon + \udelta_o + \udelta_m + \delta_o + \delta_m \leq \delta$, then the constant mass aspect function foliation $\{ \bSigma_u \}$ exists for $u\in [0, +\infty)$, and there exist constants $c(n,p,\cdot)$ such that the estimates in assumption \ref{assum 8.1} hold.
\end{theorem}
\begin{remark}\label{rem 8.3}
Comparing with the assumption \ref{assum 7.1} of the estimates on the geometry of the initial leaf $\bSigma_{u=0}$ in section \ref{sec 7}, the assumption of theorem \ref{thm 8.2} on $\udelta_m$ is stronger which requires $\udelta_m \leq \delta$, while assumption \ref{assum 7.1} requires $\epsilon \udelta_m \leq \delta$.
\end{remark}

We shall prove theorem \ref{thm 8.2} by the bootstrap argument. We briefly explain the strategy of the proof in the following.
\begin{enumerate}[label=\alph*.]
\item
Define the interval $I_a$ of $u_a$ where bootstrap assumption \ref{assum 8.1} holds, i.e.
\begin{align*}
	I_a
	=
	\{ u _a: \text{ bootstrap assumption \ref{assum 8.1} holds on $[0,u_a]$} \}.
\end{align*}
The goal is to prove that there exist $\delta>0$ and constants $c(n,p,\cdot)$ such that $I_a= [0, +\infty)$, equivalently $I_a$ is nonempty, closed and open.

\item
Nonemptiness of $I_a$. By choosing $\delta>0$ such that estimates on the geometry of the initial leaf in section \ref{sec 7} hold and $c(n,p,\cdot)$ greater than the constants in these estimates, we can guarantee that $0\in I_a$, thus $I_a$ is nonempty.

\item
Closedness of $I_a$. Suppose that $u_a^{\sup} = \sup I_a < +\infty$, we need to show that $u_a^{\sup} \in I_a$. This is straightforward by letting $u$ passing to the limit $u_a^{\sup}$ and conclude that all the estimates in bootstrap assumption still hold, since all the inequalities in the estimates are closed.

\item
Openness of $I_a$. This is the essential difficult step in the bootstrap argument. The proof of the openness is done in two steps. We shall find appropriate $\delta>0$ and constants $c(n,p,\cdot)$ in the process such that the following two assertions can be proved.
\begin{enumerate}[label=\roman*.]
\item
Improvement of the estimates in the bootstrap assumption to strict inequalities. Assume that $[0,u_a] \subset I_a$, we can make use of the bootstrap assumption \ref{assum 8.1}, the basic equations \eqref{eqn 6.3}-\eqref{eqn 6.13} of the foliation and the elliptic equation \eqref{eqn 8.1} of $\flu$ to improve the estimates in $[0,u_a]$ to strict inequalities.

\item
Extension of the interval $[0, u_a] \subset I_a$. Given the first assertion proved, the estimates of the bootstrap assumption is improved to strict inequalities at the end point $u_a$, we can show that the constant mass aspect function foliation can be extended further beyond $u_a$ with the estimates in the bootstrap assumption hold in the extension. This assertion implies that $I_a$ must be open.
\end{enumerate}
\end{enumerate}

\subsection{Notation of constants in bootstrap arguments}

In the following, we improve the estimates in the bootstrap assumption to strict inequalities, i.e. step i. in part d. of the strategy to prove theorem \ref{thm 8.2}. Assume that the bootstrap assumption holds on the interval $[0,u_a]$. 

For the estimates in this section, we introduce the following notation for constants depending on $c(n,p,\cdot)$ in the bootstrap assumption in a polynomial way:
\begin{align*}
	C(n,p,a_1,a_2, \cdots, a_k)
	=
	P(c(n,p,a_1), c(n,p,a_2), \cdots, c(n,p,a_k)),
\end{align*}
where $P$ is a polynomial of $k$ variables. In the following, we abuse this notation to denote constants depending on $c(n,p, \cdots)$ with different polynomials by the same $C(n,p,\cdots)$, since main concerns in the application are which constants $c(n,p,\cdot)$ enter the dependence instead of the precise form of the polynomials. However we shall  that emphasise that all $C(n,p,\cdots)$ in the estimates are effective, whose expressions can be written down explicitly. 

We shall introduce the following basic condition on the constant $\delta$ in theorem \ref{thm 8.2} for the improvement of the estimates in the bootstrap assumption throughout this section that for all constants $C(n,p,\cdot)$ appearing in this section,
\begin{align}
	\delta \cdot C(n,p,\cdot) \leq \frac{1}{2}.
	\label{eqn 8.2}
\end{align}

For the sake of the simplification, we introduce the notation $d_o$ to denote the bound of $\bslashd \flu$ in the bootstrap assumption, i.e.
\begin{align}
	d_o = c(n,p, \bslashd \flu) ( \epsilon + \delta_o + \ud_{o,\uh} ),
	\quad
	\Vert \bslashd \flu \Vert^{n+1,p} \leq d_o \br_u.
	\label{eqn 8.3}
\end{align}

\subsection{Estimates of parameterisation functions}\label{sec 8.4}

\subsubsection{Estimates of parameterisation functions $\uflu$ and $\dslashd \uh$, $\circnabla^2 \uh$}
Their estimates follow from propositions \ref{prop 3.3}, \ref{prop 5.2}, \ref{prop 5.5}. By the bootstrap assumptions of $\flu$ and $\br_u$, we have that 
\begin{align*}
	\overline{\flu}^{\circg} 
	&=
	\flu_S + \overline{\flu - \flu_S}^{\circg}
\\
	&\geq
	\flu_S - c(n,p,\flu) (\epsilon + \delta_o + \delta_m + \ud_{o,\uh}) \br_u
\\
	&\geq
	[1- c(n,p,\flu) (\epsilon + \delta_o + \delta_m + \ud_{o,\uh})] \br_u - 1.1 r_0
\\
	&\geq
	[1- c(n,p,\flu) (\epsilon + \delta_o + \delta_m + \ud_{o,\uh})] [ \br_{u,S} - c(n,p,\br_u) (\epsilon + \delta_o + \delta_m) r_0]  
	- 1.1 r_0,
\end{align*}
and
\begin{align*}
	\frac{\Vert \bslashd \flu \Vert^{n+1,p}}{r_0 + \overline{\flu}^{\circg}}
	\leq
	\frac{c(n,p,\bslashd \flu)(\epsilon + \delta_o + \ud_{o,\uh}) \br_u}{[1- c(n,p,\flu) (\epsilon + \delta_o + \delta_m + \ud_{o,\uh})] \br_u - 0.1 r_0}
\end{align*}

Introduce the constant $\bar{\delta}$ by
\begin{align}
	\bar{\delta} = \min \{ \delta_{\ref{prop 3.3}}, \delta_{\ref{prop 5.2}}, \delta_{\ref{prop 5.5}} \}
	\label{eqn 8.4}
\end{align}
where $\delta_{\ref{prop 3.3}}, \delta_{\ref{prop 5.2}}, \delta_{\ref{prop 5.5}}$ denote the constants $\delta$ in propositions \ref{prop 3.3}, \ref{prop 5.2}, \ref{prop 5.5} respectively.

In order to improve the estimates of $\uflu$, $(\dslashd \uh)|_{\bSigma_u}$, $(\circnabla^2 \uh)|_{\bSigma_u}$ to strict inequalities, we simply require that
\begin{align*}
	\begin{aligned}
		&
		[1- c(n,p,\flu) (\epsilon + \delta_o + \delta_m + \ud_{o,\uh})] \br_u - 1.1 r_0
		> - \frac{r_0}{2},
	\\
		&
		\frac{c(n,p,\bslashd \flu)(\epsilon + \delta_o + \ud_{o,\uh}) \br_u}{[1- c(n,p,\flu) (\epsilon +\delta_o + \delta_m + \ud_{o,\uh})] \br_u - 0.1 r_0}
		<
		\bar{\delta},
	\end{aligned}
\end{align*}
where a sufficient condition is the following
\begin{align}
	\left\{
	\begin{aligned}
		&
		1- c(n,p,\flu) (\epsilon + \delta_o + \delta_m + \ud_{o,\uh}) >0.9,
	\\
		&
		1 - c(n,p,\br_u) (\epsilon + \delta_o + \delta_m) > 0.9,
	\\
		&
		\frac{c(n,p,\bslashd \flu)(\epsilon + \delta_o + \ud_{o,\uh})}{0.9 - 1/9}
		<
		\bar{\delta},
	\end{aligned}
	\right.
	\Leftrightarrow
	\left\{
	\begin{aligned}
		&
		c(n,p,\flu) (\epsilon + \delta_o + \delta_m + \ud_{o,\uh}) < 0.1,
	\\
		&
		c(n,p,\br_u) (\epsilon + \delta_o + \delta_m) < 0.1,
	\\
		&
		c(n,p,\bslashd \flu)(\epsilon + \delta_o + \ud_{o,\uh})
		<
		(0.9 - 1/9)\bar{\delta}.
	\end{aligned}
	\right.
	\label{eqn 8.5}
\end{align}
Thus we conclude the following lemma for the improvements on the estimates of $\uflu$, $(\dslashd \uh)|_{\bSigma_u}$, $(\circnabla^2 \uh)|_{\bSigma_u}$.
\begin{lemma}\label{lem 8.4}
Given the conditions \eqref{eqn 8.5}, the estimates \eqref{eqn 5.15} of $\uflu$, $(\dslashd \uh)|_{\bSigma_u}$, $(\circnabla^2 \uh)|_{\bSigma_u}$ can be improved to strict inequalities.
\end{lemma}

\subsubsection{Estimate of the parameterisation function $\flu - \flu_S$}
\begin{lemma}\label{lem 8.5}
By equation \eqref{eqn 6.3},
\begin{align*}
	\vert \flu - \flu_S \vert
	\leq
	C(n,p,\balu) (\epsilon + \delta_o + \delta_m + \ud_{o,\uh}  + \uslashd_m ) \br_u.
\end{align*}
\end{lemma}
\begin{proof}
We derive the propagation equation of $\flu-\flu_S$ from equation \eqref{eqn 6.3},
\begin{align*}
	\uL (\flu - \flu_S) 
	= 
	\balu - \balu_S.
\end{align*}
Integrating the above equation, we obtain that
\begin{align*}
	\vert \flu - \flu_S \vert
	&
	\leq
	\vert \fl{u=0} - \fl{u=0}_S \vert
	+
	\int_0^u \vert \bal{u'} - \bal{u'}_S \vert \ed u
\\
	&
	\leq
	c(n,p) ( \delta_o + \delta_m) r_0
	+
	c(n,p,\balu) (\epsilon + \delta_o + \delta_m + \ud_{o,\uh}  + \uslashd_m) u
\\
	&
	\leq
	C(n,p,\balu) (\epsilon + \delta_o + \delta_m + \ud_{o,\uh}  + \uslashd_m) \br_u.
\end{align*}
\end{proof}

\subsubsection{Estimate of parameterisation function $\bslashd \flu$}

As explained in section \ref{sec 8.1.2}, we shall obtain the estimate of $\bslashd \flu$ from the elliptic equation \eqref{eqn 8.1}
\begin{align}
	\balu^{-1} \btr \bchilu' 
	= 
	\ddtr \ddchilu' 
	= 
	- 2 (r_S|_{\bSigma_u})^{-2} \circDelta \flu 
	+ \tr \chi'_S|_{\bSigma_u}
	+ \hi{\ddtr \ddchilu'}.
	\tag{\ref{eqn 8.1}\ensuremath{'}}
\end{align}
In the case of the spherically symmetric foliation in the Schwarzschild spacetime, equation \eqref{eqn 8.1} becomes
\begin{align}
	\balu_S^{-1} (\btr \bchilu')_S
	= 
	(\ddtr \ddchilu')_S
	= 
	- 2 \br_{u,S}^{-2} \circDelta \flu_S
	+ (\tr \chi'_S)|_{\bSigma_{u,S}}.
	\tag{\ref{eqn 8.1}$'_S$}
\end{align}
Introduce the surface $\bSigma_u'$ which is viewed as the average of $\bSigma_u$, i.e.
\begin{align*}
	\bSigma_u' = \Sigma (\us=\overline{\uflu}^{\circg}, s= \overline{\flu}^{\circg}).
\end{align*}
Then we derive the equation of $\flu - \overline{\flu}^{\circg}$ from equation \eqref{eqn 8.1},
\begin{align*}
	&\phantom{=\ }
	- 2 (r_S|_{\bSigma_u'})^{-2} \circDelta (\flu -  \overline{\flu}^{\circg})
	+ \tr \chi'_S|_{\bSigma_u} - \tr \chi'_S|_{\bSigma_u'}
\\	
	&=
	2 [(r_S|_{\bSigma_u})^{-2} - (r_S|_{\bSigma_u'})^{-2}] \circDelta \flu
	+ \balu_S^{-1} (\btr \bchilu')_S - \tr \chi'_S|_{\bSigma_u'}
	-\hi{\ddtr \ddchilu'}.
\end{align*}
Decompose the left hand side of the above equation as a linear operator on $\flu - \overline{\flu}^{\circg}$ and a remainder term and rewrite the equation as follows
\begin{align}
	\begin{aligned}
	&\phantom{= \ }
	- 2 (r_S|_{\bSigma_u'})^{-2} \circDelta ( \flu - \overline{\flu}^{\circg} )
	+ (\partial_s \tr \chi'_S)|_{\bSigma_u'} (\flu - \overline{\flu}^{\circg} )
\\
	&= 
	- [ \tr \chi'_S|_{\bSigma_u} - (\tr \chi'_S)|_{\bSigma_u'} - (\partial_s \tr \chi'_S)|_{\bSigma_u'} (\flu - \overline{\flu}^{\circg})]
	+ 2 [(r_S|_{\bSigma_u})^{-2} - (r_S|_{\bSigma_u'})^{-2}] \circDelta \flu 
\\
	&\phantom{=\ }
	+ \balu^{-1} \btr \bchilu' - \tr \chi'_S|_{\bSigma_u'}
	- \hi{\ddtr \ddchilu'}.
	\end{aligned}
	\label{eqn 8.6}
\end{align}
Since $\partial_s \tr \chi'_S = -\frac{2(r-r_0)}{r^3} + \frac{2 r_0}{r(r_0+s)^3}$, the linear operator in the left hand side of equation \eqref{eqn 8.6} is
\begin{align*}
	- 2 (r_S|_{\bSigma_u'})^{-2} \circDelta + (\partial_s \tr \chi'_S)|_{\bSigma_u'}
	=
	2 (r_S|_{\bSigma_u'})^{-2} (- \circDelta - 1) + \frac{2 r_0}{r_S|_{\bSigma_u'} ( r_0 + \overline{\flu}^{\circg})^2} + \frac{2r_0}{(r_S|_{\bSigma_u'})^3}.
\end{align*}
Introduce the operator $\Lambda_u$ as
\begin{align}
	\Lambda_{\bSigma_u'}
	= 
	(- \circDelta - 1) 
	+ \frac{ r_0 r_S|_{\bSigma_u'} }{( r_0 + \overline{\flu}^{\circg})^2} 
	+ \frac{r_0}{r_S|_{\bSigma_u'}}.
	\label{eqn 8.7}
\end{align}
This operator is bounded and invertible from the Sobolev space $\dot{\mathrm{W}}^{n+2,p}(\mathbb{S}^2,\circg)$ to $\dot{\mathrm{W}}^{n,p}(\mathbb{S}^2,\circg)$ where
\begin{align*}
	\dot{\mathrm{W}}^{m,p}(\mathbb{S}^2,\circg)
	= 
	\{ f \in \mathrm{W}^{n,p}(\mathbb{S}^2,\circg): \int_{\mathbb{S}^2} f \dvol_{\circg} =0 \}.
\end{align*}
There exists a constant $c(n,p)$ depending on $n,p$, such that for all $u\geq 0$, we have
\begin{align}
	\Vert \slashd f \Vert^{n+1,p} \leq c(n,p) \Vert  \Lambda_{\bSigma_u'} (f-\overline{f}^{\circg}) \Vert^{n,p}.
	\label{eqn 8.8}
\end{align}
Introduce the notation $\osc{\cdot}$ to denote a quantity subtracting its mean value with respect to $\circg$, i.e.
\begin{align*}
	\osc{f} = f - \overline{f}^{\circg}.
\end{align*}
Then by the above estimate and following equation derived from equation \eqref{eqn 8.6},
\begin{align}
	\begin{aligned}
	&\phantom{= \ }
	- 2 (r_S|_{\bSigma_u'})^{-2} \Lambda_{\bSigma_u'} ( \osc{\flu} )
\\
	&= 
	- \osc{ \tr \chi'_S|_{\bSigma_u} - (\tr \chi'_S)|_{\bSigma_u'} - (\partial_s \tr \chi'_S)|_{\bSigma_u'} \cdot \osc{\flu}}
\\
	&\phantom{=\ }
	+ 2 \osc{[(r_S|_{\bSigma_u})^{-2} - (r_S|_{\bSigma_u'})^{-2}] \circDelta \flu }
	+ \osc{\balu^{-1} \btr \bchilu'}
	- \osc{\hi{\ddtr \ddchilu'}},
	\end{aligned}
	\tag{\ref{eqn 8.6}\ensuremath{'}}
	\label{eqn 8.6'}
\end{align}
we can obtain the estimate of $\bslashd \flu$.

\begin{lemma}\label{lem 8.6}
By equation \eqref{eqn 8.6}, we obtain that
\begin{align*}
	\Vert \bslashd \flu \Vert^{n+1,p}
	\leq
	C(n,p, \bslashd \balu, \bslashd \btr \bchilu' ) (\epsilon + \delta_o + \uh_{o,\uh}) \br_u.
\end{align*}
\end{lemma}
\begin{proof}
We estimate each term on the right hand side of equation \eqref{eqn 8.6'}.
\begin{enumerate}[label=\textit{\textbf{\alph*.}}]
\item
$\pmb{ \osc{[(r_S|_{\bSigma_u})^{-2} - (r_S|_{\bSigma_u'})^{-2}] \circDelta \flu }}$. By the Sobolev inequality, we have that
\begin{align*}
	&\phantom{\leq \ }
	\Vert \osc{[(r_S|_{\bSigma_u})^{-2} - (r_S|_{\bSigma_u'})^{-2}] \circDelta \flu } \Vert^{n,p}
\\
	&\leq
	c(n,p) \Vert \bslashd \{ [(r_S|_{\bSigma_u})^{-2} - (r_S|_{\bSigma_u'})^{-2}] \circDelta \flu \} \Vert^{n-1,p}
\\
	&\leq
	c(n,p) \Vert \bslashd [(r_S|_{\bSigma_u})^{-2}] \Vert^{n-1,p} \cdot \Vert \circDelta \flu \Vert^{n-1,p}
\\
	&\phantom{\leq \ }
	+ c(n,p)  \Vert (r_S|_{\bSigma_u})^{-2} - (r_S|_{\bSigma_u'})^{-2} \Vert^{n-1,p} \cdot \Vert \bslashd \circDelta \flu \Vert^{n-1,p}.
\end{align*}
We need to estimate each term on the right hand side.
\begin{enumerate}[label=\textit{\textbf{a.\arabic*.}}]
\item
$\pmb{\bslashd [(r_S|_{\bSigma_u})^{-2}]}$.
For this term, we have that
\begin{align*}
\bslashd [(r_S|_{\bSigma_u})^{-2}] 
= 
-2 (r_S|_{\bSigma_u})^{-3} ( \partial_s r_S \cdot  \bslashd \flu +  \partial_{\us} r_S \cdot \bslashd \uflu ).
\end{align*}
By the estimates of $\bslashd \flu$ and $\bslashd \uflu$ in the bootstrap assumption,
\begin{align*}
	\Vert \bslashd [(r_S|_{\bSigma_u})^{-2}]  \Vert^{n,p}
	\leq
	\frac{C(n,p,\bslashd \flu)(\epsilon + \delta_o + \ud_{o,\uh})}{\br_u^2}.
\end{align*}

\item
$\pmb{(r_S|_{\bSigma_u})^{-2} - (r_S|_{\bSigma_u'})^{-2}}$. For $r_S|_{\bSigma_u} - r_S|_{\bSigma_u'}$, we have that
\begin{align*}
	r_S|_{\bSigma_u} - r_S|_{\bSigma_u'}
	&=
	r_S(\uflu, \flu) - r_S(\overline{\uflu}^{\circg}, \overline{\flu}^{\circg})
\\
	&=
	\int_0^1 \partial_s r_S(t \uflu + (1-t) \overline{\uflu}^{\circg}, t \flu + (1-t) \overline{\flu}^{\circg}) \ed t \cdot (\flu- \overline{\flu}^{\circg})
\\
	&\phantom{=\ }
	+
	\int_0^1 \partial_{\us} r_S(t \uflu + (1-t) \overline{\uflu}^{\circg}, t \flu + (1-t) \overline{\flu}^{\circg}) \ed t \cdot (\uflu- \overline{\uflu}^{\circg}),
\end{align*}
then
\begin{align*}
	\vert r_S|_{\bSigma_u} -  r_S|_{\bSigma_u'} \vert
	\leq
	C(n,p,\bslashd \flu)(\epsilon + \delta_o + \ud_{o,\uh}) \br_u.
\end{align*}
Thus
\begin{align*}
	\Vert (r_S|_{\bSigma_u})^{-2} - (r_S|_{\bSigma_u'})^{-2} \Vert^{n+1,p}
	&\leq
	\frac{C(n,p,\bslashd \flu)(\epsilon + \delta_o + \ud_{o,\uh})}{\br_u^2}.
\end{align*}
\end{enumerate}
Combining the above estimates of $\bslashd [(r_S|_{\bSigma_u})^{-2}]$, $(r_S|_{\bSigma_u})^{-2} - (r_S|_{\bSigma_u'})^{-2}$ and the estimate of $\flu$ in the bootstrap assumption, we obtain that
\begin{align*}
	\Vert \osc{[(r_S|_{\bSigma_u})^{-2} - (r_S|_{\bSigma_u'})^{-2}] \circDelta \flu } \Vert^{n,p}
	& \leq
	\frac{C(n,p, \bslashd \flu)(\epsilon + \delta_o + \ud_{o,\uh})^2}{\br_u}.
\end{align*}

\item
$\pmb{\osc{[ \tr \chi'_S|_{\bSigma_u} - (\tr \chi'_S)|_{\bSigma_u'} - (\partial_s \tr \chi'_S)|_{\bSigma_u'} \cdot \osc{\flu}]}}$.
We have the following formula for $ \tr \chi'_S|_{\bSigma_u} - (\tr \chi'_S)|_{\bSigma_u'} - (\partial_s \tr \chi'_S)|_{\bSigma_u'} \cdot \osc{\flu} $ that
\begin{align*}
	&\phantom{=\ }
	\tr \chi'_S|_{\bSigma_u} - (\tr \chi'_S)|_{\bSigma_u'} - (\partial_s \tr \chi'_S)|_{\bSigma_u'} \cdot \osc{\flu} 
	 \\
	 &=
	 \underbrace{\int_0^1 [ \partial_s \tr \chi'_S(t \uflu + (1-t) \overline{\uflu}^{\circg}, t \flu + (1-t) \overline{\flu}^{\circg}) - \partial_s \tr \chi'_S( \overline{\uflu}^{\circg}, \overline{\flu}^{\circg})  ] \ed t  \cdot (\flu- \overline{\flu}^{\circg})}_{\pmb{b.1}}
\\
	&\phantom{=\ }
	+
	\underbrace{\int_0^1 \partial_{\us} \tr \chi'_S(t \uflu + (1-t) \overline{\uflu}^{\circg}, t \flu + (1-t) \overline{\flu}^{\circg}) \ed t  \cdot (\uflu- \overline{\uflu}^{\circg})}_{\pmb{b.2}}
\end{align*}
Similar to the estimates of $(r_S|_{\bSigma_u})^{-2} - (r_S|_{\bSigma_u'})^{-2}$,
\begin{align*}
	&\phantom{\leq\ }
	\Vert  \partial_s \tr \chi'_S(t \uflu + (1-t) \overline{\uflu}^{\circg}, t \flu + (1-t) \overline{\flu}^{\circg}) - \partial_s \tr \chi'_S( \overline{\uflu}^{\circg}, \overline{\flu}^{\circg}) \Vert^{n+1,p}
\\
	&\leq
	t\cdot \frac{C(n,p,\bslashd \flu)(\epsilon + \delta_o + \ud_{o,\uh})}{\br_u^2},
\\
	&\phantom{\leq\ }
	\Vert \bslashd [ \partial_s \tr \chi'_S(t \uflu + (1-t) \overline{\uflu}^{\circg}, t \flu + (1-t) \overline{\flu}^{\circg}) - \partial_s \tr \chi'_S( \overline{\uflu}^{\circg}, \overline{\flu}^{\circg}) ] \Vert^{n,p}
\\
	&\leq
	t\cdot \frac{C(n,p,\bslashd \flu)(\epsilon + \delta_o + \ud_{o,\uh})}{\br_u^2},
\end{align*}
thus we have that
\begin{align*}
	\Vert \osc{\pmb{b.1}} \Vert^{n+1,p}
	& \leq
	\frac{C(n,p,\bslashd \flu)(\epsilon + \delta_o + \ud_{o,\uh})^2}{\br_u}.
\end{align*}
For $\pmb{b.2}$, we have that
\begin{align*}
	\Vert \osc{\pmb{b.2}} \Vert^{n+1,p}
	\leq
	\frac{c(n,p) \ud_o r_0}{\br_u^2}.
\end{align*}
Combining the estimates of $\pmb{b.1}$ and $\pmb{b.2}$, we obtain the estimate of $\Vert \osc{[ \tr \chi'_S|_{\bSigma_u} - (\tr \chi'_S)|_{\bSigma_u'} - (\partial_s \tr \chi'_S)|_{\bSigma_u'} \cdot \osc{\flu}]} \Vert^{n+1,p}$.

\item
$\pmb{\osc{\balu^{-1} \btr \bchilu'}}$.
By the Sobolev inequality, we have that
\begin{align*}
	\Vert \osc{\balu^{-1} \btr \bchilu'} \Vert^{n,p}
	&\leq
	c(n,p) [\Vert \bslashd (\balu^{-1}) \Vert^{n-1,p} \Vert \btr \bchilu' \Vert^{n-1,p}
	+ \Vert \balu^{-1} \Vert ^{n-1,p} \Vert \bslashd \btr \bchilu' \Vert^{n-1,p}]
\\
	&\leq
	\frac{C(n,p,\bslashd \balu, \bslashd \btr \bchilu')(\epsilon + \delta_o + \ud_{o,\uh})}{\br_u}.
\end{align*}

\item
$\pmb{\osc{\hi{ \ddtr \ddchilu'}}}$. By proposition \ref{prop 5.10}, we have
\begin{align*}
	\Vert \osc{\hi{ \ddtr \ddchilu'}} \Vert^{n,p}
	\leq
	\frac{c(n,p) \epsilon}{\br_u} + \frac{c(n,p) \delta_o ( \delta_o + \ud_{\uh}) r_0}{\br_u^2}.
\end{align*}
\end{enumerate}
Combining the above estimates of each term in the right hand side of equation \eqref{eqn 8.6'}, we obtain that
\begin{align*}
	\Vert \bslashd \flu \Vert^{n+1,p}
	&\leq
	c(n,p) \Vert \Lambda_{\bSigma_u'} ( \osc{\flu} ) \Vert^{n,p}
	\leq
	C(n,p,\bslashd \balu, \bslashd \btr \bchilu')(\epsilon + \delta_o + \ud_{o,\uh}) \br_u.
\end{align*}
\end{proof}

\subsection{Estimates of metric components}\label{sec 8.5}

\subsubsection{Estimate of lapse function $\balu - \balu_S$}
We integrate equation \eqref{eqn 6.4} to obtain the estimate of $\balu - \balu_S$.
\begin{lemma}\label{lem 8.7}
By equation \eqref{eqn 6.4},
\begin{align*}
	\vert \balu - \balu_S \vert
	\leq
	c(n,p) ( \delta_m + \uslashd_m )
	+ C(n,p,\bslashd \balu) ( \epsilon + \delta_o + \ud_{o,\uh} ).
\end{align*}
\end{lemma}
\begin{proof}
Recall equation \eqref{eqn 6.4},
\begin{align}
	\buL \log \balu = 2 \buomegalu - 2 \balu \cdot \duomegalu.
	\tag{\ref{eqn 6.4}\ensuremath{'}}
\end{align}
In the case of the spherically symmetric foliation in Schwarzschild spacetime, equation \eqref{eqn 6.4} reads
\begin{align*}
	\buL \log \balu_S 
	= 2 \buomegalu_S - 2 \balu_S \cdot \duomegalu_S 
	= - 2 \balu_S \cdot \duomegalu_S.
\end{align*}
We derive the equation for the mean value $\overline{\balu}^{u}$ from equation \eqref{eqn 6.4}
\begin{align*}
	\buL \overline{\log \balu}^{u}
	=
	2 \overline{\buomegalu}^{u}
	- 2 \overline{\balu \cdot \duomegalu}^{u}
	- \overline{(\btr \buchilu - \overline{\btr \buchilu}^{u})(\log \balu - \overline{\log \balu}^{u})}^{u}.
\end{align*}
Taking the difference of the two equations, we have
\begin{align}
	\begin{aligned}
		\buL (\overline{\log \balu}^{u} - \log \balu_S)
		&=
		2 \overline{\buomegalu}^{u} 
		- 2 (\overline{\balu \cdot \duomegalu}^{u} - \balu_S \cdot \duomegalu_S)
	\\
		&\phantom{= }
		- \overline{(\btr \buchilu - \overline{\btr \buchilu}^{u})(\log \balu - \overline{\log \balu}^{u})}^{u}.
	\end{aligned}
	\label{eqn 8.9}
\end{align}

We shall estimate $\overline{\log \balu}^{u} - \log \balu_S$ by integrating the above equation \eqref{eqn 8.9}, then by the Sobolev inequality
\begin{align*}
	\vert \log \balu - \overline{\log \balu}^{u} \vert
	\leq
	c(n,p) \Vert \frac{\bslashd \balu}{\balu} \Vert^{n,p}
	\leq
	C(n,p, \bslashd \balu) (\epsilon + \delta_o + \ud_{o,\uh}),
\end{align*}
we can obtain the estimate of $\log \balu - \log \balu_S$. To derive the estimate of $\balu - \balu_S$, we use the following formula
\begin{align*}
	\balu - \balu_S
	=
	\balu_S [ \exp (\log \balu - \log \balu_S) - 1].
\end{align*}
Following the above procedure, we estimate $\balu - \balu_S$.
\begin{enumerate}[label=\textbf{\textit{\alph*.}}]
\item
$\pmb{\overline{\log \balu}^{u} - \log \balu_S}$.
We estimate each term on the right hand side of equation \eqref{eqn 8.9}.
\begin{enumerate}[label*=\textbf{\textit{\arabic*.}}]
\item
$\pmb{\overline{\buomegalu}^{u}}$.
By the bootstrap assumption,
\begin{align*}
	\vert \overline{\buomegalu}^{u} \vert
	\leq
	\frac{c(n,p,\overline{\buomegalu}^{u})(\epsilon + \delta_o + \ud_{o,\uh})^2 r_0^2}{\br_u^3}.
\end{align*}

\item
$\pmb{\duomegalu - \duomegalu_S}$. By the decomposition \eqref{eqn 5.8d} of $\duomegalu$, we have that
\begin{align*}
	\duomegalu - \duomegalu_S
	=
	\uomega_S|_{\bSigma_u} - \uomega_S|_{\bSigma_{u,S}}
	+ \hi{\duomegalu}.
\end{align*}
For $\uomega_S|_{\bSigma_u} - \uomega_S|_{\bSigma_{u,S}}$, we have that
\begin{align*}
	\uomega_S|_{\bSigma_u} - \uomega_S|_{\bSigma_{u,S}}
	&=
	\int_0^1 \partial_s \uomega_S(t \uflu + (1-t) \uflu_S, t \flu + (1-t) \flu_S) \ed t \cdot (\flu- \flu_S)
\\
	&\phantom{=\ }
	+
	\int_0^1 \partial_{\us} \uomega_S(t \uflu + (1-t) \uflu_S, t \flu + (1-t) \flu_S) \ed t \cdot (\uflu- \uflu_S).
\end{align*}
Then by the estimates of $\flu, \uflu$ in the bootstrap assumption and the estimate of $\hi{\duomegalu}$ in proposition \ref{prop 5.9}, we have that
\begin{align*}
	\vert \uomega_S|_{\bSigma_u} - \uomega_S|_{\bSigma_{u,S}} \vert
	&\leq
	\frac{C(n,p,\flu)(\udelta_m + \uslashd_m) (\epsilon + \delta_o + \delta_m + \ud_{o,\uh}) r_0^2}{\br_u^3}
	+ \frac{c(n,p) \uslashd_m r_0^2}{\br_u^3},
\\
	\vert \hi{\duomegalu} \vert
	&\leq
	\frac{c(n,p)(\epsilon + \ud_{\uh}^2) r_0^2}{\br_u^3},
\end{align*}
which imply that
\begin{align*}
	\vert \duomegalu - \duomegalu_S \vert
	&\leq
	\frac{c(n,p)(\epsilon + \ud_{\uh}^2 + \uslashd_m) r_0^2}{\br_u^3}
	+ \frac{C(n,p,\flu)\udelta_m ( \delta_o + \delta_m + \ud_{o,\uh}) r_0^2}{\br_u^3}.
\end{align*}

\item
$\pmb{\overline{\balu \cdot \duomegalu}^{u} - \balu_S \cdot \duomegalu_S}$.
We have that
\begin{align*}
	\overline{\balu \cdot \duomegalu}^{u} - \balu_S \cdot \duomegalu_S
	=
	\overline{(\balu - \balu_S) \duomegalu}^{u}
	+ \balu_S ( \overline{\duomegalu}^{u} - \duomegalu_S)
\end{align*}
By the estimates of $\duomegalu$ in propositions \ref{prop 5.7} and \ref{prop 5.9}, and the estimate of $\duomegalu - \duomegalu_S$ above, we obtain that
\begin{align*}
	\vert \overline{\balu \cdot \duomegalu}^{u} - \balu_S \cdot \duomegalu_S \vert
	&\leq
	\frac{c(n,p)(\epsilon + \ud_o + \ud_m + \ud_h^2) r_0^2}{\br_u^3} \Vert \balu - \balu_S \Vert^{\infty}
\\
	&\phantom{\leq\ }
	+ \frac{c(n,p)(\epsilon + \ud_{\uh}^2 + \uslashd_m) r_0^2}{\br_u^3}
	+ \frac{C(n,p,\flu)\udelta_m ( \delta_o + \delta_m + \ud_{o,\uh}) r_0^2}{\br_u^3}.
\end{align*}

\item
$\pmb{\overline{(\btr \buchilu - \overline{\btr \buchilu}^{u})(\log \balu - \overline{\log \balu}^{u})}^{u}}$.
By the estimates of $\bslashd \btr \buchilu$ in the bootstrap assumption and the previous estimate of $\log \balu - \overline{\log \balu}^{u}$ from $\bslashd \balu$, we have
\begin{align*}
	\vert \overline{(\btr \buchilu - \overline{\btr \buchilu}^{u})(\log \balu - \overline{\log \balu}^{u})}^{u} \vert
	\leq
	\frac{C(n,p,\bslashd \btr \buchilu, \bslashd \balu) ( \epsilon + \delta_o + \ud_{o,\uh})^2 r_0}{\br_u^2}.
\end{align*}
\end{enumerate}

Combining the above estimates of the right hand side of equation \eqref{eqn 8.9}, we have
\begin{align*}
	\vert \buL (\overline{\log \balu}^{u} - \log \balu_S) \vert
	&\leq
	\frac{c(n,p)(\epsilon + \ud_o + \ud_m + \ud_h^2) r_0^2}{\br_u^3} \Vert \balu - \balu_S \Vert^{\infty}
\\
	&\phantom{\leq\ }
	+ \frac{c(n,p)(\epsilon + \ud_{\uh}^2 + \uslashd_m) r_0^2}{\br_u^3}
	+ \frac{C(n,p,\flu)\udelta_m ( \delta_o + \delta_m + \ud_{o,\uh}) r_0^2}{\br_u^3}
\\
	&\phantom{\leq\ }
	+ \frac{C(n,p,\bslashd \btr \buchilu, \bslashd \balu) ( \epsilon + \delta_o + \ud_{o,\uh})^2 r_0}{\br_u^2}.
\end{align*}
Thus by integrating the above inequality, we obtain that
\begin{align*}
	\vert \overline{\log \balu}^{u} - \log \balu_S \vert
	&\leq
	\int_0^u \frac{c(n,p)(\epsilon + \ud_o + \ud_m + \ud_h^2) r_0^2}{\br_u^3} \Vert \balu - \balu_S \Vert^{\infty} \ed u
\\
	&\phantom{\leq\ }
	+ c(n,p)(\epsilon + \delta_o + \delta_m + \ud_{o, \uh} + \uslashd_m).
\end{align*}
Substituting the estimate of $\balu - \balu_S$ in the bootstrap assumption, we have
\begin{align*}
	\vert \overline{\log \balu}^{u} - \log \balu_S \vert
	\leq
	C(n,p,\balu) (\epsilon + \delta_o + \delta_m + \ud_{o,\uh} + \uslashd_m)
	\leq 
	1.	
\end{align*}

\item
$\pmb{\log \balu - \log \balu_S}$.
We have that
\begin{align*}
	\vert \log \balu - \log \balu_S \vert
	&\leq
	\vert \log \balu - \overline{\log \balu}^{u} \vert
	+ \vert \overline{\log \balu}^{u} - \log \balu_S \vert
\\
	&\leq
	\int_0^u \frac{c(n,p)(\epsilon + \ud_o + \ud_m + \ud_h^2) r_0^2}{\br_u^3} \Vert \balu - \balu_S \Vert^{\infty} \ed u
\\
	&\phantom{\leq\ }
	+ c(n,p) (\delta_m + \uslashd_m)
	+ c(n,p,\bslashd \balu )(\epsilon + \delta_o + \ud_{o,\uh}).
\end{align*}

\item
$\pmb{\balu - \balu_S}$.
By the formula $\balu - \balu_S = \balu_S [ \exp (\log \balu - \log \balu_S) - 1]$, we have that
\begin{align*}
	\Vert \balu - \balu_S \Vert^{\infty}
	&\leq
	\int_0^u \frac{c(n,p)(\epsilon + \ud_o + \ud_m + \ud_h^2) r_0^2}{\br_u^3} \Vert \balu - \balu_S \Vert^{\infty} \ed u
\\
	&\phantom{\leq\ }
	+ c(n,p) (\delta_m + \uslashd_m)
	+ C(n,p,\bslashd \balu )(\epsilon + \delta_o + \ud_{o,\uh}),
\end{align*}
thus by Gronwall's inequality, we obtain that
\begin{align*}
	\Vert \balu - \balu_S \Vert^{\infty}
	\leq
	c(n,p) (\delta_m + \uslashd_m)
	+ C(n,p,\bslashd \balu )(\epsilon + \delta_o + \ud_{o,\uh}).
\end{align*}
\end{enumerate}
\end{proof}

\subsubsection{Rotational vector fields $\{ \barR_i \}_{i=1,2,3}$ on $\{ \bSigma_u \}$}
We introduce the rotational vector fields $\{\barR_i \}_{i=1,2,3}$ on $\{ \bSigma_u \}$, similar to the rotational vector fields $\{ \dR_i \}_{i=1,2,3}$ on $\{ \Sigma_s \}$ in section \ref{sec 5.1.1}. 

Restrict the mapping from the Schwarzschild spacetime $(M,g_S)$ to $(\mathbb{S}^2, \circg)$ on the foliation leaf $\bSigma_u$ and pull the rotational vector fields $\{ R_i \}_{i=1,2,3}$ on $(\mathbb{S}^2, \circg)$ back to $\bSigma_u$, which give rise to the rotational vector fields $\{\barR_i \}_{i=1,2,3}$ on $\bSigma_u$. We have that
\begin{align*}
	\barR_i = \dR_i + (\barR_i \flu) \ddpartial_u,
	\quad
	\lie_{\ddpartial_u} \barR_i = [\ddpartial_u, \barR_i] =0.
\end{align*}
We shall use the notation of the rotational vector component of a tensor field on $(\bSigma_u, \circg)$ as introduced in section \ref{sec 5.1.1}.

\subsubsection{Estimate of lapse function $\bslashd \balu$}
To estimate $\bslashd \balu$, we employ the rotational vector derivatives $\{\barR_i\}_{i=1,2,3}$ on $\{\bSigma_u\}$ to the propagation equation \eqref{eqn 6.4} of $\log \balu$ to derive the propagation equation for the rotational vector derivatives of $\log \balu$. Rewrite equation \eqref{eqn 6.4} as follows
\begin{align}
	\ddpartial_u \log \balu + \balu \db^i \ddpartial_i \log \balu
	=
	2 \buomegalu - 2 \balu \duomegalu.
	\tag{\ref{eqn 6.4}\ensuremath{''}}
\end{align}
In the following, use $(\log \balu)_{\barR,k}$ to denote the rotational vector derivative $\barR_k \log \balu$ of $\log \balu$, and similar notations to other functions. Differentiating equation \eqref{eqn 6.4} with respect to $\barR_k$, we obtain the propagation equation for the rotational derivative of $\log \balu$
\begin{align}
	(\ddpartial_u + \balu \db^i \ddpartial_i) ( \log \balu)_{\barR,k}
	=
	-[\barR_k, \balu \db^i \ddpartial_i] \log \balu
	+ 2 \buomegalu_{\barR, k}
	- 2 \duomegalu \cdot \balu_{\barR,k}
	- 2 \balu \cdot \duomegalu_{\barR, k}.
	\label{eqn 8.10}
\end{align}
We shall estimate the right hand side of the propagation equation, then integrate it to derive the estimate of $\bslashd \balu$.
\begin{lemma}\label{lem 8.8}
By equation \eqref{eqn 8.10},
\begin{align*}
	\Vert \bslashd \balu \Vert^{n,p}
	\leq
	C(n,p, \bslashd \buomegalu) ( \epsilon + \delta_o + \ud_{o,\uh})
\end{align*}
\end{lemma}
\begin{proof}
We shall estimate the shifting vector field $\balu \db^i \ddpartial_i$ and the right hand side of propagation, then integrate it.
\begin{enumerate}[label=\textbf{\textit{\alph*.}}]
\item
$\pmb{\balu \db, [\barR_k, \balu \db^i \ddpartial_i]}$.
By the estimates of $\balu$ in bootstrap assumption and $\db$ in propositions \ref{prop 5.7}, \ref{prop 5.9}, we have
\begin{align*}
	\Vert \balu \db \Vert^{n+1,p}, \Vert [\barR_k, \balu \db^i \ddpartial_i] \Vert^{n,p}
	\leq
	\frac{c(n,p) \ud_{\uh} r_0}{\br_u^2}
	+ \frac{c(n,p) \epsilon ( \ud_o + \ud_m ) r_0^2}{\br_u^3}.
\end{align*}

\item
$\pmb{[\barR_k, \balu \db^i \ddpartial_i] \log \balu}$.
We have
\begin{align*}
	\Vert [\barR_k, \balu \db^i \ddpartial_i] \log \balu \Vert^{n,p}
	&\leq
	c(n,p) \Vert [\barR_k, \balu \db^i \ddpartial_i] \Vert^{n,p} \cdot \Vert \frac{\bslashd \balu}{\balu} \Vert^{n,p}
\\
	&\leq
	c(n,p) [\frac{\ud_{\uh} r_0}{\br_u^2}
	+ \frac{\epsilon ( \ud_o + \ud_m ) r_0^2}{\br_u^3}]
	\cdot \Vert \bslashd \balu \Vert^{n,p}
\\
	&\leq
	C(n,p,\bslashd \balu)
	[ \frac{ \ud_{\uh} r_0}{\br_u^2}
	+ \frac{\epsilon ( \ud_o + \ud_m )r_0^2}{\br_u^3}]
	(\epsilon + \delta_o + \ud_{o,\uh}).
\end{align*}

\item
$\pmb{\buomegalu_{\barR,k}}$.
By the estimate of $\bslashd \buomegalu$ in the bootstrap assumption, we have
\begin{align*}
	\Vert \buomegalu_{\barR,k} \Vert^{n,p}
	\leq
	\frac{C(n,p,\bslashd \buomegalu) \epsilon r_0^{\frac{3}{2}}}{\br_u^{\frac{5}{2}}}
	+
	\frac{C(n,p,\bslashd \buomegalu) (\delta_o + \ud_{o,\uh}) r_0^2}{\br_u^3}.
\end{align*}

\item
$\pmb{\duomegalu \cdot \balu_{\barR, k}}$.
By the estimates of $\bslashd \balu$ in bootstrap assumption and $\duomegalu$ in propositions \ref{prop 5.7}, \ref{prop 5.9}, we have
\begin{align*}
	\Vert \duomegalu \cdot \balu_{\barR, k} \Vert^{n,p}
	\leq
	\frac{C(n,p,\bslashd \balu) (\epsilon + \ud_o + \ud_h^2 + \ud_m)(\epsilon + \delta_o + \ud_{o,\uh}) r_0^2}{\br_u^3}.
\end{align*}

\item
$\pmb{\balu \cdot \duomegalu_{\barR,k}}$.
By the estimates of $\balu$ in bootstrap assumption and $\duomegalu$ in propositions \ref{prop 5.7}, \ref{prop 5.9}, we have
\begin{align*}
	\Vert \balu \cdot \duomegalu_{\barR, k} \Vert^{n,p}
	\leq
	\frac{c(n,p) [\epsilon + (\ud_m + \ud_o) d_o + \ud_o + \ud_h^2] r_0^2}{\br_u^3}.
\end{align*}
\end{enumerate}
Combining the above estimates, we obtain that
\begin{align*}
	\Vert (\ddpartial_u + \balu \db^i \ddpartial_i) (\log \balu)_{\barR,k} \Vert^{n,p}
	&\leq
	c(n,p) [\frac{\ud_{\uh} r_0}{\br_u^2}
	+ \frac{\epsilon ( \ud_o + \ud_m ) r_0^2}{\br_u^3}]
	\cdot \Vert \bslashd \balu \Vert^{n,p}
\\
	&\phantom{\leq\ }
	+ \frac{C(n,p,\bslashd \buomegalu) \epsilon r_0^{\frac{3}{2}}}{\br_u^{\frac{5}{2}}}
	+ \frac{C(n,p,\bslashd \buomegalu) (\delta_o + \ud_{o,\uh}) r_0^2}{\br_u^3}.
\end{align*}
By the estimate of $\balu \db^i \ddpartial_i$ and the initial estimate of $\balu$ in proposition \ref{prop 7.2}, we integrate the above estimate of $(\ddpartial_u + \balu \db^i \ddpartial_i) (\log \balu)_{\barR,k}$ to obtain that
\begin{align*}
	\Vert \bslashd \balu \Vert^{n,p}
	\leq
	\int_0^u c(n,p) [\frac{\ud_{\uh} r_0}{\br_u^2}
	+ \frac{\epsilon ( \ud_o + \ud_m ) r_0^2}{\br_u^3}]
	\cdot \Vert \bslashd \balu \Vert^{n,p} \ed u
	+ C(n,p,\bslashd \buomegalu) (\epsilon + \delta_o + \ud_{o,\uh}).
\end{align*}
By Gronwall's inequality, we obtain the estimate of $\bslashd \balu$ in the lemma.
\end{proof}

\subsubsection{Estimate of area radius $\br_u$}

\begin{lemma}\label{lem 8.9}
By the choice of the parameter by the area radius,
\begin{align*}
	\vert \br_u - \br_{u,S} \vert
	\leq
	c(n,p) (\epsilon + \delta_o + \delta_m) r_0.
\end{align*}
\end{lemma}
\begin{proof}
The estimate simply follows from $\br_u - \br_{u,S} = \br_{u=0} - \br_{u=0,S}$ and the estimate of initial area radius $\br_{u=0} - \br_{u=0,S}$ in lemma \ref{lem 7.3}.
\end{proof}

\subsubsection{Estimate of intrinsic metric $\bslashglu - \bslashglu_S$}

\begin{lemma}\label{lem 8.10}
We have
\begin{align*}
	\vert \bslashglu - \bslashglu_S \vert
	\leq
	C(n,p,\flu) (\epsilon + \delta_o + \delta_m + \ud_{o,\uh} + \uslashd_m) \br_u^2.
\end{align*}
\end{lemma}
\begin{proof}
Decompose $\bslashglu - \bslashglu_S$ as follows
\begin{align*}
	\bslashglu - \bslashglu_S
	&=
	\lo{\ddslashg|_{\bSigma_u}} - \bslashglu_S + \hi{\ddslashg|_{\bSigma_u}}
\\
	&=
	[(r_S|_{\bSigma_u})^2 - (r_S|_{\bSigma_{u,S}})^2] \circg+ \hi{\ddslashg|_{\bSigma_u}}.
\end{align*}
By the estimates of $\flu - \flu_S$ and $\uflu - \uflu_S$ in the bootstrap assumption, we have that
\begin{align*}
	\vert r_S|_{\bSigma_u} - r_S|_{\bSigma_{u,S}} \vert
	\leq
	C(n,p, \flu) ( \epsilon + \delta_o + \delta_m + \ud_{o,\uh} + \uslashd_m) \br_u.
\end{align*}
By the estimate of $\hi{\ddslashg|_{\bSigma_u}}$ in proposition \ref{prop 5.10}, we have that
\begin{align*}
	\vert \hi{\ddslashg|_{\bSigma_u}} \vert
	\leq
	c(n,p) \epsilon \br_u^2.
\end{align*}
Therefore combining the above estimates, we obtain the estimate of $\bslashglu - \bslashglu_S$ in the lemma.
\end{proof}

\subsubsection{Estimate of intrinsic metric $\bcircnabla \bslashglu$}

\begin{lemma}\label{lem 8.11}
We have
\begin{align*}
	\Vert \bcircnabla \bslashglu \Vert^{n,p}
	\leq
	C(n,p, \bslashd \flu) ( \epsilon + \delta_o + \ud_{o,\uh}) \br_u^2.
\end{align*}
\end{lemma}
\begin{proof}
Decompose $\bcircnabla \bslashglu$ as follows
\begin{align*}
	\bcircnabla \bslashglu
	=
	\bslashd (r_S^2|_{\bSigma_u}) \cdot \circg
	 + \bslashd \hi{\ddslashg|_{\bSigma_u}}.
\end{align*}
By the estimates of $\bslashd \flu$ and $\bslashd \uflu$ in the bootstrap assumption, we have
\begin{align*}
	\Vert \bslashd (r_S^2|_{\bSigma_u}) \Vert^{n,p}
	\leq
	c(n,p) ( d_o \br_u +  \ud_o r_0 ) \br_u
	\leq
	C(n,p, \bslashd \flu) ( \epsilon + \delta_o + \ud_{o,\uh}) \br_u^2.
\end{align*}
Combining with the estimate of $\hi{\ddslashg|_{\bSigma_u}}$ in proposition \ref{prop 5.10}, we obtain the estimate of $\bcircnabla \bslashglu$ in the lemma.
\end{proof}

\subsection{Estimates of connection coefficients}\label{sec 8.6}

\subsubsection{Estimate of null expansion $\overline{\btr \buchilu}^{u} - (\btr \buchilu)_S$}
\begin{lemma}\label{lem 8.12}
We have
\begin{align*}
	\vert \overline{\btr \buchilu}^{u} - (\btr \buchilu)_S \vert
	\leq
	\frac{C(n,p,\br_u)(\epsilon + \delta_o + \delta_m) r_0}{\br_u^2}.
\end{align*}
\end{lemma}
\begin{proof}
By the choice of the parameter by the area radius,
\begin{align*}
	\overline{\btr \buchilu}^{u} - (\btr \buchilu)_S
	=
	\frac{2}{\br_u} - \frac{2}{\br_{u,S}},
\end{align*}
then the estimate follows from the estimate of $\br_u - \br_{u,S}$.
\end{proof}

\subsubsection{Estimate of null expansion $\bslashd \btr \buchilu$ (coupled with acceleration $\bslashd \buomegalu$)}
We shall use the propagation equation \eqref{eqn 6.8} to obtain the estimate of $\bslashd \btr \buchilu$. Apply the rotational vector fields $\{ \barR_k \}_{k=1,2,3}$ to equation \eqref{eqn 6.8} to derive the following propagation equation for $\barR_k \btr \buchilu = (\btr \buchilu)_{\barR,k}$,
\begin{align}
	\begin{aligned}
		&\phantom{=\ }
		\buL (\btr \buchilu)_{\barR,k}
		+ \btr \buchilu (\btr \buchilu)_{\barR,k}
	\\
		&=
		- [ \barR_k, \balu \db^i \ddpartial_i] \btr \buchilu
		+ 2 \buomegalu (\btr \buchilu)_{\barR,k}
		+ 2 \btr \buchilu (\buomegalu)_{\barR, k}
		- ( \vert \hatbuchilu \vert^2)_{\barR, k}.
	\end{aligned}
	\label{eqn 8.11}
\end{align}
Integrating the above equation by propagation lemma \ref{lem c.2} gives the estimate of $\bslashd \btr \buchilu$.

A subtle term in equation \eqref{eqn 8.11} is the term $\btr \buchilu (\buomegalu)_{\barR, k}$. When integrating equation \eqref{eqn 8.11}, substituting the estimate of $\bslashd \buomegalu$ from the bootstrap assumption is not sufficient to close the bootstrap argument, since as we shall see later that the estimate of $\bslashd \buomegalu$ by the elliptic equation \eqref{eqn 6.13} involves $\bslashd \btr \buchilu$ as a linear term.

Therefore in order to close the bootstrap argument, we shall treat the propagation equation \eqref{eqn 8.11} coupled with the elliptic equation \eqref{eqn 6.13} of $\buomegalu$, and substituting the estimate of $\bslashd \buomegalu$ obtained from the elliptic equation \eqref{eqn 6.13} to integrate the propagation equation \eqref{eqn 8.11}. In this way, the term $\btr \buchilu (\buomegalu)_{\barR, k}$ is treated as a top order linear term of $\bslashd \btr \buchilu$, instead of a merely integrable term.
\begin{lemma}\label{lem 8.13}
By the propagation equation \eqref{eqn 8.11},
\begin{align*}
	\Vert \bslashd \btr \buchilu \Vert^{n-1,p}
	\leq
	\frac{C(n,p,\bubetalu) (\epsilon + \ud_{o,\uh}) r_0}{\br_u^2}
	+ \frac{c(n,p) \delta_o r_0}{\br_u^2}.
\end{align*}
\end{lemma}
\begin{proof}
We estimate each term on the right hand side of equation \eqref{eqn 8.11}.
\begin{enumerate}[label=\textbf{\textit{\alph*.}}]
\item
$\pmb{[ \barR_k, \balu \db^i \ddpartial_i] \btr \buchilu.}$
By the estimate of $[ \barR_k, \balu \db^i \ddpartial_i]$ in lemma \ref{lem 8.8}, we have that
\begin{align*}
	\Vert [ \barR_k, \balu \db^i \ddpartial_i] \btr \buchilu \Vert^{n-1,p}
	&\leq
	c(n,p) 
	[ \frac{ \ud_{\uh} r_0}{\br_u^2}
	+ \frac{\epsilon ( \ud_o + \ud_m )r_0^2}{\br_u^3}]
	\cdot
	\Vert \bslashd \btr \buchilu \Vert^{n-1,p}.
\end{align*}

\item
$\pmb{\buomegalu (\btr \buchilu)_{\barR,k}.}$ By the estimate of $\buomegalu$ in the bootstrap assumption, we have
\begin{align*}
	\Vert \buomegalu  (\btr \buchilu)_{\barR,k} \Vert^{n-1,p}
	\leq
	C(n,p,\bslashd \buomegalu) [\frac{\epsilon r_0^{\frac{3}{2}}}{\br_u^{\frac{5}{2}}}
	+ \frac{ (\delta_o + \ud_{o,\uh}) r_0^2}{\br_u^3}]
	\cdot
	\Vert \bslashd \btr \buchilu \Vert^{n-1,p}.
\end{align*}

\item
$\pmb{\btr \buchilu (\buomegalu)_{\barR, k}.}$ By the first estimate of $\bslashd \buomegalu$ in lemma \ref{lem 8.20} proved later which is independent of the proof of this lemma \ref{lem 8.13}, we have
\begin{align*}
	\Vert \btr \buchilu (\buomegalu)_{\barR, k} \Vert^{n-1,p}
	&\leq
	\frac{c(n,p) r_0}{\br_u^2} \Vert \bslashd \btr \buchilu \Vert^{n-1,p}
	+ \frac{C(n,p,\bubetalu) \epsilon r_0^{\frac{3}{2}}}{\br_u^{\frac{7}{2}}}
\\
	&\phantom{\leq\ }
	+ \frac{C(n,p,\bubetalu) \ud_{o,\uh} r_0^2}{\br_u^4}
	+ \frac{C(n,p,\btalu,\hatbuchilu) \delta_o^2 r_0^2}{\br_u^4}.
\end{align*}

\item
$\pmb{( \vert \hatbuchilu \vert^2)_{\barR, k}.}$ By the estimate of $\hatbuchilu$ in the bootstrap assumption, we have
\begin{align*}
	\Vert ( \vert \hatbuchilu \vert^2)_{\barR, k} \Vert^{n-1,p}
	\leq
	\frac{C(n,p,\hatbuchilu) ( \epsilon + \delta_o + \ud_{o,\uh})^2 r_0^2}{\br_u^4}.
\end{align*}
\end{enumerate}
Combining the above estimates, we obtain that
\begin{align*}
	\Vert  \buL (\btr \buchilu)_{\barR,k} + \btr \buchilu (\btr \buchilu)_{\barR, k} \Vert^{n-1,p}
	&\leq
	\frac{c(n,p) r_0}{\br_u^2} \Vert \bslashd \btr \buchilu \Vert^{n-1,p}
	+ \frac{C(n,p,\bubetalu) \epsilon r_0^{\frac{3}{2}}}{\br_u^{\frac{7}{2}}}
\\
	&\phantom{\leq\ }
	+ \frac{C(n,p,\bubetalu) \ud_{o,\uh} r_0^2}{\br_u^4}
	+ \frac{C(n,p,\btalu,\hatbuchilu) \delta_o^2 r_0^2}{\br_u^4}.
\end{align*}
By propagation lemma \ref{lem c.2} in appendix \ref{appen p.l.}, we obtain that
\begin{align*}
	\br_u^2 \Vert \bslashd \btr \buchilu \Vert^{n-1,p}
	&\leq
	\int_0^u \frac{c(n,p) r_0}{\br_{u'}^2} \cdot ( \br_{u'}^2 \Vert \bslashd \btr \buchil{u'} \Vert^{n-1,p}) \ed u'
	+ c(n,p) r_0^2 \Vert \bslashd \btr \buchil{u=0} \Vert^{n-1,p}
\\
	&\phantom{\leq\ }
	+ C(n,p,\bubetalu) \epsilon r_0
	+ C(n,p,\bubetalu) \ud_{o,\uh} r_0
	+ C(n,p,\btalu,\hatbuchilu) \delta_o^2 r_0.
\end{align*}
Since $\int_0^u \frac{c(n,p) r_0}{\br_{u'}^2}\ed u' \leq c(n,p)$, then by Gronwall's inequality, we have
\begin{align*}
	\br_u^2 \Vert \bslashd \btr \buchilu \Vert^{n-1,p}
	&\leq
	C(n,p,\bubetalu) (\epsilon + \ud_{o,\uh}) r_0
	+ c(n,p) \delta_o r_0.
\end{align*}
The lemma follows.
\end{proof}

\subsubsection{Estimate of shear $\hatbuchilu$}
We estimate $\hatbuchilu$ by equation \eqref{eqn 6.10}, which we rewrite as follows
\begin{align}
	\bslashdiv \hatbuchilu
	=
	\frac{1}{2} \bslashd \btr \buchilu
	+ \hatbuchilu \cdot \btalu 
	- \frac{1}{2} \btr \buchilu \, \btalu
	- \bubetalu.
	\label{eqn 6.10'}
	\tag{\ref{eqn 6.10}\ensuremath{'}}
\end{align}
\begin{lemma}\label{lem 8.14}
By equation \eqref{eqn 6.10},
\begin{align*}
	\Vert \hatbuchilu \Vert^{n,p}
	&\leq
	C(n,p,\bslashd \btr \buchilu, \btalu)(\epsilon + \delta_o + \ud_{o,\uh}) r_0
	+ \frac{C(n,p,\bubetalu) \epsilon r_0^{\frac{3}{2}}}{\br_u^{\frac{1}{2}}}
	+ \frac{C(n,p,\bubetalu) \ud_{o,\uh} r_0^2}{\br_u}.
\end{align*}
\end{lemma}
\begin{proof}
We estimate each term on the right hand side of equation \eqref{eqn 6.10'}.
\begin{enumerate}[label=\textbf{\textit{\alph*.}}]
\item
$\pmb{\bslashd \btr \buchilu.}$
By the estimate of $\bslashd \btr \buchilu$ in the bootstrap assumption, we have
\begin{align*}
	\Vert \bslashd \btr \buchilu \Vert^{n-1,p}
	\leq
	\frac{c(n,p,\bslashd \btr \buchilu)(\epsilon + \delta_o + \ud_{o,\uh}) r_0}{\br_u^2}.
\end{align*}

\item
$\pmb{\hatbuchilu \cdot \btalu.}$
By the estimate of $\btalu$ in the bootstrap assumption, we have
\begin{align*}
	\Vert \hatbuchilu \cdot \btalu \Vert^{n,p}
	\leq
	\frac{C(n,p,\btalu)(\epsilon + \delta_o + \ud_{o,\uh}) r_0}{\br_u^3}
	\Vert \hatbuchilu \Vert^{n,p}.
\end{align*}

\item
$\pmb{\btr \buchilu \, \btalu.}$
By the estimate of $\btalu$ in the bootstrap assumption, we have
\begin{align*}
	\Vert \btr \buchilu \, \btalu \Vert^{n,p}
	\leq
	\frac{C(n,p,\btalu)(\epsilon + \delta_o + \ud_{o,\uh}) r_0}{\br_u^2}.
\end{align*}

\item
$\pmb{\bubetalu.}$
By the estimate of $\bubetalu$ in the bootstrap assumption, we have
\begin{align*}
	\Vert \bubetalu \Vert^{n,p}
	\leq
	\frac{c(n,p,\bubetalu) \epsilon r_0^{\frac{3}{2}}}{\br_u^{\frac{5}{2}}}
	+ \frac{c(n,p,\bubetalu) \ud_{o,\uh} r_0^2}{\br_u^3}.
\end{align*}
\end{enumerate}
Combining the above estimates, we obtain that
\begin{align*}
	\Vert \bslashdiv \hatbuchilu \Vert^{n-1,p}
	&\leq
	\frac{C(n,p,\btalu)(\epsilon + \delta_o + \ud_{o,\uh}) r_0}{\br_u}
	\Vert \hatbuchilu \Vert^{n,p}
	+ \frac{C(n,p,\bslashd \btr \buchilu, \btalu)(\epsilon + \delta_o + \ud_{o,\uh}) r_0}{\br_u^2}
\\
	&\phantom{\leq\ }
	+ \frac{c(n,p,\bubetalu) \epsilon r_0^{\frac{3}{2}}}{\br_u^{\frac{5}{2}}}
	+ \frac{c(n,p,\bubetalu) \ud_{o,\uh} r_0^2}{\br_u^3}.
\end{align*}
By the theory of elliptic equations on the sphere and the estimate of $\bslashglu$ in the bootstrap assumption, we obtain that
\begin{align*}
	\Vert \hatbuchilu \Vert^{n,p}
	&\leq
	\frac{C(n,p,\btalu)(\epsilon + \delta_o + \ud_{o,\uh}) r_0}{\br_u}
	\Vert \hatbuchilu \Vert^{n,p}
	+ C(n,p,\bslashd \btr \buchilu, \btalu)(\epsilon + \delta_o + \ud_{o,\uh}) r_0
\\
	&\phantom{\leq\ }
	+ \frac{C(n,p,\bubetalu) \epsilon r_0^{\frac{3}{2}}}{\br_u^{\frac{1}{2}}}
	+ \frac{C(n,p,\bubetalu) \ud_{o,\uh} r_0^2}{\br_u}.
\end{align*}
Since $C(n,p,\btalu)(\epsilon + \delta_o + \ud_{o,\uh}) \leq C(n,p,\btalu) \delta \leq  \frac{1}{2}$, the lemma follows.
\end{proof}

\subsubsection{Estimate of null expansion $\overline{\btr \bchilu'}^{u} - (\btr \bchilu')_S$}
We derive the following propagation equation for $\overline{\btr \bchilu'}$ from equation \eqref{eqn 6.8} of $\btr \bchilu'$
\begin{align*}
	\begin{aligned}
		\frac{\ed}{\ed u} \overline{\btr \bchilu'}^{u}
		&=
		-2 \overline{ \buomegalu}^{u}\, \overline{\btr \bchilu' }^{u}
		- 2 \overline{ (\buomegalu -  \overline{ \buomegalu}^{u}) (\btr \bchilu' -  \overline{\btr \bchilu' }^{u})}^{u}
		- \frac{1}{2} \overline{\btr \buchilu}^{u} \overline{\btr \bchilu' }^{u}
	\\
		&\phantom{=\ }
		- 2 \overline{| \btalu |^2}^{u}
		+ 2 \bmulu
		+ \frac{1}{2} \overline{(\btr \buchilu - \overline{\btr \buchilu}^{u})(\btr \bchilu' - \overline{\btr \bchilu'}^{u})}^{u}.
	\end{aligned}
\end{align*}
In the Schwarzschild spacetime, we have
\begin{align*}
	\frac{\ed}{\ed u} (\btr \bchilu')_S
	=
	- \frac{1}{2} (\btr \buchilu)_S \, (\btr \bchilu')_S
	+ 2 \bmulu_S.
\end{align*}
Thus taking the difference of the above equations, we derive the propagation equation of $\overline{\btr \bchilu'}^{u} - (\btr \bchilu')_S$ that
\begin{align}
	\begin{aligned}
		\frac{\ed}{\ed u} ( \overline{\btr \bchilu'}^{u} - (\btr \bchilu')_S)
		&=
		- \frac{1}{2} \overline{\btr \buchilu}^{u} [\overline{\btr \bchilu' }^{u} -(\btr \bchilu')_S]
		- \frac{1}{2} [\overline{\btr \buchilu}^{u} - (\btr \buchilu)_S ] (\btr \bchilu')_S
	\\
		&\phantom{=\ }
		- 2 \overline{ \buomegalu}^{u}\, \overline{\btr \bchilu' }^{u}
		- 2 \overline{ (\buomegalu -  \overline{ \buomegalu}^{u}) (\btr \bchilu' -  \overline{\btr \bchilu' }^{u})}^{u}
	\\
		&\phantom{=\ }
		- 2 \overline{| \btalu |^2}^{u} 
		+ 2 (\bmulu - \bmulu_S)
	\\
		&\phantom{=\ }
		+ \frac{1}{2} \overline{(\btr \buchilu - \overline{\btr \buchilu}^{u})(\btr \bchilu' - \overline{\btr \bchilu'}^{u})}^{u}.
	\end{aligned}
	\label{eqn 8.12}
\end{align}
We obtain the estimate of $\overline{\btr \bchilu'}^{u} - (\btr \bchilu')_S$ by integrating equation \eqref{eqn 8.12}.
\begin{lemma}\label{lem 8.15}
By equation \eqref{eqn 8.12}, we have
\begin{align*}
	\vert \overline{\btr \bchilu'}^{u} - (\btr \bchilu')_S \vert
	\leq
	\frac{C(n,p,\bmulu, \overline{\btr \buchilu}^{u})(\epsilon + \delta_o + \delta_m)}{\br_u}
	+ \frac{C(n,p,\bmulu) \ud_{o,\uh} u}{\br_u^2}.
\end{align*}
\end{lemma}
\begin{proof}
Substituting $\overline{\btr \buchilu}^{u} = \frac{2}{\br_u}$ and rewrite equation \eqref{eqn 8.12} as
\begin{align}
	\begin{aligned}
		\frac{\ed}{\ed u} [\br_u ( \overline{\btr \bchilu'}^{u} - (\btr \bchilu')_S)]
		&=
		- \frac{\br_u}{2} [\overline{\btr \buchilu}^{u} - (\btr \buchilu)_S ] (\btr \bchilu')_S
	\\
		&\phantom{=\ }
		- 2 \br_u \overline{ \buomegalu}^{u}\, \overline{\btr \bchilu' }^{u}
		- 2 \br_u \overline{ (\buomegalu -  \overline{ \buomegalu}^{u}) (\btr \bchilu' -  \overline{\btr \bchilu' }^{u})}^{u}
	\\
		&\phantom{=\ }
		- 2 \br_u \overline{| \btalu |^2}^{u}
		+ 2 \br_u (\bmulu - \bmulu_S)
	\\
		&\phantom{=\ }
		+ \frac{\br_u}{2} \overline{(\btr \buchilu - \overline{\btr \buchilu}^{u})(\btr \bchilu' - \overline{\btr \bchilu'}^{u})}^{u}.
	\end{aligned}
	\tag{\ref{eqn 8.12}\ensuremath{'}}
	\label{eqn 8.12'}
\end{align}
We estimate each term on the right hand side of equation \eqref{eqn 8.12'}.
\begin{enumerate}[label=\textbf{\textit{\alph*.}}]
\item
$\pmb{\br_u [\overline{\btr \buchilu}^{u} - (\btr \buchilu)_S ] (\btr \bchilu')_S.}$ 
By the estimate of $\overline{\btr \buchilu}^{u} - (\btr \buchilu)_S$ in the bootstrap assumption,
\begin{align*}
	\vert \br_u [\overline{\btr \buchilu}^{u} - (\btr \buchilu)_S ] (\btr \bchilu')_S \vert
	\leq
	\frac{C(n,p,\overline{\btr \buchilu}^{u}) (\epsilon + \delta_o + \delta_m)r_0}{\br_u^2}.
\end{align*}

\item
$\pmb{\br_u \overline{ \buomegalu}^{u}\, \overline{\btr \bchilu' }^{u}.}$ 
By the estimate of $\overline{\buomegalu}^{u}$ in the bootstrap assumption,
\begin{align*}
	\vert \br_u \overline{ \buomegalu}^{u}\, \overline{\btr \bchilu' }^{u} \vert
	\leq
	\frac{C(n,p,\overline{\buomegalu}^{u})(\epsilon + \delta_o + \ud_{o,\uh})^2 r_0^2}{\br_u^3}.
\end{align*}

\item
$\pmb{\br_u \overline{ (\buomegalu -  \overline{ \buomegalu}^{u}) (\btr \bchilu' -  \overline{\btr \bchilu' }^{u})}^{u}.}$
By the estimates of $\bslashd \buomegalu$ and $\bslashd \btr \bchilu'$ in the bootstrap assumption,
\begin{align*}
	\vert \br_u \overline{ (\buomegalu -  \overline{ \buomegalu}^{u}) (\btr \bchilu' -  \overline{\btr \bchilu' }^{u})}^{u} \vert
	\leq
	\frac{C(n,p,\bslashd \buomegalu, \bslashd \btr \bchilu') ( \epsilon + \delta_o + \ud_{o,\uh})^2 r_0^{\frac{3}{2}}}{\br_u^{\frac{5}{2}}}.
\end{align*}

\item
$\pmb{\br_u \overline{| \btalu |^2}^{u}.}$
By the estimate of $\btalu$ in the bootstrap assumption,
\begin{align*}
	\vert \br_u \overline{| \btalu |^2}^{u} \vert
	\leq
	\frac{C(n,p,\btalu) ( \epsilon + \delta_o + \ud_{o,\uh})^2 r_0^2}{\br_u^3}.
\end{align*}

\item
$\pmb{\br_u (\bmulu - \bmulu_S).}$
By the estimate of $\bmulu - \bmulu_S$ in the bootstrap assumption,
\begin{align*}
	\vert \br_u (\bmulu - \bmulu_S) \vert
	\leq
	\frac{C(n,p,\bmulu)(\epsilon + \delta_o + \delta_m) r_0}{\br_u^2}
	+ \frac{C(n,p,\bmulu) \ud_{o,\uh} u r_0}{\br_u^3}.
\end{align*}

\item
$\pmb{\br_u \overline{(\btr \buchilu - \overline{\btr \buchilu}^{u})(\btr \bchilu' - \overline{\btr \bchilu'}^{u})}^{u}.}$
By the estimates of $\bslashd \btr \buchilu$ and $\bslashd \btr \bchilu'$ in the bootstrap assumption,
\begin{align*}
	\vert \br_u \overline{(\btr \buchilu - \overline{\btr \buchilu}^{u})(\btr \bchilu' - \overline{\btr \bchilu'}^{u})}^{u} \vert
	\leq
	\frac{C(n,p, \bslashd \btr \buchilu, \bslashd \btr \bchilu')(\epsilon + \delta_o + \ud_{o,\uh})^2 r_0}{\br_u^2}.
\end{align*}
\end{enumerate}

Combining the above estimates, we obtain that
\begin{align*}
	\vert \frac{\ed}{\ed u} [\br_u ( \overline{\btr \bchilu'}^{u} - (\btr \bchilu')_S)] \vert
	\leq
	\frac{C(n,p, \bmulu, \overline{\btr \buchilu}^{u})(\epsilon + \delta_o + \delta_m) r_0}{\br_u^2}
	+ \frac{C(n,p, \bmulu) \ud_{o,\uh} r_0}{\br_u^2}.
\end{align*}
Therefore by integrating the above inequality and the estimate $\btr \bchil{u=0}'$ on the initial leaf, we obtain that
\begin{align*}
	\vert \br_u ( \overline{\btr \bchilu'}^{u} - (\btr \bchilu')_S) \vert
	\leq
	C(n,p, \bmulu, \overline{\btr \buchilu}^{u})(\epsilon + \delta_o + \delta_m) 
	+ \frac{C(n,p, \bmulu) \ud_{o,\uh} u }{\br_u}.
\end{align*}
The lemma follows.
\end{proof}

\subsubsection{Estimate of null expansion $\bslashd \btr \bchilu'$}
We derive the propagation equation for the rotational derivative of $\btr \bchilu'$ from the propagation equation \eqref{eqn 6.9},
\begin{align}
	\begin{aligned}
		&\phantom{=\ }
		\buL (\btr \bchilu')_{\barR,k}
		+ \frac{1}{2} \btr \buchilu (\btr \bchilu')_{\barR,k}
	\\
		&=
		- [ \barR_k, \balu \db^i \ddpartial_i ] \btr \bchilu'
		- \frac{1}{2} \btr \bchilu' (\btr \buchilu)_{\barR, k}
		- 2 \buomegalu (\btr \bchilu')_{\barR,k}
		- 2 \btr \bchilu' (\buomegalu)_{\barR,k}
		- 2 (\vert \btalu \vert^2)_{\barR,k}.
	\end{aligned}
	\label{eqn 8.13}
\end{align}
Integrate the above equation by propagation lemma \ref{lem c.2} to obtain the estimate of $\bslashd \btr \bchilu'$.
\begin{lemma}\label{lem 8.16}
By equation \eqref{eqn 8.13},
\begin{align*}
	\Vert \bslashd \btr \bchilu' \Vert^{n-1,p}
	\leq
	\frac{C(n,p,\bslashd \btr \buchilu, \bslashd \buomegalu) (\epsilon + \delta_o + \ud_{o,\uh})}{\br_u}.
\end{align*}
\end{lemma}
\begin{proof}
We estimate each term on the right hand side of equation \eqref{eqn 8.13}.
\begin{enumerate}[label=\textbf{\textit{\alph*}}]
\item
$\pmb{[ \barR_k, \balu \db^i \ddpartial_i ] \btr \bchilu'.}$
By the estimate of $[ \barR_k, \balu \db^i \ddpartial_i]$ in lemma \ref{lem 8.8}, we have that
\begin{align*}
	\Vert [ \barR_k, \balu \db^i \ddpartial_i] \btr \bchilu' \Vert^{n-1,p}
	&\leq
	c(n,p) 
	[ \frac{ \ud_{\uh} r_0}{\br_u^2}
	+ \frac{\epsilon ( \ud_o + \ud_m )r_0^2}{\br_u^3}]
	\cdot
	\Vert \bslashd \btr \bchilu' \Vert^{n-1,p}.
\end{align*}

\item
$\pmb{\btr \bchilu' (\btr \buchilu)_{\barR, k}.}$
By the estimate of $\bslashd \btr \buchilu$ in the bootstrap assumption,
\begin{align*}
	\Vert \btr \bchilu' (\btr \buchilu)_{\barR, k} \Vert^{n-1,p}
	\leq
	\frac{C(n,p, \bslashd \btr \buchilu) (\epsilon + \delta_o + \ud_{o,\uh} ) r_0}{\br_u^3}. 
\end{align*}

\item
$\pmb{\buomegalu (\btr \bchilu')_{\barR,k}.}$
By the estimate of $\buomegalu$ in the bootstrap assumption, we have
\begin{align*}
	\Vert \buomegalu  (\btr \bchilu')_{\barR,k} \Vert^{n-1,p}
	\leq
	C(n,p,\bslashd \buomegalu) 
	[\frac{\epsilon r_0^{\frac{3}{2}}}{\br_u^{\frac{5}{2}}}
	+ \frac{ (\delta_o + \ud_{o,\uh}) r_0^2}{\br_u^3}]
	\cdot
	\Vert \bslashd \btr \bchilu' \Vert^{n-1,p}.
\end{align*}

\item
$\pmb{\btr \bchilu' (\buomegalu)_{\barR,k}.}$
 By the estimate of $\bslashd \buomegalu$ in the bootstrap assumption, we have
\begin{align*}
	\Vert \btr \bchilu' (\buomegalu)_{\barR,k} \Vert^{n-1,p}
	\leq
	C(n,p,\bslashd \buomegalu) [\frac{\epsilon r_0^{\frac{3}{2}}}{\br_u^{\frac{7}{2}}}
	+ \frac{ (\delta_o + \ud_{o,\uh}) r_0^2}{\br_u^4}].
\end{align*}

\item
$\pmb{(\vert \btalu \vert^2)_{\barR,k}.}$
By the estimate of $\btalu$ in the bootstrap assumption,
\begin{align*}
	\Vert (\vert \btalu \vert^2)_{\barR,k} \Vert^{n,p}
	\leq
	\frac{C(n,p, \btalu) (\epsilon + \delta_o + \ud_{o,\uh})^2 r_0^2}{\br_u^4}.
\end{align*}
\end{enumerate}
Combining the above estimates, we obtain that
\begin{align*}
	\Vert \buL (\btr \bchilu')_{\barR,k} + \frac{1}{2} \btr \buchilu ( \btr \bchilu' )_{\barR,k} \Vert^{n-1,p}
	&\leq
	C(n,p,\bslashd \buomegalu)
	[\frac{\epsilon r_0^{\frac{3}{2}}}{\br_u^{\frac{5}{2}}}
	+ \frac{\delta_o r_0^2}{\br_u^3}
	+ \frac{ \ud_{o,\uh} r_0}{\br_u^2}]
	\cdot
	\Vert \bslashd \btr \bchilu' \Vert^{n-1,p}
\\
	&\phantom{\leq\ }
	+ \frac{C(n,p,\bslashd \btr \buchilu, \bslashd \buomegalu) (\epsilon + \delta_o + \ud_{o,\uh}) r_0}{\br_u^3}.
\end{align*}
By propagation lemma \ref{lem c.2},
\begin{align*}
	\br_u \Vert \bslashd \btr \bchilu' \Vert^{n-1,p}
	&\leq
	\int_0^u \frac{C(n,p,\bslashd \btr \buchilu, \bslashd \buomegalu) ( \epsilon + \delta_o + \ud_{o,\uh}) r_0}{\br_{u'}^2} (\br_{u'} \Vert \bslashd \btr \bchil{u'}' \Vert^{n-1,p}) \ed u'
\\
	&\phantom{\leq\ }
	+ c(n,p) r_0 \Vert \bslashd \btr \bchil{u=0}' \Vert^{n-1,p}
	+ C(n,p,\bslashd \btr \buchilu, \bslashd \buomegalu) (\epsilon + \delta_o + \ud_{o,\uh}).
\end{align*}
Therefore by Gronwall's inequality, we obtain
\begin{align*}
	\br_u \Vert \bslashd \btr \bchilu' \Vert^{n-1,p}
	\leq
	C(n,p,\bslashd \btr \buchilu, \bslashd \buomegalu) (\epsilon + \delta_o + \ud_{o,\uh}).
\end{align*}
The lemma follows.
\end{proof}

\subsubsection{Estimate of shear $\hatbchilu'$}
We estimate $\hatbchilu'$ by equation \eqref{eqn 6.11}, which we rewrite as follows
\begin{align}
	\bslashdiv \hatbchilu'
	=
	\frac{1}{2} \bslashd \btr \bchilu'
	- \hatbchilu' \cdot \btalu
	+ \frac{1}{2} \btr \bchilu'\ \btalu
	-\bbetalu.
	\tag{\ref{eqn 6.11}\ensuremath{'}}
	\label{eqn 6.11'}
\end{align}
\begin{lemma}\label{lem 8.17}
By equation \eqref{eqn 6.11},
\begin{align*}
	\Vert \hatbchilu' \Vert^{n,p}
	&\leq
	C(n,p,\bslashd \btr \bchilu', \bbetalu) \epsilon \br_u
	+ C(n,p,\btalu) \epsilon r_0
\\
	&\phantom{\leq\ }
	+ C(n,p,\bslashd \btr \bchilu') (\delta_o + \ud_{o,\uh}) \br_u
	+ C(n,p, \btalu, \bbetalu) (\delta_o + \ud_{o,\uh}) r_0.
\end{align*}
\end{lemma}
\begin{proof}
We estimate each term on the right hand side of \eqref{eqn 6.11'}.
\begin{enumerate}[label=\textbf{\textit{\alph*.}}]
\item
$\pmb{\bslashd \btr \bchilu'.}$
By the bootstrap assumption,
\begin{align*}
	\Vert \bslashd \btr \bchilu' \Vert^{n-1,p}
	\leq
	\frac{c(n,p,\bslashd \btr \bchilu') (\epsilon + \delta_o + \ud_{o,\uh})}{\br_u}.
\end{align*}

\item
$\pmb{\hatbchilu' \cdot \btalu.}$
By the estimate of $\btalu$ in the bootstrap assumption,
\begin{align*}
	\Vert \hatbchilu' \cdot \btalu \Vert^{n-1,p}
	\leq
	\frac{C(n,p,\btalu) ( \epsilon + \delta_o + \ud_{o,\uh}) r_0}{\br_u^3} \Vert \hatbchilu' \Vert^{n-1,p}.
\end{align*}

\item
$\pmb{\btr \bchilu'\ \btalu.}$
By the estimates of $\btr \bchilu'$ and $\btalu$ in the bootstrap assumption,
\begin{align*}
	\Vert \btr \bchilu'\ \btalu \Vert^{n-1,p}
	\leq
	\frac{C(n,p,\btalu) ( \epsilon + \delta_o + \ud_{o,\uh}) r_0}{\br_u^2}.
\end{align*}

\item
$\pmb{\bbetalu.}$
By the bootstrap assumption,
\begin{align*}
	\Vert \bbetalu \Vert^{n,p}
	\leq
	\frac{c(n,p,\bbetalu) \epsilon}{\br_u}
	+ \frac{c(n,p,\bbetalu) (\delta_o + \ud_{o,\uh}) r_0}{\br_u^2}.
\end{align*}
\end{enumerate}
Combining the above estimates, we obtain that
\begin{align*}
	\Vert \bslashdiv \hatbchilu' \Vert^{n-1,p}
	&\leq
	\frac{C(n,p,\btalu) ( \epsilon + \delta_o + \ud_{o,\uh}) r_0}{\br_u^3} \Vert \hatbchilu' \Vert^{n-1,p}
\\
	&\phantom{\leq\ }
	+ \frac{C(n,p,\bslashd \btr \bchilu',\bbetalu) \epsilon}{\br_u}
	+ \frac{C(n,p,\btalu) \epsilon r_0}{\br_u^2}
\\
	&\phantom{\leq\ }
	+ \frac{C(n,p,\bslashd \btr \bchilu') (\delta_o + \ud_{o,\uh})}{\br_u}
	+ \frac{C(n,p,\bbetalu, \btalu) (\delta_o + \ud_{o,\uh}) r_0}{\br_u^2}.
\end{align*}
By the theory of elliptic equations on the sphere and the estimate of $\bslashglu$, the lemma follows.
\end{proof}

\subsubsection{Estimate of torsion $\btalu$}
We estimate $\btalu$ by equation \eqref{eqn 6.12}, which we rewrite as follows
\begin{align}
	\left\{
	\begin{aligned}
		&
		\bslashcurl \btalu
		=
		\frac{1}{2} \hatbchilu' \wedge \hatbuchilu 
		+ \bsigmalu,
	\\
		&
		\bslashdiv \btalu
		= 
		- (\brholu - \overline{\brholu}^u)
		- \frac{1}{2} [( \hatbuchilu, \hatbchilu' ) - \overline{( \hatbuchilu, \hatbchilu' )}^u].
	\end{aligned}
	\right.
	\tag{\ref{eqn 6.12}\ensuremath{'}}
	\label{eqn 6.12'}
\end{align}
\begin{lemma}\label{lem 8.18}
By equation \eqref{eqn 6.12},
\begin{align*}
	\Vert \btalu \Vert^{n+1,p}
	\leq
	\frac{C(n,p, \bslashd \brholu,\bsigmalu) \epsilon r_0}{\br_u}
	+ \frac{C(n,p,\bslashd \brholu) (\delta_o + \ud_{o,\uh}) r_0}{\br_u}.
\end{align*}
\end{lemma}
\begin{proof}
We estimate each term on the right hand side of equation \eqref{eqn 6.12'}.
\begin{enumerate}[label=\textbf{\textit{\alph*.}}]
\item
$\pmb{\hatbchilu' \wedge \hatbuchilu.}$
By the estimates of $\hatbchilu'$ and $\hatbchilu$ in the bootstrap assumption,
\begin{align*}
	\Vert \hatbchilu' \wedge \hatbuchilu \Vert^{n,p}
	\leq
	\frac{C(n,p,\hatbchilu', \hatbuchilu)(\epsilon + \delta_o + \ud_{o,\uh})^2 r_0}{\br_u^3}.
\end{align*}

\item
$\pmb{\bsigmalu.}$
By the bootstrap assumption,
\begin{align*}
	\Vert \bsigmalu \Vert^{n,p}
	\leq
	\frac{c(n,p, \bsigmalu)(\epsilon + \delta_o \ud_{o,\uh} + \ud_{o,\uh}^2) r_0}{\br_u^3}.
\end{align*}

\item
$\pmb{\brholu - \overline{\brholu}^u.}$
By the estimate of $\bslashd \brholu$ in the bootstrap assumption,
\begin{align*}
	\Vert \brholu - \overline{\brholu}^u \Vert^{n,p}
	\leq
	\frac{C(n,p,\brholu) (\epsilon + \delta_o + \ud_{o,\uh}) r_0}{\br_u^3}.
\end{align*}

\item
$\pmb{( \hatbuchilu, \hatbchilu' ) - \overline{( \hatbuchilu, \hatbchilu' )}^u.}$
By the estimates of $\hatbchilu'$ and $\hatbchilu$ in the bootstrap assumption,
\begin{align*}
	\Vert ( \hatbuchilu, \hatbchilu' ) - \overline{( \hatbuchilu, \hatbchilu' )}^u \Vert^{n,p}
	\leq
	\frac{C(n,p,\hatbchilu', \hatbuchilu)(\epsilon + \delta_o + \ud_{o,\uh})^2 r_0}{\br_u^3}.
\end{align*}
\end{enumerate}
Combining the above estimates, we obtain that
\begin{align*}
	\Vert \bslashdiv \btalu \Vert^{n,p},
	\Vert \bslashcurl \btalu \Vert^{n,p}
	\leq
	\frac{C(n,p, \bslashd \brholu,\bsigmalu) \epsilon r_0}{\br_u^3}
	+ \frac{C(n,p,\bslashd \brholu) (\delta_o + \ud_{o,\uh}) r_0}{\br_u^3}.
\end{align*}
By the theory of elliptic equations on the sphere and the estimate of $\bslashglu$ in the bootstrap assumption, the lemma follows from the above estimates.
\end{proof}

\subsubsection{Estimate of acceleration $\overline{\buomegalu}^u$}
We estimate $\overline{\buomegalu}^u$ by equation \eqref{eqn 6.13} that
\begin{align*}
	\overline{\buomegalu}^u
	&=
	-\frac{r_u}{2} \overline{(\buomegalu - \overline{\buomegalu}^u)(\btr\buchilu - \overline{\btr\buchilu}^u)}^u
	- \frac{r_u}{8} \overline{(\btr \buchilu - \overline{\btr \buchilu}^u)^2}^u
	+ \frac{r_u}{4} \overline{|\hatbuchilu|^2}^u.
\end{align*}
\begin{lemma}\label{lem 8.19}
By equation \eqref{eqn 6.13}, we have that
\begin{align*}
	\vert \overline{\buomegalu}^{u} \vert
	&\leq
	C(n,p,\bslashd \btr \buchilu, \bslashd \buomegalu) (\epsilon + \delta_o + \ud_{o,\uh}) [\frac{\epsilon r_0^{\frac{5}{2}}}{\br_u^{\frac{7}{2}}} + \frac{(\delta_o + \ud_{o,\uh}) r_0^3}{\br_u^4} ]
\\
	&\phantom{\leq\ }
	+ \frac{C(n,p,\bslashd \btr \buchilu) (\epsilon + \delta_o + \ud_{o,\uh})^2 r_0^2}{\br_u^3}
	+ \frac{C(n,p,\hatbuchilu) (\epsilon + \delta_o + \ud_{o,\uh})^2 r_0^2}{\br_u^3}
\end{align*}
\end{lemma}
\begin{proof}
The proof is straightforward by the estimates of $\bslashd \buomegalu$, $\bslashd \btr \buchilu$, $\hatbuchilu$ in the bootstrap assumption and the Sobolev inequality.
\end{proof}

\subsubsection{Estimate of acceleration $\bslashd \buomegalu$}
We shall estimate $\bslashd \buomegalu$ by the elliptic equation \eqref{eqn 6.13} which we cite as follows
\begin{align*}
	2 \bslashDelta \buomegalu
	&=
	-\frac{3}{2} \bmulu ( \btr \buchilu - \overline{\btr \buchilu}^u )
	+ \frac{1}{2} ( \btr \buchilu | \btalu |^2 - \overline{\btr \buchilu | \btalu |^2}^u )
\\
	&
	\phantom{=}
	+ \frac{1}{4} ( \btr \bchilu' | \hatbuchilu |^2 - \overline{ \btr \bchilu' | \hatbuchilu|^2}^u )
	+ 4 (\bslashdiv \hatbuchilu, \btalu )
	+ 4 ( \hatbuchilu, \bslashnabla \btalu )
	- 2 \bslashdiv \bubetalu.
\end{align*}
\begin{lemma}\label{lem 8.20}
By the elliptic equation \eqref{eqn 6.13}, we have that
\begin{align*}
	\Vert \bslashd \buomegalu \Vert^{n,p}
	\leq
	\frac{c(n,p) r_0}{\br_u} \Vert \bslashd \btr \buchilu \Vert^{n-1,p}
	+ \frac{C(n,p,\bubetalu) \epsilon r_0^{\frac{3}{2}}}{\br_u^{\frac{5}{2}}}
	+ \frac{C(n,p,\bubetalu) \ud_{o,\uh} r_0^2}{\br_u^3}
	+ \frac{C(n,p,\btalu,\hatbuchilu) \delta_o^2 r_0^2}{\br_u^3}.
\end{align*}
Substituting the estimate of $\bslashd \btr \buchilu$ from the bootstrap assumption, we have that
\begin{align*}
	\Vert \bslashd \buomegalu \Vert^{n,p}
	\leq
	\frac{C(n,p, \bubetalu) \epsilon r_0^{\frac{3}{2}}}{\br_u^{\frac{5}{2}}}
	+  \frac{C(n,p, \bslashd \btr \buchilu) (\epsilon + \delta_o) r_0^2}{\br_u^3}
	+  \frac{C(n,p, \bubetalu, \bslashd \btr \buchilu) \ud_{o,\uh} r_0^2}{\br_u^3}.
\end{align*}
\end{lemma}
\begin{proof}
We estimate each term on the right hand side of equation \eqref{eqn 6.13}.
\begin{enumerate}[label=\textbf{\textit{\alph*.}}]
\item
$\pmb{\bmulu ( \btr \buchilu - \overline{\btr \buchilu}^u ).}$ By the Sobolev inequality on the sphere, we have that
\begin{align*}
	\Vert \bmulu ( \btr \buchilu - \overline{\btr \buchilu}^u ) \Vert^{n,p}
	\leq
	\frac{c(n,p) r_0}{\br_u^3} \Vert \bslashd \btr \buchilu \Vert^{n-1,p}.
\end{align*}

\item
$\pmb{\btr \buchilu | \btalu |^2 - \overline{\btr \buchilu | \btalu |^2}^u.}$ By the Sobolev inequality on the sphere, we have that
\begin{align*}
	\Vert \btr \buchilu | \btalu |^2 - \overline{\btr \buchilu | \btalu |^2}^u \Vert^{n,p}
	\leq
	c(n,p) \Vert \bslashd (\btr \buchilu | \btalu |^2) \Vert^{n-1,p}
	\leq
	\frac{C(n,p,\btalu) (\epsilon + \delta_o + \ud_{o,\uh})^2 r_0^2}{\br_u^5} .
\end{align*}

\item
$\pmb{\btr \bchilu' | \hatbuchilu |^2 - \overline{ \btr \bchilu' | \hatbuchilu|^2}^u.}$ By the Sobolev inequality on the sphere, we have that
\begin{align*}
	\Vert \btr \bchilu' | \hatbuchilu |^2 - \overline{\btr \bchilu' | \hatbuchilu |^2}^u \Vert^{n,p}
	\leq
	c(n,p) \Vert \bslashd (\btr \bchilu' | \hatbuchilu |^2) \Vert^{n-1,p}
	\leq
	\frac{C(n,p,\hatbuchilu) (\epsilon + \delta_o + \ud_{o,\uh})^2 r_0^2}{\br_u^5}.
\end{align*}

\item
$\pmb{(\bslashdiv \hatbuchilu, \btalu ), ( \hatbuchilu, \bslashnabla \btalu ).}$ We have that
\begin{align*}
	\Vert (\bslashdiv \hatbuchilu, \btalu ) \Vert^{n-1,p},
	\Vert ( \hatbuchilu, \bslashnabla \btalu ) \Vert^{n-1,p}
	&\leq
	c(n,p) \Vert \hatbuchilu \Vert^{n,p} \cdot \Vert \btalu \Vert^{n,p}
\\
	&\leq
	\frac{C(n,p,\hatbuchilu, \btalu) (\epsilon + \delta_o + \ud_{o,\uh})^2 r_0^2}{\br_u^5}.
\end{align*}

\item
$\pmb{\bslashdiv \bubetalu.}$ We have that
\begin{align*}
	\Vert \bslashdiv \bubetalu \Vert^{n-1,p}
	\leq
	\frac{C(n,p,\bubetalu) \epsilon r_0^{\frac{3}{2}}}{\br_u^{\frac{9}{2}}}
	+ \frac{C(n,p,\bubetalu) \ud_{o,\uh} r_0^2}{\br_u^5}
\end{align*}
\end{enumerate}

Combining the above estimates, we obtain that
\begin{align*}
	\Vert \bslashDelta \buomegalu \Vert^{n-1,p}
	\leq
	\frac{c(n,p) r_0}{\br_u^3} \Vert \bslashd \btr \buchilu \Vert^{n-1,p}
	+ \frac{C(n,p,\bubetalu) \epsilon r_0^{\frac{3}{2}}}{\br_u^{\frac{9}{2}}}
	+ \frac{C(n,p,\bubetalu) \ud_{o,\uh} r_0^2}{\br_u^5}
	+ \frac{C(n,p,\btalu,\hatbuchilu) \delta_o^2 r_0^2}{\br_u^5}.
\end{align*}
Then by the theory of the elliptic equation on the sphere and the estimate of $\bslashglu$ on $\bSigma_u$, we obtain that
\begin{align*}
	\Vert \bslashd \buomegalu \Vert^{n,p}
	\leq
	\frac{c(n,p) r_0}{\br_u} \Vert \bslashd \btr \buchilu \Vert^{n-1,p}
	+ \frac{C(n,p,\bubetalu) \epsilon r_0^{\frac{3}{2}}}{\br_u^{\frac{5}{2}}}
	+ \frac{C(n,p,\bubetalu) \ud_{o,\uh} r_0^2}{\br_u^3}
	+ \frac{C(n,p,\btalu,\hatbuchilu) \delta_o^2 r_0^2}{\br_u^3}.
\end{align*}
Substituting the estimate of $\bslashd \btr \buchilu$ in the bootstrap assumption, the lemma follows.
\end{proof}

\subsection{Estimate of mass aspect function}\label{sec 8.7}
We shall estimate $\bmulu - \bmulu_S$ by the following propagation equation derived from equation \eqref{eqn 6.7} 
\begin{align*}
	\frac{\ed}{\ed u} (\bmulu - \bmulu_S)
	&=
	- \frac{3}{2} \overline{\btr \buchilu}^u (\bmulu - \bmulu_S)
	- \frac{3}{2} (\overline{\btr \buchilu}^u - (\btr \buchilu)_S) \bmulu_S
	+ \frac{1}{4} \overline{\btr \bchilu' | \hatbuchilu |^2}^u
	+ \frac{1}{2} \overline{\btr \buchilu | \btalu |^2}^u.
\end{align*}
\begin{lemma}\label{lem 8.21}
By equation \eqref{eqn 6.7}, we have that
\begin{align*}
	\vert \bmulu - \bmulu_S \vert
	\leq
	\frac{C(n,p,\overline{\btr \buchilu}^u) (\epsilon + \delta_o + \delta_m) r_0}{\br_u^3}.
\end{align*}
\end{lemma}
\begin{proof}
Since $\overline{\btr \buchilu}^u = \frac{2}{\br_u}$, we obtain that
\begin{align*}
	\frac{\ed}{\ed u}[ \br_u^3 (\bmulu - \bmulu_S) ]
	=- \frac{3 \br_u^3 }{2} (\overline{\btr \buchilu}^u - (\btr \buchilu)_S) \bmulu_S
	+ \frac{\br_u^3}{4} \overline{\btr \bchilu' | \hatbuchilu |^2}^u
	+ \frac{\br_u^3}{2} \overline{\btr \buchilu | \btalu |^2}^u.
\end{align*}
By the estimates of $\overline{\btr \buchilu}^u - (\btr \buchilu)_S$, $\hatbuchilu$, $\btalu$ in the bootstrap assumption, we have that
\begin{align*}
	\vert \frac{\ed}{\ed u}[ \br_u^3 (\bmulu - \bmulu_S) ] \vert
	\leq
	\frac{C(n,p,\overline{\btr \buchilu}^u) (\epsilon + \delta_o + \delta_m) r_0^2}{\br_u^2}.
\end{align*}
At the initial leaf $u=0$, by proposition \ref{prop 7.10}, $\vert \br_{u=0}^3 (\bmulu - \bmulu_S)  \vert \leq c(n,p)(\epsilon + \delta_o + \delta_m) r_0$. Therefore the lemma follows from integrating $\frac{\ed}{\ed u}[ \br_u^3 (\bmulu - \bmulu_S) ]$.
\end{proof}

\subsection{Estimates of curvature components}\label{sec 8.8}
We estimate the curvature components by propositions \ref{prop 5.8}, \ref{prop 5.10} and the estimates of parameterisation functions in the bootstrap assumption.
\begin{lemma}\label{lem 8.22}
For the curvature components, we have
\begin{align*}
	&
	\Vert \bubetalu \Vert^{n,p}
	\leq
	\frac{c(n,p) \epsilon r_0^{\frac{3}{2}}}{\br_u^{\frac{5}{2}}}
	+ \frac{c(n,p) \ud_{o,\uh} r_0^2}{\br_u^3}.
\\
	&
	\vert \brholu - \brholu_S \vert
	\leq
	\frac{C(n,p,\flu) (\epsilon + \delta_o + \delta_m + \ud_{o,\uh} + \uslashd_m)r_0}{\br_u^3}.
\\
	&
	\Vert \bslashd \brholu \Vert^{n-1,p}
	\leq
	\frac{C(n,p,\bslashd \flu) ( \epsilon + \delta_o + \ud_{o,\uh}) r_0}{\br_u^3}.
\\
	&
	\Vert \bsigmalu \Vert^{n,p}
	\leq
	\frac{c(n,p) \epsilon r_0}{\br_u^3}
	+ \frac{C(n,p, \bslashd \flu)(\delta_o \ud_{o,\uh} + \ud_{o,\uh}^2) r_0^2}{\br_u^4}.
\\
	&
	\Vert \bbetalu \Vert^{n,p}
	\leq
	\frac{c(n,p) \epsilon}{\br_u}
	+ \frac{C(n,p,\bslashd \flu)( \delta_o + \ud_{o,\uh}) r_0}{\br_u^2}.
\end{align*}
\end{lemma}
\begin{proof}
We estimate each curvature component in the following. Recall the notation $d_o$ introduced in equation \eqref{eqn 8.3} that $d_o = c(n,p,\bslashd \flu)(\epsilon + \delta_o + \ud_{o,\uh})$.
\begin{enumerate}
\item[$\pmb{\bubetalu.}$]
We have $\bubetalu = \balu \cdot \ddubeta|_{\bSigma_u}$. By propositions \ref{prop 5.8}, \ref{prop 5.10},
\begin{align*}
	\Vert \bubetalu \Vert^{n,p}
	&\leq
	\Vert \balu \Vert^{n,p} \Vert \lo{\bubetalu} \Vert^{n,p}
	+ \Vert \balu \Vert^{n,p} \Vert \hi{\bubetalu} \Vert^{n,p}
\\
	&\leq
	\frac{c(n,p) \epsilon r_0^{\frac{3}{2}}}{\br_u^{\frac{5}{2}}}
	+ \frac{c(n,p) \ud_{o,\uh} r_0^2}{\br_u^3}.
\end{align*}

\item[$\pmb{\brholu.}$]
By the decomposition of $\brholu$, we have
\begin{align*}
	\brholu - \brholu_S
	=
	\rho_S|_{\bSigma_u} - \rho_S|_{\bSigma_{u,S}}
	+ \hi{\ddrholu}
\end{align*}
By the estimates of $\flu$, $\uflu$ in the bootstrap assumption and proposition \ref{prop 5.10},
\begin{align*}
	\vert \brholu - \brholu_S \vert
	\leq
	\frac{C(n,p,\flu) (\epsilon + \delta_o + \delta_m + \ud_{o,\uh} + \uslashd_m)r_0}{\br_u^3}.
\end{align*}

\item[$\pmb{\bslashd \brholu.}$]
By the decomposition of $\brholu$, we have
\begin{align*}
	\bslashd \brholu
	=
	\bslashd (\rho_S|_{\bSigma_u})
	+ \bslashd \hi{\ddrholu}
	=
	\bslashd \flu \cdot (\partial_s \rho_S)|_{\bSigma_u}
	+ \bslashd \uflu \cdot (\partial_{\us} \rho_S)|_{\bSigma_u}
	+ \bslashd \hi{\ddrholu}.
\end{align*}
By the estimates of $\flu$, $\uflu$ in the bootstrap assumption and proposition \ref{prop 5.10},
\begin{align*}
	\Vert \bslashd \brholu \Vert^{n-1,p}
	\leq
	\frac{C(n,p,\bslashd \flu) (\epsilon + \delta_o + \ud_{o,\uh})r_0}{\br_u^3}.
\end{align*}

\item[$\pmb{\bsigmalu.}$]
By propositions \ref{prop 5.8}, \ref{prop 5.10},
\begin{align*}
	\Vert \bsigmalu \Vert^{n,p}
	\leq
	\frac{c(n,p) \epsilon r_0}{\br_u^3}
	+ \frac{c(n,p) d_o \ud_{\uh} r_0^2}{\br_u^4}
	\leq
	\frac{c(n,p) \epsilon r_0}{\br_u^3}
	+ \frac{C(n,p,\bslashd \flu) (\delta_o \ud_{o,\uh} + \ud_{o,\uh}^2) r_0^2}{\br_u^4}.
\end{align*}

\item[$\pmb{\bbetalu.}$]
By propositions \ref{prop 5.8}, \ref{prop 5.10},
\begin{align*}
	\Vert \bbetalu \Vert^{n,p}
	\leq
	\frac{c(n,p) d_o r_0}{\br_u^2}
	+ \frac{c(n,p) \epsilon}{\br_u}
	\leq
	\frac{c(n,p) \epsilon }{\br_u}
	+ \frac{C(n,p,\bslashd \flu) (\delta_o + \ud_{o,\uh}) r_0}{\br_u^2}.
\end{align*}
\end{enumerate}
The lemma is proved.
\end{proof}

\subsection{Improvement of estimates in bootstrap assumption: choices of constants}\label{sec 8.9}
In this section, we shall carefully choose the constants $c(n,p,\cdot)$ in the bootstrap assumption such that by the estimates obtained in lemmas \ref{lem 8.5}-\ref{lem 8.22}, the estimates in the bootstrap assumption can be improved to strict inequalities. Then we choose the constant $\delta$ in theorem \ref{thm 8.2}, such that lemmas \ref{lem 8.4}-\ref{lem 8.22} hold.

\begin{figure}[h]
\centering
\begin{tikzpicture}[scale=0.8]
\draw (0,0) rectangle (1,1); \node at (0.5,0.5) {\footnotesize $\bubetalu$};
\draw[->] (1,0.5) -- (1.8,0.5); 
\draw (1.8,0) rectangle (3.2,1); \node at (2.5,0.5) {\footnotesize $\bslashd \btr \buchilu$};
\draw[->] (1,0.2) -- (1.5,0.2) -- (1.5,-1);  \draw[->] (1.8,0.2) -- (1.5,0.2) -- (1.5,-1);
\draw (1,-1) rectangle (2,-2); \node at (1.5, -1.5) {\footnotesize $\bslashd \buomegalu$};
\draw[->] (2.5,0) -- (2.5,-1.5) -- (2.8,-1.5); \draw[->] (2,-1.5) -- (2.8,-1.5);
\draw (2.8,-2) rectangle (4.2,-1); \node at (3.5,-1.5) {\footnotesize $\bslashd \btr \bchilu'$};
\draw[->] (1.5,-2) -- (1.5,-3);
\draw (1,-3) rectangle (2,-4); \node at (1.5, -3.5) {\footnotesize $\bslashd \balu$};
\draw[->] (3.5,-2) -- (3.5,-3.5) -- (4,-3.5); \draw[->] (2,-3.5) -- (4,-3.5);
\draw (4,-4) rectangle (5,-3); \node at (4.5,-3.5) {\footnotesize $\bslashd \flu$};
\draw[->] (1.5,-4) -- (1.5,-5);
\draw (1,-5) rectangle (2,-6); \node at (1.5, -5.5) {\footnotesize $\balu$};
\draw[->] (2,-5.5) -- (3,-5.5);
\draw (3,-5) rectangle (4,-6); \node at (3.5,-5.5) {\footnotesize $\flu$};
\draw[->] (4,-5.5) -- (5,-5.5);
\draw (5,-4.5) rectangle (6,-6.5); \node at (5.5,-5) {\footnotesize $\bslashglu$}; \node at (5.5,-6) {\footnotesize $\brholu$};
\draw[->] (5,-3.5) -- (6.5,-3.5);
\draw (6.5,-1) rectangle (8.5,-5);
\node at (7.5,-1.5) {\footnotesize $\bbetalu$};
\node at (7.5,-2.5) {\footnotesize $\bslashd \brholu$};
\node at (7.5,-3.5) {\footnotesize $\bsigmalu$};
\node at (7.5,-4.5) {\footnotesize $\bcircnabla \bslashglu$};
\draw (7,-2) rectangle (8,-4);
\draw[->] (8,-3) -- (9,-3);
\draw (9,-2.5) rectangle (10,-3.5); \node at (9.5,-3) {\footnotesize $\btalu$};
\draw (9,0) rectangle (12,1); \node at(10.5,0.5) {\footnotesize $\bslashd \btr \buchilu, \bubetalu,\btalu$};
\draw[->] (12,0.5) -- (13,0.5); 
\draw (13,0) rectangle (14,1); \node at (13.5,0.5) {\footnotesize $\hatbuchilu$};
\draw (9,-2) rectangle (12,-1); \node at(10.5,-1.5) {\footnotesize $\bslashd \btr \bchilu', \bbetalu,\btalu$};
\draw[->] (12,-1.5) -- (13,-1.5); 
\draw (13,-2) rectangle (14,-1); \node at (13.5,-1.5) {\footnotesize $\hatbchilu'$};
\draw (9,-5) rectangle (12,-4); \node at(10.5,-4.5) {\footnotesize $\bslashd \buomegalu, \bslashd \btr \buchilu, \hatbuchilu$};
\draw[->] (12,-4.5) -- (13,-4.5); 
\draw (13,-5) rectangle (14,-4); \node at (13.5,-4.5) {\footnotesize $\overline{\buomegalu}^{u}$};
\draw (7,-6.5) rectangle (8,-5.5); \node at (7.5,-6) {\footnotesize $\br_u$}; 
\draw[->] (8,-6) -- (9,-6);
\draw (9,-6.5) rectangle (10,-5.5); \node at (9.5,-6) {\footnotesize $\overline{\btr \buchilu}^u$}; 
\draw[->] (10,-6) -- (11,-6);
\draw (11,-6.5) rectangle (12,-5.5); \node at (11.5,-6) {\footnotesize $\bmulu$}; 
\draw[->] (12,-6) -- (13,-6);
\draw (13,-6.5) rectangle (14,-5.5); \node at (13.5,-6) {\footnotesize $\overline{\btr \bchilu'}^u$}; 
\end{tikzpicture}
\caption{The sequence to determine the constants $c(n,p,\cdot)$}
\label{fig 4}
\end{figure}
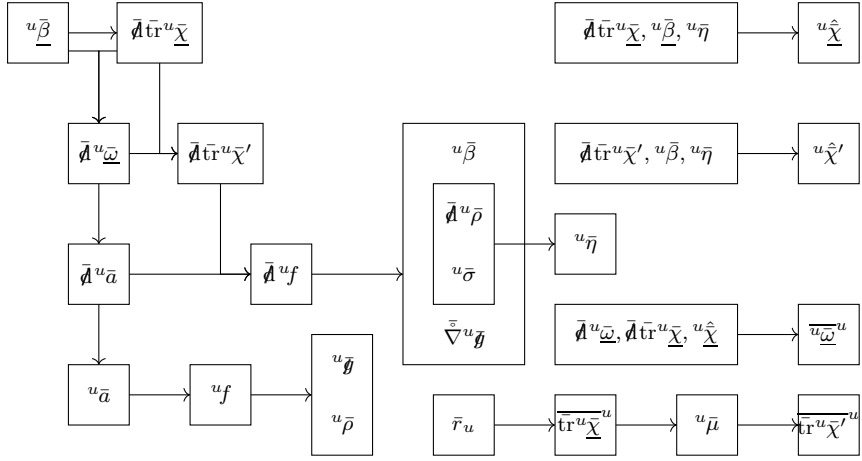

We collect the requirements on the constants $c(n,p,\cdot)$ by the estimates in lemmas \ref{lem 8.5}-\ref{lem 8.22}.
\begin{enumerate}[label=\textbullet]
\item
On the constants for the parameterisation functions:
\begin{align}
	&
	c(n,p,\flu)
	>
	C(n,p,\balu),
	\quad
	c(n,p, \bslashd \flu)
	>
	C(n,p,\bslashd \balu, \bslashd \btr \bchilu').
	\label{eqn 8.14}
\end{align}

\item
On the constants for the metric components:
\begin{align}
\begin{aligned}
	&
	c(n,p,\balu)
	>
	\max \{ c(n,p), C(n,p,\bslashd \balu) \},
	\quad
	c(n,p,\bslashd \balu)
	>
	C(n,p,\bslashd \buomegalu),
\\
	&
	c(n,p,\br_u)
	>
	c(n,p),
\\
	&
	c(n,p,\bslashglu)
	>
	C(n,p,\flu),
	\quad
	c(n,p,\bcircnabla \bslashglu)
	>
	C(n,p,\bslashd \flu).
\end{aligned}
\label{eqn 8.15}
\end{align}

\item
On the constants for the connection coefficients:
\begin{align}
\begin{aligned}
	&
	c(n,p, \overline{\btr \buchilu}^u)
	>
	C(n,p,\br_u),
	\quad
	c(n,p, \bslashd \btr \buchilu)
	>
	\max\{ c(n,p), C(n,p,\bubetalu)\},
\\
	&
	c(n,p,\hatbuchilu)
	>
	\max\{ C(n,p,\bslashd \btr \buchilu, \btalu), C(n,p,\bubetalu) \},
\\
	&
	c(n,p,\overline{\btr \bchilu'}^u)
	>
	\max\{ C(n,p,\bmulu, \overline{\btr \buchilu}^u), C(n,p,\bmulu) \},
\\
	&
	c(n,p,\bslashd \btr \bchilu')
	>
	C(n,p,\bslashd \btr \buchilu, \bslashd \buomegalu),
\\
	&
	c(n,p,\hatbchilu')
	>
	\max\{C(n,p,\bslashd \btr \bchilu', \bbetalu), C(n,p,\btalu), C(n,p,\bslashd \btr\bchilu'), C(n,p,\btalu, \bbetalu) \},
\\
	&
	c(n,p,\btalu)
	>
	\max\{ C(n,p,\bslashd \brholu, \bsigmalu), C(n,p,\bslashd \brholu) \},
\\
	&
	c(n,p,\overline{\buomegalu}^u)
	>
	\max\{ C(n,p, \bslashd \btr \buchilu, \bslashd \buomegalu), C(n,p,\bslashd \btr \buchilu), C(n,p,\hatbuchilu) \},
\\
	&
	c(n,p,\bslashd \buomegalu)
	>
	\max\{C(n,p,\bubetalu), C(n,p,\bslashd \btr \buchilu), C(n,p,\bubetalu, \bslashd \btr \buchilu) \},
\end{aligned}
\label{eqn 8.16}
\end{align}

\item
On the constant for the mass aspect function:
\begin{align}
	c(n,p,\bmulu)
	>
	C(n,p,\overline{\btr \buchilu}^u).
	\label{eqn 8.17}
\end{align}

\item
On the constants for the curvature components:
\begin{align}
\begin{aligned}
	&
	c(n,p,\bubetalu)
	>
	c(n,p),
\\
	&
	c(n,p,\brholu)
	>
	C(n,p,\flu),
\\
	&
	c(n,p,\bslashd \brholu)
	>
	C(n,p,\bslashd \flu),
\\
	&
	c(n,p,\bsigmalu)
	>
	\max \{ c(n,p), C(n,p,\bslashd \flu) \},
\\
	&
	c(n,p,\bbetalu)
	>
	\max \{ c(n,p), C(n,p,\bslashd \flu) \}.
\end{aligned}
\label{eqn 8.18}
\end{align}
\end{enumerate}
Then we can choose the constants $c(n,p,\cdot)$ by the sequence illustrated in figure \ref{fig 4}.

After the constants $c(n,p,\cdot)$ are determined, we choose $\delta$ sufficiently small such that conditions \eqref{eqn 8.2}, \eqref{eqn 8.5} hold, then lemmas \ref{lem 8.4} - \ref{lem 8.22} hold. Therefore for such chosen constants $c(n,p,\cdot)$ and $\delta$, the estimates in the bootstrap assumption are improved. We summarise the above in the following lemma.
\begin{lemma}\label{lem 8.23}
There exist constants $c(n,p,\cdot)$ and $\delta$ depending on $n,p$ satisfying conditions \eqref{eqn 8.14} - \eqref{eqn 8.18}, \eqref{eqn 8.2}, \eqref{eqn 8.5} such that the estimates in the bootstrap assumption are improved.
\end{lemma}

\subsection{Openness of the interval in bootstrap argument}\label{sec 8.10}

Recall the interval $I_a$ defined in section \ref{sec 8.2}, which is the interval of $u_a$ where the bootstrap assumption \ref{assum 8.1} holds,
\begin{align*}
	I_a
	=
	\{ u _a: \text{ bootstrap assumption \ref{assum 8.1} holds on $[0,u_a]$} \}.
\end{align*}
In section \ref{sec 8.2}, we have shown that $I_a$ is closed in step c. in the strategy of the proof of theorem \ref{thm 8.2}. In this section, we show that $I_a$ is also open: suppose $u_a \in I_a$, then $[u_a, u_a + \tau] \subset I_a$ for some positive number $\tau$.

\subsubsection{Local existence of the constant mass aspect function foliation}
Before proving the existence of $[u_a, u_a+\tau) \subset I_a$, we state the local existence theorem of the constant mass aspect function foliation referring to \cite{S2008}. For the sake of self-containedness, we sketch the proof in the appendix \ref{appen thm 8.24}.
\begin{theorem}\label{thm 8.24}
Let $\ucalH$ be an incoming null hypersurface with a background foliation $\{\Sigma_s\}$ in a vacuum spacetime, where $s$ increases in the past-pointed direction. Suppose that $\{s,\vartheta\}$ be a coordinate system of $\ucalH$, where $\vartheta$ is the coordinate system on the sphere $\mathbb{S}^2$. Associated with this coordinate system, the intrinsic metric of the null hypersurface takes the form
\begin{align*}
	g|_{\ucalH}
	=
	\slashg_{ab} ( \ed \theta^a - b^a \ed s) \otimes (\ed \theta^b - b^b \ed s).
\end{align*}

Let $\uL$ be the null tangential vector field on $\ucalH$ with $\uL s= 1$. Let $L'$ be the conjugate null normal vector on $\Sigma_s$ that $g(\uL, L')=2$. Associated with this conjugate null frame $\{\uL, L'\}$, we have the corresponding connection coefficients $\uchi, \chi', \eta, \uomega$, and curvature components $\ualpha, \ubeta, \rho, \sigma, \beta, \alpha$. Assume that
\begin{align*}
	&
	\slashg, b \in C^{n+2},
\\
	&
	\uchi, \eta, \uomega \in C^{n+2},
\\
	&
	\chi', \ualpha, \ubeta, \rho \in C^{n+1},
\end{align*}

Let $\bSigma_{u_a}$ be a spacelike surface in $\ucalH$ with the parameterisation function $\fl{u_a}$ that
\begin{align*}
	\bSigma_{u_a} 
	= 
	\{ (s,\vartheta): s= \fl{u_a}(\vartheta) \}.
\end{align*}
If
\begin{align*}
	\fl{u_a} \in \mathrm{W}^{n+2,p}(\mathbb{S}^2),
\end{align*}
where $n\geq n_p$, $p >1$, then there exists a positive constant $\tau$ such that locally near $\bSigma_{u_a}$, we can construct a constant mass aspect function foliation $\{ \bSigma_u \}_{u\in (u_a - \tau, u_a + \tau)}$ emanating from $\bSigma_{u_a}$ and parametrised by area radius $r_u - r_{u_a} = u - u_a$. Let $\flu$ be the paramterisation function of the leaf $\bSigma_u$, we have that $\flu \in \mathrm{W}^{n+2,p}(\mathbb{S}^2)$ and the Sobolev norm $\Vert \flu \Vert^{n+2,p}$ is continuous w.r.t. $u$. Moreover, $\{ \flu \}$ is a continuous family in $\mathrm{W}^{n+1,p}(\mathbb{S}^2)$.
\end{theorem}

As a corollary of the above theorem, we have the following local existence result of the constant mass aspect function foliation.
\begin{proposition}\label{prop 8.25}
Let $\ucalH$ be an incoming null hypersurface in $(M,g)$ with the background foliation $\{ \Sigma_s = \ucalH \cap C_s\}$. Suppose that associated with the double null coordinate system $\{ \us, s, \vartheta \}$, we have
\begin{align*}
	&
	\slashg, b \in C^{n+3},
\\
	&
	\chi', \uchi, \eta, \uomega \in C^{n+2},
\\
	&
	\ualpha, \ubeta, \rho, \sigma, \beta, \alpha \in C^{n+1}.
\end{align*}
Let $\uh$ be the parameterisation function of $\ucalH$ in the double null coordinate system,
\begin{align*}
	\ucalH
	=
	\{ (\us, s, \vartheta): \us = \uh(s,\vartheta) \}.
\end{align*}
Assume that for $s\in ( s_a(\vartheta), s_b(\vartheta) )$,
\begin{align*}
	\uh(s,\cdot) 
	\in
	\mathrm{W}^{n+5+n_p,p}_{loc}(\mathbb{S}^2),
	\quad
	n_p= 	
	\left\{
	\begin{aligned}
		&
		1, 
		&&
		p>2,
	\\
		&
		2, 
		&&
		2\geq p >1,
	\end{aligned}
	\right.
\end{align*}
where $n\geq n_p$, $p >1$.

Let $\bSigma_{u_a}$ be a spacelike surface in $\ucalH$ with the parameterisation function $\fl{u_a}$ in the background coordinate system $\{ s, \vartheta\}$ of $\ucalH$. Suppose that
\begin{align*}
	\fl{u_a}(\vartheta) \in (s_a(\vartheta), s_b(\vartheta))
	\quad
	\fl{u_a} \in \mathrm{W}^{n+2,p}(\mathbb{S}^2),
\end{align*}
then there exists a constant mass aspect function foliation parameterised by area radius $\{ \bSigma_u \}_{u \in (u_a - \tau, u_a + \tau)}$, of which the parameterisation function $\{\flu \}_{u \in (u_a - \tau, u_a + \tau)}$ is a family in $\mathrm{W}^{n+2,p}(\mathbb{S}^2)$ with the Sobolev norm $\Vert \flu \Vert^{n+2,p}$ varying continuously. Moreover $\{ \flu \}$ is a continuous family in $\mathrm{W}^{n+1,p}(\mathbb{S}^2)$.
\end{proposition}
\begin{proof}
Recall the equation for $\uh$
\begin{align}
	\partial_s \uh = -b^i \partial_i \uh + \Omega^2 (\slashg^{-1})^{ij} \partial_i \uh \partial_j \uh.
	\tag{\ref{eqn 3.1}\ensuremath{'}}.
\end{align}
Then the assumption on $\uh(s,\cdot)$ implies that
\begin{align*}
	\uh 
	\in
	\mathrm{W}^{n+5+n_p,p}(\{(s,\vartheta): s\in (s_a(\vartheta),s_b(\vartheta))\})
	\hookrightarrow
	C^{n+4},
	\quad
	n\geq n_p, p>1.
\end{align*}
Therefore on $\ucalH$, the metric components, connection coefficients and curvature components associated with the background foliation $s$ satisfy that
\begin{align*}
	&
	\dslashg, \db \in C^{n+3},
\\
	&
	\duchi, \dchi', \deta, \duomega \in C^{n+2},
\\
	&
	\dualpha, \dubeta, \drho, \dsigma, \dbeta, \dalpha \in C^{n+1}.
\end{align*}
Then the local existence of the constant mass aspect function foliation follows from theorem \ref{thm 8.24}.
\end{proof}

The following lemma replaces the assumption on $\uh(s,\cdot)$ in proposition \ref{prop 8.25} by the assumption on $\ufl{s=0}$.
\begin{lemma}\label{lem 8.26}
Let $\ucalH$ be an incoming null hypersurface in $(M,g)$ with the background foliation $\{ \Sigma_s = \ucalH \cap C_s\}$. Let $\uh$ be the parameterisation function of $\ucalH$ in the double null coordinate system $\{ \us, s, \vartheta \}$,
\begin{align*}
	\ucalH
	=
	\{ (\us, s, \vartheta): \us = \uh(s,\vartheta) \}.
\end{align*}
Assume that associated with the double null coordinate system,
\begin{align*}
	&
	\Omega, \slashg, b \in C^{n+5+n_p},
\end{align*}
where $n\geq n_p$, $p>1$. Suppose that the parametrisation function $\ufl{s=0}$ of $\Sigma_{s=0}$ satisfies the following estimate
\begin{align*}
	\Vert \slashd \ufl{s=0} \Vert^{n+1,p}
	\leq
	\udelta_o r_0,
	\quad
	\vert \overline{\ufl{s=0}}^{\circg} \vert
	\leq
	\udelta_m r_0,
\end{align*}
and moreover $\ufl{s=0} \in \mathrm{W}^{n+5+n_p,p}(\mathbb{S}^2)$. Let $\delta$ be the positive constant in proposition \ref{prop 3.1}. If $\epsilon, \udelta_o, \udelta_m$ are suitably bounded that $\epsilon, \udelta_o, \epsilon \udelta_m \leq \delta$, then $\uh(s, \cdot) \in \mathrm{W}^{n+5+n_p,p}(\mathbb{S}^2)$. Moreover there exists a constant $C$ depending on $n,p$ such that
\begin{align*}
	&
	\Vert \uh(s,\cdot) \Vert^{n+3+n_p,p}
	\leq
	e^{C [\delta + \underline{B}_{\Omega^2 \subslashg^{-1}, b}^{n+3+n_p}(s)] \cdot \vert s \vert} 
	\cdot \Vert \ufl{s=0} \Vert^{n+3+n_p,p},
\end{align*}
and
\begin{align*}
	\Vert \uh(s,\cdot) \Vert^{n+5+n_p,p}
	&\leq
	e^{C [\delta + \underline{B}_{\Omega^2 \subslashg^{-1}, b}^{n+5+n_p}(s)]\cdot \vert s \vert} 
\\
	&
	\phantom{ \leq}
	\cdot 
	\{
		 \Vert \ufl{s=0} \Vert^{n+5+n_p,p}
		+ 
		C \underline{B}_{\Omega, \subslashg, b}^{n+3+n_p}(s) 
		[\Vert \uh(s,\cdot) \Vert^{n+3+n_p}]^2 
		\cdot 
		|s| 
	\},
\end{align*}
where
\begin{align*}
	\underline{B}_{\Omega^2 \subslashg^{-1}, b}^{m} (s) 
	= 
	\sup_{s' \in [0,s] (\text{or }[s,0]) } 
	\{ \vert \circnabla^k \partial_{\us}^l (\Omega^2 \slashg^{-1}) \vert,  \vert \circnabla^k \partial_{\us}^l b \vert \}_{ k+l \leq m}.
\end{align*}
\end{lemma}
\begin{proof}
$\Vert \uh(s,\cdot) \Vert^{n+2,p}$ satisfies the estimate in proposition \ref{prop 3.1}.

For the estimate of $\Vert \uh(s,\cdot) \Vert^{n+3+n_p,p}$, we take the differential of equation \eqref{eqn 3.1} and integrate the equation to obtain the estimate. Note that by the choice of $n$ and $n_p$, we have
\begin{align*}
	(n+3 + n_p) + 2 \leq 2(n+2 + 1).
\end{align*}
Suppose that $m \leq n+3+n_p$. Therefore in the $m$-th differential of equation \eqref{eqn 3.1}, there is no term involving the product of the high order differentials of $\uh$ in which both orders are greater than $n+2$. Thus the high order differentials of $\uh$ with the order $\geq n+2$ appear linearly in the the $m$-th differential of equation \eqref{eqn 3.1}. Hence the $m$-th differential of equation \eqref{eqn 3.1} takes the form
\begin{align*}
	&\phantom{=}
	\partial_s \uh_{R, i_1 \cdots i_m} 
	+
	\Xl{s}^i \partial_i  \uh_{R, i_1 \cdots i_m} 
\\
	&=
	\sum_{\substack{k+l\leq n+3+n_p \\ m_0 \leq n+3+n_p}} 
	\circnabla^k \partial_{\us}^l \{ b, (\Omega^2 \slashg^{-1}) \} 
	\cdot
	(\prod_{r,m_r \leq n+2} \uh_{R, i_{r,1} \cdots i_{r,m_r}} )
	\cdot
	\uh_{R, i_{0,1} \cdots i_{0,m_0}},
\end{align*}
where $\Xl{s}^i = b^i \partial_i - 2\Omega^2(\slashg^{-1})^{ij} \partial_j \uh$. Then the estimate of $\Vert \uh \Vert^{n+3+n_p,p}$ follows from integrating the propagation equation.

The estimate of $\Vert \uh \Vert^{n+5+n_p,p}$ is similar. Suppose that $m\leq n+5+n_p$. Note that the $m$-th differential of equation \eqref{eqn 3.1} takes the form
\begin{align*}
	&\phantom{=}
	\partial_s \uh_{R, i_1 \cdots i_m} 
	+
	\Xl{s}^i \partial_i  \uh_{R, i_1 \cdots i_m} 
\\
	&=
	\sum_{\substack{k+l\leq n+5+n_p \\ m_0 \leq n+5+n_p}} 
	\circnabla^k \partial_{\us}^l \{ b, (\Omega^2 \slashg^{-1}) \}
	\cdot
	(\prod_{r,m_r \leq n+2} \uh_{R, i_{r,1} \cdots i_{r,m_r}} )
	\cdot
	\uh_{R, i_{0,1} \cdots i_{0,m_0}}
\\
	&\phantom{=}
	+
	\sum_{\substack{k+l\leq n+3+n_p \\ m_0, m'_0 \leq n+3+n_p}} 
	\circnabla^k \partial_{\us}^l \{ b, (\Omega^2 \slashg^{-1}) \}
	\cdot
	(\prod_{r,m_r \leq n+2} \uh_{R, i_{r,1} \cdots i_{r,m_r}} )
	\cdot
	\uh_{R, i_{0,1} \cdots i_{0,m_0}}
	\cdot
	\uh_{R, i'_{0,1} \cdots i'_{0,m'_0}},
\end{align*}
then the last term in the above equation is estimated in terms of $[\Vert \uh \Vert^{n+3+n_p,p}]^2$ and the estimate of $\Vert \uh \Vert^{n+5+n_p,p}$ follows from integrating the above propagation equation.
\end{proof}

In the local existence theorem \ref{thm 8.24}, given a uniform bound of $\Vert \flu \Vert^{n+1,p}$ for $u \in (u_a - \tau, u_a + \tau)$, we can obtain an estimate of $\Vert \flu \Vert^{n+2,p}$ which is linear in $\vert u- u_a \vert$ and $\Vert \fl{u=u_a} \Vert^{n+2,p}$.
\begin{lemma}\label{lem 8.27}
Under the assumption of theorem \ref{thm 8.24}, suppose that $\{ \flu \}_{u\in(u_a-\tau, u_a+ \tau)}$ satisfies the following estimate in $\mathrm{W}^{n+1,p}(\mathbb{S}^2)$,
\begin{align*}
	\Vert \dslashd \flu \Vert^{n,p} \leq \delta_o (r_0+\os_u),
	\quad
	\os_u = \overline{\flu}^{\circg},
	\quad
	-\frac{r_0}{2} \leq \flu \leq s_0,
\end{align*}
then we have the following estimate for $\{ \flu \}_{u\in(u_a-\tau, u_a+ \tau)}$ in $\mathrm{W}^{n+2,p}(\mathbb{S}^2)$,
\begin{align*}
	\Vert \dslashd \flu \Vert^{n+1,p}
	\leq
	\exp[C \vert u-u_a \vert] \cdot 
	(\Vert \dslashd \fl{u=u_a} \Vert^{n+1,p} + C\vert u-u_a \vert ),
\end{align*}
where $C$ is a constant depending on $n,p,\delta_o, s_0$ and the following bounds
\begin{align*}
	&
	\sup_{s\in[-\frac{r_0}{2}, s_0]} \{ \circnabla^k \partial_s^l \slashg, \circnabla^k \partial_s^l b\}_{k+l\leq n+1}, 
\\
	&
	\sup_{s\in[-\frac{r_0}{2}, s_0]} \{ \circnabla^k \partial_s^l \uchi, \circnabla^k \partial_s^l \eta, \circnabla^k \partial_s^l \uomega\}_{k+l\leq n+1},
\\
	&
	\sup_{s\in[-\frac{r_0}{2}, s_0]} \{ \circnabla^k \partial_s^l \chi',  \circnabla^k \partial_s^l \rho, \circnabla^k \partial_s^l \ubeta, \circnabla^k \partial_s^l \ualpha \}_{k+l \leq n}.
\end{align*}
\end{lemma}
\begin{proof}
Note that the constants $C_{\xi}$ and $B_{\xi}$ in the proof of theorem \ref{thm 8.24} in appendix \ref{appen thm 8.24} depend on $n,p$, $\delta_o$, $s_0$ and the bounds in the lemma. Then the estimate follows from the proof of equation \eqref{eqn d.2} in appendix \ref{appen thm 8.24}.
\end{proof}

\subsubsection{Regularised approximation of the foliation}
In order to apply the local existence proposition \ref{prop 8.25} to the constant mass aspect function foliation $\bSigma_{u_a}$, we shall approximate the parameterisation function $(\ufl{s=0}, \ufl{u_a})$ by more regular functions in $\mathrm{W}^{n+6+n_p,p}(\mathbb{S}^2) \times \mathrm{W}^{n+3,p}(\mathbb{S}^2)$ to construct the regularised approximating foliation, and prove the convergence of the approximation. This approximation approach relies on the following lemma.

The lemma is on the local estimate of the parameterisation function $\flu$ of the constant mass aspect function foliation, given the first kind parameterisation function $(\ufl{s=0}, \fl{u=u_a})$ in $\mathrm{W}^{n+2,p}(\mathbb{S}^2) \times \mathrm{W}^{n+2,p}(\mathbb{S}^2)$ and the associated geometric quantities at the leaf $\bSigma_{u=u_a}$ satisfying the following bootstrap assumption stronger than assumption \ref{assum 8.1}.
\begin{assumption}\label{assum 8.28}
Let $\epsilon'$ be a positive number. On the leaf $\bSigma_u$ of the constant mass aspect function foliation, the first kind parameterisation function $(\ufl{s=0}, \flu)$ and the associated geometric quantities satisfy the stronger estimates by $\epsilon'$, i.e. for any quantity $A$ and the corresponding norm $N(A)$ with the estimate $N(A) \leq B$ in the bootstrap assumption \ref{assum 8.1}, we have that
\begin{align*}
	N(A|_{\bSigma_u}) \leq B - \epsilon'.
\end{align*}
\end{assumption}

\begin{lemma}\label{lem 8.29}
Let $\{ \bSigma_u \}$ be a constant mass aspect function foliation on an incoming null hypersurface $\ucalH$ in $(M,g)$. Let $\epsilon'$ be a positive number. Suppose that assumption \ref{assum 8.28} holds on the leaf $\bSigma_{u=u_a}$. There exists a small positive constant $\delta$ depending on $n, p$, such that if $\epsilon, \udelta_o, \udelta_m, \delta_o$ are suitably bounded that $\epsilon + \udelta_o + \udelta_m + \delta_o \leq \delta$, then there exists a positive constant $\tau_{\epsilon'}$ depending on $n, p, r_0$ and $\epsilon'$, such that assumption \ref{assum 8.1} holds on the foliation $\{\bSigma_u\}$ for $u\in [u_a - \tau_{\epsilon'}, u_a + \tau_{\epsilon'}]$.
\end{lemma}
We leave the proof of lemma \ref{lem 8.29} in appendix \ref{appen lem 8.29}. Now we are ready to construct the regularised approximation of the constant mass aspect function foliation $\{ \bSigma_u \}$ near $\bSigma_{u_a}$ as follows.
\begin{construction}[Regularised approximation]\label{constr 8.30}
Let $\bSigma_{u_a}$ be the leaf in the constant mass aspect function foliation $\{ \bSigma_u \}$ of $\ucalH$ in $(M,g)$. Let $( \ufl{s=0}, \fl{u_a})$ be the first kind parametersation of $\bSigma_{u_a}$. Assume that assumption \ref{assum 8.28} with the positive constant $\epsilon'$ holds on $\bSigma_{u_a}$. Let $n\geq n_p$, $p>1$. Suppose that
\begin{align*}
	&
	\Omega, \slashg, b \in C^{n+6+n_p},
\\
	&
	\chi', \uchi, \eta, \uomega \in C^{n+3},
\\
	&
	\ualpha, \ubeta, \rho, \sigma, \beta, \alpha \in C^{n+2}.
\end{align*}
Approximate the parameterisation $(\ufl{s=0}, \fl{u_a})$ in the space $\mathrm{W}^{n+2,p}(\mathbb{S}^2) \times \mathrm{W}^{n+2,p}(\mathbb{S}^2)$ by more regular functions in $\mathrm{W}^{n+6+n_p,p}(\mathbb{S}^2) \times \mathrm{W}^{n+3,p}(\mathbb{S}^2)$.
\begin{enumerate}[label=\textit{\roman*}.]
\item
For any positive integer $k\in \mathbb{Z}_{>0}$, choose a pair of functions $(\ufl{s=0,k}, \fl{u_a,k}) \in \mathrm{W}^{n+5+n_p,p}(\mathbb{S}^2) \times \mathrm{W}^{n+3,p}(\mathbb{S}^2)$ and
\begin{align*}
	\Vert \ufl{s=0,k} - \ufl{s=0} \Vert^{n+2,p} \leq \frac{1}{k},
	\quad
	\Vert \fl{u_a,k} - \fl{u_a} \Vert^{n+2,p} \leq \frac{1}{k}.
\end{align*}
Then introduce the null hypersurface $\ucalH_k$ with the parameterisation $\uhl{k}(s,\vartheta)$
\begin{align*}
	\uhl{k}(s=0,\cdot) = \ufl{s=0,k},
\end{align*}
and the surface $\bSigma_{k,u_a}$ in $\ucalH_k$ with the first kind parameterisation $(\ufl{s=0,k}, \fl{u_a,k})$.

\item
Construct the approximation foliation $\{\bSigma_{k,u}\}$ in $\ucalH_k$. The local existence of $\{\bSigma_{k,u}\}$ follows from proposition \ref{prop 8.25} and lemma \ref{lem 8.26}.
\end{enumerate}
\end{construction}

For the above construction, we show that there exists a uniform local existence interval for the approximation foliation $\{ \bSigma_{k,u} \}$.
\begin{lemma}\label{lem 8.31}
There exists a positive constant $\delta$ depending on $n,p$ such that if $\epsilon$, $\udelta_o$, $\udelta_m$, $\delta_o$, $\delta_m$ are suitably bounded that $\epsilon + \udelta_o + \udelta_m + \delta_o + \delta_m \leq \delta$, then there exists a positive constant $\tau$ depending on $n,p$, $r_0$, $\epsilon'$ and a sufficiently large integer $K$, such that the approximation foliation $\{\bSigma_{k,u}\}$ exists for $u\in [u_a - \tau_u, u_a + \tau_u]$ and $k \geq K$. Moreover assumption \ref{assum 8.1} holds on $\{ \bSigma_{k,u} \}_{[u_a - \tau_u, u_a + \tau_u], k \geq K}$.
\end{lemma}
\begin{proof}
The local existence proposition \ref{prop 8.25} implies the local existence of $\{\bSigma_{k,u}\}$, however it doesnot implies that the local existence interval is uniform for large $k$. In order to prove the uniform local existence interval, we shall also make use of lemma \ref{lem 8.29} on the local estimate of the foliation, where the interval of the estimate is uniform.

Choose $\delta$ in lemmas \ref{lem 8.23} and \ref{lem 8.29}. There exists $K$ sufficiently large such that assumption \ref{assum 8.28} with a smaller positive constant $\frac{\epsilon'}{2}$ holds on $\bSigma_{k,u_a}$ for $k\geq K$. Let $\tau_{\epsilon'/2}$ be the constant in lemma \ref{lem 8.29}. We shall show $\tau_{\epsilon'/2}$ is the constant $\tau_u$ required in lemma \ref{lem 8.31}.

Define $\tau_k$ as the supremum of the following set:
\begin{align*}
	I_k = \{ \tau:  \
	\begin{aligned}
		&
		\{\bSigma_{k,u}\} \text{ exists for } u\in [u_a - \tau, u_a + \tau],
	\\
		&
		\text{assumption \ref{assum 8.1} holds on } \{ \bSigma_{k,u} \}_{[u_a - \tau_u, u_a + \tau_u]}.
	\end{aligned}
	\}
\end{align*}
We show that $\tau_k \geq \tau_{\epsilon'/2}$ for $k \geq K$.

Suppose $\tau \in I_k$ and $\tau < \tau_{\epsilon'/2}$. Then by proposition \ref{prop 8.25}, lemma \ref{lem 8.26}, there exists $\tau' > \tau$, such that $\{ \bSigma_{k,u} \}$ exists on $[u_a - \tau', u_a + \tau']$ and the parameterisation function $\{ \fl{k,u} \}_{[u_a - \tau', u_a + \tau']}$ is a continuous family in $\mathrm{W}^{n+2,p}(\mathbb{S}^2)$. Let $\tau_m = \min\{ \tau', \tau_{\epsilon'/2} \} > \tau'$. Therefore by lemma \ref{lem 8.29}, assumption \ref{assum 8.1} holds on $\{ \bSigma_{k,u} \}$ for $[u_a - \tau_m, u_a + \tau_m]$. Hence $\tau_m > \tau$ and $\tau_m \in I_k$. Thus $\tau_k = \sup I_k \geq \tau_{\epsilon'/2}$.
\end{proof}

\subsubsection{Convergence of regularised approximations}
By lemma \ref{lem 8.31}, we know that the regularised approximation $\{\bSigma_{k,u}\}$ in construction \ref{constr 8.30} exists and satisfies assumption \ref{assum 8.1} on a uniform interval $[u_a-\tau_u, u_a + \tau_u]$ for sufficiently large $k$. Then we shall consider the convergence of the regularised approximations, and show that assumption \ref{assum 8.1} is preserved under the convergence.

\begin{lemma}\label{lem 8.32}
Let $\{ \bSigma_{k,u} \}_{[u_a - \tau_u, u_a + \tau_u]}$ be the regularised approximation foliation in lemma \ref{lem 8.31}. There exists a subsequence $\{ k' \}$ of $\{k\}$ that the following convergences hold for $\{\bSigma_{k',u}\}$ as $k' \rightarrow +\infty$: convergence in the weaker norms in assumption \ref{assum 8.1} and weak convergence in the norms in assumption \ref{assum 8.1}.
\begin{enumerate}[label=\roman*.]
	\item
	Convergence of the parameterisation functions.
	\begin{enumerate}[label=\raisebox{0.1ex}{\scriptsize$\bullet$}]
		\item
		$\fl{u,k'}$ converges in $\mathrm{W}^{n+1,p}(\mathbb{S}^2)$ and weakly converges in $\mathrm{W}^{n+2,p}(\mathbb{S}^2)$.

		\item
		$\ufl{u,k'}$ converges in $\mathrm{W}^{n,p}(\mathbb{S}^2)$ and weakly converges in $\mathrm{W}^{n+1,p}(\mathbb{S}^2)$.

		\item
		$(\dslashd \uhl{k'})|_{\bSigma_{k',u}}$ converges in $\mathrm{W}^{n,p}(\mathbb{S}^2)$ and weakly converges in $\mathrm{W}^{n+1,p}(\mathbb{S}^2)$.

		\item
		$(\circnabla^2 \uhl{k'})|_{\bSigma_{k',u}}$ converges in $\mathrm{W}^{n-1,p}(\mathbb{S}^2)$ and weakly converges in $\mathrm{W}^{n,p}(\mathbb{S}^2)$.

	\end{enumerate}

	\item
	Convergence of the metric components.
	\begin{enumerate}[label=\raisebox{0.1ex}{\scriptsize$\bullet$}]
		\item
		$\bslashgl{u,k'}$ converges in $\mathrm{W}^{n,p}(\mathbb{S}^2)$ and weakly converges in $\mathrm{W}^{n+1,p}(\mathbb{S}^2)$.

		\item
		$\br_{k',u}$ converges.
		
		\item
		$\bal{u,k'}$ converges in $\mathrm{W}^{n,p}(\mathbb{S}^2)$ and weakly converges in $\mathrm{W}^{n+1,p}(\mathbb{S}^2)$.

	\end{enumerate}

	\item
	Convergence of the connection coefficients.
	\begin{enumerate}[label=\raisebox{0.1ex}{\scriptsize$\bullet$}]
		\item
		$\btr \buchil{u,k'}$, $\btr \bchil{u,k'}'$ converge in $\mathrm{W}^{n-1,p}(\mathbb{S}^2)$ and weakly converge in $\mathrm{W}^{n,p}(\mathbb{S}^2)$.

		\item
		$\hatbuchil{u,k'}$, $\hatbchil{u,k'}'$ converge in $\mathrm{W}^{n-1,p}(\mathbb{S}^2)$ and weakly converge in $\mathrm{W}^{n,p}(\mathbb{S}^2)$.
	
		\item
		$\btal{u,k'}$ converges in $\mathrm{W}^{n,p}(\mathbb{S}^2)$ and weakly converges in $\mathrm{W}^{n+1,p}(\mathbb{S}^2)$.

		\item
		$\buomegal{u,k'}$ converges in $\mathrm{W}^{n,p}(\mathbb{S}^2)$ and weakly converges in $\mathrm{W}^{n+1,p}(\mathbb{S}^2)$.

	\end{enumerate}
	
	\item
	Convergence of the mass aspect function. $\bmul{u,k'}$ converges.
	
	\item
	Convergence of the curvature components.
	\begin{enumerate}[label=\raisebox{0.1ex}{\scriptsize$\bullet$}]
		\item
		$\bubetal{u,k'}$ converge in $\mathrm{W}^{n,p}(\mathbb{S}^2)$.

		\item
		$\brhol{u,k'}$ converge in $\mathrm{W}^{n,p}(\mathbb{S}^2)$.

		\item
		$\bsigmal{u,k'}$ converges in $\mathrm{W}^{n,p}(\mathbb{S}^2)$.

		\item
		$\bbetal{u,k'}$ converges in $\mathrm{W}^{n,p}(\mathbb{S}^2)$.

	\end{enumerate}

\end{enumerate}
\end{lemma}
\begin{proof}[Proof of lemma \ref{lem 8.32}]
Note the parameterisation function $\{\fl{u,k}\}_{[u_a-\tau_u, u_a + \tau_u]}$ is a uniformly bounded family in $L^{\infty}([u_a-\tau_u, u_a + \tau_u], \mathrm{W}^{n+2,p}(\mathbb{S}^2))$ and $\{ \partial_u \fl{u,k}\}_{[u_a-\tau_u, u_a + \tau_u]}$ is a uniformly bounded family in $L^{\infty}([u_a-\tau_u, u_a + \tau_u], \mathrm{W}^{n+1,p}(\mathbb{S}^2))$. Therefore by the Sobolev embedding $\mathrm{W}^{n+2,p}(\mathbb{S}^2)) \hookrightarrow \mathrm{W}^{n+1,p}(\mathbb{S}^2))$ and the Arzel\`{a}-Ascoli theorem, there exists a subsequence $\{k'\} \subset \{k\}$ such that $\{\fl{u,k'}\}_{[u_a-\tau_u, u_a + \tau_u]}$ converges in $L^{\infty}([u_a-\tau_u, u_a + \tau_u], \mathrm{W}^{n+1,p}(\mathbb{S}^2))$. In the following, we conclude that all the other convergences follows from the convergence of $\fl{u,k'}$ in $\mathrm{W}^{n+1,p}(\mathbb{S}^2)$.

\begin{enumerate}[label=\textit{\roman*}.]
	\item
	\begin{enumerate}[label=\raisebox{0.1ex}{\scriptsize$\bullet$}]
		\item
		$\pmb{\fl{u,k'}}$: convergence in $\mathrm{W}^{n+1,p}(\mathbb{S}^2)$ and the boundedness in $\mathrm{W}^{n+2,p}(\mathbb{S}^2)$ imply the weak convergence in $\mathrm{W}^{n+2,p}(\mathbb{S}^2)$. Similarly the weak convergence of other terms follow from the convergence in the corresponding norm $\mathrm{W}^{m,p}(\mathbb{S}^2)$ and the boundedness in $\mathrm{W}^{m+1,p}(\mathbb{S}^2)$.
		
		\item
		$\pmb{\ufl{u,k'}}$: note that $\ufl{u,k'}$ is obtained by solving equation \eqref{eqn 3.3} cited here
		\begin{align}
			\begin{aligned}
			&
			\partial_t \uflt 
			= 
			F( f,\ t b^i f_i,\ t e^i f_i,\ t \ue^i f_i,\ \uvarepsilon,\  b^i (\uflt)_i,\  e^i (\uflt)_i,\  \ue^i (\uflt)_i  ),
		\\
			&
			F= f \cdot [ 1-tb^i f_i - t\ue^i f_i - t e^i f_i \cdot \uvarepsilon ]^{-1} \cdot [ \uvarepsilon - b^i (\uflt)_i - \ue^i (\uflt)_i - e^i (\uflt)_i \cdot \uvarepsilon ].
		\end{aligned}
		\tag{\ref{eqn 3.3}\ensuremath{'}}
		\end{align}
		For the estimates of $\ufl{u,k'}$, we integrate equation \eqref{eqn 3.4} cited here
		\begin{align}
			\tag{\ref{eqn 3.4}\ensuremath{'}}
			\partial_t ( \ddcircDelta \uflt )
			=
			\Xlt^i \partial_i ( \ddcircDelta \uflt ) + \relt.
		\end{align}
		We have that for $\uflt$, $f$ satisfying the estimates in assumption \ref{assum 8.1}, $F$ is continuously differentiable as a map
		\begin{align*}
			&
			(f, tf, \uflt)
			\quad
			\mapsto
			\quad
			F (f, tf, \uflt)
			=
			f \cdot [ 1-tb^i f_i - t\ue^i f_i - t e^i f_i \cdot \uvarepsilon ]^{-1}  
		\\
			&\phantom{(f, \uflt, tf) \quad \mapsto \quad F (f, \uflt, tf)=f }
			\cdot
			[ \uvarepsilon - b^i (\uflt)_i - \ue^i (\uflt)_i - e^i (\uflt)_i \cdot \uvarepsilon ]
		\end{align*}
		from $\mathrm{W}^{n+1,p}(\mathbb{S}^2) \times \mathrm{W}^{n+1,p}(\mathbb{S}^2) \times \mathrm{W}^{n,p}(\mathbb{S}^2)$ to $L^{\infty}(\mathbb{S}^2)$.
		
		For the terms in equation \ref{eqn 3.4}, we have that for $\uflt$, $f$ satisfying the estimates in assumption \ref{assum 8.1}, $\Xlt$ is continuously differentiable as a map
		\begin{align*}
			(f, tf, \uflt)
			\quad
			\mapsto
			\quad
			\Xlt (f, tf, \uflt),
		\end{align*}
		from 
		$\mathrm{W}^{n+1,p}(\mathbb{S}^2) \times \mathrm{W}^{n+1,p}(\mathbb{S}^2) \times \mathrm{W}^{n,p}(\mathbb{S}^2)$ 
		to $\mathrm{W}^{n-1,p}(\mathbb{S}^2)$, 
		and $\relt$ is continuously differentiable as a map
		\begin{align*}
			(f, tf, \uflt)
			\quad
			\mapsto
			\quad
			\relt (f, tf, \uflt),
		\end{align*}
		from 
		$\mathrm{W}^{n+1,p}(\mathbb{S}^2) \times \mathrm{W}^{n+1,p}(\mathbb{S}^2) \times \mathrm{W}^{n,p}(\mathbb{S}^2)$ 
		to 
		$\mathrm{W}^{n-2,p}(\mathbb{S}^2)$.
		
		Considering the differences of equations \eqref{eqn 3.3}, \eqref{eqn 3.4} for $\fl{u,k'_1}$ and $\fl{u,k'_2}$, we have
		\begin{align*}
			&
			\partial_t (\ufl{t,\{u,k'_2\}} - \ufl{t,\{u,k'_2\}})
			=
			F(\fl{u,k'_2}, t \fl{u,k'_2}, \ufl{t,\{u,k'_2\}})
			-
			F(\fl{u,k'_1}, t \fl{u,k'_1}, \ufl{t,\{u,k'_1\}}),
		\end{align*}
		and
		\begin{align*}
			&\phantom{=}
			\partial_t [\ddcircDelta(\ufl{t,\{u,k'_2\}} - \ufl{t,\{u,k'_1\}})]
			-
			\Xl{t,\{u,k'_1\}}^i \partial_i [\ddcircDelta(\ufl{t,\{u,k'_2\}} - \ufl{t,\{u,k'_1\}})]
		\\
			&
			=
			\rel{t,\{u,k'_2\}} - \rel{t,\{u,k'_1\}}
			+
			(\Xl{t,\{u,k'_2\}} - \Xl{t,\{u,k'_1\}})^i \partial_i (\ddcircDelta \ufl{t,\{u,k'_2\}}),
		\end{align*}
		where
		\begin{align*}
			\Xl{t,\{u,k'\}}
			=
			\Xlt( \fl{u,k'}, t \cdot \fl{u,k'}, \ufl{t,\{u,k'\}}),
		\\	
			\rel{t,\{u,k'\}}
			=
			\relt( \fl{u,k'}, t \cdot \fl{u,k'}, \ufl{t,\{u,k'\}}).
		\end{align*}
		Therefore by 
		\begin{enumerate}[label=\alph*.]
		\item
		the continuous differentiabilities of $F$, $\Xlt$, $\relt$, 
		\item
		the convergence of $\fl{u,k'}$ in $\mathrm{W}^{n+1,p}$,
		\item
		the uniform boundedness of $\ufl{t,\{u,k'\}}$ in $\mathrm{W}^{n+1,p}$ for $t\in[0,1]$,
		\end{enumerate}
		integrating the above equations implies the convergence of $\{ \ufl{t,\{u,k'\}} \}$ in $\mathrm{W}^{n,p}(\mathbb{S}^2)$ for any fixed $t\in [0,1]$. In particular, taking $t=1$, we obtain the convergence of $\{ \ufl{t,\{u,k'\}} \}$ in $\mathrm{W}^{n,p}(\mathbb{S}^2)$.
		
		\item
		$\pmb{(\dslashd \uhl{k'})|_{\bSigma_{k',u}}}$:
		we obtain the estimate of $(\dslashd \uhl{k'})|_{\bSigma_{k',u}}$ by equation \eqref{eqn 5.3} cited here
		\begin{align*}
			\partial_t\, \uhlt_{R,k}
			=
			\Xlt_{\uh}^i \ddpartial_i\, \uhlt_{\dR,k} + \relt_{\uh,k}.
		\end{align*}
		We have that for $\uflt$, $f$ satisfying the estimates in assumption \ref{assum 8.1}, $\Xlt_{\uh}$ is continuously differentiable as a map
		\begin{align*}
			(f, tf, \uflt, \uhlt_{R,i})
			\quad
			\mapsto
			\quad
			\Xlt_{\uh} (f, tf, \uflt, \uhlt_{R,i}),
		\end{align*}
		from 
		$\mathrm{W}^{n+1,p}(\mathbb{S}^2) 
			\times 
			\mathrm{W}^{n+1,p}(\mathbb{S}^2) 
			\times 
			\mathrm{W}^{n,p}(\mathbb{S}^2) 
			\times 
			\mathrm{W}^{n,p}(\mathbb{S}^2)$
		to 
		$\mathrm{W}^{n,p}(\mathbb{S}^2)$, and $\relt_{\uh,k}$ is continuously differentiable as a map
		\begin{align*}
			(f, tf, \uflt, \uhlt_{R,i})
			\quad
			\mapsto
			\quad
			\relt_{\uh,k} (f, tf, \uflt, \uhlt_{R,i}),
		\end{align*}
		from 
		$\mathrm{W}^{n+1,p}(\mathbb{S}^2) 
			\times 
			\mathrm{W}^{n+1,p}(\mathbb{S}^2) 
			\times 
			\mathrm{W}^{n,p}(\mathbb{S}^2) 
			\times 
			\mathrm{W}^{n,p}(\mathbb{S}^2)$
		to $\mathrm{W}^{n,p}(\mathbb{S}^2)$.
		Considering the difference of equation \eqref{eqn 5.3} for $\uhl{t,k'_1}_{R,k}$ and $\uhl{t,k'_2}_{R,k}$, we have
		\begin{align*}
			&\phantom{=}
			\partial_t ( \uhl{t,\{u,k'_2\}}_{R,k} - \uhl{t,\{u,k'_1\}}_{R,k}  )
			-
			\Xl{t,\{u,k'_1\}}_{\uh}^i \ddpartial_i ( \uhl{t,\{u,k'_2\}}_{R,k} - \uhl{t,\{u,k'_1\}}_{R,k} )
		\\
			&
			=
			\rel{t,\{u,k'_2\}}_{\uh,k} - \rel{t,\{u,k'_1\}}_{\uh,k}
			+
			( \Xl{t,\{u,k'_2\}}_{\uh}  - \Xl{t,\{u,k'_1\}}_{\uh})^i \ddpartial_i (\uhl{t,\{u,k'_2\}}_{R,k})
		\end{align*}
		Therefore by 
		\begin{enumerate}[label=\alph*.]
		\item
		the continuous differentiabilities of $\Xlt_{\uh}$, $\relt_{\uh,k}$, 
		\item
		the convergence of $\fl{u,k'}$ in $\mathrm{W}^{n+1,p}(\mathbb{S}^2)$, $\ufl{u,k'}$ in $\mathrm{W}^{n,p}(\mathbb{S}^2)$,
		\item
		the uniform boundedness of $\uhl{t,\{u,k'_1\}}_{\dR,k}$ in $\mathrm{W}^{n+1,p}(\mathbb{S}^2)$,
		\end{enumerate}
		integrating the above equation implies the convergence of $\{ \uhl{t,\{u,k'_1\}}_{\dR,k} \}$ in $\mathrm{W}^{n,p}(\mathbb{S}^2)$ for any fixed $t\in[0,1]$. In particular, taking $t=1$, we obtain the convergence of $\{ (\dslashd \uhl{k'})|_{\bSigma_{k',u}}\}$ in $\mathrm{W}^{n,p}(\mathbb{S}^2)$.
		
		\item
		$\pmb{(\circnabla^2 \uhl{k'})|_{\bSigma_{k',u}}}$: the proof of the convergence of $\{ (\circnabla^2 \uhl{k'})|_{\bSigma_{k',u}} \}$ in $\mathrm{W}^{n-1,p}(\mathbb{S}^2)$ by equation \eqref{eqn 5.3} follows the same strategy as the above proof for $(\dslashd \uhl{k'})|_{\bSigma_{k',u}}$.
	
	\end{enumerate}
		
	\item
	\begin{enumerate}[label=\raisebox{0.1ex}{\scriptsize$\bullet$}]
		\item
		$\pmb{\bslashgl{u,k'}}$: the convergence of $\bslashgl{u,k'}$ in $\mathrm{W}^{n,g}$ follows from the convergences of $\fl{u,k'}$ in $\mathrm{W}^{n+1,p}(\mathbb{S}^2)$, $\ufl{u,k'}$ in $\mathrm{W}^{n,p}(\mathbb{S}^2)$, $(\dslashd \uhl{k'})|_{\bSigma_{k',u}}$ in $\mathrm{W}^{n,p}(\mathbb{S}^2)$ and the formula of $\ddslashg$ in section \ref{sec 4.2}.
				
		\item
		$\pmb{\br_{u,k'}}$: the convergence of $\br_{u,k'}$ follows from the convergence of $\bslashgl{u,k'}$.
		
		\item
		$\pmb{\bal{u,k'}}$:
		considering the difference of equation \eqref{eqn 6.5} for $\bal{u,k'_1}$ and $\bal{u,k'_2}$,
		\begin{align*}
			&\phantom{=}
			\bslashDelta_{\subbslashgl{u,k'_1}} (\log \bal{u,k'_2} - \log \bal{u,k'_1} )
		\\
			&=
			(\bslashDelta_{\subbslashgl{u,k'_2}} - \bslashDelta_{\subbslashgl{u,k'_1}}) \log \bal{u,k'_2}
			- (\ddrhol{u,k'_2} - \ddrhol{u,k'_1}) + (\overline{\ddrhol{u,k'_2}}^{u,k'_2} -\overline{\ddrhol{u,k'_1}}^{u,k'_1}) 
		\\
			&\phantom{=}
			- \frac{1}{2} [ (\hatdduchil{u,k'_2}, \hatddchil{u,k'_2}')_{\subbslashgl{u,k'_2}} -  (\hatdduchil{u,k'_1}, \hatddchil{u,k'_1}')_{\subbslashgl{u,k'_1}} ]
		\\
			&\phantom{=}
			+ \frac{1}{2} [ \overline{(\hatdduchil{u,k'_2}, \hatddchil{u,k'_2}')_{\subbslashgl{u,k'_2}}}^{u,k'_2} -  \overline{(\hatdduchil{u,k'_1}, \hatddchil{u,k'_1}')_{\subbslashgl{u,k'_1}}}^{u,k'_1} ]
		\\
			&\phantom{=}
			- \bslashdiv_{\subbslashgl{u,k'_1}} (\ddetal{u,k'_2} - \ddetal{u,k'_1})
			- ( \bslashdiv_{\subbslashgl{u,k'_2}} - \bslashdiv_{\subbslashgl{u,k'_1}})  \ddetal{u,k'_2},
		\end{align*}
		where
		\begin{align*}
			&\phantom{=}
			(\bslashDelta_{\subbslashgl{u,k'_2}} - \bslashDelta_{\subbslashgl{u,k'_1}}) \log \bal{u,k'_2}
		\\
			&=
			[(\bslashgl{u,k'_2}^{-1})^{ij} - (\bslashgl{u,k'_2}^{-1})^{ij}] \cdot \bslashnabla_{\subbslashgl{u,k'_1},ij}^2 \log \bal{u,k'_2}
		\\
			&\phantom{=}
			+
			(\bslashgl{u,k'_2}^{-1})^{ij} \cdot (\bslashgl{u,k'_2}^{-1})^{kl} \cdot \Gamma( \bslashnabla_{\subbslashgl{u,k'_1}} (\bslashgl{u,k'_2} - \bslashgl{u,k'_1}) )_{ijl} \cdot (\dslashd \log \bal{u,k'_2})_k,
		\\
			&\phantom{=}
			(\bslashdiv_{\subbslashgl{u,k'_2}} - \bslashdiv_{\subbslashgl{u,k'_1}}) \ddetal{u,k'_2}
		\\
			&=
			[(\bslashgl{u,k'_2}^{-1})^{ij} - (\bslashgl{u,k'_2}^{-1})^{ij}] \cdot \bslashnabla_{\subbslashgl{u,k'_1},i} \ddetal{u,k'_1}_j 
		\\
			&\phantom{=}
			+
			(\bslashgl{u,k'_2}^{-1})^{ij} \cdot (\bslashgl{u,k'_2}^{-1})^{kl} \cdot \Gamma( \bslashnabla_{\subbslashgl{u,k'_1}} (\bslashgl{u,k'_2} - \bslashgl{u,k'_1}) )_{ijl} \cdot \ddetal{u,k'_1}_k,
		\end{align*}
		where $\Gamma(T)_{ijl} = T_{ilj} + T_{jli} - T_{lij}$. Therefore the convergence of $\bslashd \bal{u,k'}$ in $\mathrm{W}^{n-1,p}(\mathbb{S}^2)$ follows from
		\begin{enumerate}[label=\alph*.]
		\item
		the convergences of $\fl{u,k'}$ in $\mathrm{W}^{n+1,p}(\mathbb{S}^2)$, $\ufl{u,k'}$ in $\mathrm{W}^{n,p}(\mathbb{S}^2)$, $(\dslashd \uhl{k'})|_{\bSigma_{k',u}}$ in $\mathrm{W}^{n,p}(\mathbb{S}^2)$, $(\circnabla^2 \uhl{k'})|_{\bSigma_{k',u}}$ in $\mathrm{W}^{n-1,p}(\mathbb{S}^2)$, $\bslashgl{u,k'}$ in $\mathrm{W}^{n,p}(\mathbb{S}^2)$,
		
		\item
		the convergences of $\ddrhol{u,k'}$ in $\mathrm{W}^{n,p}(\mathbb{S}^2)$, $\hatdduchil{u,k'}$, $\hatddchil{u,k'}'$, $\ddetal{u,k'}$ in $\mathrm{W}^{n-1,p}(\mathbb{S}^2)$ by the formulae of $\ddrho$, $\hatdduchi$, $\hatddchi'$, $\ddeta$ in section \ref{sec 4.3},
		
		\item
		the uniform boundedness of $\bslashd \bal{u,k'}$ in $\mathrm{W}^{n,p}(\mathbb{S}^2)$.
		\end{enumerate}
		The convergence of $\bal{u,k'}$ in $L^{\infty}(\mathbb{S}^2)$ follows from the convergence of $\br_{u,k'}$, the convergence of $\ddtr \dduchil{u,k'}$ in $\mathrm{W}^{n-1,p}(\mathbb{S}^2)$, $\bslashd \bal{u,k'}$ in $\mathrm{W}^{n-1,p}(\mathbb{S}^2)$ and equation \eqref{eqn 6.5}.

	\end{enumerate}
	
	\item
	\begin{enumerate}[label=\raisebox{0.1ex}{\scriptsize$\bullet$}]
		\item
		$\pmb{\btr \buchil{u,k'}}$, $\pmb{\btr \bchil{u,k'}'}$:
		the convergences of $\btr \buchil{u,k'}, \btr \bchil{u,k'}'$ in $\mathrm{W}^{n-1,p}(\mathbb{S}^2)$ follow from the convergences of $\ddtr \dduchil{u,k'}, \ddtr \ddchil{u,k'}'$ in $\mathrm{W}^{n-1,p}(\mathbb{S}^2)$ and $\bal{u,k'}$ in $\mathrm{W}^{n,p}(\mathbb{S}^2)$.

		\item
		$\pmb{\hatbuchil{u,k'}}$, $\pmb{\hatbchil{u,k'}'}$:
		the convergences of $\hatbuchil{u,k'}, \hatbchil{u,k'}'$ in $\mathrm{W}^{n-1,p}(\mathbb{S}^2)$ follow from the convergences of $\hatdduchil{u,k'}, \hatddchil{u,k'}'$ in $\mathrm{W}^{n-1,p}(\mathbb{S}^2)$ and $\bal{u,k'}$ in $\mathrm{W}^{n,p}(\mathbb{S}^2)$.
	
		\item
		$\pmb{\btal{u,k'}}$:
		considering the difference of equation \eqref{eqn 6.12} for $\btal{u,k'_1}$ and $\btal{u,k'_2}$,
		\begin{align*}
			&\phantom{=}
			\bslashcurl_{\subbslashgl{u,k'_1}} (\btal{u,k'_2} - \btal{u,k'_1})
		\\
			&=
			-(\bslashcurl_{\subbslashgl{u,k'_2}} - \bslashcurl_{\subbslashgl{u,k'_1}} )\btal{u,k'_2}
			+ \ddsigmal{u,k'_2} - \ddsigmal{u,k'_1}
		\\	
			&\phantom{=}
			+ \frac{1}{2} ( \hatbchil{u,k'_2} \wedge_{\bslashgl{u,k'_2}} \hatbuchil{u,k'_2} - \hatbchil{u,k'_1} \wedge_{\bslashgl{u,k'_1}} \hatbuchil{u,k'_1} )
		\\	
			&\phantom{=}
			\bslashdiv_{\subbslashgl{u,k'_1}} (\btal{u,k'_2} - \btal{u,k'_1})
		\\
			&=
			-(\bslashdiv_{\subbslashgl{u,k'_2}} - \bslashdiv_{\subbslashgl{u,k'_1}} )\btal{u,k'_2}
			- (\ddrhol{u,k'_2} - \ddrhol{u,k'_1} )
			- ( \bmul{u,k'_2} - \bmul{u,k'_1})
		\\	
			&\phantom{=} 
			-\frac{1}{2} [ (\hatbuchil{u,k'_2}, \hatbchil{u,k'_2}')_{\subbslashgl{u,k'_2}} -  (\hatbuchil{u,k'_1}, \hatbchil{u,k'_1}')_{\subbslashgl{u,k'_1}} ]
		\end{align*}
		where
		\begin{align*}
			&\phantom{=}
			(\bslashcurl_{\subbslashgl{u,k'_2}} - \bslashcurl_{\subbslashgl{u,k'_1}} )\btal{u,k'_2}
		\\
			&=
			[(\epsilon_{\subbslashgl{u,k'_2}})^{ij} - (\epsilon_{\subbslashgl{u,k'_2}})^{ij}] \cdot (\bslashnabla_{\subbslashgl{u,k'_1}} \btal{u,k'_2})_{ij}
		\\
			&\phantom{=}
			+
			(\epsilon_{\subbslashgl{u,k'_2}})^{ij}  \cdot (\bslashgl{u,k'_2}^{-1})^{kl} \cdot \Gamma( \bslashnabla_{\subbslashgl{u,k'_1}} (\bslashgl{u,k'_2} - \bslashgl{u,k'_1}) )_{ijl} \cdot \btal{u,k'_2}_k.
		\end{align*}
		Therefore the convergence of $\btal{u,k'}$ in $\mathrm{W}^{n,p}(\mathbb{S}^2)$ follows from
		\begin{enumerate}[label=\alph*.]
		\item
		the convergences of $\bslashgl{u,k'}$ in $\mathrm{W}^{n,p}(\mathbb{S}^2)$, $\hatbuchil{u,k'}$, $\hatbchil{u,k'}'$ in $\mathrm{W}^{n-1,p}(\mathbb{S}^2)$,
		
		\item
		the convergences of $\fl{u,k'}$ in $\mathrm{W}^{n+1,p}(\mathbb{S}^2)$, $\ufl{u,k'}$ in $\mathrm{W}^{n,p}(\mathbb{S}^2)$, $(\dslashd \uhl{k'})|_{\bSigma_{k',u}}$ in $\mathrm{W}^{n,p}(\mathbb{S}^2)$,
		
		\item
		the convergence of $\brhol{u,k'}$, $\bsigmal{u,k'}$ in $\mathrm{W}^{n,p}(\mathbb{S}^2)$ by formulae in section \ref{sec 4.3},
		
		\item
		the convergence of $\bmul{u,k'}$ proved later in \textit{iv}..
		
		\end{enumerate}

		\item
		$\pmb{\buomegal{u,k'}}$:
		considering the difference of equation \eqref{eqn 6.13} for $\buomegal{u,k'_1}$ and $\buomegal{u,k'_2}$,
		\begin{align*}
			&\phantom{=}
			\bslashDelta_{\subbslashgl{u,k'_1}} (\buomegal{u,k'_2} - \buomegal{u,k'_1} )
		\\
			&=
			(\bslashDelta_{\subbslashgl{u,k'_2}} - \bslashDelta_{\subbslashgl{u,k'_1}}) \buomegal{u,k'_2}
			+
			\text{difference of r.h.s. of equation \eqref{eqn 6.13}}.
		\end{align*}
		Similarly to the proof of the convergence of $\bslashd \bal{u,k'}$, the convergence of $\bslashd \buomegal{u,k'}$ in $\mathrm{W}^{n-1,p}(\mathbb{S}^2)$ follows from
		\begin{enumerate}[label=\alph*.]
		\item
		the convergence of $\bslashgl{u,k'}$ in $\mathrm{W}^{n,p}(\mathbb{S}^2)$,
		
		\item
		the convergences of $\btr \buchil{u,k'}$, $\btr \bchil{u,k'}'$, $\hatbuchil{u,k'}$, $\hatbchil{u,k'}'$ in $\mathrm{W}^{n-1,p}(\mathbb{S}^2)$, $\btal{u,k'}$ in $\mathrm{W}^{n,p}(\mathbb{S}^2)$,
		
		\item
		the convergence of $\bubetal{u,k'}$ in $\mathrm{W}^{n-1,p}(\mathbb{S}^2)$ proved later in \textit{v}.
		
		\end{enumerate}
		The convergence of $\buomegal{u,k'}$ in $\mathrm{W}^{n,p}$ follows from the convergence of $\bslashd \buomegal{u,k'}$ in $\mathrm{W}^{n-1,p}(\mathbb{S}^2)$, the convergence of $\overline{\buomegal{u,k'}}^{u,k'}$ by formula \eqref{eqn 6.13}.

	\end{enumerate}
	
	\item
	$\pmb{\bmul{u,k'}}$:
	the convergence of $\bmul{u,k'}$ follows from the formula
	\begin{align*}
		\bmul{u,k'} = 
		- \overline{\brhol{u,k'}}^{u,k'} 
		- \frac{1}{2} \overline{(\hatbuchil{u,k'}, \hatbchil{u,k'}')_{\subbslashgl{u,k'}}}^{u,k'}
	\end{align*}
	and the convergences of the terms on the right hand side.
	
	\item
	$\pmb{\bubetal{u,k'}}$,
	$\pmb{\brhol{u,k'}}$,
	$\pmb{\bsigmal{u,k'}}$,
	$\pmb{\bbetal{u,k'}}$:
	the convergences in $\mathrm{W}^{n,p}(\mathbb{S}^2)$ follow from
	\begin{enumerate}[label=\alph*.]
		\item
		the convergences of $\fl{u,k'}$ in $\mathrm{W}^{n+1,p}(\mathbb{S}^2)$, $\ufl{u,k'}$ in $\mathrm{W}^{n,p}(\mathbb{S}^2)$, $(\dslashd \uhl{k'})|_{\bSigma_{k',u}}$ in $\mathrm{W}^{n,p}(\mathbb{S}^2)$,
		
		\item
		the convergence of $\bal{u,k'}$ in $\mathrm{W}^{n,p}(\mathbb{S}^2)$,
		
		\item
		formulae of $\bubetalu$, $\brholu$, $\bsigmalu$, $\bbetalu$ in section \ref{sec 4.3}.
	\end{enumerate}
\end{enumerate}
\end{proof}

\subsubsection{Openness of the interval}

We apply lemma \ref{lem 8.32} to prove the openness of the interval $I_a$ in the bootstrap argument step d.ii. in section \ref{sec 8.2}.
\begin{lemma}\label{lem 8.33}
Let $\{ \bSigma_u \}$ be a constant mass aspect function foliation in an incoming null hypersurface $\ucalH$ in $(M,g)$. Let $I_a$ be the interval of $u_a$ where the bootstrap assumption \ref{assum 8.1} holds,
\begin{align*}
	I_a
	=
	\{ u _a: \text{ bootstrap assumption \ref{assum 8.1} holds on $[0,u_a]$} \}.
\end{align*}
Then there exists a positive constant $\delta$ depending on $n,p$ such that if $u_a \in I_a$, then $[u_a, u_a + \tau] \subset I_a$ for some positive number $\tau$.
\end{lemma}
\begin{proof}
Choose $\delta$ as in lemma \ref{lem 8.31}. Then construct the regularised approximation foliation $\{\bSigma_{k,u}\}$ for $\bSigma_{u_a}$. By lemma \ref{lem 8.32}, there exists a subsequence $\{ k' \}$ and a constant $\tau_u$, that $\{\bSigma_{k',u}\}_{[u_a-\tau_u, u_a + \tau_u]}$ satisfies assumption \ref{assum 8.1} and converges to a limit foliation $\{\bSigma_{\infty,u}\}_{[u_a - \tau_u, u_a + \tau_u]}$, where $\bSigma_{\infty,u_a} = \bSigma_{u_a}$. Lemma \ref{lem 8.32} implies that the limit foliation $\{\bSigma_{\infty,u}\}_{[u_a-\tau_u, u_a + \tau_u]}$ is a constant mass aspect function foliation and satisfies assumption \ref{assum 8.1}. Thus $[u_a , u_a + \tau_u] \subset I_a$.
\end{proof}

Finally we complete the bootstrap argument to prove theorem \ref{thm 8.2} by lemmas \ref{lem 8.23} and \ref{lem 8.33}.

\section{Asymptotic geometry of the foliation}\label{sec 9}
As mentioned in the introduction, the asymptotic geometry of the constant mass aspect function foliation is essential for the application of the foliation to study the null Penrose inequality. In this section, we study the asymptotic behaviour of the geometry of the foliation at past null infinity.

\subsection{Derivatives of parameterisation functions along the foliation}
We derive the formulae for the derivatives of the parameterisation functions $\flu$, $\uflu$ and the differential $(\dslashd \uh)|_{\bSigma_u}$, the Hessian $(\circnabla \uh)|_{\bSigma_u}$ along the foliation.
\begin{lemma}\label{lem 9.1}
Along a constant mass aspect function foliation $\bSigma_u$ on an incoming null hypersurface in $(M,g)$, we have the following formulae for the derivatives of parameterisation functions.
\begin{enumerate}[label=\alph*.]
\item
For $\flu$, 
\begin{align*}
	\ddpartial \flu 
	= 
	\balu - \balu \db^i (\ddpartial_i \flu).
\end{align*}

\item
For $\uflu$,
\begin{align}
	\begin{aligned}
		&
		\ddpartial_u \uflu 
		= 
		F( \ddpartial_u \flu,\ b^i (\flu)_i,\ e^i (\flu)_i,\ \ue^i (\flu)_i,\ \uvarepsilon,\  b^i (\uflu)_i,\  e^i (\uflu)_i,\  \ue^i (\uflu)_i  ),
	\\
		&
		F= \ddpartial_u \flu \cdot [ 1- b^i (\flu)_i - \ue^i (\flu)_i - e^i (\flu)_i \cdot \uvarepsilon ]^{-1} \cdot [ \uvarepsilon - b^i (\uflu)_i - \ue^i (\uflu)_i - e^i (\uflu)_i \cdot \uvarepsilon ],
	\end{aligned}
	\tag{\ref{eqn 3.3}\ensuremath{'}}
\end{align}
where
\begin{align*}
	\begin{aligned}
		&
		\underline{\varepsilon} = \frac{ -|\underline{e}|^2}{(2\Omega^2 + e\cdot \underline{e}) + \sqrt{(2\Omega^2 + e\cdot \underline{e})^2 -|e|^2 |\underline{e}|^2}},
	\\
		&
		|e|^2 = \slashg_{ij}e^ie^j,
		\quad
		|\underline{e}|^2 = \slashg_{ij} \underline{e}^i \underline{e}^j, 
		\quad
		e\cdot \underline{e} =\slashg_{ij} e^i \underline{e}^j,
	\\
		&
		e^k =-2\Omega^2 (\flu)_i (B^{-1})_j^i (\slashg^{-1})^{jk}, 
		\quad
		\underline{e}^k = -2\Omega^2 (\uflu)_i (B^{-1})_j^i (\slashg^{-1})^{jk},
	\\
		&
		B_i^j= \delta_i^j - \flu_i b^j.
	\end{aligned}
\end{align*}

\item
For $\uhl{u}_{\dR,k} = \uh_{\dR,k} |_{\bSigma_u}$, and $\uhl{u}_{\dR,lk} = \uh_{\dR,lk} |_{\bSigma_u}$
\begin{align}
	\begin{aligned}
		&
		\ddpartial_u\, \uhl{u}_{\dR,k}
		=
		\Xl{u}_{\uh}^i \ddpartial_i\, \uhl{u}_{\dR,k} + \rel{u}_{\uh,k},
	\\
		&
		\ddpartial_u \uhl{u}_{\dR,lk}
		=
		\Xl{u}_{\uh}^i \ddpartial_i \uhl{u}_{\dR,lk}
		+ \rel{u}_{\uh,lk},
	\end{aligned}
	\tag{\ref{eqn 5.3}\ensuremath{'}}
\end{align}
where
\begin{align*}
	&
	\Xl{u}_{\uh} 
	=
	(\ddpartial_u \flu) [ 1- b^m (\flu)_m - 2 \Omega^2 ( \slashg^{-1} )^{R,\overline{m}j} (\flu)_m (\uhl{u}_{\dR,j})  ]^{-1}
\\
	&\phantom{\Xl{u}_{\uh} = (\ddpartial_u \flu)}
	\cdot
	[ - b^i + 2 \Omega^2 ( \slashg^{-1} )^{R,\overline{i}j} (\uhl{u}_{\dR,j})  ]
	\ddpartial_i,
\\
	&
	\rel{u}_{\uh,k}
	=
	(\ddpartial_u \flu) [ 1- b^m (\flu)_m - 2 \Omega^2 ( \slashg^{-1} )^{R,\overline{m}j} (\flu)_m (\uhl{u}_{\dR,j}) ]^{-1}
	\cdot
	Q_k|_{\bSigma_u},
\\
	&
	\rel{u}_{\uh,lk}
	=
	(\ddpartial_u \flu) [ 1- b^m (\flu)_m - 2 \Omega^2 ( \slashg^{-1} )^{R,\overline{m}j} (\flu)_m (\uhl{u}_{\dR,j}) ]^{-1}
	\cdot
	Q_{lk}|_{\bSigma_u},
\end{align*}
and $Q_k$, $Q_{lk}$ are given in equations \eqref{eqn 5.1}, \eqref{eqn 5.2}.
\end{enumerate}
\end{lemma}
\begin{proof}
The formula for $\flu$ follows from equation \eqref{eqn 6.3} and $\buL = \ddpartial_u + \balu \db^i \ddpartial_i$. The formulae of $\uflu$ and $\uhl{u}_{\dR,k}$, $\uhl{u}_{\dR,lk}$ follow the same derivations for equation \eqref{eqn 3.3}, \eqref{eqn 5.3}.
\end{proof}

We estimate the derivatives of the parameterisation functions along the foliation by the formulae in the above lemma.
\begin{lemma}\label{lem 9.2}
Along the constant mass aspect function foliation $\{\bSigma_u\}$ in theorem \ref{thm 8.2}, the following estimates hold.
\begin{enumerate}[label=\alph*.]
\item
For $\ddpartial_u \flu$,
\begin{align*}
	&
	\vert \ddpartial_u \flu - \partial_u \flu_S \vert
	\leq
	c(n,p)(\epsilon + \delta_o + \delta_m + \ud_{o,\uh} + \uslashd_m),
\\
	&
	\Vert \bslashd ( \ddpartial_u \flu) \Vert^{n,p}
	\leq
	c(n,p) (\epsilon + \delta_o + \ud_{o,\uh}).
\end{align*}

\item
For $\ddpartial_u \uflu$,
\begin{align*}
	\Vert \ddpartial_u \uflu \Vert^{n,p}
	\leq
	\frac{c(n,p) r_0^2}{\br_u^2}[\ud_o^2 + \frac{\epsilon \ud_o \ud_m r_0}{\br_u}].
\end{align*}

\item
For $\ddpartial_u [(\dslashd \uh)|_{\bSigma_u}]$, $\ddpartial_u [(\circnabla^2 \uh)|_{\bSigma_u}]$,
\begin{align*}
	&\phantom{\leq}
	\Vert \ddpartial_u [(\dslashd \uh)|_{\bSigma_u}] \Vert^{n,p},
	\Vert \ddpartial_u [(\circnabla^2 \uh)|_{\bSigma_u}] \Vert^{n-1,p}
\\
	&
	\leq
	\frac{(c(n,p)r_0^2}{\br_u^2} 
	[ \ud_{\uh}^2 + \frac{\epsilon (\ud_o + \ud_m) \cdot \ud_{\uh} r_0}{\br_u} + \frac{\ud_{\uh}^3 r_0}{\br_u} ].
\end{align*}
\end{enumerate}
\end{lemma}
\begin{proof}
The estimate of $\ddpartial_u \flu$ follows from assumption \ref{assum 8.1}. The other estimates follow the same routes of estimating $F$ in \eqref{eqn 3.6} and $\Xlu_{\uh}$, $\relu_{\uh,k}$, $\relu_{\uh,lk}$ in the proofs of propositions \ref{prop 5.2}, \ref{prop 5.5}.
\end{proof}

\subsection{Renormalised metric and Gauss curvature}
In section \ref{sec 8}, we construct the global constant mass aspect function foliation on a nearly spherically symmetric null hypersurface, emanating from a spacelike surface near $C_{s=0}$ to the past null infinity. Moreover, the geometry of the foliation satisfies the estimates in assumption \ref{assum 8.1}. Then we introduce the following renormalised metric and renormalised Gauss curvature of the foliation.
\begin{definition}\label{def 9.3}
Let $\{ \bSigma_u \}$ be a constant mass aspect function foliation on an incoming null hypersurface $\ucalH$ in theorem \ref{thm 8.2}. Define the renormalised metric and the renormalised Gauss curvature for the foliation as
\begin{align*}
	\bslashgl{u,r} = \br_u^{-2} \cdot \bslashglu,
	\quad
	\bKl{u,r} = \br_u^2 \cdot \bKlu.
\end{align*}
Define the limit renormalised metric and Gauss curvature at the past null infinity as
\begin{align*}
	\bslashgl{\infty,r}
	=
	\lim_{u\rightarrow \infty} \bslashgl{u,r},
	\quad
	\bKl{\infty,r}
	=
	\lim_{u\rightarrow \infty} \bKl{u,r}.
\end{align*}
\end{definition}
The limit renormalised metric and Gauss curvature at the past null infinity in definition \ref{def 9.3} are well defined.
\begin{proposition}\label{prop 9.4}
Along the constant mass aspect function foliation $\{\bSigma_u\}$ in theorem \ref{thm 8.2}, we have
\begin{enumerate}[label=\alph*.]
\item
the renormalised metric $\bslashgl{u,r}$ converges to the limit renormalised metric $\bslashgl{\infty,r}$ in $\mathrm{W}^{n,p}(\mathbb{S}^2)$, weakly in $\mathrm{W}^{n+1,p}(\mathbb{S}^2)$, and
\begin{align*}
	\Vert \bslashgl{\infty,r} - \circg \Vert^{n+1,p}
	\leq
	c(n,p) ( \epsilon + \delta_o + \ud_{o,\uh} ).
\end{align*}

\item
the renormalised Gauss curvature $\bKl{u,r}$ converges to the limit renormalised Gauss curvature $\bKl{\infty,r}$ in $\mathrm{W}^{n-1,p}(\mathbb{S}^2)$, weakly in $\mathrm{W}^{n,p}(\mathbb{S}^2)$, and
\begin{align*}
	\Vert \bKl{\infty,r} - 1 \Vert^{n,p}
	\leq
	c(n,p) ( \epsilon + \delta_o + \ud_{o,\uh} ).
\end{align*}
\end{enumerate}
\end{proposition}
\begin{proof}
\begin{enumerate}[label=\alph*.]
\item
Note $\ddpartial_u = \buL - \balu \db^i \ddpartial_i$, then
\begin{align*}
	\ddpartial_u \bslashgl{u,r}_{ij} 
	&= 
	\lie_{\buL} \bslashgl{u,r} - \lie_{\balu \db^i \ddpartial_i} \bslashgl{u,r}
\\
	&=
	\frac{2}{\br_u^2}
	\{ \frac{1}{2} ( \btr \buchilu - \frac{2}{\br_u}) \bslashgl{u}_{ij}
	+ \hatbuchil{u}_{ij}
	- \sym (\bslashd \balu \otimes \db)_{ij}
	- \balu \cdot \sym (\bslashnabla \db)_{ij} \}
\end{align*}
By estimates in assumption \ref{assum 8.1}, we have
\begin{align*}
	\Vert \ddpartial_u \bslashgl{u,r} \Vert^{n,p}
	\leq
	\frac{c(n,p)(\epsilon + \delta_o + \ud_{o,\uh})r_0}{\br_u^2}.
\end{align*}
Since $\frac{r_0}{\br_u^2}$ is integrable, then $\bslashgl{u,r}$ converges to $\bslashgl{\infty,r}$ in $\mathrm{W}^{n,p}(\mathbb{S}^2)$.

On the other hand, by estimates in assumption \ref{assum 8.1}, we have
\begin{align*}
	\Vert \circnabla \bslashgl{u,r} \Vert^{n,p}
	\leq
	c(n,p) (\epsilon + \delta_o + \ud_{o,\uh}),
\end{align*}
which implies that
\begin{align*}
	\Vert \bslashgl{u,r} - \circg \Vert^{n+1,p}
	\leq
	c(n,p) (\epsilon + \delta_o + \ud_{o,\uh}).
\end{align*}
Therefore the convergence in $\mathrm{W}^{n,p}(\mathbb{S}^2)$ and the above boundedness in $\mathrm{W}^{n+1,p}(\mathbb{S}^2)$ of $\bslashgl{u,r}$ implies the weak convergence in $\mathrm{W}^{n+1,p}(\mathbb{S}^2)$ of $\bslashgl{u,r}$.

\item
Consider the derivative of $\bKl{u,r}$
\begin{align*}
	\ddpartial_u \bKl{u,r}
	&=
	\buL \bKl{u,r} - \balu \db^i \ddpartial_i \bKl{u,r}
\\
	&=
	2 \br_u \cdot \bKlu + \br_u^2 \cdot \buL \bKlu - \br_u^2 \cdot \balu \db^i (\bslashd \bKlu)_i.
\end{align*}
We have that
\begin{align*}
	\buL \bKlu
	=
	\bslashdiv \bslashdiv \hatbuchilu
	- \frac{1}{2} \bslashDelta \btr \buchilu
	- \bKlu \btr \buchilu.
\end{align*}
Substituting equation \eqref{eqn 6.10} above
\begin{align*}
	\bslashdiv \hatbuchilu - \frac{1}{2} \bslashd \btr \buchilu
	&=
	\hatbuchilu \cdot \btalu  - \frac{1}{2} \btr \buchilu \btalu
	- \bubetalu,
\end{align*}
we obtain that
\begin{align*}
	\buL \bKlu
	&=
	\hatbuchilu (\btalu, \btalu)  - \frac{1}{2} \btr \buchilu \vert \btalu \vert_u^2
	 - \bubetalu \cdot \btalu
\\
	&\phantom{=}
	+ \hatbuchilu \cdot \bslashnabla \btalu  - \frac{1}{2} \btr \buchilu \cdot \bslashdiv \btalu
	- \bslashdiv \bubetalu
	- \bKlu \btr \buchilu.
\end{align*}
therefore
\begin{align*}
	\ddpartial_u \bKl{u,r}
	&=
	\br_u^2 
	\{ \hatbuchilu (\btalu, \btalu)  - \frac{1}{2} \btr \buchilu\ \vert \btalu \vert_u^2
	 - \bubetalu \cdot \btalu
\\
	&\phantom{= \br_u^2 \{ }
	+ \hatbuchilu \cdot \bslashnabla \btalu  - \frac{1}{2} \btr \buchilu\ \bslashdiv \btalu
	- \bslashdiv \bubetalu\}
\\ 
	&\phantom{=}
	- \br_u^2 ( \btr \buchilu - \frac{2}{\br_u}) \bKlu
	- \br_u^2 \cdot \balu \db^i (\bslashd \bKlu)_i.
\end{align*}
By estimates in assumption \ref{assum 8.1},
\begin{align*}
	\Vert \ddpartial_u \bKl{u,r} \Vert^{n-1,p}
	\leq
	\frac{c(n,p) (\epsilon + \delta_o + \ud_{o,\uh})r_0}{\br_u^2},
\end{align*}
then $\bKl{u,r}$ converges to $\bKl{\infty,r}$ in $\mathrm{W}^{n-1,p}(\mathbb{S}^2)$.

On the other hand, $\bKl{u,r}$ is bounded in $\mathrm{W}^{n,p}(\mathbb{S}^2)$ that
\begin{align*}
	\Vert \bslashd \bKl{u,r} \Vert^{n-1,p}
	\leq
	c(n,p) ( \epsilon + \delta_o + \ud_{o,\uh} ),
\end{align*}
thus it implies the weak convergence of $\bKl{u,r}$ in $\mathrm{W}^{n,p}(\mathbb{S}^2)$.
\end{enumerate}
\end{proof}

\subsection{Limit renormalised curvature component $\bP$ and shears $\itbSigma$, $\itbXi$}
We construct the limit renormalised curvature component $\bP$ for $\brholu$ and the renormalised shears $\itbSigma$, $\itbXi$ for $\hatbuchilu$, $\hatbchilu'$ at the past null infinity.
\begin{definition}\label{def 9.5}
Let $\{ \bSigma_u \}$ be a constant mass aspect function foliation on an incoming null hypersurface $\ucalH$ in theorem \ref{thm 8.2}. Define the renormalised curvature component $\brhol{u,r}$ and renormalised shear $\hatbchil{u,r}'$ at the past null infinity for the foliation as
\begin{align*}
	\brhol{u,r} = \br_u^3 \cdot \brholu,
	\quad
	\hatbchil{u,r}' = \br_u^{-1} \cdot \hatbchilu'.
\end{align*}
Define the limit renormalised curvature component $\bP$ and limit renormalised shears $\itbSigma$, $\itbXi$ at the past null infinity as
\begin{align*}
	\bP
	=
	\lim_{u\rightarrow \infty} \brhol{u,r},
	\quad
	\itbSigma
	=
	\lim_{u\rightarrow \infty} \hatbuchil{u},
	\quad
	\itbXi
	=
	\lim_{u\rightarrow \infty} \hatbchil{u,r}'.
\end{align*}
\end{definition}
Similarly as proposition \ref{prop 9.4}, we show the limit renormalised curvature component $\bP$ and shears $\itbSigma$, $\itbXi$ in definition \ref{def 9.5} are well defined.
\begin{proposition}\label{prop 9.6}
Along the constant mass aspect function foliation $\{\bSigma_u\}$ in theorem \ref{thm 8.2}, we have
\begin{enumerate}[label=\alph*.]
\item
the renormalised curvature component $\brhol{u,r}$ converges to the limit renormalised curvature component $\bP$ in $\mathrm{W}^{n-1,p}(\mathbb{S}^2)$, weakly in $\mathrm{W}^{n,p}(\mathbb{S}^2)$, and
\begin{align*}
	&
	\vert \bP - \bP_S \vert
	\leq
	c(n,p) ( \epsilon + \delta_o + \delta_m + \ud_{o,\uh} + \uslashd_m ) r_0,
	\quad
	\bP_S
	=
	-r_0,
\\
	&
	\Vert \slashd \bP \Vert^{n-1,p}
	\leq
	c(n,p) ( \epsilon + \delta_o + \ud_{o,\uh} ) r_0.
\end{align*}

\item
the shear $\hatbuchilu$ converges to the limit renormalised shear $\itbSigma$ in $\mathrm{W}^{n-1,p}(\mathbb{S}^2)$, weakly in $\mathrm{W}^{n,p}(\mathbb{S}^2)$, and
\begin{align*}
	\Vert \itbSigma \Vert^{n,p}
	\leq
	c(n,p) ( \epsilon + \delta_o + \ud_{o,\uh} ) r_0.
\end{align*}

\item
the renormalised shear $\hatbchil{u,r}'$ converges to the limit renormalised shear $\itbXi$ in $\mathrm{W}^{n-1,p}(\mathbb{S}^2)$, weakly in $\mathrm{W}^{n,p}(\mathbb{S}^2)$, and
\begin{align*}
	\Vert \itbXi \Vert^{n,p}
	\leq
	c(n,p) ( \epsilon + \delta_o + \ud_{o,\uh} ).
\end{align*}
\end{enumerate}
\end{proposition}
\begin{proof}
\begin{enumerate}[label=\alph*.]
\item
The derivative of $\brholu$ along $\{\bSigma_u\}$ is given by the following equation derived from the Bianchi equations, see proposition 1.2 in \cite{C2009}
\begin{align*}
	\buL (\brholu)
	+
	\frac{3}{2} \btr \buchilu \cdot \brholu
	-
	\bslashdiv \bubetalu
	-
	\btalu \cdot \bubetalu
	-
	\frac{1}{2} ( \hatbchilu', \bualphalu)_{\subbslashglu}
	=
	0.
\end{align*}
Then we derive the propagation equation of $\brhol{u,r}$,
\begin{align*}
	\buL (\brhol{u,r})
	=
	- \frac{3}{2} ( \btr \buchilu - \frac{2}{\br_u}) \brhol{u,r}
	+
	\br_u^3 [ \bslashdiv \bubetalu
	+
	\btalu \cdot \bubetalu
	+
	\frac{1}{2} ( \hatbchilu', \bualphalu)_{\subbslashglu} ].
\end{align*}
By assumption \ref{assum 8.1}, we have
\begin{align*}
	\Vert \buL (\brhol{u,r}) \Vert^{n-1,p}
	\leq
	\frac{c(n,p)(\epsilon +\delta_o + \ud_{o,\uh}) r_0^{\frac{3}{2}}}{\br_u^{\frac{3}{2}}}.
\end{align*}
Since $\ddpartial_u (\brhol{u,r}) = \buL (\brhol{u,r}) - \balu \db^i \ddpartial_i (\brhol{u,r})$, therefore
\begin{align*}
	\Vert \ddpartial_u (\brhol{u,r}) \Vert^{n-1,p}
	\leq
	\frac{c(n,p)(\epsilon +\delta_o + \ud_{o,\uh}) r_0^{\frac{3}{2}}}{\br_u^{\frac{3}{2}}}.
\end{align*}
Hence $\brhol{u,r}$ converges in $\mathrm{W}^{n,p}(\mathbb{S}^2)$. By the boundedness of $\brhol{u,r}$ in $\mathrm{W}^{n+1,p}(\mathbb{S}^2)$ from assumption \ref{assum 8.1},
\begin{align*}
	&
	\vert \brhol{u,r} + r_0 \vert
	\leq
	c(n,p) ( \epsilon + \delta_o + \delta_m + \ud_{o,\uh} + \uslashd_m ) r_0,
\\
	&
	\Vert \bslashd \brhol{u,r} \Vert^{n-1,p}
	\leq
	c(n,p) ( \epsilon + \delta_o + \ud_{o,\uh} ) r_0,
\end{align*}
we have that $\brhol{u,r}$ converges weakly in $\mathrm{W}^{n+1,p}(\mathbb{S}^2)$ and the limit satisfies the estimates in the proposition.

\item
$\hatbuchilu$ satisfies the following propagation equation
\begin{align*}
	\buL \hatbuchilu
	=
	2\buomegalu \hatbuchilu
	+ (\hatbuchilu, \hatbuchilu)_{\subbslashglu} \bslashglu
	- \bualphalu.
\end{align*}
By assumption \ref{assum 8.1}, we have
\begin{align*}
	\Vert \buL \hatbuchilu \Vert^{n,p}
	\leq
	\frac{c(n,p)(\epsilon + \delta_o + \ud_{o,\uh}) r_0^{\frac{3}{2}}}{\br_u^{\frac{3}{2}}}.
\end{align*}
Since $\ddpartial_u \hatbuchilu = \buL \hatbuchilu - \lie_{\balu \db^i \ddpartial_i} \hatbuchilu$, we have
\begin{align*}
	\Vert \ddpartial_u \hatbuchilu \Vert^{n-1,p}
	\leq
	\frac{c(n,p)(\epsilon + \delta_o + \ud_{o,\uh}) r_0^{\frac{3}{2}}}{\br_u^{\frac{3}{2}}}.
\end{align*}
Hence $\hatbuchilu$ converges in $\mathrm{W}^{n-1,p}(\mathbb{S}^2)$. By the boundedness of $\hatbuchilu$ in $\mathrm{W}^{n,p}(\mathbb{S}^2)$ from assumption \ref{assum 8.1},
\begin{align*}
	\Vert \hatbuchilu \Vert^{n,p}
	\leq
	c(n,p)(\epsilon + \delta_o + \ud_{o,\uh}) r_0,
\end{align*}
we have that $\hatbuchilu$ converges weakly in $\mathrm{W}^{n,p}(\mathbb{S}^2)$ and the limit satisfies the estimates in the proposition.

\item
$\hatbchilu'$ satisfies the following propagation equation
\begin{align*}
	\buL (\hatbchilu')
	&=
	-2 \buomegalu \hatbchilu'
	+
	\frac{1}{2} \btr \buchilu \cdot \hatbchilu' 
	-
	\frac{1}{2} \btr \bchilu' \cdot \hatbuchilu
	+
	(\hatbuchilu, \hatbchilu')_{\subbslashglu} \bslashglu
\\
	&\phantom{=}
	-
	2 (\widehat{\sym}\{\bslashnabla \btalu \} 
		+ \widehat{\btalu \otimes \btalu}).
\end{align*}
Then we derive the propagation equation of $\hatbchil{u,r}$,
\begin{align*}
	\buL (\hatbchil{u,r}')
	&=
	-2 \buomegalu \cdot \hatbchil{u,r}'
	+
	\frac{1}{2} (\btr \buchilu - \frac{2}{\br_u})  \cdot \hatbchil{u,r}' 
	-
	\frac{1}{2 \br_u} \btr \bchilu' \cdot \hatbuchilu
	+
	\frac{1}{\br_u} (\hatbuchilu, \hatbchilu')_{\subbslashglu} \bslashglu
\\
	&\phantom{=}
	-
	\frac{2}{\br_u} [\widehat{\sym}\{\bslashnabla \btalu \} 
		+ \widehat{\btalu \otimes \btalu}].
\end{align*}

By assumption \ref{assum 8.1}, we have
\begin{align*}
	\Vert \buL (\hatbchil{u,r}') \Vert^{n,p}
	\leq
	\frac{c(n,p)(\epsilon + \delta_o + \ud_{o,\uh}) r_0}{\br_u^2}.
\end{align*}
Since $\ddpartial_u (\hatbchil{u,r}') = \buL (\hatbchil{u,r}') - \lie_{\balu \db^i \ddpartial_i} \hatbchil{u,r}'$, we have
\begin{align*}
	\Vert \ddpartial_u (\hatbchil{u,r}') \Vert^{n-1,p}
	\leq
	\frac{c(n,p)(\epsilon + \delta_o + \ud_{o,\uh}) r_0}{\br_u^2}.
\end{align*}
Hence $\hatbchil{u,r}'$ converges in $\mathrm{W}^{n-1,p}(\mathbb{S}^2)$. By the boundedness of $\hatbuchilu$ in $\mathrm{W}^{n,p}(\mathbb{S}^2)$ from assumption \ref{assum 8.1},
\begin{align*}
	\Vert \hatbchil{u,r}' \Vert^{n,p}
	\leq
	c(n,p)(\epsilon + \delta_o + \ud_{o,\uh}),
\end{align*}
we have that $\hatbchil{u,r}'$ converges weakly in $\mathrm{W}^{n,p}(\mathbb{S}^2)$ and the limit satisfies the estimates in the proposition.
\end{enumerate}
\end{proof}

\subsection{Asymptotic reference frame and energy-momentum vector}
In this section, we introduce the asymptotic reference frame defined by an asymptotically round foliation and the formula of the corresponding energy-momentum vector.

\begin{definition}\label{def 9.7}
Let $\{\bSigma_u\}_{(u_a, +\infty)}$ be a foliation (not necessary being a constant mass aspect function foliation) on a null hypersurface $\ucalH$ near the past null infinity in a vacuum spacetime $(M,g)$. Let
\begin{align*}
	\phi:
	\quad
	(u_a ,+\infty) \times \mathbb{S}^2
	\rightarrow
	\bigcup_{(u_a, +\infty)} \{\bSigma_u\} \subset \ucalH
\end{align*}
be a parameterisation of $\{\bSigma_u\}_{(u_a, +\infty)}$. Consider the limit of the renormalised metric $\bslashgl{u,r} = \br_u^{-2} \cdot \bslashglu$ by identifying it with its pull back on $\mathbb{S}^2$ via the map $\phi|_u$. Suppose the limit renormalised metric 
\begin{align*}
	\bslashgl{\infty,r}
	=
	\lim_{u\rightarrow +\infty} \bslashgl{u,r}.
\end{align*}
We call the foliation $\{\bSigma_u\}_{(u_a, +\infty)}$ defining an asymptotic reference frame on $\ucalH$ at the past null infinity if the limit renormalised metric is round, i.e. $\{\bSigma_u\}_{(u_a, +\infty)}$ is asymptotically round.
\end{definition}

\begin{definition}
Let $\{\bSigma_u \}_{(u_a, +\infty)}$ be an asymptotically round foliation which defines an asymptotic reference frame $\gamma_{\{\bSigma_u\}}$ at the past null infinity. Let $\{Y_{1,i}\}_{i=1,2,3}$ be a $L^2(\mathbb{S}^2, \bslashgl{\infty,r})$ unit orthogonal basis of the spherical harmonics of degree one. Suppose that the limit renormalised curvature component $\bP$ and limit renormalised shears $\itbSigma$, $\itbXi$ defined as in definition \ref{def 9.5} exist. Then define the limit renormalised mass aspect function $\bN$ by
\begin{align*}
	\bN
	=
	- \bP
	- \frac{1}{2}(\itbSigma, \itbXi)_{\subbslashgl{\infty,r}}.
\end{align*}
The Bondi energy $E^{\gamma_{\{\bSigma_u\}}}$ corresponding to the asymptotic reference frame $\gamma_{\{\bSigma_u\}}$ is given by
\begin{align*}
	E^{\gamma_{\{\bSigma_u\}}}
	=
	\lim_{u \rightarrow +\infty} m_{H}(\bSigma_u)
	=
	\frac{1}{8\pi} \int_{\mathbb{S}^2} \bN \dvol_{\subbslashgl{\infty,r}}.
\end{align*}
The asymptotic linear momentum $\vec{P}^{\gamma_{\{\bSigma_u\}}} = (P_i^{\gamma_{\{\bSigma_u\}}})_{i=1,2,3}$ corresponding to the asymptotic reference frame $\gamma_{\{\bSigma_u\}}$ and the basis $\{Y_{1,i}\}_{i=1,2,3}$ is given by
\begin{align*}
	P_i^{\gamma_{\{\bSigma_u\}}}
	=
	\frac{1}{8\pi} \int_{\mathbb{S}^2} \bN \cdot (\sqrt{\frac{4\pi}{3}} Y_{1,i} ) \dvol_{\subbslashgl{\infty,r}}.
\end{align*}
If $\vec{P}^{\gamma_{\{\bSigma_u\}}}=0$, then the asymptotic reference frame $\gamma_{\{\bSigma_u\}}$ is an asymptotic center-of-mass frame, and the corresponding Bondi energy $E^{\gamma_{\{\bSigma_u\}}}$ is the Bondi mass $m_B(\ucalH)$ on $\ucalH$, which is the minimum of all the Bondi energies on $\ucalH$,
\begin{align*}
	E^{\gamma_{\{\bSigma_u\}}}
	=
	m_B(\ucalH)
	=
	\min \{ E^{\gamma} : \text{all the asymptotic reference frame $\gamma$ on $\ucalH$} \}.
\end{align*}
\end{definition}

\section*{Acknowledgements}
\addcontentsline{toc}{section}{Acknowledgements}
The author is supported by the National Natural Science Foundation of China under Grant No. 12201338. This paper refines the result in the author's thesis \cite{L2018} on the global existence of the constant mass aspect function foliation in a vacuum perturbed Schwarzschild spacetime. The author is grateful to Demetrios Christodoulou for his constant encouragement and generous guidance. The author also thanks Alessandro Carlotto and Lydia Bieri for many helps on the refinements of the manuscript.

\appendix
\section{Proof of proposition \ref{prop 3.3}}\label{appen prop 3.3}
It is sufficient to prove estimates \eqref{eqn 3.6} in order to complete the proof of proposition \ref{prop 3.3}.
\begin{proof}[Proof of estimates \eqref{eqn 3.6}]
By the bootstrap assumption, we have that $\uflt$ satisfies the following estimates in $[0,t_a]$
\begin{align*}
\Vert \slashd \uflt \Vert^{n,p} \leq \ud_o r_0,
\quad
\vert \overline{\uflt}^{\circg} \vert \leq \ud_m r_0.
\end{align*}
Estimates \eqref{eqn 3.6} shall follow from the estimates of $f$, $t b^i f_i$, $t e^i f_i$, $t \ue^i f_i$, $\uvarepsilon$, $b^i (\uflt)_i$, $e^i (\uflt)_i$, $\ue^i (\uflt)_i$. Then we shall estimate the background quantities, including $\vec{b}, \Omega, \slashg$ as well as their derivatives with respect to $\circnabla, \partial_s, \partial_{\us}$, restricted on any surface $S_t$. The following inequality is useful: let $a$ be such a background quantity, then for any $w\leq n$
\begin{align*}
&
\ddcircnabla_i a = \circnabla_i a + \uflt_i \cdot \partial_{\us} a + t f_i \cdot \partial_s a,
\\
&
\Vert \ddcircnabla a \Vert^{w,p} 
\leq
\Vert \circnabla_i a \Vert^{w,p}_{S_t}
+
c(n,p) \Vert \slashd \uflt \Vert^{n,p} \cdot \Vert \partial_{\us} a \Vert^{w,p}_{S_t}
+ 
c(n,p) t \Vert \slashd f \Vert^{n,p} \cdot \Vert \partial_s a \Vert^{w,p}_{S_t}
\\
&\phantom{\Vert \ddcircnabla a \Vert^{w,p} }
\leq
\Vert \circnabla_i a \Vert^{w,p}_{S_t}
+
c(n,p) \ud_o r_0 \cdot \Vert \partial_{\us} a \Vert^{w,p}_{S_t}
+ 
c(n,p) t \delta_o (r_0+\os) \cdot \Vert \partial_s a \Vert^{w,p}_{S_t}.
\end{align*}
Then we have the following estimates for background quantities: for any $w\leq n$, $k,l,m\geq 0$
\begin{gather*}
\begin{aligned}
& 
\Vert \circnabla^k \partial_{s}^l \vec{b} \Vert^{w+1,p}_{S_t} 
\leq 
\frac{c(n,p,k,l)\epsilon r_0^2}{(r_0+t\os)^{3+l}} ( \ud_o +  \ud_m ),
\quad
\Vert \circnabla^k \partial_s^l \partial_{\us}^m \vec{b} \Vert^{w+1,p}_{S_t} \leq \frac{c(n,p,k,l,m)\epsilon}{(r_0+t\os)^{3+l} r_0^{m-2}},
\end{aligned}
\\
\begin{aligned}
&
\Vert \circnabla^k \partial_s^l \partial_{\us}^m (\Omega-\Omega_S) \Vert^{w+1,p}_{S_t} 
\underset{m=0,1}{\leq} 
\frac{c(n,p,k,l,m) \epsilon r_0}{(r_0+t\os)^{1+l+m}},
\\
&
\Vert \circnabla^k \partial_s^l \partial_{\us}^m (\Omega-\Omega_S) \Vert^{w+1,p}_{S_t} 
\underset{m\geq2}{\leq }
\frac{c(n,p,k,l,m) \epsilon r_0}{(r_0+t\os)^{2+l}},
\end{aligned}
\\
\begin{aligned}
\Vert \circnabla^k (\slashg-\slashg_S) \Vert^{w+1,p}_{S_t} 
\leq 
c(n,p,k,l,m) \epsilon (r_0+t \os)^2,
\quad
\Vert \circnabla^k \partial_s^l \partial_{\us}^m (\slashg-\slashg_S) \Vert^{w+1,p}_{S_t} 
\leq 
\frac{c(n,p,k,l,m) \epsilon}{(r_0+t \os)^{l-1}r_0^{m-1}}.
\end{aligned}
\end{gather*}
and for the Schwarzschild background quantities $\Omega_S$, $\slashg_S$ and $l,m\geq 1$
\begin{gather*}
\begin{aligned}
&
\Vert \Omega_S \Vert^{w+1,p}_{S_t} 
\leq
c(n,p),
\quad
\Vert \partial_s^l \Omega_S \Vert^{w+1,p}_{S_t} 
\leq
\frac{c(n,p,l)r_0^2}{(r_0+t\os)^{2+m}}(\ud_o+\ud_m),
\\
&
\Vert \partial_s^l \partial_{\us}^m \Omega_S \Vert^{w+1,p}_{S_t} 
\leq
\frac{c(n,p,l,m) \epsilon r_0}{(r_0+t\os)^{1+l+m}},
\end{aligned}
\\
\begin{aligned}
&
\Vert \slashg_S \Vert^{w+1,p}_{S_t} 
\leq
c(n,p) (r_0+ t\os)^2,
\quad
\Vert \partial_{\us}^m \slashg_S \Vert^{w+1,p}_{S_t} 
\leq
c(n,p,l,m)t(\delta_o + \os/r_0)(r_0+t\os)^{2-m},
\\
&
\Vert \partial_s^l \slashg_S \Vert^{w+1,p}_{S_t} 
\leq
c(n,p,l) (r_0+ t\os)^{2-l},
\quad
\Vert \partial_s^l \partial_{\us}^m \slashg_S \Vert^{w+1,p}_{S_t} 
\leq
c(n,p,l,m)(r_0+t\os)^{2-l-m}.
\end{aligned}
\end{gather*}
Note that to obtain the norm $\Vert \circnabla^k \partial_s^l \partial_{\us}^m a \Vert^{w+1,p}$ of a background quantity $a$, it requires the $L^{\infty}$ bounds of $\circnabla^{k+k'} \partial_s^{l+l'} \partial_{\us}^{m+m'} a$ for $k'+l'+m'\leq m+1$.

We can use the above estimates of $\uflt$ and the background quantities to obtain the estimates of $f$, $t b^i f_i$, $t e^i f_i$, $t \ue^i f_i$, $\uvarepsilon$, $b^i (\uflt)_i$, $e^i (\uflt)_i$, $\ue^i (\uflt)_i$ listed as follows,
\begin{align*}
&
\Vert f \Vert^{n+2,p} \leq c r_0 (\delta_o + |\os|/r_0), 
\\
&
\Vert t b^i f_i \Vert^{n+1,p}
\leq
\frac{c(n,p) r_0^2}{(r_0+t\os)^2} \epsilon ( \ud_o +  \ud_m ) \delta_o  t,
\\
&
\Vert t e^i f_i \Vert^{n+1,p}
\leq
c(n,p) \delta_o^2 t^2,
\\
&
\Vert t \ue^i f_i \Vert^{n,p}, \Vert e^i\, \uflt_i \Vert^{n,p}
\leq
\frac{c(n,p) r_0}{r_0+t\os} \delta_o \ud_o t,
\\
&
\Vert \uvarepsilon \Vert^{n,p}, \Vert \ue^i\, \uflt_i \Vert^{n,p}
\leq
\frac{c(n,p) r_0^2}{(r_0+ t\os)^2} \ud_o^2,
\\
&
\Vert b^i\, \uflt_i \Vert^{n,p}
\leq
\frac{c(n,p)r_0^3}{(r_0+t\os)^3} \epsilon  ( \ud_o + \ud_m   ) \ud_o.
\end{align*}
Then the above gives the estimates of $\vert F \vert$, $\Xlt$ and $\relt$ in \eqref{eqn 3.6}. The derivations are straightforward, thus we mainly remark several points in the derivations.
\begin{enumerate}[label=\alph*.]
\item 
The top order terms determine the regularities of $\Xlt$ and $\relt$ are the following terms: the terms involving $\slashd \uflt$ for $\Xlt$, and the terms involving $\ddcircnabla^2 \uflt$ and $\ddcircnabla^3 f$ for $\relt$.

\item 
We use the following inequality to estimate the product of two terms:
\begin{align*}
\Vert h_1 h_2 \Vert^{m,p} \leq c(m,p) \Vert h_1 \Vert^{m,p} \Vert h_2 \Vert^{m,p},
\quad
m\geq 1, p>2 \text{ or } m\geq 2, 2\geq p >1.
\end{align*}
Thus $(n-1, p)$ lies in the above range where the above inequality holds.

\item
In the estimate of $\Vert \Xlt \Vert^{n,p}$, it is required the $L^{\infty}$ bounds of $\circnabla^k \partial_s^l \partial_{\us}^m a$, where $a=\vec{b}, \Omega, \slashg$ and $k+l+m \leq n$, which comes from the terms $\vec{b}, \Omega, \slashg$ in $\Xlt$.

\item
In the estimate of $\Vert \relt \Vert^{n-1,p}$, it is required the $L^{\infty}$ bounds of $\circnabla^k \partial_s^l \partial_{\us}^m a$, where $a=\vec{b}, \Omega, \slashg$ and $k+l+m \leq n+1$, which comes from the second order derivatives $\mathrm{D}^2 a$, $a=\vec{b}, \Omega, \slashg$, $D=\circnabla, \partial_s, \partial_{\us}$ in $\relt$.
\end{enumerate}
The above implies that proposition requires the $L^{\infty}$ bounds of the metric components $\vec{b}, \Omega, \slashg$ up to $n+1$-th order derivatives.
\end{proof}

\section{Some properties of rotational vector component and derivative}\label{appen r.v.}

We list the properties of the rotational vector component and the rotational vector derivative in the following.
\begin{enumerate}[label=\arabic*.]
\item
Contraction and inner product via rotational vector components:
\begin{align*}
&
\tr T = T_i^i = T_{R,i}^{R,i},
\quad
\omega(v) = \omega_i v^i = \omega_{R,i} v^{R,i},
\\
&
\circg(v, w) = \circg_{R,ij} v^{R,i} w^{R,j} = \sum_{i=1,2,3} v^{R,i} w^{R,j},
\\
&
\circg^{-1} (\omega, \sigma) = (\circg^{-1})^{R,ij} \omega_{R,i} \sigma_{R,j} = \sum_{i=1,2,3} \omega_{R,i} \sigma_{R,j}.
\end{align*}

\item
Lie bracket and covariant derivative of the rotational vector field:
\begin{align*}
[R_i, R_j] = - \sum_{k=1,2,3} \epsilon_{ijk} R_k,
\quad
\circnabla_{R_i} R_j = - \sum_{k,l=1,2,3} \epsilon_{ikl} x_j x_k R_l.
\end{align*}

\item
Lie derivative with respect to the rotational vector field in terms of the rotational vector component. For a vector field $v$ and a one-form $\omega$,
\begin{align*}
&
(\lie_{R_i} v)^{R,j} 
= 
[R_i, v]^{R,j} 
= 
R_i(v^{R,j}) + \sum_{k=1,2,3} \epsilon_{ijk} v^{R,k},
\\
&
(\lie_{R_i} \omega)_{R,j} 
=
R_i(\omega_{R,j})
+
\sum_{k=1,2,3} \epsilon_{ijk} \omega_{R,k}.
\end{align*}
For a general tensor field $T$,
\begin{align*}
( \lie_{R_i} T )_{R,j_1 \cdots j_p}^{R, k_1 \cdots k_q}
=
R_i ( T_{R, j_1 \cdots j_p}^{R, k_1 \cdots k_q}  )
+
\sum_{\substack{s=1,\cdots,p \\ t=1,2,3}} \epsilon_{i j_s t} T_{R,j_1 \cdots \underset{\hat{j_s}}{t} \cdots j_p}^{R, k_1 \cdots k_q}
+
\sum_{\substack{l=1,\cdots,q \\ r=1,2,3}} \epsilon_{i k_l r} T_{R,j_1 \cdots j_p}^{R, k_1 \cdots \overset{\hat{k_l}}{r} \cdots k_q}.
\end{align*}

\item
Covariant derivative with respect to the rotational vector field in terms of the rotational vector component. For a vector field $v$ and a one-form $\omega$,
\begin{align*}
&
( \circnabla_{R_i} v )^{R,j}
= R_i ( v^{R,j} ) + \sum_{k,l=1,2,3} \epsilon_{ikl} x_j x_k v^{R,l},
\\
&
( \circnabla_{R_i} \omega )_{R,j}
= R_i ( \omega_{R,j} ) + \sum_{k,l=1,2,3} \epsilon_{imn} x_j x_k \omega_{R,l},
\end{align*}
and for a general tensor field $T$,
\begin{align*}
( \circnabla_{R_i} T )_{R,j_1 \cdots j_p}^{R, k_1 \cdots k_q}
&=
R_i ( T_{R, j_1 \cdots j_p}^{R, k_1 \cdots k_q}  )
+
\sum_{\substack{s=1,\cdots,p \\ m,n=1,2,3}} \epsilon_{i m t} x_{j_s} x_m T_{R,j_1 \cdots \underset{\hat{j_s}}{n} \cdots j_p}^{R, k_1 \cdots k_q}
\\
&\phantom{=}
+
\sum_{\substack{l=1,\cdots,q \\ m, n=1,2,3}} \epsilon_{i m r} x_{k_l} x_m T_{R,j_1 \cdots j_p}^{R, k_1 \cdots \overset{\hat{k_l}}{n} \cdots k_q}.
\end{align*}
By induction arguments, for the higher order covariant derivatives of $T$,
\begin{align*}
	(\circnabla^l_{R,i_1\cdots i_l} T)_{R,j_1\cdots j_p}^{R,k_1,\cdots k_q}
	&=
	R_{i_1} \cdots R_{i_l} (T_{R,j_1\cdots j_p}^{R,k_1,\cdots k_q})
\\
	&\phantom{=}
	+ 
	\sum_{
		\substack{\{i_1, \cdots, i_s \} \subset \{1,\cdots, l \}
			\\
			n_1,\cdots,m_1,\cdots = 1,2,3
		}
	} P_{2(l-s)} (x_1,x_2,x_3) \cdot R_{i_{r_1}} \cdots R_{i_{r_s}} (T_{R,n_1\cdots n_p}^{R,m_1,\cdots m_q}),
\end{align*}
where $P_{2(l-s)}(x_1,x_2,x_3)$ is a polynomial of $x_1,x_2, x_3$ of degree $2(l-s)$.

\item
Hessian of a function in terms of the rotational vector field,
\begin{align*}
( \circnabla^2 f )_{R,ij}
=
\circnabla^2 f (R_i, R_j)
=
R_i R_j f + \sum_{k,l=1,2,3} \epsilon_{ikl} x_j x_k R_l f.
\end{align*}
Higher order covariant derivatives of a function in terms of the rotational vector field,
\begin{align*}
	(\circnabla^l f)_{R,i_1\cdots i_l}
	=
	R_{i_1} \cdots R_{i_l} f
	+ 
	\sum_{\{i_1, \cdots, i_s \} \subset \{1,\cdots, l \}} 
	P_{2(l-s)} (x_1,x_2,x_3) \cdot R_{i_{r_1}} \cdots R_{i_{r_s}} f.
\end{align*}

\end{enumerate}

We introduce the Sobolev norm in terms of the rotational vector field derivative. Let $f$ be a function on $\Sigma$, then define the norm
\begin{align*}
\Vert f \Vert^{n,p}_{R} 
= 
\sum_{\substack{k=0,\cdots, n\\ i_1, \cdots ,i_k =1,2,3}}
\big( 
\int_{\Sigma} 
\vert R_{i_1} \cdots R_{i_k} f \vert^{p} \dvol_{\circg} 
\big)^{1/p},
\quad
p\geq 1.
\end{align*}
For a general tensor field $T$, define the norm
\begin{align*}
\Vert T \Vert^{n,p}_{R} 
= 
\sum_{\substack{l=0,\cdots, n\\ i_1, \cdots ,j_1 \cdots, k_1, \cdots=1,2,3}}
\big( 
\int_{\Sigma} 
\vert R_{i_1} \cdots R_{i_l} ( T_{R, j_1 \cdots j_s}^{R, k_1 \cdots k_r}  ) \vert^{p} \dvol_{\circg} 
\big)^{1/p},
\quad
p\geq 1.
\end{align*}
There exists a constant $c(n,p,r,s)$ such that for any $(r,s)$-type tensor field $T$
\begin{align*}
c(n,p,r,s)^{-1}\Vert T \Vert_R^{n,p} 
\leq 
\Vert T \Vert^{n,p} 
\leq
c(n,p,r,s) \Vert T \Vert_R^{n,p}.
\end{align*}

\section{Propagation lemma}\label{appen p.l.}
We consider the following propagation equation on $[0,+\infty) \times \mathbb{S}^2$
\begin{align}
	\partial_u \flu + \Xlu^i \partial_i \flu +  \lambda \Mlu \flu = \relu,
	\label{eqn c.1}
\end{align}
where $\lambda$ is a constant, $\Xlu^i \partial_i$ is a tangental vector field of $\{ u \} \times \mathbb{S}^2$ and $\Mlu$, $\relu$ are functions on $\{ u \} \times \mathbb{S}^2$. We shall obtain the estimate of the solution $\flu$ under the propagation of equation \eqref{eqn c.1}. In the following, use $r_u = r_0 + u$ to simplify the notation.

We state a lemma on the propagation of the volume form invariant under the flow generated by the vector field $\partial_u + \Xlu^i \partial_i$.
\begin{lemma}\label{lem c.1}
Define the following tangential volume form $\dvol_u$ on $[0,+\infty) \times \mathbb{S}^2$ invariant under the flow generated by the vector field $\partial_u + \Xlu^i \partial_i$,
\begin{align*}
	\dvol_{u=0} = \dvol_{\circg},
	\quad
	\lie_{\partial_u + \Xlu^i \partial_i} \dvol_u = 0.
\end{align*}
There exists a function $\philu$ on $[0,+\infty) \times \mathbb{S}^2$ such that
\begin{align*}
	\dvol_u = \philu \cdot \dvol_{\circg}.
\end{align*}
$\philu$ satisfies the equation
\begin{align*}
	(\partial_u + \Xlu^i \partial_i) \log\philu 
	= 
	- \circdiv \Xlu.
\end{align*}
Suppose that
\begin{align*}
	\vert \circdiv \Xlu \vert
	\leq
	\frac{\delta r_0}{r_u^2},
\end{align*}
then
\begin{align*}
	\vert \log \philu \vert 
	\leq
	\frac{\delta u}{r_u}.
\end{align*}
Define the norm with respect to the volume form $\dvol_u$,  i.e. 
$\Vert f \Vert^p_{\dvol_u} = ( \int_{\mathbb{S}^2} f^p \dvol_u )^{\frac{1}{p}}$.
Then
\begin{align*}
	e^{-\frac{\delta u}{p r_u}} \Vert f \Vert^p_{\dvol_{\circg}}
	\leq
	\Vert f \Vert^p_{\dvol_u}
	\leq
	e^{\frac{\delta u}{p r_u}} \Vert f \Vert^p_{\dvol_{\circg}}
\end{align*}.
\end{lemma}

The propagation lemma states the estimate of $\flu$ via integrating the above equation.
\begin{lemma}\label{lem c.2}
Suppose that
\begin{align*}
	\Vert \Xlu \Vert^{n,p}
	\leq
	\frac{\delta r_0}{r_u^2},
	\quad
	\Vert \slashd \Mlu \Vert^{n-1,p}
	\leq
	\frac{\delta r_0}{r_u^2},
	\quad
	\vert \overline{M}^{\circg} -\frac{2}{r_u} \vert
	\leq
	\frac{\delta r_0}{r_u^2},
\end{align*}
where $n\geq 1+n_p$, $p >1$. Then there exists a constant $c(n,p)$ such that the solution $\flu$ of equation \eqref{eqn c.1} satisfies the following differential inequality
\begin{align*}
	\Big\vert \frac{\ed}{\ed u} (r_u^{2\lambda} \Vert \flu \Vert^{n,p}_{R,\dvol_u}) \Big\vert
	\leq
	\frac{c(n,p)(1 + \lambda) \delta r_0}{r_u^2} (r_u^{2\lambda} \Vert \flu \Vert^{n,p}_{R,\dvol_u}) 
	+ r_u^{2\lambda}\Vert \relu \Vert^{n,p}_{R,\dvol_u}.
\end{align*}
Thus
\begin{align*}
	r_u^{2\lambda} \Vert \flu \Vert^{n,p}_{R,\dvol_u}
	&\leq
	\int_0^u \frac{c(n,p)(1 + \lambda) \delta r_0}{r_u^2} (r_{u'}^{2\lambda}\Vert \fl{u'} \Vert^{n,p}_{R,\dvol_{u'}}) \ed u'
\\
	&\phantom{=\ }
	+ r_0^{2\lambda} \Vert \fl{u=0} \Vert^{n,p}_R
	+ \int_0^u r_{u'}^{2\lambda}\Vert \rel{u'} \Vert^{n,p}_{R,\dvol_{u'}} \ed u'.
\end{align*}
By Gronwall's inequality,
\begin{align*}
	r_u^{2\lambda} \Vert \flu \Vert^{n,p}_{R,\dvol_u}
	&\leq
	e^{\frac{c(n,p)(1 + \lambda)\delta u}{r_u}}
	(r_0^{2\lambda} \Vert \fl{u=0} \Vert^{n,p}_R
	+ \int_0^u r_{u'}^{2\lambda}\Vert \rel{u'} \Vert^{n,p}_{R,\dvol_{u'}} \ed u').
\end{align*}
Since there exists a constant $c$ that $\vert \circdiv \Xlu \vert \leq c \Vert \Xlu \Vert^{n,p} \leq \frac{c\delta r_0}{r_u^2}$, then by lemma \ref{lem c.1} we have
\begin{align*}
	r_u^{2\lambda} \Vert \flu \Vert^{n,p}_R
	&\leq
	e^{\frac{c(n,p)(1 + \lambda)\delta u}{r_u}} \cdot e^{\frac{2c\delta u}{p r_u}}
	(r_0^{2\lambda} \Vert \fl{u=0} \Vert^{n,p}_R
	+ e^{\frac{c\delta u}{p r_u}} \int_0^u  r_{u'}^{2\lambda}\Vert \rel{u'} \Vert^{n,p}_R \ed u').
\end{align*}
\end{lemma}
The proof of the above lemma relies on the following equation for the rotational derivatives of $\flu$ derived from \eqref{eqn c.1}
\begin{align}
	\begin{aligned}
		&\phantom{=\ }
		(\partial_u + \Xlu^i \partial_i) \flu_{R,i_1, \cdots, i_n}
		+ \lambda \Mlu \flu_{R,i_1, \cdots, i_n}
	\\
		&=
		\underbrace{
		-\mkern-20mu \sum_{
		\substack{
			\{ l_1, \cdots, l_s \} \cup \{k_1, \cdots, k_{n-s} \} = \{1,\cdots, n\}
		\\
			l_1 < \cdots < l_s, k_1 < \cdots < k_{n-s}
		\\
			n-s\leq n-1
		}}\mkern-80mu
		\Big\{
		[ R_{i_{l_1}}, \cdots R_{i_{l_s}}, \Xlu] (\flu_{R,i_{k_1}, \cdots, i_{k_{n-s}}})
		+
		\lambda \Mlu_{R,i_{l_1}, \cdots, i_{l_s}} \cdot \flu_{R,i_{k_1}, \cdots, i_{k_{n-s}}}
		\Big\}}_{Q_{R,i_1, \cdots, i_n}}
	\\
		&\phantom{=\ }
		+ \relu_{R,i_1, \cdots, i_n}.
	\end{aligned}
	\label{eqn c.2}
\end{align}

\begin{proof}
Rewrite equation \eqref{eqn c.2} as
\begin{align*}
	&\phantom{=\ }
	(\partial_u + \Xlu^i \partial_i) \flu_{R,i_1, \cdots, i_n}
	+ \frac{2\lambda}{r_u}\flu_{R,i_1, \cdots, i_n}
\\
	&=
	\lambda (\frac{2}{r_u} - \Mlu) \flu_{R,i_1, \cdots, i_n}
	+ Q_{R,i_1, \cdots, i_n}
	+  \relu_{R,i_1, \cdots, i_n}.
\end{align*}
Then calculate the Lie derivative of $(r_u^{2\lambda } \flu_{R,i_1,\cdots, i_n})^p \dvol_u$ w.r.t. $\partial_u + \Xlu^i \partial_i$,
\begin{align*}
	&\phantom{=\ }
	\lie_{\partial_u + \Xlu^i \partial_i} [ (r_u^{2\lambda }\cdot \flu_{R,i_1,\cdots, i_n})^p \dvol_u]
\\
	&=
	p (r_u^{2\lambda }\cdot \flu_{R,i_1,\cdots, i_n})^{p-1} \dvol_u 
	\cdot 
	r_u^{2\lambda}[ (\partial_u + \Xlu^i \partial_i) \flu_{R,i_1, \cdots, i_n}
	+ \frac{2\lambda}{r_u}\flu_{R,i_1, \cdots, i_n} ]
\\
	&=
	p (r_u^{2\lambda }\cdot \flu_{R,i_1,\cdots, i_n})^{p-1} 
	\cdot 
	r_u^{2\lambda}[ \lambda (\frac{2}{r_u} - \Mlu) \flu_{R,i_1, \cdots, i_n}
	+ Q_{R,i_1, \cdots, i_n}
	+  \relu_{R,i_1, \cdots, i_n}]
	\dvol_u.
\end{align*}
Therefore
\begin{align*}
	&\phantom{=\ }
	\Big\vert \frac{\ed}{\ed u} (r_u^{2\lambda} \Vert \flu \Vert^{n,p}_{R,\dvol_u})^p \Big\vert 
\\
	&\leq
	p (r_u^{2\lambda} \Vert \flu \Vert^{n,p}_{R,\dvol_u})^{p-1}
	\cdot
	[\frac{c(n,p) \lambda \delta r_0}{r_u^2} \cdot ( r_u^{2\lambda} \Vert \flu \Vert^{n,p}_{R,\dvol_u})
	 + \sum_{\{i_1, \cdots ,i_k\}_{k\leq n}}  r_u^{2\lambda} \Vert Q_{R,i_1,\cdots, i_k} \Vert^p_{R,\dvol_u}
\\
	&\phantom{\leq p (r_u^{2\lambda} \Vert \flu \Vert^{n,p}_{R,\dvol_u})^{p-1} \cdot [\ }
	 + r_u^{2\lambda} \Vert \relu_{R,i_1,\cdots, i_k} \Vert^{n,p}_{R,\dvol_u}]
\end{align*}
By Sobolev's inequality,
\begin{align*}
	\sum_{\{i_1, \cdots ,i_k\}_{k\leq n}}  r_u^{2\lambda} \Vert Q_{R,i_1,\cdots, i_k} \Vert^p_{R,\dvol_u}
	\leq
	\frac{c(n,p) \delta r_0}{r_u^2} \cdot (r_u^{2\lambda} \Vert \flu \Vert^{n,p}_{R,\dvol_u}).
\end{align*}
Hence we obtain
\begin{align*}
	\Big\vert \frac{\ed}{\ed u} ( r_u^{2\lambda} \Vert \flu \Vert^{n,p}_{R,\dvol_u}) \Big\vert
	\leq
	\frac{c(n,p)(1 + \lambda) \delta r_0}{r_u^2} \cdot ( r_u^{2\lambda} \Vert \flu \Vert^{n,p}_{R,\dvol_u})
	 + r_u^{2\lambda} \Vert \relu_{R,i_1,\cdots, i_k} \Vert^{n,p}_{R,\dvol_u}.
\end{align*}
The rest of the lemma follows.
\end{proof}

Lemma \ref{lem c.2} concerns about the Sobolev norm in terms of the rotational vector field derivatives. We can also work with the covariant derivatives. Recall that
\begin{align*}
	\circnabla^l_{R,i_1\cdots i_l} f
	=
	R_{i_1} \cdots R_{i_l} f
	+ 
	\sum_{\{i_1, \cdots, i_s \} \subsetneq \{1,\cdots, l \}} 
	P_{2(l-s)} (x_1,x_2,x_3) \cdot R_{i_{r_1}} \cdots R_{i_{r_s}} f,
\end{align*}
then
\begin{align*}
	&\phantom{=}
	(\partial_u + \Xlu^i \partial_i) (\circnabla^l_{R,i_1\cdots i_l} f) + \lambda \Mlu (\circnabla^l_{R,i_1\cdots i_l} f)
\\
	&=
	(\partial_u + \Xlu^i \partial_i) f_{R,i_1 \cdots i_l}  + \lambda \Mlu \cdot f_{R,i_1 \cdots i_l} 
\\
	&\phantom{=}
	+
	\sum_{\{i_1, \cdots, i_s \} \subset \{1,\cdots, l \}} 
	P_{2(l-s)} (x_1,x_2,x_3) \cdot [(\partial_u + \Xlu^i \partial_i) + \lambda \Mlu] f_{R,i_{r_1} \cdots i_{r_s}} 
\\
	&\phantom{=}
	+
	\sum_{\{i_1, \cdots, i_s \} \subset \{1,\cdots, l \}} 
	\Xlu^i \partial_i [P_{2(l-s)} (x_1,x_2,x_3)] \cdot f_{R,i_{r_1} \cdots i_{r_s}} 
\\
	&=
	Q_{R,i_1 \cdots i_l} + \relu_{r,i_1 \cdots i_l}
\\
	&\phantom{=}
	+
	\sum_{\{i_1, \cdots, i_s \} \subset \{1,\cdots, l \}} 
	P_{2(l-s)} (x_1,x_2,x_3) \cdot [Q_{R,i_{r_1} \cdots i_{r_s}}  + \relu_{R,i_{r_1} \cdots i_{r_s}}]
\\
	&\phantom{=}
	+
	\sum_{\{i_1, \cdots, i_s \} \subset \{1,\cdots, l \}} 
	\Xlu^i \partial_i [P_{2(l-s)} (x_1,x_2,x_3)] \cdot f_{R,i_{r_1} \cdots i_{r_s}}.
\end{align*}
Since
\begin{align*}
	\vert \circnabla^l f \vert^2
	=
	\sum_{i_1,\cdots,i_l=1,2,3} (\circnabla^l_{R,i_1\cdots i_l} f)^2,
\end{align*}
we have
\begin{align}
\begin{aligned}
	&\phantom{=}
	(\partial_u + \Xlu^i \partial_i) (\vert \circnabla^l f \vert^2) + 2\lambda \Mlu (\vert \circnabla^l f \vert^2)
\\
	&=
	2 [Q_{R,i_1 \cdots i_l} + \relu_{R,i_1 \cdots i_l}] \cdot \circnabla^l_{R,i_1\cdots i_l} f
\\
	&\phantom{=}
	+
	\sum_{\{i_1, \cdots, i_s \} \subsetneq \{1,\cdots, l \}} 
	2P_{2(l-s)} (x_1,x_2,x_3) \cdot [Q_{R,i_{r_1} \cdots i_{r_s}}  + \relu_{R,i_{r_1} \cdots i_{r_s}}] \cdot \circnabla^l_{R,i_1\cdots i_l} f
\\
	&\phantom{=}
	+
	\sum_{\{i_1, \cdots, i_s \} \subsetneq \{1,\cdots, l \}} 
	2 \Xlu^i \partial_i [P_{2(l-s)} (x_1,x_2,x_3)] \cdot f_{R,i_{r_1} \cdots i_{r_s}} \cdot \circnabla^l_{R,i_1\cdots i_l} f.
\end{aligned}
\label{eqn c.3}
\end{align}
Therefore similar to lemma \ref{lem c.2}, we otain the following estimate for $\Vert \flu \Vert^{n,p}$.
\begin{lemma}\label{lem c.3}
Suppose that
\begin{align*}
	\Vert \Xlu \Vert^{n,p}
	\leq
	\frac{\delta r_0}{r_u^2},
	\quad
	\Vert \slashd \Mlu \Vert^{n-1,p}
	\leq
	\frac{\delta r_0}{r_u^2},
	\quad
	\vert \overline{M}^{\circg} -\frac{2}{r_u} \vert
	\leq
	\frac{\delta r_0}{r_u^2},
\end{align*}
where $n\geq 1+n_p$, $p >1$. Then there exists a constant $c(n,p)$ such that the solution $\flu$ of equation \eqref{eqn c.1} satisfies the following differential inequality
\begin{align*}
	\Big\vert \frac{\ed}{\ed u} (r_u^{2\lambda} \Vert \flu \Vert^{n,p}_{\dvol_u}) \Big\vert
	\leq
	\frac{c(n,p)(1 + \lambda) \delta r_0}{r_u^2} (r_u^{2\lambda} \Vert \flu \Vert^{n,p}_{\dvol_u}) 
	+ c(n,p) r_u^{2\lambda}\Vert \relu \Vert^{n,p}_{\dvol_u}.
\end{align*}
Thus
\begin{align*}
	r_u^{2\lambda} \Vert \flu \Vert^{n,p}_{\dvol_u}
	&\leq
	\int_0^u \frac{c(n,p)(1 + \lambda) \delta r_0}{r_u^2} (r_{u'}^{2\lambda}\Vert \fl{u'} \Vert^{n,p}_{\dvol_{u'}}) \ed u'
\\
	&\phantom{=\ }
	+ r_0^{2\lambda} \Vert \fl{u=0} \Vert^{n,p}
	+ c(n,p) \int_0^u r_{u'}^{2\lambda}\Vert \rel{u'} \Vert^{n,p}_{\dvol_{u'}} \ed u'.
\end{align*}
By Gronwall's inequality,
\begin{align*}
	r_u^{2\lambda} \Vert \flu \Vert^{n,p}_{\dvol_u}
	&\leq
	e^{\frac{c(n,p)(1 + \lambda)\delta u}{r_u}}
	(r_0^{2\lambda} \Vert \fl{u=0} \Vert^{n,p}
	+ c(n,p) \int_0^u r_{u'}^{2\lambda}\Vert \rel{u'} \Vert^{n,p}_{\dvol_{u'}} \ed u').
\end{align*}
Since there exists a constant $c$ that $\vert \circdiv \Xlu \vert \leq c \Vert \Xlu \Vert^{n,p} \leq \frac{c\delta r_0}{r_u^2}$, then by lemma \ref{lem c.1} we have
\begin{align*}
	r_u^{2\lambda} \Vert \flu \Vert^{n,p}
	&\leq
	e^{\frac{c(n,p)(1 + \lambda)\delta u}{r_u}} \cdot e^{\frac{2c\delta u}{p r_u}}
	(r_0^{2\lambda} \Vert \fl{u=0} \Vert^{n,p}
	+ c(n,p) e^{\frac{c\delta u}{p r_u}} \int_0^u  r_{u'}^{2\lambda}\Vert \rel{u'} \Vert^{n,p} \ed u').
\end{align*}
\end{lemma}

\section{Proof of local existence theorem \ref{thm 8.24}}\label{appen thm 8.24}
The proof is similar to the one of the local existence of the constant mass aspect function foliation in \cite{S2008}. 
\begin{proof}
Define the inverse lapse operator $\calA$ for the constant mass aspect function foliation by $\calA[f] = a_f$ where $\bSigma_f$ be the surface parameterised by $f$ in the $\{s,\vartheta\}$ coordinate system that $\bSigma_f = \{ (s,\vartheta): s= f(\vartheta)\}$,
\begin{align*}
	\left\{
	\begin{aligned}
		&
		\bslashDelta \log a_f
		=
		\dmu_{f} - \overline{\dmu_{f}}^{\subbslashg_f}
		=
		\drho_f - \overline{\drho_f}^{\subbslashg_f}
		-
		\frac{1}{2} [(\hatduchi_f, \hatdchi_f') -  \overline{(\hatduchi_f, \hatdchi_f')}^{\subbslashg_f} ]
		-
		\bslashdiv \deta_f,
	\\
		&
		\overline{a_f \cdot \dtr \duchi_f}^{\subbslashg_f}
		=
		\overline{\btr \buchi}^{\subbslashg_f}
		=
		\frac{2}{\br_{\bSigma_f}}.
	\end{aligned}
	\right.
\end{align*}
The following maps
\begin{align*}
	&
	\bslashg_f:
	\quad
	\mathrm{W}^{n+2,p}(\mathbb{S}^2)
	\rightarrow
	\mathrm{W}^{n+1,p}(\mathbb{S}^2),
	\quad
	f
	\mapsto
	\bslashg_f,
\\
	&
	\drho_f:
	\quad
	\mathrm{W}^{n+2,p}(\mathbb{S}^2)
	\rightarrow
	\mathrm{W}^{n,p}(\mathbb{S}^2),
	\quad
	f
	\mapsto
	\drho_f,
\\
	&
	(\hatduchi_f, \hatdchi_f'):
	\quad
	\mathrm{W}^{n+2,p}(\mathbb{S}^2)
	\rightarrow
	\mathrm{W}^{n,p}(\mathbb{S}^2),
	\quad
	f
	\mapsto
	(\hatduchi_f, \hatdchi_f'),
\\
	&
	\bslashdiv \deta_f:
	\quad
	\mathrm{W}^{n+2,p}(\mathbb{S}^2)
	\rightarrow
	\mathrm{W}^{n,p}(\mathbb{S}^2),
	\quad
	f
	\mapsto
	\bslashdiv \deta_f,
\\
	&
	\dtr \duchi_f:
	\quad
	\mathrm{W}^{n+2,p}(\mathbb{S}^2)
	\rightarrow
	\mathrm{W}^{n,p}(\mathbb{S}^2),
	\quad
	f
	\mapsto
	\dtr \duchi_f,
\end{align*}
are $C^1$. Therefore $\calA$ is $C^1$ as
\begin{align*}
	\calA:
	\quad
	\mathrm{W}^{n+2,p}
	\rightarrow
	\mathrm{W}^{n+2,p},
	\quad
	f
	\mapsto
	\calA[f] = a_f.
\end{align*}

We shall solve the following equation to construct the constant mass aspect function foliation,
\begin{align*}
	\partial_u \flu + a_{\flu} b^i \partial_i \flu
	=
	a_{\flu}.
\end{align*}
Constructing the following sequence $\{\flu_k\}_{k \in \mathbb{N}}$ to approximate the solution:
\begin{enumerate}[label=\roman*.]
\item
$\flu_a = \fl{u=u_a} + u-u_a$,
\item
Solving the following equation to construct $\flu_k$ successively,
\begin{align}
	\partial_u \flu_k + a_{\flu_{k-1}} b^i|_{\bSigma_{\flu_{k-1}}} \partial_i \flu_k
	=
	a_{\flu_{k-1}}.
	\label{eqn d.1}
\end{align}
\end{enumerate}
Suppose that $\fl{u=u_a}$ is bounded in $\mathrm{W}^{n+2,p}(\mathbb{S}^2)$ by
\begin{align*}
	\Vert \fl{u=u_a} \Vert^{n+2,p}
	\leq
	N.
\end{align*}
Then there exists a positive number $\tau$, which we shall determine later, such that $\{\flu_k(u,\cdot) \}$ weakly converges to a solution in $(u_a - \tau, u_a + \tau) \times \mathrm{W}^{n+2,p} (\mathbb{S}^2)$.

Introduce the notation
\begin{align*}
	\mathrm{W}^{n+2,p}_{\fl{u=u_a},\xi}(\mathbb{S}^2)
	=
	\{ f \in \mathrm{W}^{n+2,p}(\mathbb{S}^2): 
	\begin{aligned}
		&
		 \Vert f \Vert^{n+2,p} \leq \Vert \fl{u=u_a} \Vert^{n+2,p} +\xi,
	\\
		&
		\vert f- \fl{u=u_a} \vert \leq \xi.
	\end{aligned}
	\}.
\end{align*}
We have $\calA$ is bounded in $\mathrm{W}^{n+2,p}_{\fl{u=u_a},\xi}(\mathbb{S}^2)$. View the differential $\delta \calA$ as a map from $\mathrm{W}^{n+2,p}(\mathbb{S}^2)$ to the space of linear maps between $\mathrm{W}^{n+2,p}(\mathbb{S}^2)$. Then we also have $\delta \calA$ is bounded in $\mathrm{W}^{n+2,p}_{\fl{u=u_a},\xi}(\mathbb{S}^2)$. Suppose that
\begin{align*}
	\Vert \calA[f] \Vert^{n+2,p}
	\leq
	B_{\xi},
	\quad
	\Vvert \delta \calA \Vvert^{n+2,p},
	\Vvert \delta \calA \Vvert^{n+1,p}
	\leq
	B_{\xi},
	\quad
	f\in \mathrm{W}^{n+2,p}_{\fl{u=u_a},\xi}(\mathbb{S}^2),
\end{align*}
where $\Vvert \cdot \Vvert^{m,p}$ is the operator norm from $\mathrm{W}^{m,p}(\mathbb{S}^2)$ to itself.

To prove the convergence of $\{\flu_k \}$ on $(u_a-\tau, u_a-\tau)$, we show that 
\begin{enumerate}[label=\alph*.]
\item
$\{\flu_k \}$ is bounded in $(u_a - \tau, u_a + \tau) \times \mathrm{W}^{n+2,p} (\mathbb{S}^2)$,

\item
$\{\flu_k \}$ converges in $(u_a - \tau, u_a + \tau) \times \mathrm{W}^{n+1,p} (\mathbb{S}^2)$.
\end{enumerate}
We prove the above two assertions in the following. Use $C_{\xi}$ to denote a constant depending on $n,p$ and $\xi$ which is not necessary to be the same throughout the argument.
\begin{enumerate}[label=\alph*.]
\item
Choose a positive number $\xi$, we shall show that there exists $\tau>0$ such that $\flu_k \in \mathrm{W}^{n+2,p}_{\fl{u=u_a},\xi}(\mathbb{S}^2)$ for $u\in (u_a - \tau, u_a + \tau)$. This is proved by induction arguments. Suppose that $\flu_k \in \mathrm{W}^{n+2,p}_{\fl{u=u_a},\xi}(\mathbb{S}^2)$, then
\begin{align*}
	\Vert a_{\flu_k} \Vert^{n+2,p}
	\leq
	B_{\xi},
\end{align*}
and
\begin{align*}
	\Vert b^i|_{\bSigma_{\flu_k}} \Vert^{n+2,p}
	\leq
	C_{\xi},
\end{align*}
Therefore the vector field $a_{\flu_k} b^i|_{\bSigma_{\flu_k}} \partial_i$ is bounded in $\mathrm{W}^{n+2,p}$ that
\begin{align*}
	\Vert a_{\flu_k} b^i|_{\bSigma_{\flu_k}} \Vert^{n+2,p}
	\leq
	C_{\xi} B_{\xi}.
\end{align*}
Then integrating equation \eqref{eqn d.1} of $\flu_{k+1}$, we obtain that
\begin{align*}
	&\phantom{\leq}
	\Vert \flu_{k+1} \Vert^{n+2,p}
\\	
	&\leq
	\exp (C_{\xi} B_{\xi} |u-u_a| )
\\
	&\phantom{\leq}
	\cdot
	\big[\Vert \fl{u=u_a} \Vert^{n+2,p}
		+ C_{\xi} \int_{u_a}^u \Vert a_{\fl{u=u'}_k} \Vert^{n+2,p} \ed u'
		+ \int_{u_a}^u C_{\xi} B_{\xi} \Vert \fl{u=u'}_{k+1} \Vert^{n+2,p} \ed u'
	\big]
\end{align*}
Therefore by Gronwall's inequality, we obtain that
\begin{align}
\begin{aligned}
	\Vert \flu_{k+1} \Vert^{n+2,p}
	&\leq
	\exp (C_{\xi} B_{\xi} |u-u_a| )
	\cdot
	\big[\Vert \fl{u=u_a} \Vert^{n+2,p}
		+ C_{\xi} \int_{u_a}^u \Vert a_{\fl{u=u'}_k} \Vert^{n+2,p} \ed u'
	\big]
\\
	&\leq
	\exp (C_{\xi} B_{\xi} |u-u_a| )
	\cdot
	[\Vert \fl{u=u_a} \Vert^{n+2,p}
		+ C_{\xi} B_{\xi} |u-u_a|].
\end{aligned}
\label{eqn d.2}
\end{align}
Thus for $\tau$ sufficiently small depending on $n,p$, $\xi$ and $N$, we have
\begin{align*}
	\Vert \flu_{k+1} \Vert^{n+2,p}
	\leq
	\Vert \fl{u=u_a} \Vert^{n+2,p}
	+ \xi.
\end{align*} 
For $\vert \flu_{k+1} - \fl{u=0} \vert$, we have that
\begin{align*}
	\vert \flu_{k+1} - \fl{u=u_a} \vert
	&\leq
	\int_{u_a}^u \vert a_{\fl{u=u'}_k} + a_{\fl{u=u'}_k} b^i|_{\bSigma_{\fl{u=u'}_k}} \partial_i \fl{u=u'}_{k+1} \vert \ed u'
\\
	&\leq
	C_{\xi} B_{\xi} ( 1+ N + \xi) \cdot | u - u_a |.
\end{align*}
Therefore for $\tau$ sufficiently small depending on $n, p$, $\xi$ and $N$, we have
\begin{align*}
	\vert \flu_{k+1} - \fl{u=u_a} \vert
	\leq
	\xi,
\end{align*}
which implies that $\flu_{k+1} \in \mathrm{W}_{\fl{u=u_a}, \xi}^{n+2,p}(\mathbb{S}^2)$.

\item
Denote $A_{k+1} - A_k$ by $\triangle_{k} (A)$.
By equation \eqref{eqn d.1}, we obtain that
\begin{align*}
	\partial_u \triangle_{k+1} (\flu)
	+
	a_{\flu_{k+1}} b^i|_{\bSigma_{\flu_{k+1}}} \partial_i \triangle_{k+1} ( \flu)
	=
	- \triangle_k (a_{\flu} b^i|_{\bSigma_{\flu}}) \partial_i \flu_k
	+ \triangle_k ( a_{\flu}).
\end{align*}
For the right hand side, we have that
\begin{align*}
	&
	\Vert \triangle_k (a_{\flu} b^i|_{\bSigma_{\flu}}) \Vert^{n+1,p}
	\leq
	C_{\xi}(1+B_{\xi}) \Vert \triangle_k ( \flu) \Vert^{n+1,p},
\\
	&
	\Vert \triangle_k (a_{\flu} b^i|_{\bSigma_{\flu}}) \partial_i \flu_k \Vert^{n+1,p}
	\leq
	C_{\xi}(1+B_{\xi}) (N+\xi) \Vert \triangle_k ( \flu) \Vert^{n+1,p},
\\
	&
	\Vert \triangle_k (a_{\flu}) \Vert^{n+1,p}
	\leq
	B_{\xi} \Vert \triangle_k ( \flu) \Vert^{n+1,p},
\end{align*}
therefore
\begin{align*}
	&\phantom{\leq}
	\Vert
		\partial_u \triangle_{k+1} (\flu)
		+ a_{\flu_{k+1}} b^i|_{\bSigma_{\flu_{k+1}}} \partial_i \triangle_{k+1} ( \flu)
	\Vert^{n+1,p}
\\
	&\leq
	[ C_{\xi}(1+B_{\xi}) (N+\xi) + B_{\xi} ] \Vert \triangle_k ( \flu) \Vert^{n+1,p}.
\end{align*}
Integrating the equation of $\triangle_{k+1} ( \flu)$, we obtain that
\begin{align*}
	\sup_{|u - u_a| < \tau} \Vert \triangle_{k+1} ( \flu) \Vert^{n+1}
	&\leq
	\exp (C_{\xi} B_{\xi} \tau)
\\
	&\phantom{\leq}
	\cdot
	C_{\xi} [ C_{\xi}(1+B_{\xi}) (N+\xi) + B_{\xi} ]
	\sup_{|u - u_a| < \tau} \{ \Vert \triangle_k ( \flu) \Vert^{n+1,p} \}
	\cdot \tau.
\end{align*}
Therefore choose $\tau$ sufficiently small such that
\begin{align*}
	\exp (C_{\xi} B_{\xi} \tau)
	\cdot
	C_{\xi} [ C_{\xi}(1+B_{\xi}) (N+\xi) + B_{\xi} ] \cdot \tau
	<1,
\end{align*}
then $\{\flu_k\}_{k\in \mathbb{N}}$ converges in $(u_a - \tau, u_a + \tau) \times \mathrm{W}^{n+1,p} (\mathbb{S}^2)$.
\end{enumerate}

Finally the Sobolev norm $\Vert \flu \Vert^{n+2,p}$ is continuous w.r.t. $u$ and $\{ \flu \}$ is a continuous family in $\mathrm{W}^{n+1,p}(\mathbb{S}^2)$ follow from integrating equation \eqref{eqn d.1} and the bounds obtained above.
\end{proof}

\section{Proof of lemma \ref{lem 8.29}}\label{appen lem 8.29}

\begin{proof}
This lemma is also proved by a bootstrap argument, which relies on the estimates obtained in sections \ref{sec 8.4} - \ref{sec 8.9}. Firstly, by the continuity of $\flu$ in $\mathrm{W}^{n+2,p}(\mathbb{S}^2)$, all the norms appearing in assumption \ref{assum 8.1} are continuous.

Introduce the interval $I$ of $\tau$ where assumption \ref{assum 8.1} holds, i.e.
\begin{align*}
	I = \{ \tau: \text{ assumption \ref{assum 8.1} holds on } [u_a- \tau, u_a + \tau] \}.
\end{align*}
We shall show that there exist $\delta$ depending on $n,p$ and $\tau_{\epsilon'}$ depending on $n,p$, $s_0$, $\epsilon'$ such that $\tau_{\epsilon'} \in I$. 

We shall choose $\delta$ and the constants $c(n,p,\cdot)$ in lemma \ref{lem 8.23}. Then $\tau_{\epsilon'}$ is determined by imposing the condition that the following bootstrap argument holds on $[0,\tau_{\epsilon}']$. The bootstrap argument consists of two parts.
\begin{enumerate}[label=\alph*.]
\item
Closedness of $I \cap [0, \tau_{\epsilon'}]$, i.e. suppose that $\tau = \sup\{ I \cap [0, \tau_{\epsilon'}]\}$, then $\tau \in I \cap [0, \tau_{\epsilon'}]$. This is straightforward by passing the estimates in assumption \ref{assum 8.1} to the limit $\tau$.

\item
Openness of $I \cap [0,\tau_{\epsilon'}]$ in $[0, \tau_{\epsilon'}]$.
	\begin{enumerate}[label=\roman*.]
	\item
	Improvement of the estimates similar to step d.i described in the bootstrap argument in section \ref{sec 8.2}. If $\tau \in I\cap [0,\tau_{\epsilon'}]$, then the estimates in assumption \ref{assum 8.1} can be improved to the strict inequalities.
	
	\item
	Extension of the interval $[0,\tau] \subset I \cap [0,\tau_{\epsilon'}]$. Suppose $[0,\tau] \in  I \cap [0,\tau_{\epsilon'}]$, then there exists $\tau' > \tau$ such that $I \cap [0,\tau'] \subset I \cap [0,\tau_{\epsilon'}]$.
	\end{enumerate}
\end{enumerate}
The bootstrap argument implies that $[0,\tau_{\epsilon'}] \subset I$. The key step in the above bootstrap argument is step b.i.. Step d.ii. follows from the strict inequalities obtained in step d.i. and the continuity of the norms in assumption \ref{assum 8.1} w.r.t. $u$.

In the following, we find $\tau_{\epsilon'}$ such that step d.i. holds. Examining the proof of the estimates obtained in sections \ref{sec 8.4} - \ref{sec 8.8}, except the ones obtained by integrating the corresponding propagation equations, we see that all the estimates still hold. To be more precise, the following lemmas hold on $[u_a-\tau, u_a +\tau]$.
\begin{enumerate}[label={\Roman*}.]
\item
	\begin{enumerate}[label=\raisebox{0.1ex}{\scriptsize$\bullet$}]
	\item
	Lemma \ref{lem 8.4} for $\flu$, $\dslashd \uh$, $\circnabla^2 \uh$.
	\item
	Lemma \ref{lem 8.6} for $\bslashd \flu$.
	\end{enumerate}

\item
	\begin{enumerate}[label=\raisebox{0.1ex}{\scriptsize$\bullet$}]
	\item
	Lemma \ref{lem 8.9} for $\br_u$.
	\item
	Lemmas \ref{lem 8.10}, \ref{lem 8.11} for $\bslashglu$ and $\bcircnabla \bslashglu$.
	\end{enumerate}	
	
\item
	\begin{enumerate}[label=\raisebox{0.1ex}{\scriptsize$\bullet$}]
	\item
	Lemma \ref{lem 8.12} for $\overline{\btr \buchilu}^u$.
	\item
	Lemma \ref{lem 8.14} for $\hat \buchilu$.
	\item
	Lemma \ref{lem 8.17} for $\hatbchilu'$.
	\item
	Lemma \ref{lem 8.18} for $\btalu$.
	\item
	Lemma \ref{lem 8.19}, \ref{lem 8.20} for $\overline{\buomegalu}^u$ and $\bslashd \buomegalu$.
	\end{enumerate}

\item
Lemma \ref{lem 8.22} for $\bubetalu$, $\bslashd \brholu$, $\bsigmalu$
	
\end{enumerate}
For the the terms whose estimates are obtained through integrating the propagation equations, we also derive their estimates by integrating the equations from $u=u_a$.
\begin{enumerate}[label={\Roman*}.]
\item
$\pmb{\flu}$: as in the proof of lemma \ref{lem 8.5},
\begin{align*}
	\vert \flu - \flu_S \vert
	&\leq
	\vert \fl{u_a} - \fl{u_a}_S \vert
	+ \int_{u_a}^u \vert \bal{u'} - \bal{u'}_S \vert \ed u'
\\
	&\leq
	c(n,p, \flu) ( \epsilon + \delta_o + \delta_m + \ud_{o,\uh} + \uslashd_m ) \br_u 
	- \epsilon'
	+ c(n,p ) \vert u - u_a \vert.
\end{align*}

\item
	\begin{enumerate}[label=\raisebox{0.1ex}{\scriptsize$\bullet$}]
	\item
	$\pmb{\balu}$: we have the propagation equation
	\begin{align*}
		\buL (\balu - \balu_S)
		=
		2 \balu \buomegalu - 2 (\balu^2 \cdot \duomegalu - \balu_S^2 \cdot \duomegalu_S).
	\end{align*}
	By the bootstrap assumption \ref{assum 8.1} and the estimates in the proof of \ref{lem 8.7}, we have that
	\begin{align*}
		\vert \buL (\balu - \balu_S) \vert
		\leq
		c(n,p).
	\end{align*}
	Then as in the proof of lemma \ref{lem 8.7},
	\begin{align*}
		\vert \balu - \balu_S \vert
		&\leq
		\vert \bal{u=u-0} - \bal{u=u_a}_S \vert
		+ \int_{u_a}^u \vert \buL ( \bal{u'} - \bal{u'}_S) \vert \ed u'
	\\
		&\leq
		c(n,p,\balu) (\epsilon + \delta_o + \delta_m + \ud_{o,\uh} + \uslashd_m)
		- \epsilon'
		+ c(n,p) |u - u_a|.
	\end{align*}
	
	\item
	$\pmb{\bslashd \balu}$: by equation \eqref{eqn 6.4}, we have
	\begin{align*}
		(\ddpartial_u + \balu \db^i \ddpartial_i) \balu
		=
		\rel{\balu}
		=
		2\balu \buomegalu - 2 \balu^2 \cdot \duomegalu,
	\end{align*}
	thus by equation \eqref{eqn c.3}, we obtain that
	\begin{align*}
		&\phantom{=}
		(\ddpartial_u + \balu \db^i \ddpartial_i) (\vert \bcircnabla^l \balu \vert^2)
		\\
		&=
		2 [Q_{\barR,i_1 \cdots i_l} + \rel{\balu}_{\barR,i_1 \cdots i_l}] \cdot \circnabla^l_{\barR,i_1\cdots i_l} \balu
	\\
		&\phantom{=}
		+
		\sum_{\{i_1, \cdots, i_s \} \subsetneq \{1,\cdots, l \}} 
		2P_{2(l-s)} (x_1,x_2,x_3) \cdot [Q_{\barR,i_{r_1} \cdots i_{r_s}}  + \rel{\balu}_{\barR,i_{r_1} \cdots i_{r_s}}] \cdot \circnabla^l_{\barR,i_1\cdots i_l} \balu
	\\
		&\phantom{=}
		+
		\sum_{\{i_1, \cdots, i_s \} \subsetneq \{1,\cdots, l \}} 
		2 \balu \db^i \ddpartial_i [P_{2(l-s)} (x_1,x_2,x_3)] \cdot \balu_{\barR,i_{r_1} \cdots i_{r_s}} \cdot \circnabla^l_{\barR,i_1\cdots i_l} \balu,
	\end{align*}
	where
	\begin{align*}
		Q_{\barR,i_1 \cdots i_n}
		=
		-\mkern-20mu \sum_{
		\substack{
			\{ l_1, \cdots, l_s \} \cup \{k_1, \cdots, k_{n-s} \} = \{1,\cdots, n\}
		\\
			l_1 < \cdots < l_s, k_1 < \cdots < k_{n-s}
		\\
			n-s\leq n-1
		}}\mkern-80mu
		[ R_{i_{l_1}}, \cdots R_{i_{l_s}}, \balu \db^i] (\balu_{\barR,i_{k_1}, \cdots, i_{k_{n-s}}}).
	\end{align*}
	Therefore by the bootstrap assumption \ref{assum 8.1} and the estimates in the proof of lemma \ref{lem 8.8}, we obtain that
	\begin{align*}
		\frac{\ed}{\ed u} \Vert \bslashd \balu  \Vert^{n,p}
		\leq
		c(n,p) \Vert \bslashd \balu \Vert^{n,p} + c(n,p)
		\leq
		c(n,p),
	\end{align*}
	which implies that
	\begin{align*}
		\Vert \bslashd \balu  \Vert^{n,p}
		&\leq
		\Vert \bslashd \bal{u=u_a} \Vert^{n,p}
		+ c(n,p) |u- u_a| 
	\\
		&\leq
		c(n,p,\bslashd \balu) ( \epsilon + \delta_o + \ud_{o,\uh}) - \epsilon'
		+ c(n,p) |u- u_a| .
	\end{align*}
	\end{enumerate}
	
\item
	\begin{enumerate}[label=\raisebox{0.1ex}{\scriptsize$\bullet$}]
	\item
	$\pmb{\bslashd \btr \buchilu}$: we derive the following from propagation equation \eqref{eqn 6.8} as the derivation of equation \eqref{eqn c.3},
	\begin{align*}
		&\phantom{=}
		(\ddpartial_u + \balu \db^i \ddpartial_i) (\vert \bcircnabla^l \btr \buchilu \vert^2)
		\\
		&=
		2 [Q_{\barR,i_1 \cdots i_l} + \rel{\btr \buchilu}_{\barR,i_1 \cdots i_l}] \cdot \circnabla^l_{\barR,i_1\cdots i_l} \btr \buchilu
	\\
		&\phantom{=}
		+
		\sum_{\{i_1, \cdots, i_s \} \subsetneq \{1,\cdots, l \}} 
		2P_{2(l-s)} (x_1,x_2,x_3) 
	\\
		&\phantom{=+ \sum_{\{i_1, \cdots, i_s \} \subsetneq \{1,\cdots, l \}} }
		\cdot [Q_{\barR,i_{r_1} \cdots i_{r_s}}  + \rel{\btr \buchilu}_{\barR,i_{r_1} \cdots i_{r_s}}] \cdot \circnabla^l_{\barR,i_1\cdots i_l} \btr \buchilu
	\\
		&\phantom{=}
		+
		\sum_{\{i_1, \cdots, i_s \} \subsetneq \{1,\cdots, l \}} 
		2 \balu \db^i \ddpartial_i [P_{2(l-s)} (x_1,x_2,x_3)] \cdot (\btr \buchilu)_{\barR,i_{r_1} \cdots i_{r_s}} \cdot \circnabla^l_{\barR,i_1\cdots i_l} \btr \buchilu,
	\end{align*}
	where
	\begin{align*}
		&
		\rel{\btr \buchilu} 
		= 
		2 \buomegalu \btr \buchilu - \vert \hatbuchilu \vert^2 - \frac{1}{2} (\btr \buchilu)^2,
	\\
		&
		Q_{\barR,i_1 \cdots i_n}
		=
		-\mkern-20mu \sum_{
		\substack{
			\{ l_1, \cdots, l_s \} \cup \{k_1, \cdots, k_{n-s} \} = \{1,\cdots, n\}
		\\
			l_1 < \cdots < l_s, k_1 < \cdots < k_{n-s}
		\\
			n-s\leq n-1
		}}\mkern-80mu
		[ R_{i_{l_1}}, \cdots R_{i_{l_s}}, \balu \db^i] ((\btr \buchilu)_{\barR,i_{k_1}, \cdots, i_{k_{n-s}}}).
	\end{align*}
	Therefore by the bootstrap assumption \ref{assum 8.1} and the estimates in the proof of lemma \ref{lem 8.13}, we obtain that
	\begin{align*}
		\frac{\ed}{\ed u} \Vert \bslashd \btr \buchilu  \Vert^{n-1,p}
		\leq
		c(n,p) \Vert \bslashd \btr \buchilu \Vert^{n-1,p} + c(n,p)
		\leq
		c(n,p),
	\end{align*}
	which implies that
	\begin{align*}
		\Vert \bslashd \btr \buchilu  \Vert^{n-1,p}
		&\leq
		\Vert \bslashd \btr \buchil{u=u_a} \Vert^{n-1,p}
		+ c(n,p) |u-u_a|
	\\
		&\leq
		\frac{c(n,p,\bslashd \btr \buchilu) ( \epsilon + \delta_o + \ud_{o,\uh}) r_0}{\br_{u_a}^2} - \epsilon'
		+ c(n,p) |u-u_a|
	\\
		&\leq
		\frac{c(n,p,\bslashd \btr \buchilu) ( \epsilon + \delta_o + \ud_{o,\uh}) r_0}{\br_u^2} - \epsilon'
		+ c(n,p,r_0) |u-u_a|
	\end{align*}
	
	\item
	$\pmb{\overline{\btr \bchilu'}^u}$: as in the proof of lemma \ref{lem 8.15}, we have
	\begin{align*}
		\vert \frac{\ed}{\ed u} [\br_u ( \overline{\btr \bchilu'}^{u} - (\btr \bchilu')_S)] \vert
		&\leq
		\frac{C(n,p, \bmulu, \overline{\btr \buchilu}^{u})(\epsilon + \delta_o + \delta_m) r_0}{\br_u^2}
	\\
		&\phantom{\leq}
		+ \frac{C(n,p, \bmulu) \ud_{o,\uh} r_0}{\br_u^2}
	\\
		&\leq
		c(n,p).
	\end{align*}
	Therefore
	\begin{align*}
		\vert \br_u (\overline{\btr \bchilu'}^u -  (\btr \bchilu')_S) \vert
		&\leq
		\vert \br_{u_a} (\overline{\btr \bchil{u=u_a}'}^{u_a} -  (\btr \bchil{u=u_a}')_S) \vert
		+ c(n,p) |u-u_a|
	\\
		&\leq
		c(n,p,\overline{\btr \bchilu'}^{u}) ( \epsilon + \delta_o + \delta_m)
		+ \frac{c(n,p,\overline{\btr \bchilu'}^{u}) \ud_{o,\uh} u_a}{\br_{u_a}}
	\\
		&\phantom{\leq}
		- \epsilon' \br_{u_a}
		+ c(n,p) |u-u_a|
	\\
		&=
		c(n,p,\overline{\btr \bchilu'}^{u}) ( \epsilon + \delta_o + \delta_m)
		+ \frac{c(n,p,\overline{\btr \bchilu'}^{u}) \ud_{o,\uh} u}{\br_u}
	\\
		&\phantom{\leq}
		- \epsilon' \br_{u_a}
		+ c(n,p,r_0) |u-u_a|
	\end{align*}
	
	\item
	$\pmb{\bslashd \btr \bchilu'}$: we derive the following from propagation equation \eqref{eqn 6.9} as the derivation of equation \eqref{eqn c.3},
	\begin{align*}
		&\phantom{=}
		(\ddpartial_u + \balu \db^i \ddpartial_i) (\vert \bcircnabla^l \btr \bchilu' \vert^2)
		\\
		&=
		2 [Q_{\barR,i_1 \cdots i_l} + \rel{\btr \bchilu'}_{\barR,i_1 \cdots i_l}] \cdot \circnabla^l_{\barR,i_1\cdots i_l} \btr \bchilu'
	\\
		&\phantom{=}
		+
		\sum_{\{i_1, \cdots, i_s \} \subsetneq \{1,\cdots, l \}} 
		2P_{2(l-s)} (x_1,x_2,x_3) 
	\\
		&\phantom{=+ \sum_{\{i_1, \cdots, i_s \} \subsetneq \{1,\cdots, l \}} }
		\cdot [Q_{\barR,i_{r_1} \cdots i_{r_s}}  + \rel{\btr \bchilu'}_{\barR,i_{r_1} \cdots i_{r_s}}] \cdot \circnabla^l_{\barR,i_1\cdots i_l} \btr \bchilu'
	\\
		&\phantom{=}
		+
		\sum_{\{i_1, \cdots, i_s \} \subsetneq \{1,\cdots, l \}} 
		2 \balu \db^i \ddpartial_i [P_{2(l-s)} (x_1,x_2,x_3)] \cdot (\btr \bchilu')_{\barR,i_{r_1} \cdots i_{r_s}} \cdot \circnabla^l_{\barR,i_1\cdots i_l} \btr \bchilu',
	\end{align*}
	where
	\begin{align*}
		&
		\rel{\btr \bchilu'} 
		= 
		- 2 \buomegalu \btr \bchilu' 
		- \frac{1}{2} \btr \buchilu \cdot \btr \bchilu' 
		- 2 \vert \btalu \vert^2
		+ 2 \bmulu,
	\\
		&
		Q_{\barR,i_1 \cdots i_n}
		=
		-\mkern-20mu \sum_{
		\substack{
			\{ l_1, \cdots, l_s \} \cup \{k_1, \cdots, k_{n-s} \} = \{1,\cdots, n\}
		\\
			l_1 < \cdots < l_s, k_1 < \cdots < k_{n-s}
		\\
			n-s\leq n-1
		}}\mkern-80mu
		[ R_{i_{l_1}}, \cdots R_{i_{l_s}}, \balu \db^i] ((\btr \bchilu')_{\barR,i_{k_1}, \cdots, i_{k_{n-s}}}).
	\end{align*}
	Therefore by the bootstrap assumption \ref{assum 8.1} and the estimates in the proof of lemma \ref{lem 8.16}, we obtain that
	\begin{align*}
		\frac{\ed}{\ed u} \Vert \bslashd \btr \bchilu'  \Vert^{n-1,p}
		\leq
		c(n,p) \Vert \bslashd \btr \bchilu' \Vert^{n-1,p} + c(n,p),
	\end{align*}
	which implies that
	\begin{align*}
		\Vert \bslashd \btr \buchilu  \Vert^{n-1,p}
		&\leq
		\Vert \bslashd \btr \buchil{u=u_a} \Vert^{n-1,p}
		+ c(n,p) |u- u_a|
	\\
		&\leq
		\frac{c(n,p,\bslashd \btr \bchilu') ( \epsilon + \delta_o + \ud_{o,\uh}) r_0}{\br_u^2} - \epsilon'
		+ c(n,p,r_0) |u - u_a|.
	\end{align*}
	\end{enumerate}
\end{enumerate}
We can determine the constant $\tau_{\epsilon'}$ to ensure that step d.i. holds. It only requires to improve the estimates for the following terms.
\begin{enumerate}[label=\Roman*.]
\item
$\vert \flu - \flu_S \vert$:
$c(n,p) \tau_{\epsilon'} - \epsilon' <0.$

\item
	\begin{enumerate}[label=\raisebox{0.1ex}{\scriptsize$\bullet$}]
	\item
	$\vert \balu - \balu_S \vert$:
	$c(n,p) \tau_{\epsilon'} - \epsilon' <0.$
	
	\item
	$\Vert \bslashd \balu \Vert^{n+1,p}$:
	$c(n,p) \tau_{\epsilon'} - \epsilon' <0.$
	
	\end{enumerate}

\item
	\begin{enumerate}[label=\raisebox{0.1ex}{\scriptsize$\bullet$}]
	\item
	$\Vert \bslashd \balu \Vert^{n-1,p}$:
	$c(n,p,r_0) \tau_{\epsilon'} - \epsilon' <0.$
	
	\item
	$\vert \overline{\btr \bchilu'}^u - (\btr \bchilu')_S \vert$:
	$c(n,p,r_0) \tau_{\epsilon'} - \epsilon' \br_{u_a} <0.$
	
	\item
	$\Vert \bslashd \btr \bchilu' \Vert^{n-1,p}$:
	$c(n,p,r_0) \tau_{\epsilon'} - \epsilon' <0.$
	
	\end{enumerate}

\end{enumerate}
Clearly $\tau_{\epsilon'} <  \epsilon'  \min\{ c(n,p)^{-1}, c(n,p,r_0)^{-1}, \br_{u_a} c(n,p,r_0)^{-1} \}$ solves that above inequalities. Thus there exists $\tau_{\epsilon'}$ depending on $n,p,r_0, \epsilon'$ such that the bootstrap argument holds, then lemma \ref{lem 8.29} follows.
\end{proof}

\Address


\begin{thebibliography}{100}

\bibitem[Br01]{Br2001}
Bray, H.
Proof of the Riemannian Penrose inequality using the positive mass theorem,
\emph{Journal of Differential Geometry} \textbf{59} 2001, no. 2, 177-267.

\bibitem[C03]{C2003}
Christodoulou, D.
\emph{Mathematical Problems of General Relativity II (unpublished lecture notes)}.
Lectures at ETH Z\"urich during the winter semester 2003/2004.

\bibitem[C08]{C2008}
Christodoulou, D.
\emph{Mathematical Problems of General Relativity I}.
Z\"urich Lectures in Advanced Mathematics.
European Mathematical Society, Z\"urich, 2008.

\bibitem[C09]{C2009}
Christodoulou, D.
\emph{The formation of black holes in general relativity.}
EMS Monographs in Mathematics.
European Mathematical Society, Zürich, 2009.

\bibitem[CK93]{CK1993}
Christodoulou, D.; Klainerman S.
\emph{The Global Nonlinear Stability of the Minkowski Space},
Princeton Mathematical Series 41, Princetion University Press, 1993.

\bibitem[G73]{G1973}
Geroch, R.
Energy extraction,
\emph{Ann. New York Acad. Sci} \textbf{224} 1973, 108-117.

\bibitem[H68]{H1968}
Hawking, S. W.
Gravitational Radiation in an Expanding Universe,
\emph{J. Math. Phys.} \textbf{9} 1968, 598-604.

\bibitem[HI01]{HI2001}
Huisken, G.; Ilmanen, T.
The inverse mean curvature flow and the Riemannian Penrose inequality,
\emph{Journal of Differential Geometry} \textbf{59} 2001, no. 3, 353-437.

\bibitem[JW77]{JW1977}
Jang, P. S.; Wald, R. M.
 The positive energy conjecture and the cosmic censor hypothesis,
 \emph{J. Math. Phys.} \textbf{18} 1977, 41-44.
 
\bibitem[L18]{L2018}
Le, P.
\emph{The perturbation theory of null hypersurfaces and the weak null Penrose inequality},
DISS. ETH Nr. 25387.
\url{https://doi.org/10.3929/ethz-b-000334917}

\bibitem[L20]{L2020}
Le, P.
Marginally trapped surfaces in a perturbed Schwarzschild spacetime,
arXiv: 2007.06170v2 [math.DG].
https://doi.org/10.48550/arXiv.2007.06170
 
 \bibitem[L22]{L2022}
Le, P.
Global regular null hypersurfaces in a perturbed Schwarzschild black hole exterior,
\emph{Ann. PDE} \textbf{8} (2022), no. 2, Paper no. 13, 33 pp. 

\bibitem[L23]{L2023}
Le, P.
Linearised perturbation of constant mass aspect function foliation in Schwarzschild black hole spacetime.
\emph{Commun. Math. Phys.} (2023), 55 pp.

\bibitem[W81]{W1981}
Witten, E. 
A new proof of the positive energy theorem, 
\emph{Comm. Math. Phys.} \textbf{80}, 1981, 381-402.

\bibitem[S08]{S2008}
Sauter, J.
\emph{Foliations of Null Hypersurfaces and the Penrose Inequality},
Diss. ETH No.17842, 2008.
\url{https://doi.org/10.3929/ethz-a-005713669}

\bibitem[SY79]{SY1979}
Schoen, R.; Yau, S. T.
On the proof of the positive mass conjecture in general relativity,
\emph{Comm. Math. Phys.} \textbf{65} 1979, 45-76.

\end{thebibliography}
\end{document}